\documentstyle[twoside,12pt]{article}
\newtheorem{prop}{Proposition}[subsection]
\newtheorem{opr}[prop]{Definition}
\newtheorem{theo}[prop]{Theorem}

\newtheorem{conj}[prop]{Conjecture}
\newtheorem{rem}[prop]{Remark}
\newtheorem{coro}[prop]{Corollary}
\newtheorem{exam}[prop]{Example}
\newtheorem{ques}[prop]{Question}

\newtheorem{lem}[prop]{Lemma}

\newcommand{\proof}{{\rm\bf Proof. }}

\title{Habilitationsschrift\\
\ \\
Universit\"at-GHS-Essen\\
\ \\
\ \\
Hodge Theory of Hypersurfaces in Toric Varieties and Recent 
Developments in Quantum Physics \\  \ \\
Part I: \\ \ \\
Dual Polyhedra and Mirror Symmetry for Calabi-Yau Hypersurfaces in 
Toric Varieties}

\author{Victor V. Batyrev }

\begin{document}

\date{July 11, 1993}

\maketitle

\footnote{Supported by DFG, Forschungsschwerpunkt  Komplexe 
Mannigfaltigkeiten.}

\thispagestyle{empty}

\newpage

\thispagestyle{empty}

This Habilitationsschrift consists of two parts: 

\bigskip

\bigskip

\bigskip

\bigskip
\bigskip
\bigskip
{\bf Part I: Dual Polyhedra and Mirror Symmetry for Calabi-Yau Hypersurfaces in 
Toric Varieties (p. 3-65)}
\bigskip

\begin{center}

\bigskip

{\large  {\bf Abstract}}
\end{center}

We  consider   families ${\cal F}(\Delta)$ 
of complex $(n-1)$-dimensional projective algebraic 
varieties consisting of 
compactifications of $\Delta$-regular affine  
hypersurfaces $Z_f$ defined by Laurent polynomials $f$ with a    
fixed $n$-dimensional Newton polyhedron $\Delta$ in  
$n$-dimensional algebraic torus ${\bf T} = ({\bf C}^*)^n$. 
If a  Newton polyhedron  $\Delta$ gives rise to  
the  family  ${\cal F}(\Delta)$ of  
$(n-1)$-dimensional Calabi-Yau varieties then 
the dual polyhedron $\Delta^*$ in the dual space also gives rise to a 
family ${\cal F}(\Delta^*)$ of Calabi-Yau varieties, so that we 
obtain the remarkable duality between two different families of 
Calabi-Yau varieties.  
It is shown  that for $n -1 = 3$ the properties of this duality 
coincide with the properties of 
the {\em Mirror Symmetry}  discovered by physicists 
for Calabi-Yau $3$-folds. Moreover, our method allows to construct 
many new examples of Calabi-Yau $3$-folds and new candidates for their 
mirrors which  were previously unknown for physicists.  

We conjecture that there exists an isomorphism between 
two conformal field theories corresponding  to Calabi-Yau varieties 
from two families ${\cal F}(\Delta)$ and ${\cal F}(\Delta^*)$.  
\bigskip

\thispagestyle{empty}

\newpage

\tableofcontents  

\newpage

\section{Introduction}

\bigskip

\subsection{Motivation}

\hspace*{\parindent}

Calabi-Yau  $3$-folds   caught much attention from 
theoretical physics because of their relation  with the superstring 
theory \cite{witt}. Physicists discovered a duality  for  
Calabi-Yau $3$-folds, so called {\em Mirror Symmetry}
 \cite{aspin1}-\cite{aspin3}, \cite{cand.sch}, \cite{gree}. This duality 
defines a  correspondence between two topologically different 
Calabi-Yau $3$-folds 
$V$ and $V'$ such that the Hodge numbers of $V$ and $V'$ satisfy the 
relations
\begin{equation}
  h^{1,1}(V) = h^{2,1}(V'), \;\; h^{1,1}(V') = h^{2,1}(V). 
\label{main.relation}  
\end{equation}  
 
In the recent paper  of P. Candelas, X.C. de la Ossa, P.S. Green 
and L. Parkes \cite{cand2} the mirror symmetry  was applied to 
give predictions for the number of rational curves of various 
degrees on  general   quintic $3$-folds in ${\bf P}_4$. For small degrees  
these predictions 
were confirmed by algebraic geometers (see \cite{eil.st},\cite{s.kat}).  

In \cite{mor1} Morrison has presented 
a mathematical review of the calculation of 
physicists. Applying the same method based on 
a consideration of the Picard-Fuchs equation, he has found in  \cite{mor2} 
similar predictions for the number of rational curves on general  members of 
some other one-parameter families of Calabi-Yau $3$-folds. It 
turns out, for instance, that these predictions give 
the correct number of lines on  general double coverings  
of ${\bf P}_3$ ramified along smooth  hypersurfaces  
of degree $8$. The method of 
P. Candelas et al.  was also  applied to one-parameter families of 
Calabi-Yau $3$-folds by A. Font \cite{font}, A. Klemm and S.Theisen 
\cite{klemm2}.  
\bigskip

Let us give a summary of the method of  P. Candelas et al. in \cite{cand2} 
(see also \cite{mor1}, \cite{mor2}). 
\medskip

{\bf Step 1.} Instead of the $101$-dimensional 
family ${\cal F}$ of smooth  quintic $3$-folds, P. Candelas et al. 
considered the mirror  one-parameter family ${\cal F}'$ of Calabi-Yau 
$3$-folds with the Picard number $101$. The family ${\cal F}'$ is 
constructed as follows. 

Take a special one-parameter family ${\cal Q}_{\psi}$ of smooth 
quintic $3$-folds 
defined in ${\bf P}_4$ by the homogeneous equation equation 
\[ P_{\psi}(X) = \sum_{i =1}^5 X_i^5 - 5\psi \prod_{i =1}^5 =0 .\]
Let $\mu_5$ be the multiplicative group of $5$-th complex roots of unity. 
Denote by 
$G_0$  the subgroup of $(\mu_5)^5$ consisting of elements 
$(\alpha_1, \ldots, 
\alpha_5) \in (\mu_5)^5$ such that $\prod_{i=1}^5 \alpha_i =1$. 
Define $G$ as 
the quotient of $G_0$  by the  diagonal 
subgroup  $\mu_5  \subset G_0 \subset (\mu_5)^5$. The group $G$ is 
isomorphic to $(\mu_5)^3$. We obtain the natural action of $G$ on 
${\cal Q}_{\psi}$. The action of $G$ is free 
outside the union of $10$ curves 
\[ C_{ij} = \{ P_{\psi}(X) = X_i = X_j = 0, \; i \neq j  \}. \]
One should remark that the Euler characteristic  of 
\[ {\cal Q}_{\psi} \setminus \bigcup_{i \neq j}  C_{ij} \]
is zero. The quotient space ${\cal W}_{\psi} = 
{\cal Q}_{\psi} / G$ has only Gorenstein canonical 
singularities. There exists a resolution  
$\tilde{{\cal W}}_{\psi} \rightarrow  {\cal W}_{\psi}$ of 
singularities constructed  by  toric method 
\cite{mark} such that for any $\alpha \in \mu_5$ the isomorphism 
between ${\cal W}_{\psi}$ and ${\cal W}_{\alpha \psi}$ lifts to an 
isomorphism between the resolutions ${\tilde {\cal W}}_{\psi}$ 
and ${\tilde {\cal W}}_{\alpha \psi}$. The one-parameter 
family ${\cal F}' = {\tilde {\cal W}}_{\psi}$ of Calabi-Yau $3$-folds  
consists of the varieties with the Picard number $101$.      
\medskip

{\bf Step 2.} Periods of the 
holomorphic $3$-form $\omega_{\psi}$ on the one-parameter 
family ${\tilde {\cal W}}_{\psi}$ of Calabi-Yau $3$-folds satisfy  
the Picard-Fuchs equation of order $4$.  
The singular points of the Picard-Fuchs equations  are $\psi = 1 , \infty$. 
An explicit calculation of the coefficients of the 
Picard-Fuchs differential equation shows that the monodromy operator is 
maximally unipotent at $\psi = \infty$. 
This cusp point on the $1$-dimensional moduli space of ${\cal F}'$ is called 
the "large complex structure limit" of the family ${\tilde {\cal W}}_{\psi}$. 
\medskip 

{\bf Step 3.} A basis of solutions of the Picard-Fuchs equation 
can be represented 
by generalized hypergeometric functions. The Yukawa coupling $k_{ttt}$ 
is a special gauge-normalized solution  of the Picard-Fuchs equation. 
There exists  the natural choice of the parameter $t$ 
on the universal covering of the punctured 
disc around $\psi = \infty$ such that 
the $q$-expansion  $(q = e^{2\pi it})$  of the 
gauge-normalized Yukawa coupling  
\[ k_{ttt} = 5 + \sum_{j =1}^{\infty} \frac{n_j \jmath^3 q^j}{1 - q^j} \] 
has {\em integral coefficients} $n_j$. 
The numbers $n_j$ $(j \geq 1)$ are conjectured to be the numbers of rational 
curves of degree $j$ on general quintic $3$-folds in ${\bf P}_4$. 
\bigskip

The present paper is devoted  to the  mathematical nature of  Step 1. 
\bigskip

String theory associates a conformal field theory ${\rm K}(V)$ 
to a $d$-dimensional Calabi-Yau 
manifold $V$. Moreover, many  basic mathematical properties 
of $V$, e.g., the Hodge spaces $H^{p,q}$ in cohomology ring of $V$, 
find direct expressions in the associated conformal field theory ${\rm K}(V)$ 
\cite{greene}. In this 
case the conformal field theory ${\rm K}(V)$ is said to have a 
{\em geometric interpretation}. There exists a natural involution 
${\cal MIR}$ acting on  ${\rm K}(V)$ reversing all charges. 
If we assume for a moment that any  conformal 
field theory has a geometric interpretation (see \cite{gepn}), then 
geometrical interpretation of the 
involution ${\cal MIR}$ gives rise to a new type of duality 
for  Calabi-Yau manifolds. Thus we come to the {\em mirror conjecture} 
for Calabi-Yau manifolds. A version of this conjecture claims that 
there exists {\em another} Calabi-Yau manifold $V'$ such that 
\[ H^{p,q}(V') = H^{d-p,q}(V). \]

Although  a  general mathematical definition of the  mirror $V'$   
for arbitrary Calabi-Yau variety $V$  is not known (moreover, even 
the mirror conjecture lacks for a precise mathematical formulation), 
there are special classes of Calabi-Yau $3$-folds for which mirrors 
admit  explicit constructions stemming from special conformal field 
theories, e.g. so called  Gepner's minimal models   
(see \cite{berglund,gree}). All known examples of 
explicit mirror constructions deal with Calabi-Yau $3$-folds obtained 
from $3$-dimensional hypersurfaces in $4$-dimensional weighted 
projective spaces 
${\bf P}(\omega_0, \ldots, \omega_4)$ (see \cite{schimm}). 
There  have been found 
more than 7000 different possibilities for weights  
$(\omega_0, \ldots, \omega_4)$ 
defining  Calabi-Yau $3$-folds \cite{cand.sch,klemm1}. 
The classification 
of weighted hypersurfaces shows a striking symmetry of the set of all pairs  
$(h^{1,1}, h^{2,1})$ relative to the transposition interchanging two 
components $h^{1,1}$ and $h^{2,1}$.  This fact gives an empirical evidence 
in favor of the mirror conjecture.  
\medskip 

We make  several remarks. 

\begin{rem} {\rm 
It is well-known 
that any  $n$-dimensional algebraic manifold $V$ is always birationally 
isomorphic to a hypersurface $H$ in $(n+1)$-dimensional projective space. 
In this case, the hypersurface $H$ has usually very bad singularities 
and $V$ can be considered as a desingularization of $H$. On the other hand, 
Calabi-Yau $3$-folds are examples of minimal algebraic varieties 
in sence of Minimal Model Program of Mori \cite{kol1}. It follows from 
the result of Kollar \cite{kol2} that if two Calabi-Yau $3$-folds $V_1$ 
and $V_2$ are birationally isomorphic, then $V_1$ and $V_2$ must have 
{\em the same} Hodge numbers  $h^{1,1}$ and $h^{2,1}$. Therefore, 
if a singular hypersurface $H$ admits a Calabi-Yau desingularization, then 
the Hodge numbers $h^{1,1}$ and $h^{2,1})$ of this desingularization 
are uniquely defined. So resolutions of singular hypersurfaces give  a  
common way for construction  of  Calabi-Yau 3-folds and  classification 
of their Hodge numbers.}
\end{rem}

\begin{rem} 
{\rm A mathematical  definition of  the  mirror operation seems to exist  
for a more large category than the category of smooth algebraic 
varieties.  The example of the family $\{ {\cal W}_{\psi} \}$ of singular 
Calabi-Yau varieties shows that 
the class of smooth Calabi-Yau varieties is not  sufficiently large to 
be invariant under the mirror operation.  Even the 
existence of a good resolution ${\tilde {\cal W}}_{\psi}$ does not 
significantly improve the situation, since such a resolution is not unique.  
The Minimal Model Program for higher-dimensional algebraic varieties  
provides some classes of mild singularities, e.g.,  
canonical and terminal singularities   
\cite{kol1}.  So we come to the following definition.}
\end{rem}

\begin{opr}
{\rm A complex normal irreducible $n$-dimensional projective 
algebraic variety $W$ with only 
Gorenstein canonical singularities we call a {\em Calabi-Yau variety} 
if $W$ has trivial canonical 
bundle and 
\[ H^i (W, {\cal O}_W) = 0 \]
for $0 < i < n$.}
\label{calabi-yau}
\end{opr}
\medskip

\begin{rem} {\rm The set of all pairs of Hodge numbers $(h^{1,1}, h^{2,1})$ 
of  Calabi-Yau $3$-folds obtained from weighted 3-dimensional hypersurfaces is 
not {\em completely symmetric} with respect to the involution 
interchanging their components. So we can conclude that this class of 
Calabi-Yau 3-folds is not mirror-invariant. It is natural to assume 
that there exists a larger  
mirror-invariant class of  Calabi-Yau varieties containing it. 

We note also that all examples of Calabi-Yau $3$-folds 
constructed from  weighted hypersurfaces have an  open 
subset $Z$  isomorphic to a smooth hypersurface  in $4$-dimensional 
algebraic torus ${\bf T} \cong ({\bf C}^*)^4 \subset {\bf C}^4$. 
On the other hand, the hypersurface $Z$ is defined by a 
Laurent polynomial $f$. The Newton polyhedron $\Delta = \Delta(f)$ 
defines a standard projective compactification of ${\bf T}$, 
a toric variety ${\bf P}_{\Delta}$. If we assume that 
the coefficients of $f$ are sufficiently general, then the closure 
$\overline{Z}$ of $Z$ in ${\bf P}_{\Delta}$ is a projective variety 
with mild  singularities. These singularities are called toroidal and 
they can be resolved by combinatorial methods due to  Khovanski\^i \cite{hov}.
Thus, hypersurfaces in toric varieties form a larger class of algebraic 
varieities which contains weighted hypersurfaces in weighted 
projective spaces.}
\end{rem}

\begin{rem} 
{\rm From mathematical point of view, the construction of physicists 
of the mirror family ${\cal F}'$ as  a quotient of some subfamily of 
the original family of Calabi-Yau varieties ${\cal F}$ (see Step 1) is 
not satisfactory. We can not  apply it  twice in order to  get the original 
family of Calabi-Yau varieties ${\cal F} = ({\cal F}')'$. }
\end{rem} 
\bigskip

In this paper we consider families ${\cal F}(\Delta)$ of 
Calabi-Yau hypersurfaces  which are compactifications 
in $n$-dimensional projective toric 
varieties ${\bf P}_{\Delta}$ of smooth affine hypersurfaces  
whose equations have a fixed Newton polyhedron $\Delta$ and 
sufficiently general coefficients. 
Singularities of these Calabi-Yau varieties 
are toroidal, i.e. they are analytically isomorphic to toric singularities.

The main purpose  is to give a mathematical method for construction of 
candidates for mirrors of Calabi-Yau  hypersurfaces in toric varieties. 
This  method shows that Calabi-Yau hypersurfaces in toric varieties 
form a mirror-invariant  class. On the other hand,  we 
obtain  a simple {\em  precise mathematical interpretation}  of the mirror 
involution  in terms of the classical involutive duality for convex sets.    
\bigskip

\subsection{Outline of paper}

\hspace*{\parindent}

Section  2 is devoted to two  basic terminologies and some  well-known 
results on toric varieties.  We use  two  definition of toric 
varieties. This definitions  correspond to two kind of data, contravariant and 
covariant ones. 
First definition bases on the construction of an 
$n$-dimensional projective toric variety ${\bf P}_{\Delta}$ associated with 
an $n$-dimensional integral polyhedron  
$\Delta \subset M_{\bf Q}=M \otimes {\bf Q}$, where $M$ is a free 
abelian group of rank $n$. Second one describes an abstract 
complete toric variety ${\bf P}_{\Sigma}$ by means of a complete fan 
$\Sigma \subset N_{\bf Q} = N \otimes {\bf Q}$ of rational 
polyhedral cones, where $N = {\rm Hom}\, (M, {\bf Z})$ is the 
dual abelian group. Toric Fano varieties form a subcass of all toric varieties 
for which both definitions become dual to each other. 
\medskip

In section 3 we consider general properties of hypersurfaces 
in toric varieties of both types ${\bf P}_{\Delta}$  and ${\bf P}_{\Sigma}$. 
We will always assume that these hypersurfaces  satisfy some regularity 
conditions which imply that their singularities are induced only by 
singularities of embient toric varieties. Thus, a resolution of singularities 
of ${\bf P}_{\Delta}$ authomatically gives  rise to a resolution of 
singularities of {\em all} $\Delta$-regular hypersurfaces, i.e., we 
obtain a {\em simultanious resolution}  of all members of the family 
${\cal F}(\Delta)$. 
\medskip

In section 4 we give a  simple combinatorial 
criterion for the following question.

\begin{ques}
Which lattice  polyhedra  $\Delta$ give rise to families of 
Calabi-Yau hypersurfaces in ${\bf P}_{\Delta}$$?$ 
\label{q.delta}
\end{ques}

The polyhedra $\Delta$ which give rise to Calabi-Yau families 
we call {\em reflexive polyhedra}. An integral polyhedron $\Delta$ 
is reflexive if and only if ${\bf P}_{\Delta}$ is  a 
toric Fano variety with only Gorenstein singularities.  

There exists the following crusial observation.
\bigskip

{\em It follows from the combinatorial characterization of the reflexive 
polyhedron $\Delta \subset M_{\bf Q}$ 
the corresponding dual polyhedron $\Delta^*$ in dual space $N_{\bf Q}$ 
is also reflexive.}
\bigskip

So the set of all $n$-dimensional reflexive polyhedra admits 
the involution ${MIR}\,:\, \Delta \rightarrow \Delta^*$  which induces the 
involution between families of Calabi-Yau varieties:  
\[ {MIR}\, :\, {\cal F}({\Delta}) \rightarrow {\cal F}({\Delta}^*). \]
We  use the notaion  $MIR$  for the involution acting on 
families of Calabi-Yau hypersurfaces corresponding 
to reflexive polyhedra in order to stress the following important 
conjecture.
\bigskip

\begin{conj}
The combinatorial involution $MIR$ agrees with the involution 
${\cal MIR}$ of conformal field theories 
associated to  $3$-dimensional 
Calabi-Yau varieties from ${\cal F}({\Delta})$ and  ${\cal F}({\Delta}^*)$.  
\end{conj}

The next purpose  of the paper is to give  arguments 
showing that the  Calabi-Yau families ${\cal F}(\Delta)$  and 
${\cal F}(\Delta^*)$ are good candidates to be  mirror symmetric, i.e.,  
the  involution 
$MIR : {\cal F}({\Delta}) \rightarrow {\cal F}({\Delta}^*)$ 
agrees  with  properties of the mirror operation ${\cal MIR}$ in physics.

Let ${\hat Z}$ denotes  a {\em maximal projective crepant partial} 
desingularization ({\em MPCP}-desingulari\-zation) of a projective 
Calabi-Yau hypersurface ${\overline Z}$ in ${\bf P}_{\Delta}$. 
First, using results of Danilov and Khovanski\^i,  we show  that for 
$n \leq 3$ the Hodge number $h^{n-2,1}$  of {\em MPCP}-desingularizations  
of Calabi-Yau hypersurfaces  in  the family ${\cal F}(\Delta)$ equals 
the Picard number  $h^{1,1}$ of  {\em MPCP}-desingularizations   
of  Calabi-Yau hypersurfaces  in the family ${\cal F}(\Delta^*)$. 
As a corollary, we obtain the relation (\ref{main.relation}) predicted by 
physicists for mirror symmetric Calabi-Yau $3$-folds. 
Then, we prove  that 
the nonsingular part ${\hat {Z}}$ of $\overline{Z}$ consisting of the 
union of $Z_f$ with all $(n-2)$-dimensional affine strata corresponding 
$(n-1)$-dimensional edges of $\Delta$ has {\em always} the Euler 
characteristic 
zero. Finally, we prove a simple combinatorial formula for 
the Euler characteristic of Calabi-Yau $3$-folds. This formula immediately 
implies that for any pair of dual $4$-dimensional 
reflexive polyhedra $\Delta$ and $\Delta^*$, Calabi-Yau 3-folds obtained 
as {\em MPCP}-desingularizations of Calabi-Yau hypersurfaces 
in ${\cal F}(\Delta)$ and in ${\cal F}(\Delta^*)$ have opposite Euler 
characteristics

Every $4$-dimensional reflexive polyhedron $\Delta$ gives rise to a 
family ${\cal F}(\Delta)$ of Calabi-Yau $3$-folds. There exist 
many non-isomorphic $4$-dimensional reflexive polyhedra. As a result, 
we obtain many new examples of Calabi-Yau $3$-folds and their mirrors. 
We give a simplest illustration of our method in the construction 
of mirror candidats for smooth Calabi-Yau hypersurfaces of bidegree 
$(3,3)$ in ${\bf P}_2 \times {\bf P}_2$ (see \ref{two.proj}). 
This 
example shows that the  construction of smooth Calabi-Yau 3-folds 
sometimes needs not only resolutions of 3-dimensional quotient 
singularities defined by lattice triangles (see \cite{roan0}), but 
also more complicated toroidal singularities defined by 
$2$-dimensional lattice polygons.  
\medskip

Section  5 is devoted to relations between our method and other 
already known methods of explicit constructions of mirrors. 

First, we calculate the  one-parameter  mirror family for the family 
of $(n-1)$-dimensional hypersurfaces of degree $n+1$ in 
$n$-dimensional projective space and show that our result for 
quintic therefolds coincides the well-know construction of physicists. 

Next,  we investigate the category ${\cal C}_n$ 
of reflexive pairs $(\Delta, M)$, 
where $M$ is an integral lattice of rank $n$ and $\Delta$ is an 
$n$-dimensional  reflexive 
polyhedron with vertices in $M$. Morphisms in the category ${\cal C}_n$ 
give rise to relations between families of Calabi-Yau hypersurfaces in 
toric varieties.  Namely, the existence of a morphism from  
$(\Delta_1, M_1)$ to $(\Delta_2, M_2)$ implies that 
Calabi-Yau hypersurfaces in ${\cal F}(\Delta_1)$ consist of 
quotients by action of a finite abelian group of Calabi-Yau hypersurfaces 
in ${\cal F}(\Delta_2)$. 
If a reflexive polyhedron $\Delta$ is {\em selfdual} (see \ref{selfdual}), 
then  we obtain a well-known interpretation the relation between 
families ${\cal F}(\Delta)$ and 
${\cal F}(\Delta^*)$ based on quotients by actions of finite abelian groups. 

We prove that if a reflexive polyhedron 
$\Delta$ is an $n$-dimensional simplex, then $\Delta$ is selfdual. 
This leads us to some smaller invariant under the mirror operation $MIR$ 
subclasses  of all possible Calabi-Yau families ${\cal F}(\Delta)$.  
We show that Calabi-Yau families 
${\cal F}(\Delta)$ corresponding to reflexive simplices  consist 
of deformations of quotiens of Fermat-type hypersurfaces in 
weighted projective spaces. 
As a result, we obtain that our method 
for  constructions  of mirror candidates  coincides with the method 
of Greene and Plesser \cite{gree}, so that we obtain a generalization of 
the theorem of S.-S. Roan in \cite{roan1}.  
\medskip

In section  6 we discuss the problem of the classification of 
reflexive polyhedra. There exists a general finiteness theorem for reflexive 
polyhedra $\Delta$ of fixed dimension $n$. It follows from this theorem 
that   there exist only finitely many topological types of  
$(n-1)$-dimensional Calabi-Yau varieties obtained from  
families ${\cal F}(\Delta)$. 
Therefore, we come to  the natural 
problem of classification of all reflexive polyhedra of 
any fixed dimension $n$.  
This problem is equivalent to the classification of toric Fano varieties with 
only Gorenstein  singularities. There are some  
results relating to the classification of reflexive polyhedra in 
\cite{bat0,bat11,ewald,vosk.kl,watanabe}. 
We give  an explicit algorithm for the classification of reflexive simplices 
of fixed dimension $n$. 
\medskip

In section 7 we consider Calabi-Yau hypersurfaces and our candidate 
for the  mirror operation 
from the view-point of singularity theory. We give several examples 
showing analogies between the mirror symmetry and duality phenonena 
in singularity theory 
such as the strange duality of Arnold and the duality for Tsuchihashi's cusp 
singularities.   The existence of a relationship between 
the strange duality of Arnold and the mirror symmetry has been first 
discovered by M.Lynker and R.Schimmrigk (see \cite{lyn.schimm}). However, it seems that the 
Arnold's duality can only partially explain the mirror symmetry (see 
remark \ref{orthogonal}, and also \cite{roan2}). 

If we try to extrapolate our method for construction of Calabi-Yau manifolds 
of dimension $\leq 2$, then we come to 
another remarkable relation between the mirror symmetry 
for $K3$-surfaces and special solutions of the 
Picard-Fuchs equations.  It turns out that the  
one-parameter family ${\cal F}_s$ of $K3$-surfaces relating to the 
Ap\'ery's proof of the irrationality of $\zeta (3)$ and Fermi surfaces 
for  potential zero \cite{pet} has properties which  are  analogous to 
ones of the one-parameter family $\{ {\cal W}_{\psi} \}$ of Calabi-Yau 
$3$-folds. One should remark that the number $\zeta(3)$ arises in {\em all 
calculations} of special solutions of the Picard-Fuchs 
equations for one-parameter families of Calabi-Yau $3$-folds in 
\cite{cand2,font,klemm2}. 
\bigskip

It is natural to ask whether it is possible to extend also  Step 2 and  
Step 3 of the calculation of physicists for quintic mirrors  to  
the case of Calabi-Yau hypersurfaces in arbitrary $n$-dimensional toric 
Fano varieties ${\bf P}_{\Delta}$ associated with  reflexive polyhedra  
$\Delta$.  It seems that the answer should 
be affirmative at least for reflexive polyhedra $\Delta$ of  
dimension $n=3$. For instance,  the calculation of the Gauss-Manin connection 
for affine hypersurfaces in algebraic tori (which plays the main role in 
Step 2) was given  in \cite{bat.var}. 
\bigskip

{\bf Acknowledgements.} 
The ideas presented in the paper arose  in some preliminary form during my 
visiting Japan in  April-June 1991 supported by  Japan Association 
for Mathematical Sciences and T\^ohoku University.  I am grateful to  P.M.H. Wilson whose lectures  
on Calabi-Yau $3$-folds in Sendai and Tokyo introduced me to this topic.  
Geometrical and arithmetical aspects of the theory of 
higher-dimensional Calabi-Yau varieties turned out to have deep 
connections with my  research interests due to influence of my teachers  
V.A. Iskovskih and  Yu. I. Manin. 

It is a pleasure to acknoledge helpful recommendations, questions  and 
remarks  on  preliminary versions of 
the paper from mathematicians D. Dais,  B. Hunt, Y. Kawamata, 
D. Morrison, T. Oda, S.-S. Roan, D. van Straten, J. Stienstra,    
and from physicists P. Berglund, Ph. Candelas, A. Klemm, R. Schimmrigk. 
I am also very grateful to D. Morrison, B. Greene and P. Aspinwall who  
found some serious errors in my earlier formulas for the Hodge number 
$h^{1,1}$. Making corrections of these errors, I found an easy proof 
of the formula for $h^{n-2,1}$.

I would like to thank the D.F.G. for the support, and the University of 
Essen, especially  H. Esnault and E. Viehweg, for providing 
ideal conditions for my work.    
\bigskip

\newpage

\section{The geometry of toric varieties}

\hspace*{\parindent}

Almost all statements of this section one can find in \cite{dan1}, 
\cite{oda1} and \cite{reid}.
\medskip

\subsection{Toric varieties and integral polyhedra}

\hspace*{\parindent}

Let $M$ be a free abelian group of rank $n$. Put $\overline{M} = 
{\bf Z} \oplus M$.  Define the $n$-dimensional algebraic torus 
${\bf T}$  over complex number field ${\bf C}$ as follows  
\[ {\bf T} = ({\bf C}^*)^n = \{ X = 
( X_1 , \ldots,  X_n ) \in {\bf C}^n \mid 
X_1 \cdots X_n \neq 0 \}. \]  

If we choose an isomorphism $M \cong {\bf Z}^n$, we can identify 
elements $ m \in M$ with Laurent monomials 
$X^m = X_1^{m_1} \cdots X_n^{m_n}$, where $m = 
(m_1, \ldots,  m_n) \in {\bf Z}^n$. Every element $m \in M$  
defines the algebraic character $\chi_m : {\bf T} 
\rightarrow {\bf C}^*$ :
\[ \chi_m (X_1, \ldots, X_n) = X_1^{m_1} \cdots X_n^{m_n}. \]

By the same way,  we obtain  an isomorphism between  the 
 lattice $\overline{M}$ and  the group of all algebraic characters of the  
 $(n+1)$-dimensional torus 
\[ \overline{\bf T} = ({\bf C}^*)^{n+1} = \{ \overline{X} = 
( X_0, X_1 , \ldots,  X_n ) \in {\bf C}^{n+1} \mid 
X_0 X_1 \cdots X_n \neq 0 \}. \]

Let $M_{\bf Q} = M \otimes {\bf Q}$ be the ${\bf Q}$-scalar extension 
of $M$. Take  a convex $n$-dimensional polyhedron $\Delta$ in 
$M_{\bf Q}$ with integral vertices (i.e., all vertices of $\Delta$ are 
elements of $M$). Consider the  $(n+1)$-dimensional ${\bf Q}$-linear space 
$\overline{M}_{\bf Q} = {\bf Q} \oplus M_{\bf Q}$ containing 
the integral lattice $\overline{M} = {\bf Z} \oplus M$. 
Define the convex 
$(n+1)$-dimensional cone $C_{\Delta} \subset \overline{M}_{\bf Q}$  
supporting $\Delta$ as 
\[ C_{\Delta} = 0 \cup \{ (x_0, x_1, \ldots, x_n ) \in {\bf Q} \oplus {M}_{\bf Q} 
\mid  ({x_1}/{x_0}, \ldots, {x_n}/{x_0})  \in \Delta, \; x_0 >0 \} . \]

Denote by $S_{\Delta} $ the  subring 
of ${\bf C} \lbrack X_0, X_1^{\pm 1}, \ldots , X_n^{\pm 1} \rbrack$ with 
the ${\bf C}$-basis consisting of monomials $X_0^{m_0}X_1^{m_1} \cdots X_n^{m_n}$ 
such that $(m_0, m_1, \ldots, m_n) \in C_{\Delta}$. Define the grading of 
monomials  by setting ${\rm deg}\, 
(X_0^{m_0}X_1^{m_1} \cdots X_n^{m_n}) = m_0$. Then $S_{\Delta}$ is a 
$(n+1)$-dimensional graded subring of the graded ring 
of ${\bf C} \lbrack X_0, X_1^{\pm 1}, \ldots , X_n^{\pm 1} \rbrack$.

\begin{opr}
{\rm Let 
${\bf P}_{\Delta,M} = {\rm Proj}  \, S_{\Delta}$ 
be an  $n$-dimensional 
projective algebraic variety corresponding to the graded 
ring $S_{\Delta}$ (see \cite{hartsh}, 
Ch.I, \S 2). We call ${\bf P}_{\Delta,M}$ 
{\em the  projecive  toric variety associated with 
the integral polyhedron} $\Delta$ and denote by 
${\cal O}_{\Delta}(1)$ the ample invertible sheaf on 
${\bf P}_{\Delta, M}$ corresponding to the graded $S_{\Delta}$-module 
$S_{\Delta}(-1)$. If there is no confusion, 
we put  ${\bf P}_{\Delta,M} = {\bf P}_{\Delta}$. }
\label{def.toric}
\end{opr}

\begin{opr}
{\rm Let ${\rm vol}_M\, \Delta$ be the volume of the polyhedron 
$\Delta$ relative to the integral lattice $M$. The number 
\[ d_M(\Delta) = n!{\rm vol}_M (\Delta) \]
we call the {\em degree} of the polyhedron $\Delta$ relative to $M$. 
If there is no confusion, we put simply $d_M(\Delta) = 
d(\Delta)$.}
\end{opr}

One  can easily see that $d_M(\Delta)$ coincides with the projective 
degree of ${\bf P}_{\Delta, M}$ (\cite{hartsh}, Ch.I, \S 7). 
\bigskip

Let $\Theta$ be an arbitrary $l$-dimensional polyhedral face of $\Delta$ 
$(l \leq  n)$, 
$M(\Theta)_{\bf Q}$ the minimal $l$-dimensional 
${\bf Q}$-subspace of $M_{\bf Q}$ containing all vectors $x - x'$, 
where $x, x'   \in \Theta$. One gets the   integral sublattice 
$M(\Theta) = M(\Theta)_{\bf Q} \cap M$. 
Denote by ${\bf T}_{\Theta}$ 
the $l$-dimensional algebraic torus with the group of 
characters  $M(\Theta)$ (i.e., in particular,  
${\bf T}_{\Theta} = {\bf T}$ if $\Theta = \Delta$). The embedding  $M(\Theta)  
\hookrightarrow M$ induces a surjective homomorphism 
$ {\bf T} \rightarrow {\bf T}_{\Theta}.$ 
This homomorphism and the multiplicative low action of  
${\bf T}$ on itself give rise 
to a  natural action of ${\bf T}$ on ${\bf T}_{\Theta}$. 
\medskip

\begin{opr}
{\rm The $l$-dimensional algebraic torus ${\bf T}_{\Theta}$ with the action 
of ${\bf T}$ on it is called a {\em {\bf T}-orbit associated with  
the $l$-dimensional face  $\Theta \subset \Delta$. }}
\label{orbit}
\end{opr}
\medskip

Define the  $(l+1)$-dimensional face  $C_{\Theta} \subset C_{\Delta}$   
supporting $\Theta$ as 
\[ C_{\Theta} = 0 \cup \{ (x_0, x_1, \ldots, x_n ) \in {\bf Q} \oplus {M}_{\bf Q} 
\mid  ({x_1}/{x_0}, \ldots, {x_n}/{x_0})  \in \Theta, \; x_0 >0 \} . \]
Denote by $S_{\Theta} $ the  subring of 
${\bf C} \lbrack X_0, X_1^{\pm 1}, \ldots , X_n^{\pm 1} \rbrack$ with 
the ${\bf C}$-basis consisting of monomials $X_0^{m_0}X_1^{m_1} \cdots X_n^{m_n}$ 
such that $(m_0, m_1, \ldots, m_n) \in C_{\Theta}$. 
Consider the homogeneous ideal $I_{\Theta} \subset S_{\Delta}$ whose 
${\bf C}$-basis consists  of 
monomials $X_0^{m_0}X_1^{m_1} \cdots X_n^{m_n}$ 
such that $(m_0, m_1, \ldots, m_n) \in (C_{\Delta} \setminus C_{\Theta})$.
Since $S_{\Theta} \cong S_{\Delta}/ I_{\Theta}$, the ideal $I_{\Theta}$ 
defines  an $l$-dimensional projective 
toric subvariety ${\bf P}_{\Theta} = {\rm Proj}\, S_{\Theta} \subset 
{\bf P}_{\Delta}$.

\begin{opr}
{\rm The $l$-dimensional projective subvariety ${\bf P}_{\Theta} \subset 
{\bf P}_{\Delta}$  is called the {\em toric subvariety associated with  
the $l$-dimensional face  $\Theta \subset \Delta$.} 

Put 
\[ {\bf P}^{(i)}_{\Delta} = \bigcup_{{\rm codim}\, \Theta = i} 
{\bf P}_{\Theta}.  \] }
\end{opr}
\medskip

The following theorem describes the geometry of toric varieties  
${\bf P}_{\Delta,M}$. 

\begin{theo}
A toric variety ${\bf P}_{\Delta,M}$ is a projective compactification of the 
torus ${\bf T} = {\bf T}_{\Delta}$ such that the group law action of ${\bf T}$ 
on itself extends to an algebraic  action on ${\bf P}_{\Delta,M}$ 
having the following properties 

{\rm (i)} ${\bf P}_{\Delta}$ and ${\bf P}_{\Delta}^{(i)}$  have  
decompositions into a 
disjoint unions of ${\bf T}$-orbits ${\bf T}_{\Theta}$ 
associated with faces of $\Delta$
\[ {\bf P}_{\Delta} = \bigcup_{\Theta  \subset  \Delta} {\bf T}_{\Theta
 }, \]
\[ {\bf P}_{\Delta}^{(i)} = \bigcup_{{\rm codim}\, \Theta \leq i} 
{\bf T}_{\Theta }; \]

{\rm (ii)}  for every $l$-dimensional  face  $\Theta \subset \Delta$, 
the closure of the  ${\bf T}$-orbit ${\bf T}_{\Theta}$ in 
${\bf P}_{\Delta}$ is the projective $l$-dimensional toric subvariety 
${\bf P}_{\Theta}$  associated with 
the  $l$-dimensional integral polyhedron $\Theta$$;$ 

{\rm (iii)} the action of the torus ${\bf T}$ induces 
the  one-to-one correspondence between faces $\Theta 
\subset \Delta$ and toric subvarieties ${\bf P}_{\Theta}$. Moreover,  
$ {\bf P}_{\Theta} \cap {\bf P}_{\Theta '} = {\bf P}_{\Theta \cap \Theta '}$
for any two faces $\Theta, \Theta' \subset \Delta$.  
\label{toric.poly}
\end{theo}
\bigskip

Althought, the definition of toric varieties ${\bf P}_{\Delta}$ associated 
with integral polyhedra $\Delta$ is very simple, it is not always convenient.  
Note that a polyhedron $\Delta$  defines not only the corresponding toric  
variety ${\bf P}_{\Delta}$, but also the  choice of an  
ample invertible sheaf ${\cal O}_{\Delta}(1)$ 
on ${\bf P}_{\Delta}$ and an  embedding of ${\bf P}_{\Delta}$ 
into a  projective space. In general, 
there exist infinitely many different ample sheaves on ${\bf P}_{\Delta}$. 
As a result, there are infinitely many different integral  polyhedra defining 
isomorphic toric varieties (one can take, for example,  multiples of 
$\Delta$ :  $k \Delta$, $k = 1,2, \ldots )$. If we want to get a 
one-to-one 
correspondence  between toric varieties and some combinatorial data,   
we need another definition of toric varieties in terms of fans of rational 
polyhedral cones. 
This approach to  toric varieties  gives more possibilities, 
it allows, for example, to construct affine 
and quasi-projective toric varieties as well as  complete  
toric varieties which are not quasi-projective.  
\bigskip

\subsection{Toric varieties and rational polyhedral fans} 

\hspace*{\parindent}

Let $N = {\rm Hom}\, (M, {\bf Z})$ be the dual to $M$ lattice, $N_{\bf Q}$ 
the ${\bf Q}$-scalar extension of $N$. Denote by  
$ \langle *, * \rangle \; : \; M_{\bf Q} \times N_{\bf Q} \rightarrow 
{\bf Q} $ the nondegenerate pairing between the $n$-dimensional 
${\bf Q}$-spaces $M_{\bf Q}$ and $N_{\bf Q}$.

\begin{opr}
{\rm Let $\sigma \subset  N_{\bf Q}$ be a {\em convex rational 
polyhedral cone}, i.e., an 
intersection of finitely  many half-spaces defined by equations 
\[  \sigma = \{ y \in N_{\bf Q} \mid \langle m_1, y \rangle \geq 0, 
\ldots,  \langle m_r, y \rangle \geq 0 \} , \]
where $m_1, \ldots , m_r$ are elements of $M$. 
Define {\em the dual cone} as 
\[ \check {\sigma} = \{ {\bf Q}_{\geq 0} m_{1} + \cdots + {\bf Q}_{\geq 0} 
m_{r} \} \subset  M_{\bf Q}. \]
Denote by $M_{{\sigma}}$ the finitely generated semigroup 
$\check {\sigma} \cap M$. Let $S_{\sigma} = {\bf C} \lbrack X^m \rbrack\; (
m \in  M_{\sigma})$ the corresponding semigroup algebra over ${\bf C}$. 
The affine algebraic variety 
\[ {\bf A}_{\sigma, N} = {\rm Spec}\, S_{\sigma} \]
is called {\em the affine toric variety associated with the cone $\sigma$}.

If ${\rm dim}\, \sigma = {\rm dim}\, N_{\bf Q}$, then there exists the 
unique ${\bf T}$-invariant maximal ideal in $S_{\sigma}$ generated by 
non-constant monomials $X^m$, $m \in M_{\sigma}$.  We denote by $p_{\sigma}$ 
the corresponding closed point in ${\bf A}_{\sigma, N}$. }
\label{affine}
\end{opr}
\medskip

\begin{rem}
{\rm Let  $\sigma \subset N_{\bf Q}$ be  a  convex rational 
polyhedral cone of dimension $s$ $( s \leq n)$. Assume that 
$\sigma$ does not 
contain any ${\bf Q}$-subspace of positive dimension. 
Let $N(\sigma)_{\bf Q}$ be the minimal $s$-dimensional ${\bf Q}$-subspace 
in $N_{\bf Q}$ containing $\sigma$, $N(\sigma) = N(\sigma)_{\bf Q} \cap N$   
free abelian group of rank $s$. 
Then the dual cone 
$\check {\sigma} \subset M_{\bf Q}$ has always dimension $n$ and the 
corresponding $n$-dimensional affine toric variety ${\bf A}_{\sigma, N}$ 
is isomorphic to the product 
\[ ({\bf C}^*)^{n-s} \times {\bf A}_{\sigma, N(\sigma)}. \]
In particular, one gets 
\[ {\bf A}_{\sigma, N} \cong ({\bf C}^*)^n, \;\; {\rm  if}  
\;\; {\rm dim}\, \sigma = 0, \]
\[ {\bf A}_{\sigma, N} \cong ({\bf C}^*)^{n-1} \times {\bf C},  \;\;
{\rm if} \;\;  {\rm dim}\, \sigma = 1.\] }
\label{codim1,2}
\end{rem}

\begin{opr}
{\rm {\em A rational polyhedral fan $\Sigma$}  of convex  cones 
 is a finite collection $\{ \sigma_i \}$ of convex rational polyhedral 
 cones $\sigma_i \in N_{\bf Q}$ satisfying the conditions:  

(i) every cone $\sigma_i \in \Sigma$ does not contain  ${\bf Q}$-linear 
subspaces of positive  dimension; 

(ii) if $\tau$ is a face of a cone $\sigma \in \Sigma$, then $\tau \in 
\Sigma$; 

(iii) if $\sigma_i$,  $\sigma_j \in \Sigma$, then $\sigma_i \cap \sigma_j$ is 
a common face both of $\sigma_i$ and of $\sigma_j$.  \\
The set of all $i$-dimensional cones in $\Sigma$ will be denoted by 
$\Sigma^{(i)}$. }
\end{opr}
\bigskip

By gluing together $n$-dimensional affine toric varieties 
${\bf A}_{\sigma, N}$ arising from cones $\sigma$ in a rational polyhedral 
fan $\Sigma$, we get an {\em abstract  toric variety ${\bf P}_{\Sigma}$ 
associated 
with  a rational polyhedral fan $\Sigma$}.  
The following theorem describes the geometry of ${\bf P}_{\Sigma}$ in terms 
of properties of $\Sigma$. 
\medskip

\begin{theo}
The toric variety ${\bf P}_{\Sigma,N}$ is a partial  compactification of the 
torus ${\bf T} = {\bf A}_{0,N}$ such that the group law action of ${\bf T}$ 
on itself extends to an algebraic  action on ${\bf P}_{\Sigma,N}$ 
having the following properties 

{\rm (i)} ${\bf P}_{\Sigma}$ has an open  covering by affine toric varieties 
${\bf A}_{\sigma, N}$ 
\[ {\bf P}_{\Sigma} = \bigcup_{\sigma  \subset  \Sigma} {\bf A}_{\sigma , N} 
; \]

{\rm (ii)} the action of the torus ${\bf T}$ induces 
the  one-to-one correspondence between cones  $\sigma \in \Sigma$ 
and ${\bf T}$-invariant affine open subsets ${\bf A}_{\sigma,N}$ 
in  ${\bf P}_{\Sigma}$. Moreover, ${\bf A}_{\sigma,N}  \cap 
{\bf A}_{\sigma',N} = {\bf A}_{\sigma \cap \sigma',N}$ 
for any two cones  $\sigma, \sigma' \in  \Sigma$.

{\rm (iii)} ${\bf P}_{\Sigma}$ is smooth if and only if every cone $\sigma_i 
\in \Sigma$ can be represented as  
\[ \sigma_i = \{ {\bf Q}_{\geq 0} n_{i_1} + \cdots + {\bf Q}_{\geq 0} 
n_{i_s}   \in N_{\bf Q} \}, \]
where  $\{n_{i_1}, \ldots,  n_{i_s}\} $ is  a part of a 
${\bf Z}$-basis  of $N$;

{\rm (iv)} ${\bf P}_{\Sigma}$ is compact if and only if  
the union $\bigcup_{\sigma \in \Sigma}$ 
coincides  with the whole space $N_{\bf Q}$. In this case  $\Sigma$ 
is said to be a complete fan. 
\label{def.tor2}
\end{theo}
\bigskip

\begin{opr}
{\rm Denote by  $\Sigma^{[i]}$ the subfan of $\Sigma$ 
consisting of all cones $\sigma \in \Sigma$ such that ${\rm dim}\, \sigma 
\leq i$, i.e., 
\[ \Sigma^{[i]} = \bigcup_{{\rm dim}\, \sigma \leq i} \Sigma^{(i)}. \] 
The open toric subvariety in ${\bf P}_{\Sigma}$ 
corresponding to $\Sigma^{[i]}$ we denote by 
${\bf P}^{[i]}_{\Sigma}$.}
\end{opr}

\subsection{Relations between two definitions of toric varieties}

\hspace*{\parindent}

Inspite of the fact that the definition of toric varieties via rational 
polyhedral fans is more general than via integral polyhedra, we will use both  
definitions. The choice of a definition in the 
sequel will depend on  questions  we are interested in. 
In one situation, it will be more convenient to describe properties 
of toric varieties and their subvarieties in terms of integral polyhedra. In  
another situation, it will be more convenient to use the language of 
rational polyhedral fans.  
So it is important to know how one can construct a fan $\Sigma(\Delta)$ from 
an integral polyhedron $\Delta$, and how one can construct an integral 
polyhedron $\Delta(\Sigma)$ from a rational polyhedral fan $\Sigma$. 

The first way $\Delta \Rightarrow \Sigma(\Delta)$ is rather simple:

\begin{prop}
For every  $l$-dimensional face  $\Theta \subset \Delta$, define the convex 
$n$-dimensional cone ${\check {\sigma}}(\Theta) \subset M_{\bf Q}$ 
 consisting of all vectors 
$\lambda (p - p')$,  where $\lambda \in{\bf Q}_{\geq 0}$,  
$p \in \Delta$, $p' \in \Theta$. Let $\sigma( \Theta) \subset N_{\bf Q}$ 
be the $(n-l)$-dimensional dual cone 
relative to  ${\check {\sigma}}(\Theta) \subset M_{\bf Q}$. 

 The set $\Sigma(\Delta)$ of all  cones $\sigma(\Theta)$, where 
$\Theta$ runs over all faces of $\Delta$, determines the complete 
rational polyhedral fan defining the toric variety ${\bf P}_{\Delta}$.
\label{def.fan}
\end{prop}
\bigskip

We can now describe relations between two theorems   \ref{toric.poly} and 
\ref{def.tor2}. For example, the decomposition  of 
${\bf P}_{\Delta}$ into a disjoint union of 
${\bf T}$-orbits ${\bf T}_{\Theta}$ can be reformulated via cones 
$\sigma(\Theta)$ 
in the fan $\Sigma(\Delta) = \Sigma $ as follows.

\begin{prop}
Let $\Delta$ be an $n$-dimensional $M$-integral polyhedron in $M_{\bf Q}$, 
$\Sigma = \Sigma(\Delta)$  the corresponding  
complete rational  polyhedral fan 
in $N_{\bf Q}$.   Then

{\rm (i)} For any face $\Theta \subset \Delta$,  
the affine toric variety ${\bf A}_{{\sigma}(\Theta),N}$ is the 
minimal ${\bf T}$-invariant 
affine open subset in ${\bf P}_{\Sigma}$ containing  ${\bf T}$-orbit 
${\bf T}_{\Theta}$ {\rm ($see$  \ref{orbit})}. 

{\rm (ii)} Put ${\bf T}_{\sigma}:= {\bf T}_{\sigma(\Theta)}$. 
There exists a one-to-one correspondence between 
$s$-dimensional cones $\sigma \in \Sigma$ and $(n-s)$-dimensional 
${\bf T}$-orbits ${\bf T}_{\sigma}$ such that ${\bf T}_{\sigma'}$ is 
contained in the closure of ${\bf T}_{\sigma}$ if and only if 
$\sigma$ is a face of $\sigma'$. 

{\rm (iii)} 
\[ {\bf P}_{\Sigma}^{[i]} = \bigcup_{{\rm dim}\,  \sigma \leq i} 
{\bf T}_{\sigma}  \] 
is an open ${\bf T}$-invariant subvariety 
in ${\bf P}_{\Sigma} = {\bf P}_{\Delta}$, and 
\[ {\bf P}_{\Sigma} \setminus {\bf P}_{\Sigma}^{[i]} = 
{\bf P}^{(i)}_{\Delta}. \]
\label{orbits.fan}
\end{prop}
\bigskip

Let us now consider another direction $\Sigma \Rightarrow \Delta(\Sigma)$. 
In this case, in order to construct $\Delta(\Sigma)$ 
it is not sufficient to know only a complete fan $\Sigma$. 
We need a strictly upper convex support function $h\; :\;  N_{\bf Q} 
\rightarrow {\bf Q}$. 

\begin{opr}
{\rm Let $\Sigma$ be a rational polyhedral fan and let  
$h\; :\;  N_{\bf Q} \rightarrow {\bf Q}$ be a function such that 
$h$ is linear on any cone $\sigma \subset \Sigma$. In this situation 
$h$ is called a {\em support function} for the fan $\Sigma$.  
We call $h$ {\em integral} if $h(N) \subset {\bf Z}$.  We call $h$ 
{\em upper convex} if  $h(x + x') \leq h(x) + h(x')$ for  all 
$x, x' \in N_{\bf Q}$. Finally,  $h$ is called {\em strictly upper convex} if 
$h$ is upper convex and 
for any two distinct $n$-dimensional cones $\sigma$ and $\sigma'$ in 
$\Sigma$ the restrictions $h_{\sigma}$ and $h_{\sigma'}$ of $h$ on $\sigma$ 
and $\sigma'$ are different  linear functions. }
\end{opr}

\begin{rem}
{\rm By general theory of toric varieties \cite{dan1,oda1}, support functions 
one-to-one correspond to  ${\bf T}$-invariant 
${\bf Q}$-Cartier divisors $D_h$ 
on  ${\bf P}_{\Sigma}$. A divisor $D_h$ is a ${\bf T}$-invariant 
Cartier divisor (or ${\bf T}$-linearized invertible sheaf) if and only if 
$h$ is integral. $D_h$ is numerically effective if and only if $h$ is upper 
convex.  Strictly upper convex support functions $h$ on a fan $\Sigma$ 
one-to-one correspond to ${\bf T}$-linearized {\em ample} invertible 
sheaves ${\cal O}(D_h)$  over ${\bf P}_{\Sigma}$.}
\end{rem} 

The next proposition describes the construction of a polyhedron $\Delta$  from 
a fan $\Sigma$ supplied with a strictly convex support function $h$. 

\begin{prop}
Let the convex polyhedron $\Delta = \Delta(\Sigma, h)$ to be defined as follows
\[ \Delta(\Sigma,h) = \bigcap_{\sigma \in \Sigma^{(n)}} ( -h_{\sigma} + 
{\check \sigma} ), \]
where integral linear functions $h_{\sigma}: N \rightarrow {\bf Z}$ 
 are considered as elements of 
the lattice $M$. Then one has ${\bf P}_{\Delta} \cong {\bf P}_{\Sigma}$ and  
${\cal O}_{\Delta}(1) \cong {\cal O}(D_h)$. 
\end{prop}

\subsection{Singularities of toric varieties}

\hspace*{\parindent}

From \ref{codim1,2} and \ref{def.tor2}(iii), one immediately 
obtains the following. 

\begin{prop}
The open ${\bf T}$-invariant subvariety ${\bf P}_{\Sigma}^{[1]}$ 
in ${\bf P}_{\Sigma}$ is smooth. 
\label{codim1}
\end{prop}

There is a combinatorial description of singularities of 
toric varieties ${\bf P}_{\Sigma} = {\bf P}_{\Sigma, N}$. 
\bigskip

Assume that that the dimension of  $\sigma \in \Sigma$ is $s$. 
For every two points $p, p'$ in the $(n-s)$-dimensional ${\bf T}$-orbit 
${\bf T}_{\sigma}$, the action of ${\bf T}$ induces an  isomorphism    
between local neighbourhoods of $p$ and $p'$. By \ref{orbits.fan}(i), 
${\bf T}_{\sigma}$ is contained in  the 
${\bf T}$-invariant affine open subset  ${\bf A}_{\sigma,N}$. 
By \ref{codim1,2}, ${\bf A}_{\sigma,N}$ splits into the product 
\[ ({\bf C}^*)^{n-s} \times {\bf A}_{\sigma, N(\sigma)}. \]
Therefore, we come to the following. 

\begin{prop}
Small analytical neighbourhoods in ${\bf P}_{\Sigma}$ of any point $p \in 
{\bf T}_{\sigma}$ are isomorphic to  products  of open $(n-s)$-dimensional 
ball  and a small analytical neighbourhood of the closed point $p_{\sigma} 
\in  {\bf A}_{\sigma, N(\sigma)}$ {\rm ($cf$. \ref{affine})}.  
\end{prop}
 
So the description of singularities of ${\bf P}_{\Sigma}$ reduces 
to the descriptions of  singularities at the closed points $p_{\sigma}$ 
of  affine toric varieties ${\bf A}_{\sigma, N(\sigma)}$, where 
$\sigma$ runs over all cones of $\Sigma$.

\begin{opr}
{\rm A normal variety $W$ is said to have only {\em Gorenstein} (resp.
 {\em ${\bf Q}$-Gorenstein}) singularities if the canonical divisor $K_W$ is a 
Cartier (resp.  ${\bf Q}$-Cartier) divisor. }
\end{opr}

The following two propositions are contained in the paper of M. Reid 
\cite{reid}. 

\begin{prop}
Let $n_1, \ldots, n_r \in N$ $(r \geq s)$ be primitive $N$-integral 
generators of all $1$-dimensional faces of an $s$-dimensional 
cone $\sigma$. 

{\rm (i)} the point $p_{\sigma} \in {\bf A}_{\sigma, N(\sigma)}$ is 
${\bf Q}$-factorial 
$($or quasi-smooth$)$ if ond only if the cone $\sigma$ is simplicial, 
i.e., $r = s$; 

{\rm (ii)} the point $p_{\sigma} \in {\bf A}_{\sigma, N(\sigma)}$ is 
${\bf Q}$-Gorenstein 
if ond only if the elements $n_1, \ldots, n_r$ are contained in 
an affine hyperplane 
\[ H_{\sigma} \; : \; \{ y \in N_{\bf Q} \mid \langle  
k_{\sigma}, y \rangle = 1 \}, \]
for some $k_{\sigma} \in M_{\bf Q}$ $($note that when  
${\rm dim}\, \sigma = {\rm dim}\, N$, the element 
$k_{\sigma}$ is unique if it exists$)$. Moreover, 
${\bf A}_{\sigma, N(\sigma)}$ is Gorenstein if and only if 
$k_{\sigma} \in M$. 
\label{sing1}
\end{prop}

\begin{rem}
{\rm If ${\bf P}_{\Sigma}$ is a ${\bf Q}$-Gorenstein toric variety, the 
elements $k_{\sigma} \in M_{\bf Q}$ ($\sigma \in \Sigma^{(n)}$) define 
together  
the support function $h_K$ on $N_{\bf Q}$ such that the restriction of 
$h_K$ on $\sigma \in \Sigma^{(n)}$ is $k_{\sigma}$. The support function 
$h_K$ corresponds to the anticanonical divisor on ${\bf P}_{\Sigma}$.}
\label{anti.des}
\end{rem}  

\begin{prop}
Assume that ${\bf A}_{\sigma, N(\sigma)}$ is ${\bf Q}$-Gorenstein {\rm 
($see$  \ref{sing1}(ii))}, then 

{\rm (i)} ${\bf A}_{\sigma, N(\sigma)}$ has at the point $p_{\sigma}$ 
at most { terminal} singularity  
if and only if 
\[ N \cap \sigma \cap \{y \in N_{\bf Q} \mid \langle  k_{\sigma}, 
y \rangle \leq  1 \} 
= \{ 0, n_1, \ldots, n_r \} ; \]

{\rm (ii)} ${\bf A}_{\sigma, N(\sigma)}$ has at the point $p_{\sigma}$ at most 
{ canonical } singularity  
if and only if 
\[ N \cap \sigma \cap \{y \in N_{\bf Q} 
\mid \langle  k_{\sigma} , y \rangle <  1 \} 
= \{ 0 \}. \]  
\label{sing2}
\end{prop}

Using \ref{sing1} and \ref{sing2}, we obtain: 

\begin{coro}
Any Gorenstein toric singularity is canonical. 
\label{gor.can}
\end{coro}

\begin{opr} 
{\rm Let $S$ be a $k$-dimensional simplex in ${\bf Q}^n$ $(k \leq n)$ with 
vertices in ${\bf Z}^n$, 
$A(S)$ the minimal $k$-dimensional affine ${\bf Q}$-subspace containing 
$S$. Denote by ${\bf Z}(S)$ the $k$-dimensional lattice 
$A(S) \cap {\bf Z}^n$. We call $P$ {\em elementary} if 
$S \cap {\bf Z}(S)$ contains only vertices of $S$. We call $S$ 
{\em regular}  if the degree of $S$ relative to ${\bf Z}(S)$ is $1$.}
\end{opr}

It is clear that every regular simplex is elementary. The converse is not true 
in general. However, there exists the following easy lemma. 

\begin{lem}
Every elementary simplex of dimension $\leq 2$ is 
regular. 
\label{el.reg}
\end{lem}

By \ref{def.tor2}(iii), \ref{sing1}, \ref{sing2}, we obtain:

\begin{prop}
Let  ${\bf P}_{\Sigma}$  be a toric variety with only ${\bf Q}$-Gorenstein 
singularities, i.e., for any cone $\sigma \in \Sigma$ let the corresponding 
element $k_{\sigma} \in M_{\bf Q}$ be well-defined {\rm (\ref{sing1}(ii))}. 
Then 

{\rm (i)}  ${\bf P}_{\Sigma}$ has only ${\bf Q}$-factorial terminal 
singularities if and only if for every cone $\sigma \in \Sigma$ the 
polyhedron 
\[ P_{\sigma} = \sigma \cap \{ y \in N_{\bf Q} \mid \langle  
k_{\sigma}, y \rangle \leq  1 \} \]
is an elementary simplex. 

{\rm (ii)} ${\bf P}_{\Sigma}$ is smooth if and only if for every cone 
$\sigma \in \Sigma$ the 
polyhedron 
\[ P_{\sigma} = \sigma \cap \{ y \in N_{\bf Q} \mid \langle  
k_{\sigma}, y \rangle \leq  1 \} \]
is a regular simplex. 
\label{simp}
\end{prop}

\begin{theo}
Let ${\bf P}_{\Sigma}$ be a toric variety. 

{\rm (i)} If ${\bf P}_{\Sigma}$ has only 
terminal singularities, then the open toric subvariety  
${\bf P}_{\Sigma}^{[2]}$ is smooth. 

{\rm (ii)} If ${\bf P}_{\Sigma}$ has only 
Gorenstein ${\bf Q}$-factorial 
terminal singularities, then the open toric subvariety  
${\bf P}_{\Sigma}^{[3]}$ is smooth. 
\label{codim2,3}
\end{theo}

The statement the theorem \ref{codim2,3} follows from the criterion 
\ref{def.tor2}(iii) and lemma \ref{el.reg}
\bigskip

\subsection{Proper and finite morphisms}
  
 \hspace*{\parindent} 
  
Recall a combinatorial definition of toric  morphisms  
between toric varieties. 
\bigskip

Let $\phi \; : \; N' \rightarrow N$ be a morphism of lattices, $\Sigma$ a 
fan in $N_{\bf Q}$,  $\Sigma'$ a fan in $N'_{\bf Q}$. Suppose that for each 
$\sigma' \in \Sigma'$ we can find a $\sigma \in \Sigma$ such that 
$\varphi(\sigma') \subset \sigma$. In this situation there arises  a morphism 
of toric varieties
\[ \tilde {\phi} \; : \; {\bf P}_{\Sigma',N'} \rightarrow 
{\bf P}_{\Sigma, N}. \]

\begin{exam} {Proper birational morphisms.} {\rm Assume that 
$\varphi$ is an isomorphism of lattices and $\phi(\Sigma')$ 
is a subdivision of $\Sigma$, i.e., every cone $\sigma \in \Sigma$ is a 
union of cones of $\varphi(\Sigma')$. In this case  $\tilde {\phi}$ is 
a proper birational  morphism. Such a morphism we will use  for 
constructions of desingularizations of toric singularities. }  
\label{biration}
\end{exam}
\medskip

\begin{opr}
{\rm Let $\varphi:  W' \rightarrow W$ be a proper birational 
morphism of normal ${\bf Q}$-Gorenstein algebraic varieties. The morphism 
$\varphi$ is called {\em crepant} if $\varphi^* K_{W} = K_{W'}$ 
($K_W$ and $K_{W'}$ are canonical divisors  on $W$ and $W'$ respectively). }
\end{opr}
\medskip

Using the description in \ref{anti.des} 
of support functions corresponding to anticanonical 
divisors on ${\bf P}_{\Sigma}$ and ${\bf P}_{\Sigma'}$, we obtain.  

\begin{prop}
A proper birational  morphism of ${\bf Q}$-Gorenstein toric 
varieties 
\[ \tilde {\phi} \; : \; {\bf P}_{\Sigma'} \rightarrow 
{\bf P}_{\Sigma} \]
is crepant if and only if for every cone $\sigma \in \Sigma^{(n)}$ all 
$1$-dimensional cones $\sigma'\in  \Sigma'$ which are  
contained in $\sigma$ are 
generated by primitive integral elements from $N \cap H_{\sigma}$ 
{\rm ($see$ \ref{sing1}(ii))}. 
\label{crit.crep}
\end{prop}

\begin{opr}
{\rm Let $\varphi\; : \; W' \rightarrow W$ be a projective birational 
morphism of normal ${\bf Q}$-Gorenstein algebraic varieties. The morphism 
$\varphi$ is called a {\em maximal projective crepant partial 
desingularization} ({\em MPCP-desingularization}) of $W$ 
if $\varphi$ is crepant 
and $W'$ has only ${\bf Q}$-factorial terminal singularities.}
\end{opr}
\medskip

Our next purpose is to define some  combinatorial notions which we 
use to construct $MPCP$-desingularizations of 
Gorenstein toric varieties (see \ref{crep.fano}). 
\medskip

\begin{opr}
{\rm Let $A$ be a finite subset in $\Delta \cap {\bf Z}^n$. We call $A$ 
{\em admissible} if it  contains 
all vertices of the integral polyhedron $\Delta$. }
\end{opr}

\begin{opr}
{\rm Let $A$ be an admissible subset in $\Delta \cap {\bf Z}^n$. 
By an  {\em A-triangulation} of $\Delta \subset 
{\bf Q}^n$ we mean a finite collection ${\cal T} = \{ \theta \}$ of simplices 
with vertices in $A$ having the following properties:

(i) if $\theta'$ is a face of a simplex $\theta \in {\cal T}$, then 
$\theta' \in {\cal T}$; 

(ii) the vertices of each simplex $\theta \in {\cal T}$ lie in $\Delta \cap 
{\bf Z}^n$; 

(iii) the intersection of any two simplices $\theta_1, \theta_2 \in 
{\cal T}$ either is empty or is a common face of both; 

(iv) $\Delta = \bigcup_{\theta \in {\cal T}} \theta$;

(v) every element of  $A$ is a vertex of some simplex 
$\theta \in {\cal T}$. }
\label{def.triang}
\end{opr}

\begin{opr}
{\rm An $A$-triangulation ${\cal T}$ of an integral convex polyhedron 
$\Delta \subset 
{\bf Q}^n$ is called {\em maximal} if $A = \Delta \cap {\bf Z}^n$.   }
\end{opr}

\begin{rem}
{\rm Note that ${\cal T}$ is a maximal triangulation of $\Delta$ 
if and only if every simplex $\theta \in {\cal T}$ is elementary. Therefore, 
if ${\cal T}$ is a maximal triangulation of $\Delta$, then 
for any face $\Theta \subset \Delta$ the induced triangulation of 
$\Theta$ is also maximal.}
\label{rem.max}
\end{rem}
\medskip

Assume that  $A$ is admissible. 
 Denote by ${\bf Q}^A$  the 
${\bf Q}$-space  of functions from $A$ to ${\bf Q}$. Let ${\cal T}$ be a 
triangulation of $\Delta$. Every element $\alpha \in {\bf Q}^A$ can be 
uniquely extended to a piecewise linear function 
$\alpha({\cal T})$ on $\Delta$ 
such that the restriction of $\alpha({\cal T})$ on every simplex $\theta \in 
{\cal T}$ is an affine linear function. 

\begin{opr}
{\rm Denote by $C({\cal T})$ the convex cone in ${\bf Q}^A$ consisting of 
elements $\alpha$ such that $\alpha({\cal T})$ is an upper 
convex function. We say that a 
triangulation ${\cal T}$ of $\Delta$ is {\em projective} if the cone 
$C({\cal T})$ has a nonempty  interior. In other words, ${\cal T}$ is 
projective if and only if there exists a strictly upper convex function 
$\alpha({\cal T})$. }
\label{opr.t.proj}
\end{opr}

\begin{prop} 
{\rm \cite{gelf}} 
Let $\Delta$ be an integral polyhedron. Take an admissible subset $A$ in  
$\Delta \cap {\bf Z}^n$. Then $\Delta$ admits at least one projective 
$A$-triangulation, in particular, there exists at least one maximal 
projective triangulation of $\Delta$. 
\label{triang.ex}
\end{prop}

\begin{rem}
{\rm In the paper  Gelfand, Kapranov, and Zelevinsky \cite{gelf}, it has 
been used  the notion of a {\em regular triangulation of a polyhedron} 
which is called in this paper  
a {\em projective triangulation}. The reason for such a change of 
the terminology we  will see in \ref{crep.fano}.}    
\end{rem}

\begin{exam} {\rm 
{\em Finite morphisms.} Assume that $\phi$  is injective, 
$\phi(N')$ is a sublattice of finite index 
in $N$, and $\phi(\Sigma') = \Sigma$. Then  $\tilde {\phi}$ is a finite 
surjective  morphism of toric varieties. This  morphism induces an 
\^etale covering of open subsets 
\[ {\bf P}_{\Sigma'}^{[1]} \rightarrow {\bf P}_{\Sigma}^{[1]}  \]
with the Galois group ${\rm Coker}\,\lbrack N'\rightarrow N 
\rbrack$. }
\label{final}
\end{exam}

Recall  the following  description of the fundamental group 
of toric varieties. 

\begin{prop}
The fundamental group of a toric variety ${\bf P}_{\Sigma}$ is isomorphic 
to the quotient of $N$ by sum of all sublattices $N(\sigma)$ where $\sigma$ 
runs over all cones $\sigma \in \Sigma$. In particular, the fundamental 
group of the non-singular open toric subvariety ${\bf P}_{\Sigma}^{[1]}$ 
is isomorphic to the quotient of $N$ by sublattice spanned by all 
primitive integral generators  of $1$-dimensional cones $\sigma \in 
\Sigma^{(1)}$. 
\label{fund.group}
\end{prop}
\bigskip

\subsection{Toric Fano varieties}

 \hspace*{\parindent} 

We have seen already that the construction of an integral polyhedron $\Delta$ 
from a fan $\Sigma$ such that ${\bf P}_{\Delta} \cong {\bf P}_{\Sigma}$ is 
not unique and depends on the choice of an integral 
 strictly upper convex support function $h$. 
However, there exists  an important case when we can make 
a natural choice of $h$.

Indeed, for every 
${\bf Q}$-Gorenstein toric variety ${\bf P}_{\Sigma}$ we have the  unique 
support function $h_K$ corresponding to the anticanonical 
divisor $-K_{\Sigma} = {\bf P}_{\Sigma} \setminus {\bf T}$. 

\begin{opr}
{\rm A toric variety ${\bf P}_{\Sigma}$ is said to be a 
toric {\em ${\bf  Q}$-Fano variety} if  the anticanonical support function 
$h_K$ is strictly upper convex. A toric ${\bf  Q}$-Fano variety is called 
a {\em Gorenstein toric Fano variety} if $h_K$ is integral. }
\label{def.fano1}
\end{opr}

\begin{rem}
{\rm The above definition is a toric specialization of the notion of 
${\bf Q}$-Fano variety as a normal algebraic variety whose 
anticanonical divisor is an ample ${\bf Q}$-Cartier divisor. 
This specialization follows from the standard criterion of ampleness 
for divisors on toric varieties.}
\end{rem}

With a toric ${\bf Q}$-Fano variety, one can associate two convex 
polyhedra:
\[ \Delta({\Sigma}, h_K) = \bigcap_{\sigma \in \Sigma^{(n)}} 
( -k_{\sigma} + {\check \sigma} ) \]
and 
\[ \Delta^*(\Sigma, h_K) = \{ y \in N_{\bf Q} \mid h_K(y) \leq 1 \}. \]
The polyhedra $\Delta({\Sigma}, h_K)$ and $\Delta^*(\Sigma, h_K)$ belong  
  to $M_{\bf Q}$ and $N_{\bf Q}$ respectively. 
\medskip

Using \ref{sing1}(ii)  and \ref{simp}(ii), 
one can formulate properties of toric Fano varieties in terms of 
the properties of the polyhedra $\Delta({\Sigma}, h_K)$ and 
$\Delta^*({\Sigma}, h_K)$. 

\begin{prop}
A complete toric variety ${\bf P}_{\Sigma}$ with only ${\bf Q}$-Gorenstein 
singularities is a ${\bf Q}$-Fano variety {\rm ($see$ \ref{def.fano1})} 
if and only if $\Delta({\Sigma}, h_K)$ is an 
$n$-dimensional polyhedron with vertices $-k_{\sigma}$ 
one-to-one corresponding  
to $n$-dimensional cones $\sigma \in \Sigma$.  In this situation, 

{\rm (i)} ${\bf P}_{\Sigma}$ is a Fano variety with only Gorenstein 
singularities if and only if all vertices of $\Delta({\Sigma}, h_K)$ 
belong to $M$, in particular, 
${\bf P}_{\Sigma} = {\bf P}_{\Delta({\Sigma, h_K})}$. 

{\rm (ii)} ${\bf P}_{\Sigma}$ is a smooth Fano variety if and only if for 
every $n$-dimensional cone $\sigma \in \Sigma^{(n)}$ the polyhedron  
$\Delta^*({\Sigma}, h_K) \cap \sigma$ is a  regular simplex 
of dimension $n$. 
\label{smooth.fano}
\end{prop}

Now we come to a statement which will be very important in the sequel.

\begin{theo}
Let ${\bf P}_{\Sigma}$ be a toric Fano variety with only Gorenstein  
singularities. Then ${\bf P}_{\Sigma}$ admits at least one  
MPCP-desingulari\-zation 
\[ \tilde {\phi} \; : \; {\bf P}_{\Sigma'} \rightarrow 
{\bf P}_{\Sigma}. \]
Moreover,  MPCP-desingulari\-zations of 
${\bf P}_{\Sigma}$ are defined by maximal projective triangulations of  
the polyhedron 
$ \Delta^*(\Sigma, h_K)$, 
where $h_K$ the integral strictly upper convex support function associated 
with the anticanonical divisor ${\bf P}_{\Sigma} \setminus {\bf T}$ on 
${\bf P}_{\Sigma}$. 
\label{crep.fano}
\end{theo}

\proof  Define  a finite subset $A$ in $N$ as follows 
\[ A = \{ y \in N \mid h_K(y) \leq 1 \} = N \cap \Delta^*({\Sigma}, h_K). \] 
It is clear that $A$ is an admissible subset of $\Delta^*(\Sigma, h_K)$. 
By \ref{triang.ex}, there exists 
at  least one projective $A$-triangulation ${\cal T}$ of 
$\Delta^*(\Sigma, h_K)$. Let 
$B$ be the boundary of $\Delta^*(\Sigma, h_K)$. For every 
simplex $\theta \in B$, we construct a 
convex cone $\sigma_{\theta}$ supporting $\theta$. By definition 
\ref{def.triang}, the set of all cones $\sigma_{\theta}$ 
($\theta \in {\cal T} \cap B$) defines a fan $\Sigma'$ 
which is a subdivision of $\Sigma$. 

Since generators of $1$-dimensional 
cones of $\Sigma$ are exactly elements of $A \cap B$, the morphism 
${\bf P}_{\Sigma'} \rightarrow {\bf P}_{\Sigma}$ is crepant 
(see \ref{crit.crep}).  

By \ref{opr.t.proj}, there exists a strictly upper convex function 
$\alpha({\cal T})$. We can also assume that $\alpha({\cal T})$ has 
zero value at $0 \in N$. Then $\alpha({\cal T})$ defines a 
strictly convex support function for the fan $\Sigma'$. Thus,  
${\bf P}_{\Sigma'}$ is projective.   

By \ref{simp}(i) and \ref{rem.max}, we obtain that the morphism 
${\bf P}_{\Sigma'} \rightarrow {\bf P}_{\Sigma}$ is 
a $MPCP$-desingulari\-zation. 
	
By similar arguments, one can prove that any $MPCP$-desingulari\-zation defines 
a maximal projective triangulation of $\Delta^*(\Sigma, h_K)$. 
\bigskip

\section{Hypersurfaces in toric varieties}

\subsection{Regularity conditions for hypersurfaces}

\hspace*{\parindent}

A Laurent polynomial $f = f(X)$  is a finite linear 
combination of elements of $M$ 
\[ f(X)  = \sum c_m X^m \] 
with complex coefficients $c_m$. The Newton polyhedra $\Delta(f)$ 
of $f$ is the convex 
hull in $M_{\bf Q} = M \otimes {\bf Q}$ of all elements $m$ such 
that  $c_m \neq 0$. Every Laurent polynomial $f$ with the Newton 
polyhedron $\Delta$ defines the  affine hypersurface 
\[ Z_{f, \Delta} = \{ X \in {\bf T}  \mid f(X) = 0 \}. \]
If we work with a fixed Newton polyhedron $\Delta$, we denote 
$Z_{f, \Delta}$ simply by $Z_f$. 
\bigskip

Let $\overline{Z}_{f,\Delta}$ be the 
closure of $Z_{f, \Delta} \subset {\bf T}$ in  ${\bf P}_{\Delta}$. 
For any face $\Theta \subset \Delta$, we put 
$Z_{f, \Theta} = \overline{Z}_{f, \Delta} \cap {\bf T}_{\Theta}$. 
So we obtain  the  induced 
decomposition into the disjoint union 
\[ \overline{Z}_{f, \Delta} = \bigcup_{\Theta \subset \Delta} 
Z_{f, \Theta}. \]
\medskip

\begin{opr}
{\rm 
Let $L(\Delta)$ be the space of all Laurent polynomials  with a fixed 
Newton polyhedron $\Delta$.
A Laurent polynomial $f \in L(\Delta)$ and the corresponding 
hypersurfaces $Z_{f, \Delta} \subset {\bf T}_{\Delta}$,  
$\overline{Z}_{f, \Delta} \subset {\bf P}_{\Delta}$ are  said to be 
${\Delta}$-{\em regular} if for every  face  
$\Theta \subset \Delta$ the affine variety $Z_{f,\Theta}$ 
is empty  or a smooth subvariety of codimension 1 in ${\bf T}_{\Theta}$. 
Affine varieties $Z_{f,\Theta}$  are called  {\em the strata} on 
$\overline{Z}_{f, \Delta}$ associated with  faces $\Theta \subset 
\Delta$.} 
\label{d.nondeg}
\end{opr}

\begin{rem}
{\rm Notice that if $f$ is $\Delta$-regular,  
then ${Z}_{f, \Theta} = \emptyset$  if and only if ${\rm dim}\, \Theta = 0$, 
i.e., if and only if $\Theta$  is a vertex of $\Delta$. }
\label{zero.orb}
\end{rem}
\bigskip

Since the space $L(\Delta)$ can be identified with the space of global 
sections of the ample sheaf ${\cal O}_{\Delta}(1)$ on 
${\bf P}_{\Delta}$, using  Bertini theorem, we obtain: 

\begin{prop}
The set of $\Delta$-regular hypersurfaces is a Zariski open subset 
in $L(\Delta)$. 
\end{prop} 

We extend the notion of $\Delta$-regular hypersurfaces in 
${\bf P}_{\Delta}$ to the case of hypersurfaces in general toric varieties 
${\bf P}_{\Sigma}$ associated with rational polyhedral fans $\Sigma$. 
\bigskip

\begin{opr}
{\rm Let $\overline{Z}_{f, \Sigma}$ be the closure in ${\bf P}_{\Sigma}$ 
of an  affine hypersurface $Z_f$ defined by a Laurent polynomial $f$. Consider 
the  induced decomposition into the disjoint union 
\[ \overline{Z}_{f, \Sigma } = \bigcup_{\sigma  \in  \Sigma} 
Z_{f, \sigma},  \]
where ${Z}_{f, \sigma } = \overline{Z}_{f, \Sigma } \cap 
{\bf T}_{\sigma}$. 
A Laurent polynomial $f$ and the corresponding 
hypersurfaces $Z_{f} \subset {\bf T}$,  
$\overline{Z}_{f, \Sigma} \subset {\bf P}_{\Sigma}$ are  said to be 
${\Sigma}$-{\em regular} if for every  $s$-dimensional  cone  
$\sigma  \in  \Sigma$ the variety $Z_{f,\sigma}$ 
is empty  or  a smooth subvariety of codimension 1 in ${\bf T}_{\sigma}$. 
In other words, $\overline{Z}_{f, \Sigma}$ 
has only transversal intersections with all ${\bf T}$-orbits 
${\bf T}_{\sigma}$ ($\sigma \in \Sigma$). 
Affine varieties $Z_{f,\sigma}$  are called  {\em strata} on 
$\overline{Z}_{f, \Sigma}$ associated with the cones  $\sigma  \subset 
\Sigma$. 

Denote by $Z_{f, \Sigma}^{[i]}$ the open subvariety 
in $\overline{Z}_{f, \Sigma}$ defined as follows 
\[ Z_{f, \Sigma}^{[i]} : = \bigcup_{\sigma  \in  \Sigma^{[i]}} Z_{f, \sigma} = 
\overline{Z}_{f, \Sigma} \cap {\bf P}_{\Sigma}^{[i]}.   \]}
\end{opr}
\medskip

Let $\sigma$ be an $s$-dimensional cone in $\Sigma$. 
If we apply the  implicit function theorem and the standard 
criterion of smoothness 
to the affine hypersurface $Z_{f,\sigma} \subset {\bf T}_{\sigma}$ contained 
in the open $n$-dimensional affine toric variety 
\[ {\bf A}_{\sigma, N}   \cong ({\bf C}^*)^{n-s} \times 
{\bf A}_{\sigma, N(\sigma)},  \]
then we obtain: 
\medskip

\begin{theo}
Small analytical neighbourhoods of points   on  a    
$(n -s -1)$-dimensional stratum $Z_{f, \sigma} \subset 
\overline{Z}_{f, \Sigma}$ are analytically isomorphic 
to products  of a $(s-1)$-dimensional open ball and a small analytical 
neighbourhood  of 
the point $p_{\sigma}$ on  the $(n -s)$-dimensional affine toric 
variety ${\bf A}_{\sigma, N(\sigma)}$. 
\label{anal}
\end{theo}

For the case of $\Delta$-regular hypersurfaces in a toric variety 
${\bf P}_{\Delta}$ associated with an integral polyhedron $\Delta$, 
one gets from \ref{def.fan} the following. 

\begin{coro}
For any $l$-dimensional face $\Theta \subset \Delta$, small  
analytical neighbourhoods of points  on the   
$(l -1)$-dimensional stratum $Z_{f, \Theta} \subset 
\overline{Z}_{f, \Delta}$ are analytically isomorphic 
to products  of a $(l-1)$-dimensional open ball and a small  analytical 
neighbourhood of 
the point $p_{\sigma(\Theta)}$ on  the $(n -l)$-dimensional affine toric 
variety ${\bf A}_{\sigma(\Theta), N(\sigma(\Theta))}$.    
\label{anal2}
\end{coro}

Applying  \ref{codim1} and \ref{codim2,3}, we also conclude:  

\begin{coro}
For any $\Sigma$-regular hypersurface $\overline{Z}_{f, \Sigma} \subset 
{\bf P}_{\Sigma}$, the 
open subset $Z_{f, \Sigma}^{[1]}$ consists of smooth points of 
$\overline{Z}_{f, \Sigma}$. Moreover, 

{\rm (i)} $Z_{f, \Sigma}^{[2]}$ consists of smooth points if 
${\bf P}_{\Sigma}$ has only terminal singularities. 

{\rm (ii)} $Z_{f, \Sigma}^{[3]}$ consists of smooth points if 
${\bf P}_{\Sigma}$ has only ${\bf Q}$-factorial Gorenstein 
terminal singularities. 

{\rm (iii)} $Z_{f, \Sigma}^{[n-1]} = \overline{Z}_{f, \Sigma}$ 
is  smooth  if and only if ${\bf P}_{\Sigma}^{[n-1]}$ is smooth.  
\label{sing.hyp}
\end{coro}
\bigskip

\subsection{Birational and finite morphisms of hypersurfaces}

\hspace*{\parindent}

\begin{prop}
Let $\phi\; : \Sigma' \rightarrow \Sigma$  be a subdivision of a fan $\Sigma$, 
\[ \tilde {\phi} \; : \; {\bf P}_{\Sigma'} \rightarrow 
{\bf P}_{\Sigma} \]
the corresponding proper birational morphism. Then 
for any $\Sigma$-regular hypersurface $\overline{Z}_f \subset 
{\bf P}_{\Sigma}$ the hypersurface $\overline{Z}_{{\tilde {\phi}}^{*}f} 
\subset {\bf P}_{\Sigma'}$  is  $\Sigma'$-regular.
\label{subdiv.hyp} 
\end{prop} 

\proof  The statement follows from the fact that for any cone $\sigma' \in 
\Sigma'$ such that $\phi(\sigma') \subset \sigma \in \Sigma$, one has 
\[ Z_{{\tilde {\phi}}^* f, \sigma'} \cong Z_{f, \sigma} \times 
({\bf C}^*)^{{\rm dim}\, \sigma - {\rm dim}\, \sigma'}. \]
\bigskip

One can use \ref{sing.hyp} and \ref{subdiv.hyp} in order to construct  
partial desingularizations of hypersurfaces $\overline{Z}_{f, \Sigma}$.

\begin{prop}
Let ${\bf P}_{\Sigma}$ be a projective toric variety with only Gorenstein  
singularities. Assume that 
\[ \tilde {\phi} \; : \; {\bf P}_{\Sigma'} \rightarrow 
{\bf P}_{\Sigma} \] 
is a MPCP-desingularization of ${\bf P}_{\Sigma}$. Then 
$\overline{Z}_{{\tilde \phi}^*f, \Sigma'}$ is  a 
MPCP-desingularization of $\overline{Z}_{f, \Sigma}$. 
\label{max.crep}
\end{prop}

\proof  By \ref{sing1}, \ref{sing2}, \ref{anal}, 
$\overline{Z}_{f, \Sigma}$ has at most Gorenstein singularities and 
$\overline{Z}_{{\tilde \phi}^*f, \Sigma'}$  has at most ${\bf Q}$-factorial 
terminal singularities. It suffuces now to apply the adjunction formula. 
\medskip

\begin{prop}
Let $\phi \; : \; N' \rightarrow N$ be a surjective homomorphism 
of $n$-dimensional lattices, $\Sigma$ a 
fan in $N_{\bf Q}$,  $\Sigma'$ a fan in $N'_{\bf Q}$. Assume that 
$\phi(\Sigma') = \Sigma$. Let  \[ \tilde {\phi} \; : \; 
{\bf P}_{\Sigma',N'} \rightarrow {\bf P}_{\Sigma, N} \]
be the corresponding finite surjective morphism of toric varieties. 
Then for any $\Sigma$-regular hypersurface $\overline{Z}_f$ in 
${\bf P}_{\Sigma}$ the hypersurface 
\[ {\tilde {\phi}}^{-1} (\overline{Z}_f)  =   
\overline{Z}_{{\tilde {\phi}}^{*}f} \] 
is $\Sigma'$-regular. 
\label{galois}
\end{prop}

\proof  It is sufficient to observe that for any cone $\sigma \in 
\Sigma \cong \Sigma'$ the affine variety ${Z}_{{\tilde {\phi}}^{*}f, 
\sigma}$ is 
an \^etale covering of $Z_{f, \sigma}$ whose Galois group is  isomorphic 
to the cokernel of the homomorphism $N'(\sigma) \rightarrow N(\sigma)$ 
(see  \ref{final}). 
\bigskip
   
\subsection{The Hodge structure of hypersurfaces}
   
\hspace*{\parindent}

We are interested now in the calculation of the 
cohomology groups of $\Delta$-regular 
hypersurfaces $\overline{Z}_f$ in toric varieties ${\bf P}_{\Delta}$. 
The main difficulty in this calculation is connected with 
singularities of $\overline{Z}_f$. So we will try to avoid singularities and 
to calculate cohomology groups not only of compact compact hypersurfaces 
$\overline{Z}_f$, but also ones of some naturally defined smooth open 
subsets in $\overline{Z}_f$. For the last purpose, it is more convenient 
to use cohomology with compact supports which we denote by $H^i_c(*)$. 
\medskip

First we note that there exists the following 
analog of the  Lefschetz theorem for $\Delta$-regular hypersurfaces 
proved by Bernstein, Danilov and Khovanski\^i \cite{dan.hov}.

\begin{theo} 
For any open toric subvariety $U \subset {\bf P}_{\Delta}$, 
the Gysin homomorphism 
\[ H^i_c ( \overline{Z}_f \cap U) \rightarrow H^{i+2}_c (U) \]
is bijective for $i > n -1$ and injective for $i = n - 1$. 
\label{Lef}
\end{theo}

Using this theorem, one can often reduce the calculation  of 
cohomology groups $H^i$ (or $H^i_c$) to the "interesting"  
case $i = n-1$. We consider below several typical examples of 
such a situation.

\begin{exam} 
{\rm Let $U$ be a smooth open toric subvariety in ${\bf P}_{\Delta}$ 
(e.g., $U = {\bf T}$). Then 
$V = U \cap \overline{Z}_f$ is a smooth affine open subset in 
$\overline{Z}_f$. By general properties of Stein varieties, one has 
$H^{i}(V) = 0$ for $i > n-1$.  Since the calculation of cohomology groups of 
smooth affine toric varieties is very simple, we obtain a complete 
information about all cohomology groups except for $i = n-1$ using the 
following property. }
\end{exam}

\begin{prop}
Let  $W$ be a quasi-smooth irreducible $k$-dimensional algebraic variety. Then 
there exists the Poincare pairing 
\[ H_c^i (W) \otimes H^{2k -i} (W) \rightarrow H^{2k}_c (W) \cong {\bf C}.  \]
This pairing is compatible with Hodge structures, where $H^{2k}_c(W)$ is 
assumed to be a $1$-dimensional ${\bf C}$-space of the Hodge type $(k,k)$. 
\label{duality}
\end{prop}

If we take $U = {\bf T}$, then $V = Z_f$. The Euler characteristic  of 
$Z_f$ was calculated by 
Bernstein, Khovanski\^i and Kushnirenko \cite{kuch}, \cite{hov.genus}.

\begin{theo}
$ e(Z_f) = \sum_{i \geq 0} (-1)^i {\rm dim}\, H^i(Z_f)  = 
(-1)^{n-1} d_M (\Delta)$. 
\label{euler}
\end{theo}

In particular, we obtain: 

\begin{coro}
The dimension of $H^{n-1}(Z_f)$ is equal to $d_M(\Delta) + n -1$. 
\end{coro}

\begin{opr}
Let $P$ be a compact convex subset in $M_{\bf Q}$.  Denote by $l(P)$ 
the number of integral points in $P \cap M$, and by $l^*(P)$ 
the number of integral points in the interior of $P$. 
\end{opr}
 
There exist general formulas  for the Hodge-Deligne numbers 
$h^{p,q}(Z_f)$ of the mixed  Hodge structure in the $(n-1)$-th cohomology 
group of an arbitrary 
$\Delta$-regular affine hypersurface in ${\bf T}$ (Danilov and Khovanski\^i  
\cite{dan.hov}).
For the numbers $h^{n-2,1}(Z_f)$ and $h^{n-2,0}(Z_f)$,  we get the following 
(see \cite{dan.hov}, 5.9).

\begin{prop}
Let ${\rm dim}\, \Delta = n \geq 4$, then 
\[ h^{n-2,1} (Z_f) + h^{n-2,0} (Z_f) = l^*(2\Delta) - (n+1)l^*(\Delta), \]
\[ h^{n-2,0} (Z_f) =  \sum_{{\rm codim}\,  \Theta =1} l^*(\Theta).\]
\label{hd.dh}
\end{prop}
\bigskip

We will use in the sequel the following properties 
of the Hodge-Deligne numbers $h^{p,q}(H^k_c(Z_f)$ of affine hypersurfaces 
for cohomology with compact supports (see \cite{dan.hov}):

\begin{prop}
The Hodge-Deligne numbers $h^{p,q}(H^k_c(Z_f)$ of $\Delta$-regular 
$(n-1)$-dimensional affine hypersurfaces $Z_f$ satisfy the properties:

{\rm (i)} $h^{p,q}(H^k_c(Z_f) = 0$ for $p \neq q$ and $k > n-1$; 

{\rm (ii)} $h^{p,q}(H^k_c(Z_f)) = 0$ for  $k < n-1$;

{\rm (iii)} $h^{p,q}(H^{n-1}_c(Z_f) = 0$ for  $p+ q  > n-1$. 
\label{properties.hodge}
\end{prop}
\bigskip

Although we consider the mixed Hodge structure in the cohomology group 
of the affine hypersurface $Z_f$, we  get eventually some information about 
the Hodge numbers of compactifications of $Z_f$. From general 
properties of the mixed Hodge structures \cite{deligne}, one obtains:   

\begin{theo}
Let $\overline{Z}$ be a smooth compactification of a smooth affine 
$(n-1)$-dimensional variety $Z$ such that the complementary set  
$\overline{Z} \setminus Z$ is a normal crossing divisor. Let  
\[ j\;: \;  Z \hookrightarrow \overline{Z} \] 
be the corresponding embedding. Denote by  
\[ j^* \;:\; H^{n-1} (\overline{Z}) \rightarrow H^{n-1} (Z) \] 
the induced mapping of cohomology groups. Then 
the weight subspace ${\cal W}_{n-1}H^{n-1}(Z)$ coincides with 
the image  $j^*( H^{n-1}(\overline{Z}))$.  In particular, one has 
the following inequalities between the Hodge numbers of $\overline{Z}$ and 
the Hodge-Deligne numbers of $Z$:
\[ h^{n-1-k,k}(\overline{Z}) \geq h^{n-1-k,k}(H^{n-1}(Z))\;\; 
(0 \leq k \leq n-1). \]
\label{inequal.h}
\end{theo}
\bigskip

\section{Calabi-Yau hypersurfaces in toric varieties}
\medskip

\subsection{Reflexive polyhedra and reflexive pairs}

\hspace*{\parindent}

\begin{opr} 
{\rm If  $P$ is an arbitrary 
compact convex set in $M_{\bf Q}$ containing the  
zero vector $0 \in M_{\bf Q}$ in its interior, then we call 
\[ P^* = \{ y \in N_{\bf Q} \mid \langle x, y \rangle \geq -1, \; 
{\rm for \; all} \; x \in P  \}. \]
{\em the dual set}  relative to $P$.  
 }
\label{dual} 
\end{opr}
\bigskip

The dual set $P^*$ is  a convex compact subset in $N_{\bf Q}$ 
containing  the vector zero $0 \in N_{\bf Q}$ in its interior. 
Obviously, one has $(P^*)^* = P$. 
\bigskip

\begin{exam}
{\rm Let $E$ be a Euclidian $n$-dimensional space, $\langle *, * \rangle$ 
the corresponding scalar product,  
\[ P = \{ {x} \mid \langle { x}, { x} \rangle \leq R \} \]
the  ball of  radius $R$. Using the scalar product, we can identify 
the dual space $E^*$ with $E$. Then  the dual set $P^*$ is 
the  ball of radius $1/R$. }
\end{exam}

\begin{prop}
Let $P \subset M_{\bf Q}$ be a convex set containing $0$ in its interior, 
$C_P \subset \overline{M}_{\bf Q}$ the convex cone supporting $P$, $
\overline{N}_{\bf Q} = {\bf Q} \oplus N_{\bf Q}$ the dual space, 
$C_{P^*} \subset \overline{N}_{\bf Q}$ 
the convex cone supporting $P^* \subset N_{\bf Q}$. Then $C_{P^*}$ is the 
dual cone relative to $C_P$. 
\label{dual.cones}
\end{prop}

\proof  Let $\overline{x} = (x_0, x) \in \overline{M}_{\bf Q}$, 
$\overline{y} = (y_0, y) \in \overline{M}_{\bf Q}$. 
Since $x_0$ and $y_0$ are positive, one has 
\[ \overline{x} \in C_P \Leftrightarrow {x}/{x_0} \in P \;\; {\rm and }\;\;  
\overline{y} \in C_{P^*} \Leftrightarrow {y}/{y_0} \in P^*. \]
Therefore, 
$\langle \overline{x}, \overline{y} \rangle = x_0y_0 + 
\langle {x}, {y} \rangle \geq 0 $ 
if and only if 
$\langle {x}/{x_0}, {y}/{y_0} \rangle \geq -1$. 
\bigskip

\begin{opr}
{\rm Let 
$H$ be a rational affine hyperplane in $M_{\bf Q}$, $p \in M_{\bf Q}$ an 
arbitrary integral point. Assume that $H$ is affinely generated by 
integral points $H \cap M$, i.e., there exists  a primitive integral element 
$l \in N$  such that for some integer $c$ 
\[ H = \{ x \in M_{\bf Q} \mid \langle x, l \rangle = c \}. \]  
Then the absolute value 
$\mid c - \langle p , l \rangle \mid$ is called the {\em integral distance} 
between $H$ and $p$. }
\end{opr}

\begin{opr}
{\rm Let $M$ be an integral $n$-dimensional lattice in $M_{\bf Q}$,   
$\Delta$ a convex integral polyhedron in $M_{\bf Q}$ containing the zero 
 $0 \in M_{\bf Q}$  in its  interior.  
The pair $(\Delta, M)$ is called {\em reflexive} if 
the integral distance between 
 $0$ and all affine hyperplanes   generated by  $(n -1)$-dimensional faces   
of $\Delta$ equals 1. 

If $(\Delta, M)$ is a reflexive pair, then we call $\Delta$ a {\em 
reflexive polyhedron}.}
\label{inver.p}
\end{opr}
\bigskip

The following simple property of reflexive polyhedra is very important.
\medskip

\begin{theo}
Suppose that $(\Delta, M)$ is a reflexive pair. Then $(\Delta^*, N)$ is 
again a reflexive pair.
\end{theo}

\proof  Let $\Theta_1, \ldots , \Theta_k$ be $(n-1)$-dimensional faces of 
$\Delta$, $H_1, \ldots , H_k$ the corresponding affine hyperplanes. 
By \ref{inver.p}, there exist integral elements 
$l_1, \ldots , l_k \in N_{\bf Q}$ such that for all $1 \leq i \leq k$
\[ \Theta_i = \{ x \in \Delta \mid \langle x , l_i \rangle = 1 \},\;\; 
 H_i = \{ x \in M_{\bf Q} \mid \langle x , l_i \rangle = 1 \}. \]
Therefore,
\[ \Delta = \{ x \in M_{\bf Q}  \mid \langle x , l_i \rangle \leq 1 \;
(1 \leq i \leq k ) \}. \]
So  $\Delta^*$ is a convex hull of 
the integral points $l_1, \ldots , l_k$, i.e., $\Delta^*$ is an 
integral polyhedron.  

Let $p_1, \ldots , p_m$ be vertices of $\Delta$. By  \ref{dual}, 
for any $j$ ($1 \leq j  \leq m$) 
\[ \Xi_j = \{ y \in \Delta^*  \mid \langle p_j , y \rangle = 1 \} \]
is a $(n-1)$-dimensional face of $\Delta^*$. 
Thus, $\Delta^*$ contains $0 \in N_{\bf Q}$  in its interior, 
and the integral distance between $0$ and every  $(n-1)$-dimensional 
affine linear subspace generated by $(n-1)$-dimensional 
faces of $\Delta^*$ equals $1$. 
\bigskip

We can establish the following one-to-one correspondence between faces of 
the polyhedra $\Delta$ and $\Delta^*$. 

\begin{prop}
Let $\Theta$ be an $l$-dimensional face of an $n$-dimensional reflexive 
polyhedron $\Delta \subset M_{\bf Q}$, $p_1, \ldots, p_k$ are vertices of $\Theta$, 
$\Delta^* \in N_{\bf Q}$ the dual reflexive polyhedron. 

Define the dual to $\Theta$ $(n -l -1)$-dimensional face of $\Delta^*$ as 
\[ \Theta^* = \{ y \in \Delta^* \mid \langle p_1, y \rangle  = \cdots = 
\langle p_k, y \rangle = - 1 \}.  \]

Then one gets the one-to-one  
correspondence $\Theta \leftrightarrow \Theta^*$  between faces of 
the polyhedra $\Delta$ and $\Delta^*$ reversing the incidence relation of the 
faces.   
\label{dual.edge}
\end{prop}
    
   The next theorem describes the relationship between reflexive 
polyhedra and Calabi-Yau hypersurfaces.   
   
\begin{theo}
Let $\Delta$ be an $n$-dimensional 
integral polyhedron in $M_{\bf Q}$, ${\bf P}_{\Delta}$ the corresponding 
$n$-dimensional projective toric variety, ${\cal F}(\Delta)$ the family 
of projective $\Delta$-regular hypersurfaces $\overline{Z}_f$ in 
${\bf P}_{\Delta}$.  
Then the following conditions are equivalent 

{\rm (i)} the family ${\cal F}(\Delta)$ of $\Delta$-hypersurfaces in 
${\bf P}_{\Delta}$ consists of Calabi-Yau varieties with 
canonical singularities $($see {\rm \ref{calabi-yau}}$);$

{\rm (ii)} the ample invertible sheaf ${\cal O}_{{\Delta}}(1)$ 
on the toric variety ${\bf P}_{\Delta}$ is anticanonical, i.e., 
${\bf P}_{\Delta}$ is a toric Fano variety with Gorenstein 
singularities;

{\rm (iii)} $\Delta$ contains only one integral point $m_0$ in its 
interior, and $(\Delta-m_0, M)$ is a reflexive pair. 
\label{equiv}
\end{theo}

\proof  Since $\overline{Z}_f$ is an ample Cartier divisor on 
${\bf P}_{\Delta}$, (i)$\Rightarrow$(ii) follows from the adjunction 
formula. 

The equivalence  (ii)$\Leftrightarrow$(iii) 
follows from \ref{smooth.fano}. 

Assume that (ii) and (iii) are satisfied. Let us prove (i). 
By the adjunction formula, it follows from (ii) that every hypersurface 
$\overline{Z}_f$ has 
trivial canonical divisor. By the vanishing theorem for 
arbitrary ample divisors on toric varieties \cite{dan1}, one  gets 
\[ H^i(\overline{Z}_f, {\cal O}_{\overline{Z}_f}) = 0\]
for $ 0 < i < n-1$. 
By \ref{anal2}, singularities of  $\overline{Z}_f$ are 
analytically isomorphic to toric singularities of ${\bf P}_{\Delta}$. 
Since all  singularities of ${\bf P}_{\Delta}$ are Gorenstein, by  
\ref{gor.can}, they are also  canonical. So every $\Delta$-regular 
hypersurface satisfies \ref{calabi-yau}. 
\bigskip

The next statement follows  from definitions of the polyhedra 
$\Delta(\Sigma, h_K)$ and $\Delta^*(\Sigma, h_K)$.

\begin{prop}
Let $(\Delta,M)$ be a reflexive pair, $(\Delta^*,N)$ the dual 
reflexive pair, $\Sigma$ the rational polyhedral fan 
defining the corresponding Gorenstein toric Fano variety ${\bf P}_{\Delta}$. 
Then 
\[ \Delta(\Sigma, h_K) = \Delta, \]
\[ \Delta^*(\Sigma, h_K) = \Delta^*. \]
In particular, if $\Sigma^*$ is a rational polyhedral fan defining 
the Gorenstein toric Fano variety ${\bf P}_{\Delta^*}$, then
\[ \Delta(\Sigma, h_K) = \Delta^*(\Sigma^*, h_K) \]
and 
\[ \Delta(\Sigma^*, h_K) = \Delta^*(\Sigma, h_K). \]
\end{prop}

Thus, in order to construct a rational polyhedral fan $\Sigma(\Delta)$ 
corresponding to a reflexive polyhedron $\Delta$, one can use the following 
another way: we take the dual reflexive polyhedron $\Delta^*$ and apply 
the statement
 
\begin{coro}
Let $\Delta \subset M_{\bf Q}$ be a reflexive polyhedron, $\Delta^*$ 
the dual reflexive polyhedron in $N_{\bf Q}$. For every 
$l$-dimensional face $\Theta$ of $\Delta^*$ define the $(l+1)$-dimensional 
convex cone $\sigma[\Theta]$ supporting the face  $\Theta$ as follows 
\[ \sigma[\Theta] = \{ \lambda x \in M_{\bf Q} \mid 
x \in \Theta,\;\; \lambda \in {\bf Q}_{\geq 0} \}. \] 
Then the set $\Sigma[\Delta^*]$ of all cones $\sigma[\Theta]$ where 
$\Theta$ runs over all faces of $\Delta^*$ is the complete fan defining 
the toric Fano variety ${\bf P}_{\Delta}$ associated with $\Delta$. 
Moreover, every 
$(l+1)$-dimensional cone $\sigma[\Theta]$ coincides with 
$\sigma(\Theta^*)$ 
{\rm ($see$ \ref{def.fan})}, where $\Theta^* \subset \Delta$ is the dual 
to $\Theta$ $(n-l-1)$-dimensional face of $\Delta$ 
{\rm ($see$ \ref{dual.edge})}.   
\end{coro}
\bigskip

\subsection{Singularities and morphisms of Calabi-Yau hypersurfaces}

\hspace*{\parindent}

Let $\Delta$ be a reflexive polyhedron, $\Delta^*$ the dual reflexive 
polyhedron. Take a maximal projective triangulation ${\cal T}$ of $\Delta$. 
It follows from the proof of \ref{crep.fano} that ${\cal T}$ defines a 
$MPCP$-desingularization 
\[ \varphi_{\cal T}\; :\; {\hat{\bf P}}_{\Delta} \rightarrow 
{\bf P}_{\Delta} \]
of the Gorenstein toric variety ${\bf P}_{\Delta}$. Let $\overline{Z}_f$ 
be a $\Delta$-regular Calabi-Yau hypersurface in ${\bf P}_{\Delta}$. 
Put 
\[\hat{Z}_f = \varphi^{-1}_{\cal T}(\overline{Z}_f).\]
By \ref{max.crep}, 
\[ \varphi_{\cal T}\; : \; 
 \hat{Z}_f \rightarrow \overline{Z}_f \] 
 is a $MPCP$-desingularization of 
$\overline{Z}_f$. 

\begin{opr}
{\rm We will call 
\[ \varphi_{\cal T}\; : \;  \hat{Z}_f \rightarrow \overline{Z}_f \] 
the {\em toroidal  MPCP-desingularization of $\overline{Z}_f$ corresponding to 
a maximal projective triangulation ${\cal T}$ of $\Delta^*$}.}
\end{opr}
 
Using \ref{crep.fano} and  \ref{codim2,3}(ii), one gets the following. 
 
\begin{theo}
There exist at least one toroidal MPCP-desingularization $\hat{Z}_f$ of any 
$\Delta$-regular Calabi-Yau hypersurface in ${\bf P}_{\Delta}$. Such a  
MPCP-desingularization $\hat{Z}_f$ corresponds to any  maximal 
projective triangulation ${\cal T}$ of the dual polyhedron $\Delta^*$. 
The codimension of singularities of $\hat{Z}_f$ is always at least $4$. 
\label{smooth.c}
\end{theo}

\begin{coro}
A toroidal MPCP-desingularization of a projective Calabi-Yau hypersurace 
$\overline{Z}_f$ associated with a  reflexive polyhedron  
$\Delta$ of dimension $n \leq 4$ is always a  smooth Calabi-Yau manifold. 
\label{smooth.c1}
\end{coro}

Let ${\cal T}$ be a maximal projective triangulation of $\Delta^*$. For 
any  $l$-dimensional face $\Theta$ of $\Delta$, the restriction of 
${\cal T}$ on the dual $(n - l -1)$-dimensional face $\Theta^* \subset 
\Delta^*$ is a maximal projective triangulation ${\cal T}\mid_{\Theta^*}$ of 
$\Theta^*$. By \ref{anal}, the analytical decription  of singularities along 
a stratum $Z_{f,\Theta}$  as well as of their  $MPCP$-desingularizations 
reduces to 
the combinatorial description of a $MPCP$-desingularization of 
the toric singularity at the unique 
closed ${\bf T}_{\sigma}$-invariant point $p_{\sigma}$  on 
the $(n-l)$-dimensional affine toric variety ${\bf A}_{\sigma, N(\sigma)}$, 
 where $\sigma =  \sigma[\Theta^*]$. So we introduce the following definition. 

\begin{opr} 
{\rm  We call the face $\Theta^*$ of the polyhedron $\Delta$ 
the {\em diagram} of the toric singularity 
at $p_{\sigma} \in  {\bf A}_{\sigma, N(\sigma)}$. A maximal projective 
triangulation ${\cal T}\mid_{\Theta^*}$ of $\Theta^*$ induced by a 
maximal projective triangulation ${\cal T}$ of $\Delta^*$ 
we call {\em a triangulated diagram}  of the toric singularity at 
$p_{\sigma} \in  {\bf A}_{\sigma, N(\sigma)}$.   }
\label{diagram}
\end{opr}

We have seen already in \ref{anal2} that for any face $\Theta \subset \Delta$ and 
any $\Delta$-regular Laurent polynomial $f \in L(\Delta)$ 
local neighbourhoods of points belonging to the same stratum $Z_{f, \Theta}$ 
are analitically  isomorphic. Thus, if a stratum $Z_{f, \Theta}$ consists of 
singular points of $\overline{Z}_f$, then all these singularities are 
analitically isomorphic. Our purpose now is to describe singularities along 
$Z_{f, \Theta}$ and their $MPCP$-desingularizations in terms of triangulated 
diagrams. 
\medskip

Let ${\cal T}\mid_{\Theta^*}$ be a triangulated diagram. Then we obtain 
a subdivision $\Sigma({\cal T},\Theta)$ of the cone $\sigma = \sigma[\Theta^*]$ 
into the union of subcones supporting elementary simplices of the 
triangulation  ${\cal T}\mid_{\Theta^*}$. By \ref{biration}, one has 
the corresponding  projective birational toric morphism 
\[ \varphi_{{\cal T}, \Theta^*}\; : \;{\bf P}_{\Sigma({\cal T},\Theta)} 
\rightarrow {\bf A}_{\sigma, N(\sigma)}. \]

\begin{theo}
For any $l$-dimensional face $\Theta \subset \Delta$ and any 
closed point $p \in Z_{f, \Theta}$, the fiber $\varphi_{\cal T}^{-1}(p)$ of 
a MPCP-desingularization $\varphi_{\cal T}$ is isomorphic to the fiber  
$\varphi_{{\cal T}, \Theta^*}(p_{\sigma})$ of the projective toric morphism 
$\varphi_{{\cal T}, \Theta^*}$. 

The number of irreducible $(n-l-1)$-dimensional components of 
$\varphi_{\cal T}^{-1}(p)$ equals $l^*(\Theta^*)$, i.e., the number of 
integral points in the interior of $\Theta^*$. Moreover, 
the Euler characteristic 
of $\varphi_{\cal T}^{-1}(p)$ equals the number of elementary simplices 
in the triangulated diagram ${\cal T}\mid_{\Theta^*}$. 
\label{topol.des}
\end{theo}

\proof  Since $ \varphi_{\cal T}\; :\; {\hat{\bf P}}_{\Delta} \rightarrow 
{\bf P}_{\Delta}$ is a birational toric morphism, the ${\bf T}$-action 
induces isomorphisms of fibers of $\varphi_{\cal T}$ over closed points of a 
${\bf T}$-stratum  ${\bf T}_{\Theta}$. 
Thus, we obtain isomorphisms between fibers 
of $\varphi_{\cal T}$ over closed points of  $Z_{f, \Theta} \subset 
{\bf T}_{\Theta}$. By \ref{orbits.fan}(i), 
${\bf T}_{\Theta}$ is contained in the 
${\bf T}$-invariant open subset ${\bf A}_{\sigma,N} \cong {\bf T}_{\Theta} 
\times {\bf A}_{\sigma, N(\sigma)}$. Thus, we have a commutative 
diagram 
\[ 
\begin{tabular}{ccc}
$ \varphi_{\cal T}^{-1}({\bf A}_{\sigma,N})$ &  $\rightarrow$ & $
{\bf A}_{\sigma,N}$ \\
$\downarrow$ &   & $\downarrow$ \\
$ {\bf P}_{\Sigma({\cal T},\Theta)} 
$ & $\rightarrow$ & ${\bf A}_{\sigma, N(\sigma)}$ 
\end{tabular} \]
whose vertical maps are divisions by the action of the torus 
${\bf T}_{\Theta}$. 
So the fiber $\varphi_{\cal T}^{-1}(p)$ over any closed point $p \in 
{\bf T}_{\Theta^*}$  is isomorphic to the fiber  
$\varphi_{{\cal T}, \Theta^*}^{-1}(p_{\sigma})$ of the projective 
toric morphism $\varphi_{{\cal T}, \Theta^*}$. Therefore, 
irreducible divisors of 
${\bf P}_{\Sigma({\cal T},\Theta)}$ over $p_{\sigma}$ one-to-one correspond 
to integral points in the interior of $\Theta^*$. 

On the other hand, the toric morphism $\varphi_{{\cal T}, \Theta^*}$ has 
an action the $(n-l)$-dimensional torus ${\bf T}_{\Theta}' = 
{\rm Ker}\, \lbrack {\bf T} \rightarrow {\bf T}_{\Theta} \rbrack$. Since the 
closed point $p_{\sigma} \in {\bf A}_{\sigma, N(\sigma)}$ is 
${\bf T}_{\Theta}'$-invariant, one has a decomposition of the fiber 
$\varphi_{{\cal T}, \Theta^*}^{-1}(p_{\sigma})$ into the union of 
${\bf T}_{\Theta}'$-orbits. Thus the Euler characteristic of 
$\varphi_{{\cal T}, \Theta^*}^{-1}(p_{\sigma})$ is the number of 
zero-dimensional ${\bf T}_{\Theta}'$-orbits. The latter equals the 
number of $(n-l)$-dimensional cones of $\Sigma({\cal T},\Theta)$, i.e., 
the number of elementary simplices in the triangulated diagram 
${\cal T}\mid_{\Theta^*}$.    
\medskip

\begin{exam}
{\rm Let $\Theta$ be an $(n-2)$-dimensional face of a reflexive polyhedron 
$\Delta$, and let $\Theta^*$ be the dual to $\Theta$ $1$-dimensional face 
of the dual 
polyhedron $\Delta^*$. There exists a unique maximal projective 
triangulation of $\Theta^*$ consisting of $d(\Theta^*)$ elementary,  
in fact,  regular segments. 

In this case, small analytical 
neighbourhoods of points on $Z_{f, \Theta}$ are analytically isomorphic to 
product of $(n-3)$-dimensional open ball and a small analytical 
neighbourhood of $2$-dimensional double point singularity of type 
$A_{d(\Theta^*)-1}$. 
  
The fiber of $\varphi_{\cal T}$ over any point of $\overline{Z}_{f,\Theta}$ 
is the Hirzebruch-Jung tree of $l^*(\Theta^*) = d(\Theta^*) -1$ 
smooth rational curves having 
an action of ${\bf C}^*$. }
\label{duval}
\end{exam}

\begin{exam}
{\rm Let $\Delta$ be the $4$-dimensional reflexive polyhedron with vertices 
\[ e_1 = (1,0,0,0),\; e_2 = (0,1,0,0),\; e_3 = (0,0,1,0),\; 
e_4 = (0,0,0,1),\; \]
\[ e_5 = (-1,-1,-1,-1). \] 
We will see later in \ref{first.ex} that the family ${\cal F}(\Delta)$ of 
Calabi-Yau hypersurfaces $\overline{Z}_f$ in ${\bf P}_{\Delta}$ 
coincides  with the mirror family ${\cal W}_{\psi}$ for quintic $3$-folds 
considered by P. Candelas et al. in \cite{cand2}. The dual polyhedron 
$\Delta^*$ is $5$-times  multiple of the regular $4$-dimensional simplex. 
The polyhedron $\Delta^*$ has ten  $1$-dimensional faces which are 
$5$-times multiples of 
a regular $1$-dimensional simplex, i.e., singularities along 
$1$-dimensional strata $Z_{f, \Theta}$ are isomorphic to one-parameter 
families of double points  of 
type $A_4$. There are also ten  $2$-dimensional faces of $\Delta^*$ which 
are $5$-times multiples of a regular $2$-dimensional simplex. We have 
$10$ isolated singular points on ${\cal W}_{\psi}$ whose 
$MPCP$-desingularizations  
are defined by the following  triangulated diagram (D. Morrison \cite{mor1}):

\begin{center}
\begin{picture}(200,200)
\put(100,100){\makebox(0,0){$\bullet$}}
\put(100,120){\makebox(0,0){$\bullet$}}
\put(100,80){\makebox(0,0){$\bullet$}}
\put(100,60){\makebox(0,0){$\bullet$}}
\put(120,100){\makebox(0,0){$\bullet$}}
\put(120,80){\makebox(0,0){$\bullet$}}
\put(120,60){\makebox(0,0){$\bullet$}}
\put(140,80){\makebox(0,0){$\bullet$}}
\put(140,60){\makebox(0,0){$\bullet$}}
\put(160,60){\makebox(0,0){$\bullet$}}
\put(80,140){\makebox(0,0){$\bullet$}}
\put(80,120){\makebox(0,0){$\bullet$}}
\put(80,100){\makebox(0,0){$\bullet$}}
\put(80,80){\makebox(0,0){$\bullet$}}
\put(80,60){\makebox(0,0){$\bullet$}}
\put(60,160){\makebox(0,0){$\bullet$}}
\put(60,140){\makebox(0,0){$\bullet$}}
\put(60,120){\makebox(0,0){$\bullet$}}
\put(60,100){\makebox(0,0){$\bullet$}}
\put(60,80){\makebox(0,0){$\bullet$}}
\put(60,60){\makebox(0,0){$\bullet$}}

\put(60,160){\line(1,-1){100}}
\put(60,140){\line(1,-1){80}}
\put(60,120){\line(1,-1){60}}
\put(60,60){\line(0,1){100}}
\put(80,60){\line(0,1){80}}
\put(100,60){\line(0,1){60}}
\put(60,60){\line(1,0){100}}
\put(60,80){\line(1,0){80}}
\put(60,100){\line(1,0){60}}
\put(60,60){\line(1,1){20}}
\put(60,80){\line(1,1){20}}
\put(80,60){\line(1,1){20}}

\put(120,80){\line(2,-1){40}}
\put(100,80){\line(2,-1){40}}
\put(100,100){\line(2,-1){40}}

\put(80,100){\line(-1,2){20}}
\put(80,120){\line(-1,2){20}}
\put(100,100){\line(-1,2){20}}
\end{picture}
\end{center}
}
\end{exam}

\begin{opr}
{\rm Let $(\Delta_1, M_1)$ and $(\Delta_2, M_2)$ be two reflexive 
pairs of equal dimension. {\em A finite morphism} of reflexive pair 
\[ \phi\; :\; (\Delta_1, M_1) \rightarrow (\Delta_2, M_2) \]
is  a homomorphism of lattices $\phi\,:\,  M_1 \rightarrow M_2$ such that 
$\phi (\Delta_1) = \Delta_2$. 
}
\label{morph}
\end{opr}

By \ref{final}, we obtain: 

\begin{prop}
Let 
\[ \phi\; :\; (\Delta_1, M_1) \rightarrow (\Delta_2, M_2) \] 
be a finite morphism of reflexive pairs. Then the dual finite morphism 
\[ \phi^* \; : \; (\Delta_2^*, N_2) \rightarrow (\Delta_1^*, N_1). \]
induces a finite  surjective morphism 
\[ \tilde {\phi}^* \; : \; {\bf P}_{\Delta_2,M_2} 
\rightarrow {\bf P}_{\Delta_1,M_1}.  \]
Moreover, the restriction 
\[ \tilde {\phi}^* \; : \; {\bf P}_{\Delta_2,M_2}^{[1]}  
\rightarrow {\bf P}_{\Delta_1,M_1}^{[1]}  \]
is an \^etale morphism of degree  
\[ d_{M_2}(\Delta_2) / d_{M_1}(\Delta_1) = 
d_{N_1}(\Delta_1^*) / d_{N_2}(\Delta_2^*).\]
\label{etale}
\end{prop} 

\proof  It remains to show that  
$\tilde {\phi}^* \; : \; {\bf P}_{\Delta_2,M_2}^{[1]}  
\rightarrow {\bf P}_{\Delta_1,M_1}^{[1]}$ is  \^etale. 
Take a $\Delta_1$-regular Calabi-Yau hypersurface $\overline{Z}_{f} 
\subset {\bf P}_{\Delta_1,M_1}^{[1]}$. By \ref{subdiv.hyp}, the 
hypersurface $\overline{Z}_{\tilde{\phi}^*f} = 
\tilde{\phi}^{-1}(\overline{Z}_f)$ is $\Delta_2$-regular. 
By \ref{sing.hyp} and \ref{equiv}, two quasi-projective  varieties 
$Z^{[1]}_f$ and 
${Z}_{\tilde{\phi}^*f}^{[1]}$ are smooth and have trivial canonical 
class. Therefore, any finite morphism of these varieties must be 
\^etale.

\subsection{The Hodge  number $h^{n-2,1}(\hat{Z}_f)$}

\hspace*{\parindent}

Let us  apply the result of Danilov and Khovanski\^i (see \ref{hd.dh}) 
to  the case of reflexive polyhedra $\Delta$ of dimension $\geq 4$. 
Using the properties $l^*(2\Delta) = l(\Delta)$, $l^*(\Delta) =1$, we 
can calculate  the Hodge-Deligne number 
$h^{n-2,1}(Z_f)$ of an affine Calabi-Yau hypersurface  
$Z_f$ as follows.   

\begin{theo}
Let $\Delta$ be a reflexive $n$-dimensional polyhedron $(n\geq 4)$, then 
the  Hodge-Deligne number $h^{n-2,1}$ of the cohomology group $H^{n-1}(Z_f)$ 
of any $(n-1)$-dimensional affine $\Delta$-regular Calabi-Yau 
hypersurface $Z_f$ equals 
\[ h^{n-2,1} (Z_f) = l(\Delta)  - n - 1 -  \sum_{{\rm codim}\, \Theta =1 } 
l^*(\Theta). \]
\label{n.def}
\end{theo}
\bigskip

In fact, we can calculate the Hodge-Deligne space $H^{n-2,1}(Z_f)$ itself 
(not only the dimension $h^{n-2,1}$). 

\begin{theo} 
{\rm \cite{bat.var}} 
Let $L^*_1(\Delta)$ be the subspace in $L(\Delta)$ generated by all 
monomials $X^m$ 
such  that $m$ is an interior integral point on a face 
$\Theta \subset \Delta$ of codimension $1$. Then the ${\bf C}$-space 
$H^{n-2,1}(Z_f)$ is canonically isomorphic to the quotient of $L(\Delta)$ by 
\[ L^*_1(\Delta) + {\bf C}\langle f, f_1, \ldots, f_n \rangle,  \]
where 
\[ f_i (X) = X_i \frac{\partial }{\partial X_i} f(X), 
\;\; (1 \leq i \leq n). \]
\end{theo}

We want now to calculate for $n geq 4$ the Hodge number $h^{n-2,1}$ of a  
{\em MPCP}-desingularization ${\hat Z}_f$ of the toroidal compactification 
$\overline{Z}_f$ of $Z_f$. 

For this purpose, it is  convenient to use the notion of the  
$(p,q)$-{\em Euler characteristic} $e_c^{p,q}$ introduced in \cite{dan.hov}. 

\begin{opr}
{\rm For any complex algebraic variety $V$, $e_c^{p,q}(V)$ is defined as the 
alternated sum of Hodge-Deligne numbers 
\[ \sum_{i \geq 0} (-1)^i h^{p,q}(H^i_c(V)). \] }
\end{opr}

\begin{prop}({\rm see \cite{dan.hov}}) 
Let  $V = V' \times V''$ be a product of two 
algebraic varieties. Then one has 
\[ e_c^{p,q}(V) = \sum_{(p'+ p'',q' +q'') = (p,q)} 
e_c^{p',q'}(V') \cdot e_c^{p'',q''}(V''). \]
\label{h.product}
\end{prop}

\begin{opr}
{\rm A stratification of a compact algebraic variety $V$ is a representation 
of $V$ as a disjoint union of finitely many locally closed smooth 
subvarieties $\{ V_j \}_{j \in J}$ (which are called {\em strata}) 
such that for any $j \in J$ the closure of $V_j$ in $V$ is a union of 
the strata.}
\end{opr}

The following property follows immediately from long cohomology sequences 
(see \cite{dan.hov}). 

\begin{prop}
Let $\{ V_j \}_{j \in J}$ be a stratification of $V$. Then 
\[ e_c^{p,q}(V) = \sum_{j \in J} e_c^{p,q}(V_j). \]
\label{addit}
\end{prop}

Returning to our Calabi-Yau variety $\hat{Z}_f$, we see that $\hat{Z}_f$ 
is always quasi-smooth. Therefore the cohomology groups $H_c^i(\hat{Z}_f) 
\cong H^i(\hat{Z}_f)$ have the pure Hodge structure of weight $i$. So, 
by \ref{addit},   
it suffices to calculate the $(n-2,1)$-Euler characteristic 
\[ e_c^{n-2,1}(\hat{Z}_f) = (-1)^{n-1} h^{n-2,1}(\hat{Z}_f). \]
\medskip

First, we define a convenient stratification of $\hat{Z}_f$.

Let $\varphi_{\cal T} \,: \, \hat{Z}_f \rightarrow 
\overline{Z}_f$ be the corresponding birational morphism. Then  
$\hat{Z}_f$ can be represented as a disjoint union 
\[ \hat{Z}_f = \bigcup_{\Theta \subset \Delta} 
\varphi^{-1}_{\cal T}(Z_{f, \Theta}). \]
On the other hand, all irreducible components of fibers of $\varphi_{\cal T}$ 
over closed points of  $Z_{f, \Theta}$ are toric varieties. 
Therefore, we can define 
a stratification of $\varphi^{-1}_{\cal T}(Z_{f, \Theta})$ by smooth 
affine  algebraic varieties which are isomorphic to 
products ${Z}_{f,\Theta} \times ({\bf C}^*)^k$ some nonnegative integer $k$.   
As a result, we obtain a stratification of  $\hat{Z}_f$  
by  smooth affine varieties which are isomorphic to 
products ${Z}_{f,\Theta} \times ({\bf C}^*)^k$ for some face $\Theta \subset 
\Delta$ and for some nonnegative integer $k$.

Second, we note that $(n-2,1)$-Euler characteristic of 
${Z}_{f,\Theta} \times ({\bf C}^*)^k$ might be nonzero only in two cases: 
$\Theta = \Delta$, or ${\rm dim}\, \Theta = n -2$ and $k =1$. The latter 
follows from \ref{properties.hodge}, \ref{h.product}, and from the observation 
that the $(p,q)$-Euler characteristic of an 
algebraic torus $({\bf C}^*)^k$ is nonzero only if $p =q$ .  

We already know from \ref{n.def} that 
\[ e_c^{n-2,1}(Z_f) = (-1)^{n-1}(l(\Delta)  - n - 1 -  
\sum_{{\rm codim}\, \Theta =1 } l^*(\Theta)). \]
On the other hand, the strata which are isomorphic to  
${Z}_{f,\Theta} \times {\bf C}^*$ appear from the fibers of $\varphi_{\cal T}$ 
over $(n-3)$-dimensional singular affine locally closed subvarieites 
${Z}_{f,\Theta} \subset \overline{Z}_f$ having codimension $2$ in 
$\overline{Z}_f$. 
By \ref{duval}, a $\varphi^{-1}(Z_{f, \Theta})$ consists of $l^*(\Theta^*) = 
d(\Theta^*) -1$ irreducible components  isomorphic to ${\bf P}_1 \times 
{Z}_{f,\Theta}$. As a result, for every $(n-2)$-dimensional face $\Theta 
\subset \Delta$, one obtains  $l^*(\Theta^*)$  strata isomorphic to 
${Z}_{f,\Theta} \times {\bf C}^*$. 
On the other hand, one has 

\[ e_c^{n-2,1}({Z}_{f,\Theta} \times {\bf C}^*) = 
e_c^{n-3,0}({Z}_{f,\Theta}) \cdot e_c^{1,1}({\bf C}^*). \]

It is clear that $e_c^{1,1}({\bf C}^*) =1$. 
By results of Danilov and Khovanski\^i (see \cite{dan.hov}), one has 
\[ e_c^{n-3,0}({Z}_{f,\Theta}) = (-1)^{n-3} l^*(\Theta). \]
Thus we come to the following result.  

\begin{theo}
For $n \geq 4$, the Hodge number $h^{n-2,1}({\hat Z}_f)$ equals 
\[ l(\Delta)  - n - 1 -  
\sum_{{\rm codim}\, \Theta =1 } 
l^*(\Theta) + \sum_{{\rm codim}\, \Theta =2 } 
l^*(\Theta) \cdot l^*(\Theta^*), \]
where $\Theta$ denotes a face of a reflexive $n$-dimensional 
polyhedron $\Delta$, and $\Theta^*$ denotes the corresponding dual 
face of the dual reflexive polyhedron $\Delta^*$. 
\label{n.def.com}
\end{theo}

\begin{coro}
Assume that ${\bf P}_{\Delta}$, or ${\bf P}_{\Delta^*}$ is a smooth  
toric Fano variety of dimension $n \geq 4$ (see \ref{smooth.fano}(ii)). 
Then the Hodge-Deligne number $h^{n-2,1}$ of an affine 
$\Delta$-regular hypersurface $Z_f$ coincides with the Hodge number 
of a MPCP-desingularization $\hat{Z}_f$ of its toroidal 
compactification $\overline{Z}_f$. 
\end{coro}

\proof If ${\bf P}_{\Delta}$, or ${\bf P}_{\Delta^*}$ is a smooth 
toric Fano variety, then for any face $\Theta$ of codimension $2$, 
one has $l^*(\Theta)=0$, or $l^*(\Theta^*)=0$. 
\bigskip

\subsection{The Hodge number $h^{1,1}(\hat{Z}_f)$}

\hspace*{\parindent}

First we note that the group of principal ${\bf T}$-invariant 
divisors is isomorphic to the lattice $M$. 
Applying  \ref{crep.fano}, we obtain:

\begin{prop}
Any  MPCP-desingularization $\hat{\bf P}_{\Delta}$ of a toric Fano variety 
${\bf P}_{\Delta}$ contains exactly 
\[ l(\Delta^*) - 1 = {\rm card}\, \{ N \cap \partial \Delta^* \} \]
${\bf T}$-invariant divisors, i.e., the Picard number 
$\rho(\hat{\bf P}_{\Delta}) =  h^{1,1}(\hat{\bf P}_{\Delta})$ 
equals  
\[   l(\Delta^*) - n - 1. \] 
\label{bound.points} 
\end{prop}

\begin{theo}
Let $\hat{Z}_f$ be a MPCP-desingularization 
of a projective $\Delta$-regular Calabi-Yau hypersurface 
$\overline{Z}_f$, then for $n \geq 4$ the Hodge number 
$h^{1,1}(\hat{Z}_f)$, or the Picard number of $\hat{Z}_f$, equals 
\begin{equation}
l(\Delta^*) - n - 1 -  \sum_{{\rm codim}\, \Theta^* =1 } l^*(\Theta^*) + 
 \sum_{{\rm codim}\, \Theta =2 } l^*(\Theta^*) \cdot l^*(\Theta), 
\label{for.pic}
\end{equation}
where $\Theta^*$ denotes a face of the dual to $\Delta$ 
reflexive $n$-dimensional 
polyhedron $\Delta^*$, and $\Theta$ denotes the corresponding dual to $\Theta$ 
face of the  reflexive polyhedron $\Delta$. 
\label{n.pic}
\end{theo}

\proof  
By \ref{crep.fano},  
a $MPCP$-desingularization $\hat{Z}_f$ of a  $\Delta$-regular 
Calabi-Yau hypersurface $\overline{Z}_f$ is  induced by 
a $MPCP$-desingularization $\varphi_{\cal T}\,: \, 
\hat{\bf P}_{\Delta} \rightarrow {\bf P}_{\Delta}$ of the ambient 
toric Fano variety ${\bf P}_{\Delta}$. 
Since $\hat{Z}_f$ has only terminal 
${\bf Q}$-factorial singularities, i.e., $\hat{Z}_f$ is  
quasi-smooth \cite{dan1},  the Hodge structure in cohomology groups 
of $\hat{Z}_f$ is pure and satisfies the Poincare duality 
(see \ref{duality}).  
Therefore, the number $h^{1,1}$ equals the Hodge number $h^{n-2,n-2}$ in the 
cohomology group $H^{n-2}_c (\hat{Z}_f)$ with compact 
supports. By the Lefschetz-type theorem (see \ref{Lef}), if $n \geq 4$, then 
for $i = n-3, n-4$ the Gysin homomorphisms 
\[ H^i_c(Z_f) \rightarrow H^{i+2}_c ({\bf T}) \]
are isomorphisms  of Hodge structures with the shifting 
the Hodge type by $(1,1)$. On the other hand,  $H^{2n-1}_c ({\bf T})$ 
is an $n$-dimensional space having the Hodge 
type $(n-1,n-1)$, and the space $H^{2n-2}_c ({\bf T})$ has the Hodge 
type $(n-2,n-2)$. So  
$H^{2n-3}_c (Z_f)$ is an $n$-dimensional space having the Hodge type 
$(n-2,n-2)$, and the space $H^{2n-4}_c (Z_f)$ has the Hodge type $(n-3,n-3)$.
The complementary set $Y = \hat{Z}_f \setminus Z_f$ is a 
closed subvariety of $\hat{Z}_f$ of codimension $1$. Consider 
the corresponding exact sequence of Hodge structures
\[ \cdots \rightarrow H_c^{2n-4}(Z_f) \stackrel{\beta_1}{\rightarrow}  
H_c^{2n-4}(\hat{Z}_f) {\rightarrow}  
H_c^{2n-4}(Y) {\rightarrow}  \\ 
H_c^{2n-3}(Z_f) \stackrel{\beta_2}{\rightarrow}  
H_c^{2n-3}(\hat{Z}_f) {\rightarrow} \cdots \]
Comparing the Hodge types, we immediately get that $\beta_1$ and $\beta_2$ are 
zero mappings. 
Since the space $H_c^{2n-4}(Y)$ does not have subspaces of the Hodge type 
$(n-1, n-3)$, the Hodge humber $h^{n-2,n-2}$ of 
$H_c^{2n-4}(\hat{Z}_f)$ equals 
${\rm dim}\, H_c^{2n-4} (\hat{Z}_f)$.
Thus we get the short exact sequence of cohomology 
groups of the Hodge type $(n-2, n-2)$
\[ 0 \rightarrow  H_c^{2n-4}(\hat{Z}_f) {\rightarrow}  
H_c^{2n-4}(Y) {\rightarrow}  
H_c^{2n-3}(Z_f) {\rightarrow} 0.\]

It is easy to see that the dimension of $H_c^{2n-4}(Y)$ equals the number 
of the irreducible components of the $(n-2)$-dimensional complex subvariety 
$Y$. On the other hand, $Y$ is  
the intersection of $\hat{Z}_f \subset \hat{\bf P}_{\Delta}$ with the 
union of all irreducible 
${\bf T}$-invariant  divisors 
on the corresponding maximal partial 
crepant desingularization $\hat{\bf P}_{\Delta}$ 
of the toric Fano variety ${\bf P}_{\Delta}$.  
We have seen in 
\ref{bound.points} that these  toric 
divisors on $\hat{\bf P}_{\Delta}$ are in the one-to-one correspondence to the 
integral points $\rho \in N \cap \partial \Delta^*$,i.e., we have 
exactly $l(\Delta^*) -1$ irreducible 
toric divisors on $\hat{\bf P}_{\Delta}$. Every such a divisor $D_{\rho}$ 
is the closure of a $(n-1)$-dimensional torus ${\bf T}_{\rho}$ 
whose lattice of characters consists of elements of $M$ which are orthogonal 
to $\rho \in N$. It is important to note that $D_{\rho} \cap \hat{Z}_f$ is 
the closure in $\hat{\bf P}_{\Delta}$ of the affine hypersurface 
$\varphi_{\cal T}^{-1}(Z_{f,\Theta} \cap {\bf T}_{\rho} \subset 
{\bf T}_{\rho}$ where $\varphi_{\cal T}({\bf T}_{\rho}) = {\bf T}_{\Theta}$.

Note that since 
$\overline{Z}_f$ 
does not intersect any $0$-dimensional ${\bf T}$-orbit (\ref{zero.orb}), 
$\hat{Z}_f$ does not intersect  exeptional divisors on 
$\hat{\bf P}_{\Delta}$ lying over these points. The exceptional divisors 
lying over ${\bf T}$-invariant points of ${\bf P}_{\Delta}$ correspond 
to  integral points $\rho$ of  $N$ in the interiors of $(n-1)$-dimensional 
faces $\Theta^*$ of $\Delta^*$. So we must consider only 
\[ l(\Delta^*) -1 - \sum_{{\rm codim}\, \Theta^* =1 } 
l^*(\Theta^*) \] 
integral points of $N \cap \partial \Delta$ which are contained in faces 
of codimension $\geq 2$. 

If $D_{\rho}$ is an invariant  toric divisor on $\hat{\bf P}_{\Delta}$ 
corresponding to an integral point $\rho$ belonging to the interior of a 
$(n-2)$-dimensional face $\Theta^*$ of $\Delta^*$, then $D_{\rho} \cap 
\hat{Z}_f$ consists of $d(\Theta)$ irreducible components whose 
$\varphi_{\cal T}$-images are $d(\Theta)$ distict points of the 
zero-dimensional  stratum $Z_{f, \Theta} \subset {\bf T}_{\Theta}$. 

If $D_{\rho}$ is an invariant  toric divisor on $\hat{\bf P}_{\Delta}$ 
corresponding to an integral point belonging to the interior of a 
face $\Theta^* \subset \Delta^*$ of codimension $\geq 3$, then $D_{\rho} \cap 
\hat{Z}_f$ is  irreducible because ${\bf T}_{\rho} \cap \hat{Z}_f$ is 
an irreducible affine hypersurface in ${\bf T}_{\rho}$ isomorphic 
to $Z_{f, \Theta} \times ({\bf C}^*)^{n-1 - {\rm dim}\,\Theta}$.  

Consequently, the number of irreducible components of $Y$ is 
\[ \sum_{{\rm codim}\, \Theta^* =2 } d(\Theta)\cdot l^*(\Theta^*) + \] 
\[ + {\rm number\; of\; integral\; points\; on\; faces}\; 
\Theta^* \subset \Delta^*, \; {\rm codim}\, \Theta^* \geq 3. \]
  
Since $d(\Theta) = l^*(\Theta) +1$ for any  $1$-dimensional face 
$\Theta \Delta$, we can rewrite the number of the irreducible 
components of $Y$  as follows 
\[ {\rm dim}\, H_c^{2n-4}(Y) = l(\Delta^*) - 1 - 
\sum_{{\rm codim}\, \Theta^* =1 } 
l^*(\Theta^*) +  \sum_{{\rm codim}\, \Theta =2 } 
l^*(\Theta^*) \cdot l^*(\Theta). \] 
Since  ${\rm dim}\, H_c^{2n-3}(Z_f) = n$, we obtain $(\ref{for.pic})$. 
\bigskip

Applying  \ref{n.def.com} and \ref{n.pic}, we conclude. 

\begin{theo}
For any reflexive polyhedron $\Delta$ of dimension $n \geq 4$, 
the Hodge number $h^{n-1,1}(\hat{Z}_f)$ of a MPCP-desingularization 
of a $\Delta$-regular Calabi-Yau hypersurface $\overline{Z}_f \subset 
{\bf P}_{\Delta}$ equals the 
Picard  number $h^{1,1}(\hat{Z}_g)$  of  a MPCP-desingularization  
of a $\Delta^*$-regular projective Calabi-Yau 
hypersurface $\overline{Z}_g \subset {\bf P}_{\Delta^*}$ corresponding to 
the dual reflexive polyhedron $\Delta^*$.  
\label{num.mir}
\end{theo}

\subsection{Calabi-Yau $3$-folds}

\hspace*{\parindent}

We begin with the remark that  only if ${\rm dim}\, \Delta =4$ 
both statements \ref{smooth.c1} and  \ref{num.mir} hold. 
In this case,  we deal with  $3$-dimensional Calabi-Yau hypersurfaces 
$\overline{Z}_f \subset {\bf P}_{\Delta}$ which admit   {\em smooth} 
$MPCP$-desingularizations $\hat{Z}_f$. 

 Calabi-Yau $3$-folds $\hat{Z}_f$  are of primary interest in 
theoretical  physics. The number  
\[ \frac{1}{2} e(\hat{Z}_f) = 
( h^{1,1}(\hat{Z}_f) - h^{2,1}(\hat{Z}_f)) \]
is called {\em the number 
of generations} in superstring theory \cite{witt}. 
So it is important to have 
a simple formula for the Euler characteristic $e(\hat{Z}_f)$.

We have already calculated the Hodge numbers  $h^{1,1}(\hat{Z}_f)$ and  
$h^{n-2,1}({Z}_f)$. So we obtain

\begin{coro}
For any Calabi-Yau $3$-fold $\hat{Z}_f$ defined by a $\Delta$-regular 
Laurent polynomial $f$ whose Newton polyhedron is a reflexive 
$4$-dimensional polyhedron $\Delta$, one has the following formulas 
for the Hodge numbers 
\[ h^{1,1}(\hat{Z}_f) =  
l(\Delta^*) - 5 -  \sum_{{\rm codim}\, \Theta^* =1 } l^*(\Theta^*) + 
 \sum_{{\rm codim}\, \Theta =2 } l^*(\Theta^*) \cdot l^*(\Theta), \]
\[ _+ h^{2,1} (\hat{Z}_f) = l(\Delta)  - 5 -  \sum_{{\rm codim}\, \Theta =1 } 
l^*(\Theta) + \sum_{{\rm codim}\, \Theta =2 } l^*(\Theta) \cdot 
l^*(\Theta^*). \]
\end{coro}
  
This implies also. 

\begin{coro}
\[ e(\hat{Z}_f) = (l(\Delta)  - l(\Delta^*)) -  
(\sum_{\begin{array}{c}{\rm codim}\, \Theta =1 \\ 
\Theta \subset \Delta \end{array}} l^*(\Theta)  - 
\sum_{\begin{array}{c}{\rm codim}\, \Xi =1 \\ 
\Xi \subset \Delta^* \end{array}} l^*(\Xi)\;\; ) + \]
\[ + (\sum_{\begin{array}{c}{\rm codim}\, \Theta =2 \\ 
\Theta \subset \Delta \end{array}} l^*(\Theta) \cdot 
l^*(\Theta^*) - 
\sum_{\begin{array}{c}{\rm codim}\, \Xi =2 \\ 
\Xi \subset \Delta^* \end{array}} l^*(\Xi) \cdot l^*(\Xi^*)\; ).  \]
\end{coro}
\medskip

Now we prove a more simple another 
formula for the Euler characteristic of 
Calabi-Yau $3$-folds.

\begin{theo}
Let $\hat{Z}_f$ be a MPCP-desingularization of a $3$-dimensional 
$\Delta$-regular Calabi-Yau hypersurface associated with a $4$-dimensional 
reflexive polyhedron $\Delta$. Then 
\[ e(\hat{Z}_f) =   
\sum_{\begin{array}{c}{\rm dim}\, \Theta =1 \\ 
\Theta \subset \Delta \end{array}} d(\Theta)d(\Theta^*)  - 
\sum_{\begin{array}{c}{\rm dim}\, \Theta =2 \\ 
\Theta \subset \Delta \end{array}} d(\Theta)d(\Theta^*). \]
\label{new.euler}
\end{theo}

In our proof of theorem \ref{new.euler} we will use 
one  general property of smooth quasi-projective open 
subsets ${Z}_f^{[1]}$ in   $\overline{Z}_f$ consisting of the 
union of the affine part ${Z}_f$ and all affine strata $Z_{f,\Theta}$, where 
$\Theta$ runs over all faces of $\Delta$ of codimension 1, i.e., 
\[{Z}^{[1]}_f = {\bf P}_{\Delta}^{[1]} \cap {\overline{Z}}_f . \]

\begin{theo}
For arbitrary $n$-dimensional reflexive polyhedron $\Delta$ and 
$\Delta$-regular Laurent polynomial $f \in L(\Delta)$, the 
Euler characteristic   of the smooth quasi-projective   Calabi-Yau 
variety ${Z}_f^{[1]}$ 
is always zero.
\label{euler.zero}
\end{theo}

\proof  Since  ${Z}_f^{[1]}$ is smooth,  the Euler characteristic   
of the usual cohomology groups $H^* ({Z}_f^{[1]})$  is zero if and only if 
the Euler characteristic  of the cohomology groups with compact 
supports $H_c^* ({Z}_f^{[1]})$ is zero (see \ref{duality}). 
It follows from the long exact sequence 
of cohomology groups with compact suport that 
\[ e(H_c^* ({Z}_f^{[1]})) = e(H_c^* ({Z}_f)) + 
\sum_{{\rm codim \, \Theta} = 1} 
e(H_c^* ({Z}_{f,\Theta})). \]
By \ref{euler}, 
\[ e(H_c^* ({ Z}_f)) = (-1)^{n-1} d(\Delta), \;\;
e(H_c^* ({ Z}_{f,\Theta})) = (-1)^{n-2} d(\Theta). \]
Thus,  it is sufficient to prove
\[ d(\Delta) = \sum_{{\rm codim}\,  \Theta =1} d(\Theta). \]
The latter follows immediately from the representation of the 
$n$-dimensional polyhedron $\Delta$ as a union of $n$-dimensional 
pyramids with vertex $0$ over all 
$(n-1)$-dimensional faces $\Theta \subset \Delta$. 
\bigskip

{\bf Proof of Theorem \ref{new.euler}}. 
Let $\varphi_{\cal T}\;:\; \hat{Z}_f \rightarrow \overline{Z}_f$ be a 
$MPCP$-desingularization. For any face $\Theta \subset \Delta$, let 
$F_{\Theta}$ denotes $\varphi^{-1}_{\cal T}(Z_{f,\Theta})$. 
We know that $\varphi$ is an isomorphism over 
$Z^{[1]}_f$ which is the union of 
the strata $Z_{f, \Theta}$ (${\rm dim}\, \Theta = 3,4$). 
By  \ref{euler.zero}, $e(Z^{[1]}_f) =0$. Using additivity property of 
the Euler characteristic, we obtain
\[ e(\hat{Z}_f) = 
\sum_{\begin{array}{c}{\rm dim}\, \Theta =1 \\ 
\Theta \subset \Delta \end{array}} e(F_{\Theta}) +  
\sum_{\begin{array}{c}{\rm dim}\, \Theta =2 \\ 
\Theta \subset \Delta \end{array}} e(F_{\Theta}). \]
Now the statement follows from the following lemma.

\begin{lem}
Let $\varphi_{\cal T}$ be a MPCP-desingularization as above. Then the 
Euler characteristic $e(F_{\Theta})$ 
equals $d(\Theta)d(\Theta^*)$ if ${\rm dim}\, \Theta =1$, and 
$-d(\Theta)d(\Theta^*)$ if ${\rm dim}\, \Theta =2$.
\end{lem}

\proof  If ${\rm dim}\, \Theta =1$, then $Z_{f,\Theta}$ consists of 
$d(\Theta)$ distinct points. The fiber $\varphi^{-1}_{\cal T}(p)$ 
over any such a 
point $p \in Z_{f,\Theta}$ is a union of smooth toric surfaces. 
By \ref{topol.des}, 
the Euler characteristic of $\varphi^{-1}_{\cal T}(p)$ equals the number of 
elementary simplices in the corresponding maximal projective 
triangulation of $\Theta^*$. Since ${\rm dim}\, \Theta^* =2$, the   
number of elementary simplices equals $d(\Theta^*)$ (see \ref{el.reg}). Thus, 
$e(F_{\Theta}) = e(Z_{f,\Theta})e(\varphi^{-1}_{\cal T}(p)) = 
d(\Theta)d(\Theta^*)$. 

If ${\rm dim}\, \Theta =2$, then $Z_{f,\Theta}$ is a smooth affine algebraic 
curve. By \ref{euler}, the Euler characteristic of $Z_{f,\Theta}$ equals 
$-d(\Theta)$. By \ref{duval}, the fiber $\varphi^{-1}_{\cal T}(p)$ over  
any point $p \in 
Z_{f,\Theta}$ is the Hirzebruch-Jung tree  of $d(\Theta^*)-1$ smooth rational 
curves, i.e., the Euler characteristic of $\varphi^{-1}_{\cal T}(p)$ again 
equals $d(\Theta^*)$. Thus, 
$e(F_{\Theta}) = -d(Z_{f,\Theta})e(\varphi^{-1}_{\cal T}(p)) = 
-d(\Theta)d(\Theta^*)$.
\medskip

Since the duality between $\Delta$ and $\Delta^*$ establishes 
a one-to-one correspondence between $1$-dimensional (respectively 
$2$-dimensional) faces of $\Delta$ and $2$-dimensional (respectively 
$1$-dimensional) faces of $\Delta^*$, we again obtain: 

\begin{coro}
Let $\Delta$ be a reflexive polyhedron of dimension $4$, $\Delta^*$ the 
dual reflexive polyhedron. Let $\hat{Z}_f$ be a MPCP-desingularization of 
a $\Delta$-regular Calabi-Yau hypersurface $\overline{Z}_f$ 
in ${\bf P}_{\Delta}$, $\hat{Z}_g$ be a MPCP-desingularization of 
a $\Delta^*$-regular Calabi-Yau hypersurface $\overline{Z}_g$ 
in ${\bf P}_{\Delta^*}$. Then the Euler characteristics of two Calabi-Yau 
$3$-folds $\hat{Z}_f$ and $\hat{Z}_g$ satisfy the following relation 
\[ e(\hat{Z}_f) = - e(\hat{Z}_g). \]
\label{coro.euler}
\end{coro}
\medskip

We consider one of new  examples of the costruction of candidates 
for mirrors via dual polyhedra.  

\begin{exam}
{\rm Let $\Delta$ be the reflexive polyhedron in ${\bf Q}^4$ defined 
as the convex hull of $6$ points 
\[ A_1 = (1,0,0,0),\; A_2 = (0,1,0,0), A_3 = (-1,-1,0,0), \]
\[ A_4 = (0,0,1,0),\; A_5 = (0,0,0,1), A_6 = (0,0,-1,-1). \]
It follows from standard toric description of projective spaces that the 
toric variety ${\bf P}_{\Delta^*}$ corresponding to the dual reflexive 
polyhedron $\Delta^*$ is isomorphic to ${\bf P}_2 \times {\bf P}_2$. 
On the other hand, ${\bf P}_{\Delta}$ is a singular complete intersection 
of two cubics in ${\bf P}^6$:  
\[ {\bf P}_{\Delta} = \{ (Y_0 : ... : Y_6) \in {\bf P}_6 \mid 
Y_0^3 = Y_1 Y_2 Y_3, \; Y_0^3 = Y_4 Y_5 Y_6 \}. \]
The intersection of ${\bf P}_{\Delta} \subset {\bf P}^6$ with generic 
linear subspace of codimension $1$ is a singular $3$-dimensional 
Calabi-Yau variety $\overline{Z}_f$ defines as 
\[ \overline{Z}_f = \{ (Y_0 : ... : Y_6) \in {\bf P}_{\Delta} \mid 
a_0 Y_0 + \cdots + a_6 Y_6 = 0 \}, \]
where the coefficients $a_0, \ldots, a_6$ are chosen sufficiently general and 
\[ f(X) = a_0 + a_1 X_1 + a_2 X_2 + a_3 (X_1X_2)^{-1} + 
a_4 X_3 + a_5 X_4 + a_6 (X_3X_4)^{-1}. \]
The projective hypersurface $\overline{Z}_f$ has $18$ double-point 
singularities of type $A_2$ along rational curves 
\[ C_{ijk} = \{ (Y_0 : ... : Y_6) \in \overline{Z}_f \mid 
Y_i = Y_j = Y_k = 0 \}, \]
where $\{i,j,k \}$ runs over all subsets of $\{ 1, \ldots, 6 \}$ not equal to 
$\{ 1,2,3 \}$ or $\{ 4,5,6 \}$. 

At points  of intersections of these rational curves we obtain $15$ more 
complicated isolated Gorenstein toroidal singularities: 

$6$ quotient singularities 
whose $MPCP$-resolution can be described by the triangulated diagram 

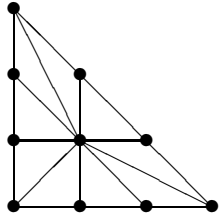
\begin{center}
\begin{picture}(200,200)
\put(100,100){\makebox(0,0){$\bullet$}}
\put(100,125){\makebox(0,0){$\bullet$}}
\put(75,150){\makebox(0,0){$\bullet$}}
\put(150,75){\makebox(0,0){$\bullet$}}
\put(125,100){\makebox(0,0){$\bullet$}}
\put(125,75){\makebox(0,0){$\bullet$}}
\put(75,125){\makebox(0,0){$\bullet$}}
\put(75,75){\makebox(0,0){$\bullet$}}
\put(100,75){\makebox(0,0){$\bullet$}}
\put(75,100){\makebox(0,0){$\bullet$}}

\put(75,75){\line(0,1){75}}
\put(75,75){\line(1,0){75}}
\put(75,75){\line(1,1){25}}

\put(75,100){\line(1,0){50}}
\put(100,75){\line(0,1){50}}
\put(75,125){\line(1,-1){50}}
\put(75,150){\line(1,-2){25}}
\put(100,100){\line(2,-1){50}}

\put(150,75){\line(-1,1){75}}
\end{picture}
\end{center}

and $9$ Gorenstein non-quotient toroidal singularities  
whose $MPCP$-resolution can be described by the triangulated diagram 

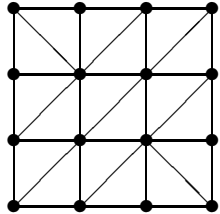
\begin{center}
\begin{picture}(200,200)
\put(100,100){\makebox(0,0){$\bullet$}}
\put(100,125){\makebox(0,0){$\bullet$}}
\put(125,125){\makebox(0,0){$\bullet$}}
\put(125,100){\makebox(0,0){$\bullet$}}
\put(125,75){\makebox(0,0){$\bullet$}}
\put(75,125){\makebox(0,0){$\bullet$}}
\put(75,75){\makebox(0,0){$\bullet$}}
\put(100,75){\makebox(0,0){$\bullet$}}
\put(75,100){\makebox(0,0){$\bullet$}}
\put(150,75){\makebox(0,0){$\bullet$}}
\put(150,100){\makebox(0,0){$\bullet$}}
\put(150,125){\makebox(0,0){$\bullet$}}
\put(150,150){\makebox(0,0){$\bullet$}}
\put(125,150){\makebox(0,0){$\bullet$}}
\put(100,150){\makebox(0,0){$\bullet$}}
\put(75,150){\makebox(0,0){$\bullet$}}

\put(75,75){\line(1,1){75}}
\put(75,100){\line(1,1){50}}
\put(100,75){\line(1,1){50}}
\put(75,75){\line(1,0){75}}
\put(75,100){\line(1,0){75}}
\put(75,125){\line(1,0){75}}
\put(75,150){\line(1,0){75}}
\put(75,75){\line(0,1){75}}
\put(100,75){\line(0,1){75}}
\put(125,75){\line(0,1){75}}
\put(150,75){\line(0,1){75}}

\put(75,150){\line(1,-1){25}}
\put(125,100){\line(1,-1){25}}
\end{picture}
\end{center}

Applying our formulas for the Hodge numbers of a $MPCP$-resolution $\hat{Z}_f$ 
of $\overline{Z}_f$, we get 
\[ h^{1,1}(\hat{Z}_f) = 83, \; h^{2,1}(\hat{Z}_f) = 2. \]
Thus, we have constructed a family ${\cal F}(\Delta^*)$ which consists of 
candidates for mirrors of 
smooth Calabi-Yau hypersurfaces of bidegree $(3,3)$ in 
${\bf P}_2 \times {\bf P}_2$ having the Hodge numbers $h^{1,1} =2$ and 
$h^{2,1} = 83$. 
  
  One can calculate the Gauss-Manin connection defined by periods of 
$\hat{Z}_f$ by general method in $\cite{bat.var}$. 
}
\end{exam}

\section{Mirror symmetry}

\subsection{Mirror candidates for hypersurfaces of degree $n+1$ in 
${\bf P}_n$}
 
\hspace*{\parindent}
Consider the  polyhedron $\Delta_n$ in $M_{\bf Q} \cong {\bf Q}^n$ defined 
by inequalities 
\[ x_1 + \ldots +  x_n \leq  1, \; x_i \geq -1 \; (1 \leq i \leq n). \]
Then $\Delta_n$ is a reflexive polyhedron, and 
${\cal F}(\Delta)$ is the family of  
all hypersurfaces of degree $n+1$ in ${\bf P}_n = {\bf P}_{\Delta_n}$. 

The polyhedron $\Delta_n$ has $n+1$  $(n-1)$-dimensional faces whose 
interiors contain exactly $n$ integral points. These $n(n+1)$ 
integral points form the root system of type $A_n$ (see \ref{root}). 
The number $l(\Delta_n)$ equals ${ 2n + 1 \choose n }$. Thus the dimension 
of the moduli space of ${\cal F}({\Delta}_n)$ equals  
\[ { 2n +1 \choose n } - (n+1)^2.\]

The dual polyhedron $\Delta^*_n$  has $n +1$ vertices
\[ u_1 =(1,0, \ldots, 0), \ldots , u_n = (0, \ldots, 0, 1), 
u_{n+1} = (-1, \ldots , -1). \]
The corresponding toric Fano variety ${\bf P}_{\Delta_n^*}$ is a singular 
toric hypersurface ${\bf H}_{n+1}$ of degree $n+1$ in ${\bf P}_{n+1}$ defined by the equation 
\[  \prod_{i =1}^{n+1} u_i = u_0^{n+1},  \]
where $(u_0: \ldots : u_{n+1})$ are homogeneous coordinates in 
${\bf P}_{n+1}$. 

Since the simplex $\Delta_n$ is $(n+1)$-times multiple of $n$-dimensional 
elementary 
simplex of degree 1, the degree $d(\Delta_n)$ equals $(n+1)^n$. On the other hand, 
the dual simplex $\Delta_n^*$ is the union of $n+1$ elementary simplices of 
degree 1, i.e., $d(\Delta_n^*) = n+1$. 

There exists a finite morphism of degree $(n+1)^{n-1} = 
d(\Delta_n)/ d(\Delta_n^*)$  of reflexive pairs 
\[ \phi \; :\; (\Delta^*_n, N) \rightarrow (\Delta_n, M), \]
where  $\phi(u_{n+1}) = (-1, \ldots , -1) \in \Delta_n$ and  

\[ \phi(u_i) = (-1, \ldots, \underbrace{n}_{i}, \ldots, -1 ) \in \Delta_n. \]

It is easy to see now that  
\[  M / \phi(N) \cong ({\bf Z}/(n+1){\bf Z})^{n-1}. \]

Let $(v_0 : v_1:  \ldots : v_n)$ be the homogeneous coordinates on 
${\bf P}_n$. The corresponding to $\phi$ \^etale  mapping of smooth 
quasi-projective toric Fano  varieties 
\[ \tilde {\phi} \;:\; {\bf P}_n^{[1]} \rightarrow 
{\bf H}_{n+1}^{[1]} \]
has the following  representation in homogeneous coordinates 
\[ (v_0 : v_1 : \dots : v_n ) \mapsto (\prod_{i= 0}^n v_i : 
v_0^{n+1} : v_1^{n+1} : \dots : v_n^{n+1} ) = 
(u_0 : u_1:  \ldots : u_{n+1} ).\]  

A Calabi-Yau hypersurface $\overline{Z}_f$  in ${\bf H}_{n+1}$ 
has an  equation 
\[ f(u)  = \sum_{i=0}^{n+1} a_i u_i = 0. \] 

Using \ref{n.def}, it is easy to show that $h^{n-2,1}(Z_f) = 1$. We can also 
describe the moduli space of the family ${\cal F}(\Delta_n^*)$. 
Since  ${\bf H}_{n+1}$  is invariant under the action 
of the $n$-dimensional torus 
\[ {\bf T} = \{ {\bf t}=(t_1, \ldots, t_{n+1}) \in ({\bf C}^*)^{n+1} \mid  
t_1  \cdots t_{n+1} = 1 \}, \] 
the equation 
\[ {\bf t}^*(f(u))  =  \sum_{i=0}^{n+1} a_i t_i u_i = 0, \] 
defines an isomorphic to $\overline{Z}_f$ hypersurface 
$\overline{Z}_{t^*(f)}$. Moreover, multiplying  
all coefficients $\{a_i\}\; (0 \leq i \leq n+1)$ 
by the same non-zero complex number $t_0$, we get also a ${\bf C}^*$-action.  
Thus, up to the action of the $(n+1)$-dimensional torus 
${\bf C}^* \times {\bf T}$ on $n+2$ coefficients $\{ a_i\}$, we get 
the one-parameter mirror family of Calabi-Yau hypersurfaces in 
${\bf H}_{n+1}$ defined by the equation 
\begin{equation}
 f_{a}(u)  = \sum_{i =1}^d u_i  + a u_0 = 0, 
 \label{one.mir}
\end{equation} 
where the number 
\[ a  =  a_0 (\prod_{j =1}^{n+1} a_j)^{\frac{-1}{n+1}} \]
is uniquely defined up to an $(n +1)$-th root of unity. 
  
Using  the homogeneous coordinates $\{ v_i \}$ on ${\bf P}_n$, we can transform  
this equation to the form 
\begin{equation}
 {\tilde {\phi}}^*f_{a} = \sum_{i =0}^n v_i^{n+1}   + 
a \prod_{i=0}^n v_i,  
\label{equ.mir}
\end{equation}
where $a^{n+1}$ is a canonical parameter of the corresponding subfamily 
in ${\cal F}(\Delta_n)$ of smooth $(n-1)$-dimensional hypersurfaces 
in ${\bf P}_n$. 
\bigskip 

\begin{theo}
The candidate for mirrors relative to the 
family of all smooth 
Calabi-Yau hypersurfaces of degree $n+1$ in ${\bf P}_n$ is the one-parameter 
family ${\cal F}(\Delta^*_n)$ of Calabi-Yau varieties consisting  of quotients   
by the action of 
the finite abelian group $({\bf Z}/ (n+1){\bf Z})^{n-1}$  of the 
hypersurfaces defined by the equation $(\ref{equ.mir})$. 
\label{first.ex}
\end{theo}

\begin{exam}
{\rm If we take  $n =4$, then the corresponding finite abelian group  
is isomorphic to  $({\bf Z}/5{\bf Z})^3$ and we come to the mirror 
family for the family of all 
$3$-dimensional quintics in ${\bf P}_4$ considered in \cite{cand2}.}
\end{exam}  
\bigskip

\subsection{A category of reflexive pairs}

\hspace*{\parindent}

The set of all reflexive pairs of dimension $n$ forms a  category ${\cal C}_n$ 
whose  morphisms are finite morphisms of reflexive pairs (\ref{morph}). 
The correspondence between dual reflexive pairs 
defines an involutive functor  
\[ {\rm Mir}\; : \; {\cal C}_n \rightarrow {\cal C}_n^* \]
\[ {\rm Mir}(\Delta, M) = (\Delta^*, N) \]
which is an isomorphism of the category ${\cal C}_n$ with the dual 
category ${\cal C}_n^*$. 
\bigskip

It is natural to describe in ${\cal C}_n$  some morphisms  satisfying 
universal properties.

\begin{opr}
{\rm Let 
\[ \phi_0 \; :\; (\Delta_0 , M_0) \rightarrow (\Delta, M) \]
be a finite morphism of reflexive pair. The morphism $\phi_0$ is said to 
be {\em minimal} if for any finite morphism 
\[ \psi \; :\; (\Delta' , M') \rightarrow (\Delta, M) \]
there exists the  unique morphism 
\[ \phi \; :\; (\Delta_0 , M_0) \rightarrow (\Delta', M') \]
such that $\phi_0 = \psi \circ \phi$.

A reflexive pair in ${\cal C}_n$ is called {\em a minimal reflexive pair}, 
if the identity morphism of this pair  is minimal.}
\end{opr}
\bigskip 
 
Note the following simple property. 
 
\begin{prop}
Assume that there exist two finite morphisms 
\[ \phi_1  \; :\; (\Delta_1, M_1) \rightarrow (\Delta_2, M_2) \]
and 
\[ \phi_2  \; :\; (\Delta_2, M_2) \rightarrow (\Delta_1, M_1). \]
Then $\phi_1$ and $\phi_2$ are isomorphisms.
\label{isom}
\end{prop}

\proof The degrees  of $\phi_1$ and $\phi_2$ are positive 
integers. On the other hand, by \ref{etale}, 
\[ d(\phi_1) = d_{M_2}(\Delta_2)/ d_{M_1}(\Delta_1), \]  
\[ d(\phi_2) = d_{M_1}(\Delta_1)/ d_{M_2}(\Delta_2). \]  
Therefore, $d_{M_1}(\Delta_1) = d_{M_2}(\Delta_2) =1$, i.e., $\phi_1$ and 
$\phi_2$ are isomorphisms of reflexive pairs.
\bigskip

\begin{coro}
Assume that 
\[ \phi_0 \; :\; (\Delta_0 , M_0) \rightarrow (\Delta, M) \]
is a minimal morphism. Then $(\Delta_0, M_0)$ is a minimal reflexive pair. 
\end{coro}

\begin{coro}
For any reflexive pair $(\Delta, M)$ there exists up to an isomorphism 
at most one minimal pair $(\Delta_0, M_0)$ with a minimal morphism
\[ \phi_0 \; : \; (\Delta_0 , M_0) \rightarrow (\Delta, M). \]
\end{coro}
\bigskip

The next proposition completely describes minimal reflexive pairs and the set 
of all finite morphisms to a fixed reflexive pair $(\Delta, M)$. 
\medskip

\begin{prop}
Let $(\Delta, M)$ be a reflexive pair. Denote by $M_{\Delta}$ the sublattice 
in $M$ generated by vertices of $\Delta$. Let $M'$ be an integral lattice 
satisfying the condition 
$M_{\Delta} \subset M' \subset M$.  Then 
$(\Delta, M')$ is also  a reflexive pair, and 
\[ (\Delta, M_{\Delta}) \rightarrow (\Delta, M' ) \]
is a minimal morphism.
\label{min.mor}
\end{prop}

The proof immediately  follows from the definition of reflexive pair 
\ref{inver.p}. 
\medskip

\begin{coro}
All reflexive pairs having a finite morphism to a fixed reflexive 
pair $(\Delta, M)$ are isomorphic to $(\Delta, M')$ for some  lattice 
$M'$ such that $M_{\Delta} \subset M' \subset M$.
\label{min.pair}
\end{coro}
\bigskip

\begin{opr}
{\rm Let 
\[ \phi^0 \; :\; (\Delta , M) \rightarrow (\Delta^0, M^0) \]
be a finite morphism of reflexive pair. The morphism $\phi^0$ is said to 
be {\em maximal} if for any finite morphism 
\[ \psi \; :\; (\Delta , M) \rightarrow (\Delta', M') \]
there exists the  unique morphism 
\[ \phi \; :\; (\Delta' , M') \rightarrow (\Delta^0, M^0) \]
such that $\phi^0 = \phi \circ \psi$.

A reflexive pair in ${\cal C}_n$ is called {\em a maximal reflexive pair}, 
if the identity morphism of this pair  is maximal.}
\end{opr}
\bigskip 

If we apply the functor ${\rm Mir}$, we get from \ref{min.mor} and 
\ref{min.pair} the following properties 
of maximal reflexive pairs.

\begin{prop}
Let $(\Delta, M)$ be a reflexive pair. Then  there exists up to an isomorphism 
the unique maximal reflexive pair $(\Delta, M^{\Delta})$ having  a maximal 
morphism
\[ \phi \; :\; (\Delta , M) \rightarrow (\Delta, M^{\Delta}). \]
Moreover, the pair  $(\Delta, M^{\Delta}) $ is dual to the minimal  pair 
$(\Delta^*, N_{\Delta^*})$ having the  morphism 
\[ \phi^*\; : \; (\Delta^* , N_{\Delta^*}) \rightarrow (\Delta^*, N) \]
as minimal.
\end{prop}
\bigskip

\begin{coro}
All reflexive pairs having finite morphisms  from a fixed reflexive 
pair $(\Delta, M)$ are isomorphic to $(\Delta, M')$ for some  lattice 
$M'$ such that $M \subset M' \subset M^{\Delta}$.
\end{coro}
\bigskip

\begin{exam}
{\rm Let $(\Delta_n, M)$ and $(\Delta^*_n, N)$ be two reflexive pairs  
from the previous section. Since the lattice $N$ is generated by 
vertices of $\Delta_n^*$, the reflexive pair  $(\Delta_n^*, N)$ is minimal. 
Therefore, $(\Delta_n, M)$ is a maximal reflexive pair. }
\end{exam}
\bigskip

\subsection{A Galois correspondence}

\hspace*{\parindent}

The existence of a finite morphism of reflexive pairs 
\[ \phi\; : \; (\Delta_1, M_1) \rightarrow (\Delta_2, M_2) \]
implies the following main geometric relation between 
Calabi-Yau hypersurfaces in toric 
Fano varieties ${\bf P}_{\Delta_1, M_1}$ and ${\bf P}_{\Delta_2, M_2}$. 
\medskip

\begin{theo}
The Calabi-Yau hypersurfaces in ${\bf P}_{\Delta_1, M_1}$ are quotients of 
some Calabi-Yau hypersurfaces in ${\bf P}_{\Delta_2, M_2}$ by the action of 
the dual to $M_2/ \phi(M_1)$ finite abelian group. 
\label{quot}
\end{theo}

\proof Consider the dual finite morphism 
\[ \phi^* \; : \; (\Delta_2^*, N_2) \rightarrow (\Delta_1^*, N_1). \]
By \ref{etale}, $\phi^*$ induces a finite  \^etale morphism 
\[ \tilde {\phi}^* \; : \;  {\bf P}_{\Delta_2,M_2}^{[1]}  
\rightarrow  {\bf P}_{\Delta_1,M_1}^{[1]}  \]
of smooth quasi-projective toric Fano varieties.  This morphism is defined 
by the surjective homomorphism of $n$-dimensional algebraic tori 
\[ {\gamma}_{\phi} \; :\; {\bf T}_{\Delta_2}  
\rightarrow {\bf T}_{\Delta_1} \] 
whose kernel is dual to the cokernel of the homomorphism of the groups of 
characters 
\[ \phi \; :\; M_1 \rightarrow  M_2.\]
Therefore, $\tilde {\phi}^*$ is the quotient by the action of the 
finite abelian group  $(M_2/ \phi(M_1))^* =  N_1 / \phi^*(N_2)$. 
The pullback of the anticanonical 
class of  ${\bf P}_{\Delta_1,M_1}^{[1]}$  
is the anticanonical class of 
${\bf P}_{\Delta_2,M_2}^{[1]}$. 
Applying \ref{galois}, we obtain that 
the smooth quasi-projective Calabi-Yau hypersurfaces 
${\hat Z}_{f,\Delta_1}$ are \^etale quotients by $N_1 / \phi^*(N_2)$  
of some smooth 
quasi-projective Calabi-Yau hypersurfaces 
${\hat Z}_{\tilde {\phi}^*f,\Delta_2}$. 
\bigskip 

\begin{coro}
The mirror mapping for families of Calabi-Yau hypersurfaces in toric 
varieties satisfies the following Galois correspondence: 

If a family ${\cal F}(\Delta_1)$ is a quotient of a family 
${\cal F}(\Delta_2 )$ by a finite abelian group ${\cal A}$, then 
the mirror family ${\cal F}(\Delta_2^*)$ is a quotient of the mirror 
family ${\cal F}(\Delta_1^*)$ by the dual finite abelian 
group ${\cal A}^*$. 
\end{coro}
\bigskip

\begin{opr}
{\rm 
Let $(\Delta, M)$ be a reflexive pair, 
$(\Delta^*,  N)$  the dual reflexive pair. 
Denote  by  $N_{\Delta^*}$ 
the sublattice  in $N$  generated by  vertices of  $\Delta^*$. 
The finite abelian  group $\pi_1 (\Delta, M) = N/ N_{\Delta^*}$ 
 is called  the {\em fundamental group of the pair } $(\Delta, M)$. }
\end{opr}
\bigskip

The fundamental group $\pi_1 (\Delta, M)$ defines  a contravariant 
functor from the category ${\cal C}_n$ to the category of finite 
abelian groups with injective homomorphisms. 
\bigskip

\begin{prop}
Let $(\Delta, M)$ be a reflexive pair, ${\bf P}_{\Delta, M}$ 
the corresponding toric Fano variety. Then the fundamental group 
$\pi_1 (\Delta, M)$ is 
isomorphic to the algebraic $($and topological$)$ fundamental group 
\[ \pi_1 ({\bf P}^{[1]}_{\Delta} ). \]
In particular, the reflexive pair $(\Delta, M)$ is maximal if and only if 
${\bf P}^{[1]}_{\Delta}$ is simply connected.  
\end{prop}

\proof The statement immediately 
follows from the description of the fundamental 
group of toric varieties in \ref{fund.group}.
\bigskip

\begin{opr}
{\rm Let $(\Delta, M)$ be a reflexive pair, $(\Delta^*, N)$ the dual 
reflexive pair, $(\Delta, M_{\Delta})$ and $(\Delta^*, N_{\Delta^*})$ are 
minimal pairs, $(\Delta, M^{\Delta})$ and $(\Delta, N^{\Delta^*})$ 
are maximal pairs. The quotients 
\[ \pi_1 (\Delta) = N^{\Delta^*} / N_{\Delta^*} \; {\rm and }\; 
\pi_1 (\Delta^*) = M^{\Delta} / M_{\Delta } \] 
 is called the {\em fundamental groups} of the reflexive polyhedra $\Delta$ 
and  $\Delta^*$ respectively. }
\end{opr}
\medskip

It is clear, $\pi_1 (\Delta)$ and $\pi_1 (\Delta^*)$ are isomorphic 
dual finite abelian  groups. 
\medskip

\begin{opr}
{\rm Assume that for a reflexive pair $(\Delta, M)$ 
there exists an isomorphism between two maximal 
reflexive pairs 
\[ \phi\; :\; (\Delta, M^{\Delta}) \rightarrow (\Delta^*, N^{\Delta^*}) .\]
Then we call  $\Delta$ {\em a selfdual reflexive polyhedron}. }
\label{selfdual}
\end{opr}
\medskip

If $\Delta$ is selfdual, then $\Delta$ and $\Delta^*$ must have the same 
combinatorial type (see \ref{dual.edge}). By \ref{quot}, we obtain.

\begin{prop}
Let $(\Delta, M)$ be a reflexive pair such that $\Delta$ is selfdual. 
Then ${\cal F}(\Delta)$ and ${\cal F}(\Delta^*)$ are  
quotients respectively by $\pi_1(\Delta, M)$ and 
$\pi_1(\Delta^*, N)$  of some  subfamilies in the family  
of Calabi-Yau hypersurfaces 
corresponding to two isomorphic maximal reflexive pairs 
$(\Delta, M^{\Delta})$ and $(\Delta, N^{\Delta^*})$. 
Moreover, the order of $\pi_1(\Delta)$ equals to the product of 
oders of $\pi_1(\Delta, M)$ and $\pi_1(\Delta^*, N)$. 
\end{prop}
\bigskip

\subsection{Reflexive simplices}

\hspace*{\parindent}

In this section we consider Calabi-Yau families ${\cal F}(\Delta)$, 
where $\Delta$ is a reflexive simplex of dimension $n$. 
Let $\{ p_0, \ldots, p_n \}$ be vertices of $\Delta$. 
There exists the  unique linear relation among $\{ p_i \}$ 
\[ \sum_{i = 0}^n w_i p_i = 0, \]
where $w_i > 0$ ($0 \leq i \leq n$) are integers  
and $g.c.d. (w_i ) = 1.$

\begin{opr}
{\rm The coefficients $w = \{ w_0, \cdots , w_n \}$ in the above linear 
relations are called  the {\em weights of the reflexive simplex} $\Delta$.}
\end{opr}
\bigskip

Let $\Delta^*$ be the dual reflexive simplex,  ${l}_0 , \ldots , {l}_n$ 
vertices of $\Delta^*$, $p_0, \dots , p_n$ vertices of $\Delta$. 
By definition 
\ref{inver.p}, we may  assume that  
\[ \langle p_i , {l}_j \rangle = -1 \; (i \neq j ), \]
i.e., that the equation $\langle x , {l}_j \rangle = -1$ defines the affine 
hyperplane in ${M}_{\bf Q}$ generated by  the $(n-1)$-dimensional face of 
$\Delta$ which does not contain  $p_j \in \Delta$. 

\begin{opr}
{\rm The  $(n+1)\times(n+1)$-matrix with integral coefficients 
\[ B(\Delta) = ( b_{ij} ) = ( \langle p_i, {l}_j \rangle ) \]
is called  {\em the matrix of the reflexive simplex $\Delta$}. }
\end{opr}
\bigskip 

\begin{theo}
Let  $\Delta$ be a reflexive $n$-dimensional simplex. Then 

{\rm (i)} the matrix $B(\Delta)$ is symmetric and its rank equals $n$;

{\rm (ii)} the diagonal coefficients $b_{ii}$ $( 0 \leq i \leq n)$ 
are positive and satisfy the equation 
\[\sum_{i =0}^n \frac{1}{b_{ii} + 1} = 1;\]

{\rm (iii)} the weights $\{ w_0, \ldots , w_n \}$ of $\Delta$ are    
the primitive integral solution of the linear  homogeneous  
system with the matrix $B(\Delta)$. Moreover, 
\[ w_i = \frac{l.c.m. (b_{ii} +1)}{b_{ii} +1}. \] 
\label{ref.mat}
\end{theo}

\proof The statement (i) follows from the fact that 
${\rm rk}\, \{ p_i \} = {\rm rk}\, \{ l_j \} = n $. One gets 
(ii) by the direct computation of the determinant of $B(\Delta)$ as a function 
on coefficients $b_{ii}$ ($0 \leq i \leq n$).  Finally, 
(iii) follows from (ii) by checking  that 
\[ \{ 1/(b_{11} +1), \dots , 1/(b_{nn} +1) \} \]
is a solution of the linear homogeneous system with the matrix $B(\Delta)$. 
\bigskip

\begin{coro}
The matrix $B(\Delta)$ depends only on the weights of $\Delta$.
\end{coro}
\bigskip

Let $\Delta$ be a reflexive simplex with weights $w = \{ w_i \}$. Using the 
vertices $\{ l_j \}$ of the dual reflexive simplex $\Delta^*$, we can 
define the homomorphism 
\[ \iota_{\Delta} \;:\;  M \rightarrow {\bf Z}^{n+1},\]
where  
\[ \iota_{\Delta} (m) = 
(\langle m, {l}_0 \rangle, \ldots , \langle m, {l}_n \rangle ). \] 
Obviously, $\iota_{\Delta}$ is injective and the image of $\iota_{\Delta}$ 
is contained in the $n$-dimensional sublattice $M(w)$ 
in ${\bf Z}^{n+1}$ defined by the equation 
\[ \sum_{i =0}^n  w_i x_i = 0.\]
Note that the image $\iota_{\Delta}(p_i)$ is the $i$-th row of $B(\Delta)$. We 
denote by $\Delta(w)$ the convex hull  of the points 
$\{ \iota_{\Delta}(p_i)\}$ in $M_{\bf Q}(w)$.
\medskip

\begin{theo}
The pair $(\Delta(w), M(w))$ is reflexive and satisfies the following 
conditions$:$ 

{\rm (i)} the corresponding to $(\Delta(w), M(w))$ 
toric Fano variety ${\bf P}_{\Delta(w)}$ 
is the weighted projective space ${\bf P}(w_0, \cdots, w_n)$;

{\rm (ii)} 
\[ \iota_{\Delta} \; : \; (\Delta, M) \rightarrow 
(\Delta(w), M(w)) \]
is a finite  morphism of reflexive pairs;

{\rm (iii)} $(\Delta(w), M(w))$ is a  maximal reflexive pair. 
\label{matrix}
\end{theo}

\proof The reflexivity of $(\Delta(w), M(w))$ and the condition (ii) 
follow immediately from 
the  definition of $\iota_{\Delta}$, since 
\[ \Delta(w) = \{ (x_0, \ldots, x_n) \in {\bf Q}^{n+1} \mid 
\sum_{i =0}^n w_i x_i =0, \; x_i \geq  -1 \; (0 \leq i \leq n) \}. \]

(i) The shifted by $(1, \ldots, 1)$ convex polyhedron 
\[ \Delta^{(1)}(w) = \Delta(w) + (1, \ldots, 1) \]
is the intersection of ${\bf Q}^{n+1}_{\geq 0}$ and the affine hyperplane 
\[ w_0 x_0 + \cdots + w_n x_n = w_0 + \cdots + w_n = 
l.c.m.\{ b_{ii} +1 \} = d. \]
Therefore, the integral points in $\Delta^{(1)}(w)$ can be identified 
with all possible monomials of degree $d$ in $n+1$ $w$-weighted 
independent variables.  

(iii) Assume that there exists a finite morphism 
\[ \phi\; :\; (\Delta(w), M(w)) \rightarrow (\Delta', M')\]
of reflexive pairs. Obviously, $\Delta'$ must be also a  simplex. 
Since  the linear mapping $\phi$    does  not change linear relations,  
$\Delta$ and $\Delta'$  must have the same weights  $w = \{w_i \}$. 
Therefore, by (ii), there exists a finite morphism 
\[ \iota_{\Delta'} \; : \; (\Delta', M') \rightarrow 
(\Delta(w), M(w)). \]
Therefore, by \ref{isom}, the reflexive pairs $(\Delta(w), M(w))$ and 
$(\Delta', M')$ are isomorphic. 
\bigskip

Since $B(\Delta) = B(\Delta^*)$ (see \ref{ref.mat} ), we obtain:  

\begin{coro}
Any reflexive simplex  $\Delta$ is selfdual.
\label{dual.simp}
\end{coro}
\bigskip

\subsection{Quotients of Calabi-Yau hypersurfaces 
in weighted projective spaces}

\hspace*{\parindent}

Theorem \ref{matrix} implies

\begin{coro}
The fundamental group $\pi_1(\Delta)$ of a  reflexive simplex 
$\Delta$ depends only on the  weights $w = \{ w_i \}$. This group is 
dual to the quotient of the lattice $M(w) \subset {\bf Z}^{n+1}$ by the 
sublattice $M_B(\Delta)$ generated by rows of the matrix $B(\Delta)$. 
\end{coro}

Now we calculate  the group  $M(w) / M_B(w)$ explicitly. 
\bigskip

Consider two integral sublattices of 
rang $n+1$ in ${\bf Z}^{n+1}$ 
\[ \tilde{M} (w) = M(w) \oplus {\bf Z}\langle (1, \ldots, 1) \rangle, \]
and 
\[ \tilde{M}_B (w) = M_B(w) \oplus {\bf Z}\langle (1, \ldots, 1) \rangle. \]
Note that $\tilde{M}(w) / \tilde{M}_B(w) \cong M(w) / M_B(w)$. 

Let $\mu_r$ denote the group of complex $r$-th roots of unity. Put   
$d_i = b_{ii} +1$ ($ 0 \leq i \leq n$), $d = l.c.m.\{ d_i \}$.  

The sublattice $\tilde{M}(w)$ is the kernel of the surjective 
homomorphism 
\[ \gamma_w \;: \; {\bf Z}^{n+1} \rightarrow \mu_d, \]
\[ \gamma_w (a_0, \ldots , a_n) = g^{w_0a_0 + \cdots + w_n a_n}, \]
where $g$ is a generator of $\mu_d$.

The sublattice $\tilde{M}_B(w)$ is generated by $(1, \ldots, 1)$ and 
\[ (d_0, 0, \ldots, 0), (0, d_1, \ldots, 0), 
\dots, (0, \ldots, 0, d_n). \]

Therefore, $\tilde{M}_B(w)$ can be represented as the sum of the infinite 
cyclic group  generated by $(1, \ldots, 1)$ and 
the kernel of the surjective homomorphism 
\[ \gamma \;: \; {\bf Z}^{n+1} \rightarrow \mu_{d_0} \times \mu_{d_1} 
\times \cdots \times \mu_{d_n} ,  \]
\[ \gamma (a_0, \ldots , a_n) = g_0^{a_0} g_1^{a_1} \cdots g_n^{a_n}, \]
where $g_i$ is a generator of $\mu_{d_i}$ ($ 0 \leq i \leq n$). 
The order of the element $(1, \ldots, 1)$ modulo ${\rm ker}\, \gamma$ 
equals  $d$. 
Thus, 
\[ {\bf Z}^{n+1} /  \tilde{M}_B(w) \cong (\mu_{d_0} \times \mu_{d_1} 
\times \cdots \times \mu_{d_n}) / \mu_d ,\]
where the subgroup $\mu_d$ is generated by 
$g_0g_1 \cdots g_n$. Finally, we obtain.

\begin{theo}
The fundamental group $\pi_1(\Delta^*) = M(w)/M_B(w)$ of the  
reflexive simplex $\Delta^*$ with weights 
$w = \{ w_i \}$ is isomorphic  to the  kernel 
of the surjective homomorphism
\[ \overline{\gamma}_w  \; : \; (\mu_{d_0} \times \mu_{d_1} 
\times \cdots \times \mu_{d_n}) / \mu_d  \rightarrow \mu_d, \]
\[ \overline{\gamma}_w ( g_0^{a_0} g_1^{a_1} \cdots g_n^{a_n} ) = 
g^{w_0a_0 + \cdots + w_n a_n}. \]
\label{fund1}
\end{theo}

By duality, we conclude. 

\begin{coro}
The fundamental group $\pi_1(\Delta)$ of a reflexive simplex $\Delta$ with 
weights $w$ is  isomorphic to the kernel of the surjective 
homomorphism 
\[ (\mu_{d_0} \times \mu_{d_1} \times \cdots \times \mu_{d_n}) / \mu_d 
\rightarrow \mu_d, \]
where the homomorphism to $\mu_d$ is the product of complex numbers 
in $\mu_{d_0}, \mu_{d_1}, \ldots , \mu_{d_n}$, and the 
embedding of $\mu_d$ in $\mu_{d_0} \times \mu_{d_1} 
\times \cdots \times \mu_{d_n}$ is defined by 
\[ g \mapsto ( g^{w_0}, \ldots, g^{w_n}). \]
\label{fund2}
\end{coro}

\begin{coro}
The order  of $\pi_1(\Delta)$ in the above theorem equals 
\[ \frac{d_0 d_1 \cdots d_n}{d^2}. \]
\end{coro}
\bigskip

\begin{exam}
{\rm Let $(d_0, d_1, \ldots , d_n )$ $(d_i > 0)$ be an 
integral solution of the equation
\begin{equation}
\label{weights.proj}
 \sum_{i =0}^{n} \frac{1}{d_i} = 1. 
\end{equation} 
Then the quasi-homogeneous equation 
\[ v_0^{d_0} + v_1^{d_1} + \cdots + v_n^{d_n} = 0\]
defines a $\Delta(w)$-regular Calabi-Yau hypersurface of Fermat-type in the 
weighted projective space 
\[ {\bf P}_{\Delta(w)} = {\bf P}(w_0, \ldots, w_n) , \]
where 
\begin{equation}
 w_i = \frac{l.c.m. (d_i)}{d_i}. 
\label{new.weights}
\end{equation}}
\end{exam}

\begin{coro}
The family ${\cal F}(\Delta(w))$ of 
Calabi-Yau hypersurfaces in the weighted projective 
space ${\bf P}_{\Delta(w)}$ consists of deformations of Fermat-type 
hypersurfaces. 
If $\Delta$ is a reflexive  simplex  with weights $w =\{ w_i \}$, then 
the corresponding family ${\cal F}(\Delta)$ in ${\bf P}_{\Delta}$ 
consists of quotients of  some subfamily in ${\cal F}(\Delta(w))$ 
by the  action of the finite abelian group $\pi_1( \Delta,M)$. 
\end{coro}

If we consider a special case $n =4$, we obtain as a corollary 
the result of Roan in \cite{roan1}. To prove this, one should  
use our general result on Calabi-Yau $3$-folds 
constructed from $4$-dimensional reflexive polyhedra (see \ref{euler.iso}) 
and apply the following simple statement.

\begin{prop} Let $(\Delta, M)$ be a reflexive pair such that 
$\Delta$ is a $4$-dimensional reflexive simplex with weights 
$(w_0, \ldots, w_4)$. Then  
the family ${\cal F}(\Delta)$ consists of quotients by $\pi_1(\Delta, M)$ 
of Calabi-Yau hypersurfaces in the weighted projective space 
${\bf P}(w_0, \ldots, w_4)$ whose equations are invariant under the 
canonical diagonal action of $\pi_1(\Delta, M)$ on 
${\bf P}(w_0, \ldots, w_4)$. 
\end{prop}
\bigskip

\section{Classification of reflexive polyhedra}

\subsection{Finiteness theorems}

\hspace*{\parindent}

The following theorem was proved in $\cite{bat01}$.

\begin{theo} 
Let $P$ be a convex $n$-dimensional polyhedron in ${\bf Q}^n$ 
with integral vertices. Assume that $P$  contains the zero point 
$0 \in {\bf Q}^n$ in its interior and  for some positive integer 
$m$ the distance between any $(n-1)$-dimensional 
face of $P$ and $0$ is not more than $m$. Then the volume of $P$ is not 
greater than some explicitly calculated 
constant $C(n,m)$ depending only on $n$ and $m$. 
\label{t.distance}
\end{theo}

If we take $m =1$, then we get immediately.

\begin{coro} 
The volume of an $n$-dimensional reflexive polyhedron $\Delta$ relative to 
the lattice $M$ is bounded by a constant $C(n)$ depening only on $n$. 

In particular, there exist (up to an isomorphism)  only finitely many 
reflexive pairs 
$(\Delta, M)$  of  fixed dimension $n$. 
\end{coro}

This corollary follows also from the  result in  \cite{hensley} and from 
the recent more strong theorem of A.A. Borisov and L.A. Borisov \cite{boris}:  

\begin{theo}
Let $P$ be a convex $n$-dimensional polyhedron in ${\bf Q}^n$ 
with integral vertices containing the origin $0$ in its interior. Assume that 
for some positive real number $\varepsilon \; (\varepsilon <1)$, one has 
\[ \varepsilon P \cap {\bf Z}^n = 0. \]
Then the volume of $P$ is bounded by some constant $C(\varepsilon, n)$ 
depending only on $\varepsilon$ and $n$.
\end{theo}
\bigskip

The complete classification of reflexive polyhedra of dimension $n \geq 3$ 
is an  unsolved combinatorial problem. This problem is 
equivalent to the classification of $n$-dimensional toric Fano varieties 
with Gorenstein  singularities. 

Since the problem of the complete classification of reflexive polyhedra 
is very hard in general, it is natural to try to classify at least 
special reflexive polyhedra satisfying  some additional conditions, 
or equivalently, {\em special} Gorenstein toric Fano varieties.  
For instance, the classification of ${\em smooth}$ 
$n$-dimensional toric Fano varieties (see \ref{smooth.fano}) was received  
for $n \leq 4$ $(\cite{bat0}, \cite{bat11}, \cite{watanabe})$.

\begin{theo}
 Let $F(n)$ be the number of different $n$-dimensional toric 
Fano varieties. Then $F(1) = 1$, $F(2) = 5$, $F(3) = 18$, $F(4) = 123$. 
\end{theo}

Classification of centrally  symmetric reflexive polyhedra $\Delta$ 
of arbitrary dimension for which ${\bf P}_{\Delta}$ is a smooth Fano 
maniflod was received  by 
Ewald \cite{ewald}, and Voskresenski\^i-Klyachko \cite{vosk.kl}. 

One can easily prove  that there exists only one  
$1$-dimensional reflexive polyhedron. 
We explain below the complete classification of $2$-dimensional 
reflexive polyhedra obtained independently by author \cite{bat.thes} and 
Koelman \cite{koelman}. This classification bases on the folowing lemma.

\begin{lem}
Asssume that  a relexive polygon $\Delta$ has $k$ integral points 
at the boundary $\partial \Delta$, then the dual polygon $\Delta^*$ has 
$12 -k$ integral points at the boundary $\partial \Delta^*$. 
\label{12}
\end{lem}

\proof Indeed, the  Del Pezzo surface ${\bf P}_{\Delta}$   
contains   only rational double singular points. The minimal  
desingularization ${\hat{\bf P}}_{\Delta}$ of ${\bf P}_{\Delta}$ 
is crepant and unique. It is easy to see that 
${\hat{\bf P}}_{\Delta}$ is a toric 
projective surface associated with the complete fan whose $1$-dimensional 
cones are generated by integral points in $\partial \Delta^* \cap N$. Therefore,  
the number $k_{\Delta^*}$ of integral  points in $\partial \Delta^* \cap N$ 
equals the 
Euler characteristic  of ${\hat{\bf P}}_{\Delta}$. On the other hand, 
the number $k_{\Delta}$  
of integral  points in $\partial \Delta \cap M$ equals $d(\Delta)$, or 
the self-intersection 
number of the anticanonical divisor on  ${\hat{\bf P}}_{\Delta}$. By 
Noether's formula for rational surfaces,  $k_{\Delta^*} + k_{\Delta} = 
12$.

\begin{theo}
There exist up to isomorphism exactly 16 different 
$2$-dimensional reflexive pairs defined by the following reflexive polygons 
$P_i$ $(P_i$ is  dual to $P_{17-i}$ for $i \leq 6):$ 

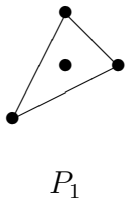
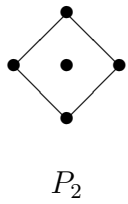
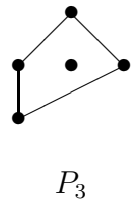
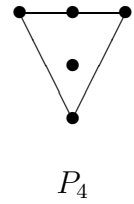
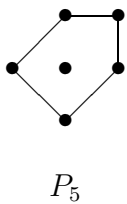
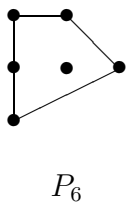
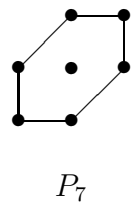
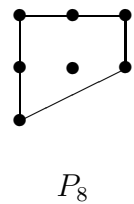
\begin{center}
\begin{picture}(450,200)

\put(50,100){\makebox(0,0){$\bullet$}}
\put(150,100){\makebox(0,0){$\bullet$}}
\put(250,100){\makebox(0,0){$\bullet$}}
\put(350,100){\makebox(0,0){$\bullet$}}

\put(50,100){\makebox(0,-90){$P_1$}}
\put(150,100){\makebox(0,-90){$P_2$}}
\put(250,100){\makebox(0,-90){$P_3$}}
\put(350,100){\makebox(0,-90){$P_4$}}

\put(50,100){\makebox(40,0){$\bullet$}}
\put(50,100){\makebox(0,40){$\bullet$}}
\put(50,100){\makebox(-40,-40){$\bullet$}}

\put(30,80){\line(2,1){40}}
\put(30,80){\line(1,2){20}}
\put(70,100){\line(-1,1){20}}

\put(150,100){\makebox(40,0){$\bullet$}}
\put(150,100){\makebox(0,40){$\bullet$}}
\put(150,100){\makebox(-40,0){$\bullet$}}
\put(150,100){\makebox(0,-40){$\bullet$}}

\put(170,100){\line(-1,1){20}}
\put(150,80){\line(-1,1){20}}
\put(150,80){\line(1,1){20}}
\put(130,100){\line(1,1){20}}

\put(250,100){\makebox(0,40){$\bullet$}}
\put(250,100){\makebox(-40,-40){$\bullet$}}
\put(250,100){\makebox(40,0){$\bullet$}}
\put(250,100){\makebox(-40,0){$\bullet$}}

\put(230,80){\line(2,1){40}}
\put(230,80){\line(0,1){20}}
\put(250,120){\line(-1,-1){20}}
\put(270,100){\line(-1,1){20}}

\put(350,100){\makebox(40,40){$\bullet$}}
\put(350,100){\makebox(0,40){$\bullet$}}
\put(350,100){\makebox(-40,40){$\bullet$}}
\put(350,100){\makebox(0,-40){$\bullet$}}

\put(350,80){\line(1,2){20}}
\put(350,80){\line(-1,2){20}}
\put(330,120){\line(1,0){40}}

\end{picture}
\end{center}

\begin{center}
\begin{picture}(450,200)

\put(50,100){\makebox(0,0){$\bullet$}}
\put(150,100){\makebox(0,0){$\bullet$}}
\put(250,100){\makebox(0,0){$\bullet$}}
\put(350,100){\makebox(0,0){$\bullet$}}

\put(50,100){\makebox(0,-90){$P_5$}}
\put(150,100){\makebox(0,-90){$P_6$}}
\put(250,100){\makebox(0,-90){$P_7$}}
\put(350,100){\makebox(0,-90){$P_8$}}

\put(50,100){\makebox(40,0){$\bullet$}}
\put(50,100){\makebox(0,40){$\bullet$}}
\put(50,100){\makebox(-40,0){$\bullet$}}
\put(50,100){\makebox(0,-40){$\bullet$}}
\put(50,100){\makebox(40,40){$\bullet$}}

\put(30,100){\line(1,1){20}}
\put(30,100){\line(1,-1){20}}
\put(50,80){\line(1,1){20}}
\put(50,120){\line(1,0){20}}
\put(70,100){\line(0,1){20}}

\put(150,100){\makebox(0,40){$\bullet$}}
\put(150,100){\makebox(-40,-40){$\bullet$}}
\put(150,100){\makebox(40,0){$\bullet$}}
\put(150,100){\makebox(-40,0){$\bullet$}}
\put(150,100){\makebox(-40,40){$\bullet$}}

\put(130,80){\line(2,1){40}}
\put(130,80){\line(0,1){40}}
\put(130,120){\line(1,0){20}}
\put(170,100){\line(-1,1){20}}

\put(250,100){\makebox(40,0){$\bullet$}}
\put(250,100){\makebox(0,40){$\bullet$}}
\put(250,100){\makebox(-40,0){$\bullet$}}
\put(250,100){\makebox(0,-40){$\bullet$}}
\put(250,100){\makebox(40,40){$\bullet$}}
\put(250,100){\makebox(-40,-40){$\bullet$}}

\put(230,80){\line(1,0){20}}
\put(230,80){\line(0,1){20}}
\put(250,80){\line(1,1){20}}
\put(230,100){\line(1,1){20}}
\put(250,120){\line(1,0){20}}
\put(270,100){\line(0,1){20}}

\put(350,100){\makebox(0,40){$\bullet$}}
\put(350,100){\makebox(-40,-40){$\bullet$}}
\put(350,100){\makebox(40,0){$\bullet$}}
\put(350,100){\makebox(-40,0){$\bullet$}}
\put(350,100){\makebox(-40,40){$\bullet$}}
\put(350,100){\makebox(40,40){$\bullet$}}

\put(330,80){\line(0,1){40}}
\put(330,80){\line(2,1){40}}
\put(330,120){\line(1,0){40}}
\put(370,100){\line(0,1){20}}
\end{picture}
\end{center}

\begin{center}
\begin{picture}(450,200)

\put(50,100){\makebox(0,0){$\bullet$}}
\put(150,100){\makebox(0,0){$\bullet$}}
\put(250,100){\makebox(0,0){$\bullet$}}
\put(350,100){\makebox(0,0){$\bullet$}}

\put(50,100){\makebox(0,-90){$P_9$}}
\put(150,100){\makebox(0,-90){$P_{10}$}}
\put(250,100){\makebox(0,-90){$P_{11}$}}
\put(350,100){\makebox(0,-90){$P_{12}$}}

\put(50,100){\makebox(40,0){$\bullet$}}
\put(50,100){\makebox(0,40){$\bullet$}}
\put(50,100){\makebox(-40,0){$\bullet$}}
\put(50,100){\makebox(0,-40){$\bullet$}}
\put(50,100){\makebox(40,40){$\bullet$}}
\put(50,100){\makebox(-40,40){$\bullet$}}

\put(50,80){\line(1,1){20}}
\put(50,80){\line(-1,1){20}}
\put(70,100){\line(0,1){20}}
\put(30,100){\line(0,1){20}}
\put(30,120){\line(1,0){40}}

\put(150,100){\makebox(0,40){$\bullet$}}
\put(150,100){\makebox(-40,-40){$\bullet$}}
\put(150,100){\makebox(40,0){$\bullet$}}
\put(150,100){\makebox(-40,0){$\bullet$}}
\put(150,100){\makebox(-40,40){$\bullet$}}
\put(150,100){\makebox(-40,80){$\bullet$}}

\put(170,100){\line(-1,1){40}}
\put(170,100){\line(-2,-1){40}}
\put(130,80){\line(0,1){60}}

\put(250,100){\makebox(40,0){$\bullet$}}
\put(250,100){\makebox(0,40){$\bullet$}}
\put(250,100){\makebox(-40,0){$\bullet$}}
\put(250,100){\makebox(0,-40){$\bullet$}}
\put(250,100){\makebox(-40,40){$\bullet$}}
\put(250,100){\makebox(-40,80){$\bullet$}}
\put(250,100){\makebox(-40,-40){$\bullet$}}

\put(270,100){\line(-1,1){40}}
\put(230,80){\line(0,1){60}}
\put(270,100){\line(-1,-1){20}}
\put(230,80){\line(1,0){20}}

\put(350,100){\makebox(40,0){$\bullet$}}
\put(350,100){\makebox(0,40){$\bullet$}}
\put(350,100){\makebox(-40,0){$\bullet$}}
\put(350,100){\makebox(0,-40){$\bullet$}}
\put(350,100){\makebox(-40,40){$\bullet$}}
\put(350,100){\makebox(40,-40){$\bullet$}}
\put(350,100){\makebox(-40,-40){$\bullet$}}

\put(330,80){\line(0,1){40}}
\put(330,80){\line(1,0){40}}
\put(330,120){\line(1,0){20}}
\put(370,80){\line(0,1){20}}
\put(370,100){\line(-1,1){20}}

\end{picture}
\end{center}

\begin{center}
\begin{picture}(450,230)

\put(50,100){\makebox(0,0){$\bullet$}}
\put(150,100){\makebox(0,0){$\bullet$}}
\put(250,100){\makebox(0,0){$\bullet$}}
\put(350,100){\makebox(0,0){$\bullet$}}

\put(50,100){\makebox(0,-90){$P_{13}$}}
\put(150,100){\makebox(0,-90){$P_{14}$}}
\put(250,100){\makebox(0,-90){$P_{15}$}}
\put(350,100){\makebox(0,-90){$P_{16}$}}

\put(50,100){\makebox(-40,120){$\bullet$}}
\put(50,100){\makebox(-40,80){$\bullet$}}
\put(50,100){\makebox(0,40){$\bullet$}}
\put(50,100){\makebox(-40,0){$\bullet$}}
\put(50,100){\makebox(0,-40){$\bullet$}}
\put(50,100){\makebox(-40,40){$\bullet$}}
\put(50,100){\makebox(40,-40){$\bullet$}}
\put(50,100){\makebox(-40,-40){$\bullet$}}

\put(30,80){\line(0,1){80}}
\put(30,80){\line(1,0){40}}
\put(70,80){\line(-1,2){40}}

\put(150,100){\makebox(40,0){$\bullet$}}
\put(150,100){\makebox(0,40){$\bullet$}}
\put(150,100){\makebox(-40,0){$\bullet$}}
\put(150,100){\makebox(0,-40){$\bullet$}}
\put(150,100){\makebox(-40,40){$\bullet$}}
\put(150,100){\makebox(40,-40){$\bullet$}}
\put(150,100){\makebox(-40,-40){$\bullet$}}
\put(150,100){\makebox(-40,80){$\bullet$}}

\put(130,80){\line(0,1){60}}
\put(130,80){\line(1,0){40}}
\put(170,80){\line(0,1){20}}
\put(170,100){\line(-1,1){40}}

\put(250,100){\makebox(40,0){$\bullet$}}
\put(250,100){\makebox(0,40){$\bullet$}}
\put(250,100){\makebox(-40,0){$\bullet$}}
\put(250,100){\makebox(0,-40){$\bullet$}}
\put(250,100){\makebox(-40,40){$\bullet$}}
\put(250,100){\makebox(40,-40){$\bullet$}}
\put(250,100){\makebox(-40,-40){$\bullet$}}
\put(250,100){\makebox(40,40){$\bullet$}}

\put(230,80){\line(0,1){40}}
\put(230,80){\line(1,0){40}}
\put(270,80){\line(0,1){40}}
\put(230,120){\line(1,0){40}}

\put(350,100){\makebox(40,0){$\bullet$}}
\put(350,100){\makebox(0,40){$\bullet$}}
\put(350,100){\makebox(-40,0){$\bullet$}}
\put(350,100){\makebox(0,-40){$\bullet$}}
\put(350,100){\makebox(-40,40){$\bullet$}}
\put(350,100){\makebox(40,-40){$\bullet$}}
\put(350,100){\makebox(-40,-40){$\bullet$}}
\put(350,100){\makebox(80,-40){$\bullet$}}
\put(350,100){\makebox(-40,80){$\bullet$}}

\put(330,80){\line(0,1){60}}
\put(330,80){\line(1,0){60}}
\put(390,80){\line(-1,1){60}}

\end{picture}
\end{center}

\end{theo} 
\bigskip

\proof Applying lemma \ref{12}, we see that it is sufficient 
to classify all smooth toric surfaces ${\hat{\bf P}}_{\Delta}$ having 
the Picard number $\leq 4$. It remains to apply the result of 
Oda on the classification of all  smooth toric surfaces with the Picard number 
$\leq 4$ (see \cite{oda1}). 
\bigskip

\subsection{Classification of reflexive simplices} 

\hspace*{\parindent}

Now we want to  explain how one can classify 
all reflexive simplices of fixed dimension $n$. 
Finiteness theorem for reflexive simplices of fixed dimension $n$ is a more 
easy statement which  follows 
from the next proposition (see for example \cite{zak}). 
\medskip 

\begin{prop}
 For any fixed positive integer $n$, there exists only finitely many  
 integral solutions $(d_0, d_1, \ldots , d_n )$ $(d_i > 0)$ of the 
equation {\rm (\ref{weights.proj})}. 
\end{prop}

Using  the formula (\ref{new.weights}), for each solution of 
(\ref{weights.proj}),  we can calculate  
the weights of for the corresponding reflexive simplex $\Delta$, i.e., the  
weights  for the weighted projective space ${\bf P}(w_0, w_1, \ldots , w_n )$. 
We can assume that $w_0 \leq \cdots \leq w_n$.  If all weights 
$w_0, \ldots, w_n$ 
are different, then there exists only one finite morphism 
\[ (\Delta, M) \rightarrow (\Delta(w), M(w)). \]
Therefore, there is the one-to-one correspondence between reflexive 
simplices $\Delta$ with weights 
$w =\{w_i\}$ and  sublattices $M \subset M(w)$ such that $M$ contains 
$M_B(w)$. So the classification reduces to the enumeration of subgroups of 
the finite abelian fundamental group $M(w)/M_B(w)$  calculated in 
\ref{fund1}.

If the  weights $w_0, \ldots, w_n$ are not different, then we must take 
into account the action of the product of symmetric groups 
\[ S(w) = S_{r_1} \times \cdots \times S_{r_l}, \]
where $l$ is the number of different weights  among $w_0, \ldots, w_n$ and 
$r_i$ is the corresponding multiplicity $( \sum_i r_i = n +1)$. 
Reflexive simplices with weights $\{w_i \}$ are in the one-to-one 
correspondence to the  orbits of the $S(w)$-action on subgroups of the 
fundamental group $\pi_1(\Delta^*)$.
\bigskip

\begin{exam}
{\rm If $n =2$, then there exist exactly $3$ different 
integral solutions of the equation  
 (\ref{weights.proj}):

{\em Case I}. $\frac{1}{3} + \frac{1}{3} + \frac{1}{3} = 1$. The  
fundamental group of the corresponding reflexive triangle is 
isomorphic to ${\bf Z}/3{\bf Z}$. There are two different subgroups 
in  ${\bf Z}/3{\bf Z}$ which define  $P_1$ and $P_{16}$. 

{\em Case II}. $\frac{1}{2} + \frac{1}{4} + \frac{1}{4} = 1$. The  
fundamental group of the corresponding reflexive triangle is 
isomorphic to ${\bf Z}/2{\bf Z}$. There are two different subgroups 
in  ${\bf Z}/2{\bf Z}$ which define  $P_4$ and $P_{13}$. 
 
{\em Case III}. $\frac{1}{2} + \frac{1}{3} + \frac{1}{6} = 1$. The  
fundamental group of the corresponding reflexive triangle is 
trivial.  So there is one trivial subgroup  defining  $P_{10}$.} 
\end{exam}

\begin{exam}
{\rm If $n =3$, then there exist exactly $14$ different integral 
solutions of the equation   (\ref{weights.proj}).  
If we take, for example, the solution $(d_1, d_2, d_3,d_4) =  (2,3,7,42)$, 
then $d = 42$, 
and the fundamental group of a simplex with the corresponding weights 
$w = (21, 14, 6, 1)$ is trivial. So there exists up to isomorphism only 
one reflexive $3$-dimensional simplex with weights $w = 
(21, 14, 6, 1)$, i.e., this simplex 
is isomorphic to its dual. }
\end{exam}   

\begin{exam}
{\rm If $n =4$, then there exist exactly $147$ different integral 
solutions of the equation   (\ref{weights.proj}). In general, one can 
obtain from every such a solution many different 
$4$-dimensional reflexive simplices 
having the same weights. 
For example, the classification of Klemm and Theisen in 
\cite{klemm2} shows 
that there exist 8 reflexive 
$4$-dimensional simplices with  weights $w = (1,1,1,1,1)$, 
28 reflexive 
$4$-dimensional simplices with  weights $w =(2,1,1,1,1)$, 
30 reflexive 
$4$-dimensional simplices with  weights $w =(4,1,1,1,1)$, and 
14 reflexive 
$4$-dimensional simplices with  weights $w= (5,2,1,1,1)$. 
For these examples the group $S(\omega)$ is isomorphic to 
$S_5$, $S_4$, or $S_3$. So the classification 
of these reflexive simplices consists of the enumeration of 
finite subgroups in $M(\omega)/M_B(\omega)$ up to the action of $S(\omega)$. 
}  
\end{exam}
\bigskip

\section{Calabi-Yau varieties in singularity theory} 

\subsection{The duality for Tsuchihashi's cusp singularities}

\hspace*{\parindent}

We consider an example of the duality due to Tsuchihashi \cite{tsuch} 
arising in cusp  singularity theory.  
This duality has several similar properties  with our construction of 
mirrors via dual polyhedra.  
\bigskip

Let $N \cong {\bf Z}^r$ be a free abelian group of rank $r$, 
${\rm GL}(r, {\bf Z})$ the group of authomorphisms of $N$. Take a 
subgroup $\Gamma \subset {\rm GL}(r, {\bf Z})$, and a convex cone 
$C \in N_{\bf R}$ which is invariant under the action of $\Gamma$. 

\begin{opr} {\rm 
\cite{oda1}  We say the pair $(C, \Gamma)$ belongs to the set ${\cal P}(N)$ 
if the pair $(C, \Gamma)$ satisfies the following  conditions:

(i) $C \subset N_{\bf R}$ is nonempty open convex cone such that for the 
closure $\overline{C}$, one has 
\[ \overline{C} \cap - \overline{C} = \{ 0 \} \in N_{\bf R}; \]

(ii) the induced action of $\Gamma$ on $D = C/{\bf R}_{> 0}$ is properly 
discontinious without fixed points and has a compact quotient $D/\Gamma$. }
\end{opr} 

The ${\bf Z}$-module $M$ dual to $N$ admits the contragredient action of 
${\rm GL}(r, {\bf Z})$ compatible with the pairing 
$\langle \;, \; \rangle : M \times N \rightarrow {\bf Z}$, i.e., 
$\langle g(m), g(n) \rangle = \langle m, n \rangle$ for all $g \in 
{\rm GL}(r, {\bf Z})$, $m \in M$ and $n \in N$. If we define the {\em open 
dual cone} of $C$ as 
\[ C' : = \{ m \in M_{\bf R} \mid \langle m, n \rangle > 0 \;{\rm for} \;{\rm 
all}\; n \in \overline{C} \setminus \{0 \} \}, \]
then $C'$ is an convex  cone invariant under the contaradient action of 
$\Gamma$, and  $(C', \Gamma)$ belongs to ${\cal P}(M)$. 
\medskip

There exists a $\Gamma$-invariant rational polyhedral 
subdivision $\Sigma$ of the cone $C$. We obtain a $\Gamma$-stabe 
open set $\tilde{U}$ of the toric variety associated with $\Sigma$. 
The set $\tilde{U}$  contains a closed analytic subset $\tilde{X}$ of 
codimension one. The action of $\Gamma$ is properly discontinious and without 
fixed points. Let 
\[ U = \tilde{U}/\Gamma \subset X = \tilde{X}/\Gamma. \]
The subvariety $X$ in $U$ can be contracted to a singular point which we 
denote by ${\rm Cusp}(C,\Gamma)$. 

\begin{opr}
{\rm The singular point ${\rm Cusp}(C,\Gamma)$ is called the {\em Tsuchihashi 
cusp singularity associated to} $(C, \Gamma)$. }
\end{opr}

\begin{opr}
{\rm Let $C'$ be the dual of $C$ cone. Then 
${\rm Cusp}(C',\Gamma)$ is called the {\em dual Tsuchihashi cusp 
singularity}.}
\end{opr}

There are the following  properties of the duality for cusp 
singularities which makes it analogous to our construction of  mirrors. 
\bigskip

(i) The duality for Tsuchihashi's cusp singularities  bases 
on the duality of convex cones in dual spaces. This is the same 
idea we used in our definition of  mirror families (cf. \ref{dual.cones}).  

(ii) If $\Gamma \subset SL(r, {\bf Z})$, then the exceptional divisor 
$X \subset U$ can be considered as a "degenerated Calabi-Yau variety", 
namely,  we can construct a differential $r$-form on $U$ which vanishes 
nowhere on $U \setminus X$, moreover, we can prove that ${\cal O}_U (-X)$ 
is the dualizing sheaf on $U$.  

(iii) The duality interchange two invariants of Tsuchihashi's cusp 
singularities ${\rm Cusp}(C, \Gamma)$, the arithmetical defect 
$\chi_{\infty}(C,\Gamma)$ and the zero value of the Ogata zeta function 
$Z(C, \Gamma,0)$, i.e., one has relations 
\[ \chi_{\infty}(C,\Gamma) = Z(C', \Gamma,0), \]
and 
\[ \chi_{\infty}(C',\Gamma) = Z(C, \Gamma,0). \]  
This property of the duality was recently considered  by Ishida in 
\cite{ishida1}, \cite{ishida2}. 
\bigskip

\subsection{The strange duality of Arnold and mirrors of $K3$-surfaces}

\hspace*{\parindent}

It was conjectured by physicists that the mirror symmetry in superconformal 
theory is related to  the strange duality of Arnold \cite{lyn.schimm}. 
We shall give 
below several results on behalf of this conjecture of physicists. 
\bigskip 
 
Recall that there exist exactly 14 exceptional quasi-homogeneous unimodular 
hypersurface 
singularities $X(p,q,r)$ connected with 14 Lobachevski triangles with angles 
$\pi/p$, $\pi/q$, $\pi/r$. The values of $p, q,r$ are 
called Dolgachev numbers of these singularities.  
Gabrielov has found the intersection forms for vanishing cycles of all 
unimodular singularities.  The lattice $E(a,b,c)$ of vanishing cycles 
are defined by diagrams $T(a,b,c)$ having the shape of a letter $T$ 
with $a,b,c$ points on its closed segments. The values of $a, b, c$ are 
called Gabrielov numbers of singularities. 

In 1974 , Arnold reported on the existence of the involution 
on 14 exceptional singularities \cite{arnold}. 
The  duality induced by the involution 
satisfies the following conditions: 

(i) The Gabrielov numbers of every exceptional singularity are the Dolgachev 
numbers of the dual one; the Gabrielov numbers of the latter are the 
Dolgachev numbers of the former. 

(ii) The sum $(p + q + r) + (p' + q' + r')$ of the Dolgachev numbers 
of two dual exceptional singularities $X(p,q,r)$ and $X(p',q',r')$ is 
always equal to $24$. 

(iii) Any two dual singularities defined by diagrams $T(a,b,c)$ and 
$T(a',b',c')$ have the same discriminants 
\[ \Delta(a,b,c) = \Delta(a',b',c'), \]
where 
\[ \Delta(x_1, x_2, x_3) = (-1)^{\sum x_i} x_1 x_2 x_3 (1 - \sum 
\frac{1}{x_i} ). \]
\medskip

Pinkham  has given in \cite{pink} an interpretation of 
the strange Arnold duality in terms of 
cycles on a $K3$-surface $W$. Let $H$ denotes 
hyperbolic lattice of rank $2$, and $E_8$ is unimodular lattice of 
rank $8$. It is well-know that the homology lattice of $2$-dimensional 
cycles $H_2(W, {\bf Z})$ is isometric to 
\[ E_8 \oplus E_8 \oplus H \oplus H \oplus H.\]
The interpretation of the strange Arnold duality bases on the 
following property

(iv) The lattices $E(a,b,c)$ and $E(a',b',c')$ of dual singularities 
have canonical embeddings as orthogonal sublattices in the quadratic 
${\bf Z}$-module 
$E_8 \oplus E_8 \oplus H \oplus H$ of rank $20$. 
\bigskip
    
Our definition of the mirror operation has  a property  
which  is analogous to the property (ii) of Arnold's duality.

\begin{theo} 
Let $\Delta$ be a $3$-dimensional reflexive polyhedron, $\Delta^*$ 
the dual $3$-dimensional polyhedron. Then 
\[ \sum_{{\rm dim}\, \Theta =1} d(\Theta) d(\Theta^*) = 24. \]
\label{24}
\end{theo}

\proof Denote by $\nu_i(\Delta)$ the number of $i$-dimensional faces of 
$\Delta$. 

Let $Z_f$ be a $\Delta$-regular affine hypersurface, then its 
projective closure $\overline{Z}_f$ is a $K3$-surface with only rational 
double points, and ${\hat{Z}}_f$ is the unique minimal 
desingularization of  $\overline{Z}_f \subset {\bf P}_{\Delta}$. Since 
${\hat{Z}}_f$ contains $Z_f$ as an open subset, we get the long 
exact sequence of cohomologies with compact supports 
\[ \cdots \rightarrow H^1_c(Y) \rightarrow H_c^{2}(Z_f) 
\stackrel{\beta_1}{\rightarrow}  
H_c^{2}(\hat{Z}_f) {\rightarrow}  
H_c^{2}(Y) {\rightarrow}  \\ 
H_c^{3}(Z_f) \stackrel{\beta_2}{\rightarrow}  
H_c^{3}(\hat{Z}_f) {\rightarrow} \cdots, \] 
where  $Y = \hat{Z}_f \setminus Z_f$. The space $H^1_c (Y)$ 
does not contain elements of the Hodge type $(1,1)$. So the homomorphism 
$\beta_1$ is injective on the subspace $H^{1,1}_c(Z_f) \subset H^2_c(Z_f)$. 
Since the cohomology group $H^3$  of a $K3$-surface is zero, $\beta_2$ is 
trivial homomorphism. Since $H^2(Y)$ consists only of  elements of the 
Hodge type 
$(1,1)$, we get the exact sequence 
\[ 0 \rightarrow H_c^{1,1}(Z_f) 
\stackrel{\beta_1}{\rightarrow}  
H_c^{1,1}(\hat{Z}_f) {\rightarrow}  
H_c^{2}(Y) {\rightarrow}  \\ 
H_c^{3}(Z_f) {\rightarrow}  0.  \]  
By the Lefschetz-type theorem (cf. \ref{Lef}), 
$H_c^{3}(Z_f) \cong H_c^5 ({\bf T})$, i.e., 
${\rm dim}\, H_c^{3}(Z_f) =3$. 

The dimension of $H_c^2(Y)$ equals 
the number of irreducible components of $Y$. Note that every $0$-dimensional 
stratum $Z_{f, \Theta}$, where $\Theta$ is a $1$-dimensional face of $\Delta$ 
consists of $d(\Theta)$ rational double  points  of type 
$A_{d(\Theta^*) -1}$. Therefore, over each of these 
points, there are $d(\Theta^*) -1$ exeptional divisors on 
${\overline{Z}}_{\rm crep}$. Thus, $Y$ consists of 
 \[ \sum_{{\rm dim}\, \Theta =1} d(\Theta) (d(\Theta^*) - 1) \]
 exeptional divisors and 
 $\nu_0(\Delta^*) = \nu_2(\Delta)$ divisors corresponding 
 to vertices of $\Delta^*$.  

For the  $K3$-surface 
$\hat{Z}_f$, one has 
\[ {\rm dim}\, H_c^{1,1}(\hat{Z}_f) = 20 .\]

By the result of Danilov and Khovanski\^i 
(\cite{dan.hov}, 5.11) for the case $n =3$
\[ {\rm dim}\, H^{1,1}(Z_f)  = \nu_0(\Delta)  + 
\sum_{{\rm dim}\, \Theta =1} (d(\Theta)  - 1) - 3. \]

It remains to apply to $\Delta$ the Euler formula $\nu_0(\Delta) 
 + \nu_2(\Delta) - \nu_1(\Delta) = 2$. 
This completes the proof. 
\bigskip

\begin{rem} {\rm 
Note that there exists much more simple proof of \ref{24}. We can calculate 
the Euler characteristic of the smooth $K3$-surface $\hat{Z}_f$ using 
the same method as in \ref{new.euler}. By \ref{euler.zero}, the 
Euler characteristic 
of ${Z}_f^{[1]}$ is zero. On the other hand, the Euler characteristic 
of fibers of the minimal resolution $\hat{Z}_f \rightarrow \overline{Z}_f$ 
over $\overline{Z}_f \setminus {Z}_f^{[1]}$ is equal to 
\[ \sum_{{\rm dim}\, \Theta =1} d(\Theta) d(\Theta^*). \]
Since $e(\hat{Z}_f) = 24$, we again obtain the statement of theorem \ref{24}. }
\end{rem}
\bigskip

\begin{rem} {\rm 
 Let $\hat{Z}_f$ and $\hat{Z}_g$ are two mirror symmetric $K3$-surfaces. 
Denote by $E(f)$ (resp. by  $E(g)$) the sublattice in $H^2(\hat{Z}_f)$ 
(resp. in $H^2(\hat{Z}_f)$ ) generated by classes of algebraic curves 
in $\hat{Z}_f \setminus Z^{[1]}_f$ (resp. in 
$\hat{Z}_g \setminus Z^{[1]}_g$). It follows from our first proof of \ref{24} 
that 
\[ {\rm rk}\, E(f) = \nu_0(\Delta) - 3 + 
\sum_{{\rm dim}\, \Theta =1} d(\Theta) (d(\Theta^*) - 1) = \nu_0(\Delta) + 
21 - \sum_{{\rm dim}\, \Theta =1} d(\Theta). \]
Applying the same arguments to the dual polyhedron $\Delta^*$, we 
obtain 
\[ {\rm rk}\, E(g) = \nu_0(\Delta^*) + 
21 - \sum_{{\rm dim}\, \Theta^* =1} d(\Theta^*). \]
By Euler formula for $2$-dimensional sphere, $\nu_0(\Delta) + 
\nu_0(\Delta^*) - 2 = \nu_1(\Delta) = \nu_1(\Delta^*)$. Thus, 
\[ {\rm rk}\, E(f) + {\rm rk}\, E(g) = 
44 - \sum_{{\rm dim}\, \Theta =1} (d(\Theta) + 
d(\Theta^*) -1). \]
Applying the inequality 
\[ \sum_{{\rm dim}\, \Theta =1} (1 - d(\Theta)) (d(\Theta^*) - 1) \leq 0, \]
we obtain
\[ {\rm rk}\, E(f) + {\rm rk}\, E(g) \leq 20. \]
Moreover, the equality holds if and only if for every $1$-dimensional face 
$\Theta \subset \Delta$ one has $d(\Theta) =1$, or   $d(\Theta^*) =1$. }
\end{rem}
\medskip

\begin{rem}
{\rm We consider below an example of two dual $3$-dimensional reflexive 
polyhedra such that two lattices $E(f)$ and $E(g)$ 
admit an interpretation  as orthogonal sublattices in the $20$-dimensional 
quadratic ${\bf Z}$-module 
$E_8 \oplus E_8 \oplus H \oplus H$, i.e., in this example we obtain the 
same property  as in (iv) for the strange duality of Arnold. We note 
that for an arbitary pair of dual $3$-dimensional polyhedra $\Delta$ and 
$\Delta^*$ the property ${\rm rk}\, E(f) + {\rm rk}\, E(g) = 20$ does not 
hold. So we obtain a similarity with Arnold's duality not for all 
families of $K3$-surfaces associated with $3$-dimensional reflexive 
polyhedra.}
\label{orthogonal}
\end{rem}

\begin{exam}
{\rm Let $\Delta$ be the reflexive polyhedron 
\[ \Delta = \{ x =(x_1, x_2, x_3 \in M_{\bf Q} \cong  
{\bf Q}^3 \mid  \; -1 \leq x_i \leq 1, \;
i =1,2,3 \} \]
This polyhedron defines the  toric Fano $3$-fold ${\bf P}_{\Delta} = 
{\bf P}_1 \times {\bf P}_1 \times {\bf P}_1$ with the ample anticanonical 
sheaf ${\cal O}(2,2,2)$. A general $K3$-surface in 
${\bf P}_{\Delta}$ has the Picard number $\rho_{\Delta} = 
3$. The N\'eron-Severi lattice ${\bf E}= E(f)$ of general 
$K3$-surface $\hat{Z}_f$ obtained from  the  $17$-dimensional 
family ${\cal F}(\Delta)$  has 
signature $(1,2)$ and  the  intersection matrix 
\begin{displaymath} \left(
\begin{array}{ccc}
 0 & 2 & 2 \\ 2 & 0 & 2 \\ 2 & 2 & 0  
\end{array} \right)
\end{displaymath}
i.e., the discriminant of the N\'eron-Severi lattice ${\bf E}$ equals $16$.

The dual polyhedron $\Delta^*$ is the convex hull of $6$ points 
$\{ \pm  e_i \}$ 
 $(i =1,2,3)$ in $N_{\bf Q} \cong {\bf Q}^3$, where $e_1, e_2, e_3$  is 
the  standard ${\bf Z}$-basis in ${\bf Z}^3 \cong N$. The polyhedron 
$\Delta^*$ 
defines the singular toric Fano variety ${\bf P}_{\Delta^*}$  
isomorphic to the complete intersection of the following three quadrics in 
${\bf P}_6$ :
\[ v_1v_4 = v_0^2,\;\;  v_2v_5 = v_0^2,\;\; v_3v_6 = v_0^2. \]
A general $K3$-surface $\overline{Z}_g$ in ${\cal F}(\Delta^*)$ has 
$12$ rational 
double singular points at the $0$-dimensional strata $Z_{f, \Theta}$ , where 
$\Theta$ are  $1$-dimensional faces of the cube $\Delta$.  The dimension of 
the  family of $K3$-surfaces in ${\cal F}(\Delta^*)$ equals $3$. 

There exists the  one-parameter subfamily ${\cal F}_s$ 
of this $3$-dimensional family with 
the affine equation  
\[ g_s (X_1, X_2, X_3) = X_1 + X_1^{-1} + X_2 + X_2^{-1} + X_3 + 
X_3^{-1} - s = 0, \]
where 
\[ X_1 = v_1/v_0 = v_0/v_4,\;\; X_2 = v_2/v_0 = v_0/v_5, \;\;
X_3 = v_3/v_0 = v_0/v_6. \]
The family ${\cal F}_s$ is famous due to its relation to Ap\'ery's proof 
of the irrationality of $\zeta(3)$ and Fermi-surfaces in physics \cite{pet}. 
 
C. Peters and J. Stienstra proved (\cite{pet}, Prop.1)  that the 
lattice ${\bf E'} = E(g)$ generated by $20$ divisors on $K3$-surfaces   
$\hat{Z}_g$ obtained from the the family ${\cal F}(\Delta^*)$ has rank $17$ 
(these $20$ divisors consist of $12$ exceptional 
rational $(-2)$-curves and $8$ curves which are 
restrictions from ${\bf P}_{\Delta^*}$ 
of toric  divisors corresponding to $8$ vertices of $\Delta$),  
signature $(1,16)$  and discriminant $16$.  Moreover, 
the orthogonal complement to ${\bf E'}$ in 
the lattice $E_8 \oplus E_8 \oplus H \oplus H \oplus H$ is isomorphic to 
$H \oplus {\bf E}$.

The lattice ${\bf E'}$ is the 
N\'eron-Severi lattice of the  minimal desingularization  of general 
$K3$-surfaces in the $3$-dimensional family ${\cal F}(\Delta)$. 
$K3$-surfaces in the one-parameter subfamily  ${\cal F}_s$ have more algebraic 
cycles forming the N\'eron-Severi lattice of rank $\geq 19$. }
\end{exam}
\medskip

The above example  shows that the Picard-Fuchs 
equation on transcendental periods of one-parameters family ${\cal F}_s$ 
of $K3$-surfaces is analogous to the Picard-Fuchs equation for the 
one-parameter mirror family ${\tilde {\cal W}}_{\psi}$ of $3$-dimensional 
Calabi-Yau varieties considered by P. Candelas et al. in \cite{cand2}. 
Indeed, in the both cases we have a reflexive polyhedron $\Delta$ and 
a  one-parameter family of Laurent polynomials 
\[ f_a(X) = \sum_{m \in \partial \Delta} X^m  + a, \]
where $a$ is a complex parameter. This complex parameter is the coefficient 
which corresponds to the unique integral point in the interior of $\Delta$. 
The value $a = \infty$ corresponds to the reducible Calabi-Yau variety 
${\bf P}_{\Delta} \setminus {\bf T}$. The point $a= \infty$ is 
the "large complex structure limit" (in the case of 
quintic mirrors ${\cal F}_{\infty}$ is the union of $5$ projective spaces).

It was proved in \cite{bat.var} the Yukawa coupling on the Hodge spaces  
$H^{2,1}$  of the family ${\cal F}_a$ degenerates to the topological 
cup-product on the mirror family of Calabi-Yau 
hypersurfaces in ${\bf P}_{\Delta^*}$. One can also check that 
the monodromy operator around 
$a =\infty$ is maximally unipotent.

\end{document}


\date{July 11, 1993}

\maketitle

\footnote{Supported by DFG, Forschungsschwerpunkt  Komplexe 
Mannigfaltigkeiten.}

\thispagestyle{empty}

\newpage

\thispagestyle{empty}
The Habilitationsschrift consists of two parts: 

\bigskip

\bigskip

\bigskip

{\bf Part II:  Variation of Mixed Hodge Structures of Affine Hypersurfaces 
in Algebraic Tori (pages 66-124)} 
\bigskip

\begin{center}

{\large  {\bf Abstract}}
\end{center}
\bigskip

Let $Z_f$ be an affine $\Delta$-regular hypersurface in ${\bf T} = 
({\bf C}^*)^n $ defined by a Laurent 
polynomial $f \in {\bf C} \lbrack X_1^{\pm 1}, \ldots, X_n^{\pm 1} \rbrack$ 
with the Newton polyhedron $\Delta \subset {\bf R}^n$. 
The paper is devoted to an explicit calculation of the 
mixed Hodge structure on cohomology of $Z_f$ by means of the Poincare 
residues  of rational differential forms on the complementary open set 
${\bf T} \setminus Z_f$. We extend the  method  of N. Katz and P. Griffiths 
for  the calculation of the Gauss-Manin connection and 
the infinitesimal period mapping for hypersurfaces in projective spaces as 
well as the method of J. Steenbrink and I. Dolgachev for hypersurfaces in 
weighted projective spaces. 

Our results are applied to the study of the mirror symmetry  
for Calabi-Yau hypersurfaces in toric varieties.

\thispagestyle{empty}


\setcounter{page}{66}

\date{}

\maketitle




%


\thispagestyle{empty} 

\newpage

\tableofcontents

\thispagestyle{empty}

\newpage

\pagestyle{myheadings}

\section{Introduction}

\hspace*{\parindent}


The study of infinitesimal variations of Hodge structure on 
cohomology of algebraic varieties developed by  J. Carlson, 
M. Green, P. Griffiths and J. Harris has been originated from the 
fundamental paper  of Griffiths on periods of complex hypersurfaces in 
projective spaces  \cite{griff1}. 
The  paper of Griffiths contains  a very fruitful appproach to the 
description of the Hodge structure 
on the primitive part $PH^{n-1}$ of  the  middle cohomology group  
of projective hypersurfaces $\overline{Z}_f \subset {\bf P}^n$ 
defined  by  homogeneous polynomials $f(Y) \in {\bf C} \lbrack Y_0,Y_1, 
 \ldots, Y_n \rbrack$. 
Griffiths has proved the existence of  an  isomorphism between 
the quotients of the Hodge filtration on  $PH^{n-1}(\overline{Z}_f)$ 
and the homogeneous components of 
an  Artinian graded ring $R_{f}$. The ring $R_f$ is called a {\em Jacobian 
ring}. It is a quotient of the polynomial 
ring ${\bf C} \lbrack Y_0,Y_1,  \ldots, Y_n \rbrack$ by the {\em Jacobian 
ideal} $J_f$  which is generated 
by a regular sequence consisting of  partial derivatives  $\partial f(Y)/ 
\partial Y_i$, $ 0(\leq i \leq n)$. 
\medskip

One should remark that the similar idea was also earlier 
discovered by Dwork in his $p$-adic proof of the Weil's  conjectures for 
hypersurfaces in projective spaces \cite{dwork1,dwork2}. 
Instead of projective hypersurfaces, Dwork has been considered only 
{\em affine parts}  $Z_f$ of  $\overline{Z}_f$ :
\[ Z_f  = \{ (Y_0 : \ldots : Y_n ) \in \overline{Z_f} \mid Y_0 \cdots Y_n 
\neq 0 \}, \]
i.e., affine hypersurfaces in $n$-dimensional algebraic torus.  
The space of solutions of $p$-adic differential equations 
investigated by Dwork 
can be represented by  the  quotient of the $d$-th Veronese subring 
\[{\bf C} \lbrack Y_0,Y_1,  \ldots, Y_n \rbrack_{(d)} \subset 
{\bf C} \lbrack Y_0,Y_1,  \ldots, Y_n \rbrack, \;\; d = {\rm deg}\, f \] 
by a regular sequence  consisting of 
{\em logarithmic partial derivatives}  $Y_i \partial f(Y)/ 
\partial Y_i$ $( 0 \leq i \leq n)$ of  $f$.   
\medskip

It was shown by N. Katz in \cite{katz} 
that the $p$-adic differential equations considered by Dwork coincide with   
the  Picard-Fuchs equations satisfied by periods of differential forms on 
the complex hypersurfaces. 

The fact that the  Picard-Fuchs equations appear 
simultaniously in $p$-adic analysis and in complex analitic 
geometry gives evidence of the unity of mathematics. The reason for it 
is the existence of the Gauss-Manin connection on cohomology 
of families of algebraic varieties. 
A simple example explaining this discovery of  Manin for the case  of 
the one-parameter family of elliptic curves  one can find in the 
beautiful introductory book of Clemens \cite{clem}. 
\medskip

The  Dwork-Griffiths method  was later 
extended by Steenbrink \cite{steen} and 
Dolgachev \cite{dolgach} to the case of  hypersurfaces in weighted projective 
spaces. 
\bigskip

The method of Dwork-Griffiths-Katz and Dolgachev-Steenbrink caught recently 
much attention from physicists because of 
the  discovered relation between Jacobian rings and chiral rings  
of primary fields in Landau-Ginzburg theories 
\cite{greene1,lerche}. On the other hand, this method 
is very useful for calculations  of periods of Calabi-Yau varieties 
\cite{cand2,font,klemm,mor2}. These periods are important due to  
the mirror symmetry \cite{aspin1,aspin3,mor1} (one attractive  
property of the mirror symmetry is connected with 
a counting of rational curves on Calabi-Yau $3$-folds). 
The polynomial representation of cohomology classes 
in the middle cohomology group of hypersurfaces allows to calculate 
explicitly the Picard-Fuchs equations, or in fact the Gauss-Manin 
connection for families of hypersurfaces. These differential equations 
seems to contain a sufficient information to describe some 
superconformal theories. 
We remark that  there exists  a general phylosophy  that 
$N =2$ supergravity and superconformal models can be identified with 
special variations of abstract mixed Hodge structures \cite{ceco}.  
\medskip

Projective spaces and weighted projective spaces are special cases of 
so called {\em toric varieties} \cite{dan2,oda1}. 
It is natural to ask whether the method of Jacobian 
rings for the calculation of the Hodge structure   can be also 
extended to the case of hypersurfaces in projective toric varieties. 

There is one more strong  reason why the calculation of the infinitesimal 
variation of the Hodge structure of hypersurfaces in {\em toric 
varieties} is especially important for physicists. 
As we mentioned, the most interesting class of variations of  
Hodge structure in physics is connected with families of Calabi-Yau 
varieties \cite{gepn}. First, it is possible to construct many examples of 
families of Calabi-Yau hypersurfaces in toric 
varieties ${\bf P}_{\Delta}$ associated with  
{\em reflexive polyhedra} $\Delta$ \cite{bat.mir}. Secondly, 
for Calabi-Yau hypersurfaces in toric varieties we can use a {\em precise 
mathematical definition of  mirrors} via the notion of dual polyhedron. 
Thus, appying the mirror symmetry 
conjecture, we can reduce the calculation of the A-variation of 
the "quantum Hodge structure" to the calculation of the 
$B$-variation of the classical Hodge structure \cite{witt}. 
This main idea can be also  
explained for families of Calabi-Yau $3$-folds as follows: 
\medskip

There are two kind of Yukawa couplings on the Hodge spaces 
$H^{2,1}(W)$ and 
$H^{1,1}(W)$ considered by physicists for 
Calabi-Yau $3$-folds $W$.  It is known that the  
Yukawa coupling on $H^{2,1}(W)$ is defined by 
the deformation of the complex structure, but the definition of the 
Yukawa coupling on $H^{1,1}(W)$ is outside the classical algebraic 
geometry. Physicists 
have discovered a trick, so called {\em mirror symmetry}, which allows to 
replace the calculation of the "misterious Yukawa coupling" on $H^{1,1}(W)$ 
by  the calculation of the Yukawa coupling on the space $H^{2,1}(W')$ 
corresponding to the deformation of the complex structure on 
{\em  another Calabi-Yau $3$-fold} $W'$. The Calabi-Yau 
variety $W'$ is called {\em mirror of} $W$. Since the mirror 
mapping for Calabi-Yau hypersurfaces in toric varieties admits a  
precise mathematical definition \cite{bat.mir}, we can compute both  
of the Yukawa couplings using only  deformations of the classical 
Hodge structures (cf. \cite{cand2}). On should remark that in general case the 
construction of mirrors for Calabi-Yau manifolds is rather difficult 
problem \cite{witt}, and Calabi-Yau hypersurfaces in toric varieties are now 
the only known cases in which this problem is solved by mathematical 
methods.  
\bigskip

Denote by ${\bf T}$ the $n$-dimensional algebraic torus $({\bf C}^*)^n$. 
The main purpose  of this  paper is to give the explicit decription of the 
mixed Hodge structure of affine $\Delta$-regular hypersurface $Z_f \subset 
{\bf T}$ defined by a Laurent polynomial $f$ with the Newton polyhedron 
$\Delta$.  Note that our definition of $\Delta$-regularity is an  
extension of the one used by Dwork in \cite{dwork1,dwork2}. 
Our calculations of 
the Hodge structure were mostly guided by the paper of N. Katz 
\cite{katz}. This paper of N. Katz 
implicitly contains the notion of mixed Hodge structure, 
e.g.,  the subspace in $H^{n-1}(Z_f)$ spanned by differential of 
the second kind  is 
the minimal nonzero weight subspace ${\cal W}_{n-1}$ of the mixed Hodge 
structure on $H^{n-1}(Z_f)$. 
\medskip

As an application of our calculation  of the infinitesimal variations 
of mixed Hodge structures, we show {\em under some restrictions}  
the following property of the mirror 
mapping for Calabi-Yau varieties (cf. theorem \ref{l.limit}): 

\begin{theo}
Topological cup-product on the cohomology ring 
\[ \bigoplus_{i=0}^{n-1}  H^{i,i}(W) \]
of smooth Calabi-Yau hypersurfaces in  toric Fano $n$-folds 
can be obtained as a "limit" of the multiplicative structure on the 
chiral ring 
\[ \bigoplus_{i=0}^n  H^{n-1-i,i}(W') \]
of the {\em mirrors} $W'$,  where the "limit" corresponds to some 
degeneration of $W'$ at a boundary point of the moduli space of complex 
structures on $W'$.
\end{theo}

This property of the mirror mapping agree with predictions of physicists who 
consider multiplicative structure on the chiral rings of mirrors $W'$ as 
topological  cup product on even-dimensional cohomology of $W$ with 
"instanton corrections". This corrections  vanish at some boundary point 
("large radius limit") of the moduli space of mirrors  $W'$ 
(cf. cite{aspin2}).  
\medskip

Let us shortly describe the content of the paper.
\medskip

In sect. 2 we give a short review of already known results from 
combinatorics relating to the Ehrhart polynomials of integral convex 
polyhedra and toric varieties. 
Our main interest will be connected with so called $\delta$-vector of 
the convex integral polyhedron $\Delta$ studied  by Stanley and 
Hibi \cite{hibi1,stan0}. 
\medskip

In sect. 3 we recall basic theorems due to Bernstein, Kouchnirenko, Danilov, 
and Khovanski\^i on cohomological  invariants of affine hypersurface $Z_f$ 
in algebraic torus ${\bf T}$ in terms of the Newton polyhedron  $\Delta$  
of the  Laurent polynomial  $f$.  Using the results of Danilov  and 
Khovanski\^i, we give the cohomological interpretation of the $\delta$-vector 
of $\Delta$ in terms of  the Hodge filtration on the primitive part of 
$H^{n-1}(Z_f)$. 
\medskip

In sect. 4 we discuss different  equivalent characterizations of 
$\Delta$-regularity of a Laurent polynomial $f$ with a fixed Newton 
polyhedron $\Delta$. The notion of $\Delta$-regularity of a polynomial $f$ has 
already appeared as a natural technical condition in the calculation  of 
Dwork of the  Weil's zeta-function \cite{dwork1,dwork2}. On the other hand, 
this nondegeneracy condition 
has been used in the theory of singularities (see, for example, 
the calculation of Varchenko \cite{varch} of the zeta-function of 
the monodromy operator). Recently 
Gelfand, Kapranov, and  Zelevinskij \cite{gel1}
have found the definition of $\Delta$-regularity  
via the notion of principal $A$-discriminant  \cite{gel0}. 
For $\Delta$-regular polynomials $f$ we define the corresponding 
{\em Jacobian ring} $R_f$. 
\medskip

In sect. 5 we compare cohomology groups of $Z_f$ and the complementary open 
subset ${\bf T} \setminus Z_f$. We  show that the primitive 
cohomology groups of affine hypersurface $Z_f$  and the complementary 
open subset ${\bf T} \setminus Z_f$ are isomorphic.  
Our main purpose in the sequel 
is to represent cohomology classes in $H^{\cdot}(Z_f)$ 
by residues of rational differential forms on ${\bf T}$ with logarithmic 
poles along $Z_f$. 
\medskip

In sect. 6 we prove that  the quotient spaces of the Hodge filtration on  
the {\em primitive part} of $H^{n}(Z_f)$ are are canonically 
represented by homogeheous components of the Jacobian ring $R_f$. 
\medskip

In sect. 7 we calculate cohomology groups 
of affine hypersurfaces $Z_f$  by method of Dwork and Katz. 
We construct two quasi-isomorphic 
Koszul complexes which are correspondingly the global de Rham complex 
of differential forms on ${\bf T} \setminus Z_f$ and the 
global logarithmic de Rham complex 
relative to a partial affine compactification $U_0$ of 
${\bf T} \setminus Z_f$. 
\medskip

In sect. 8 we  describe explicitly the space 
$H^n({\bf T} \setminus Z_f)$ with the  Hodge and weight filtrations. 
\medskip

In sect. 9 we study the pure Hodge structure arising from 
the subspace of elements of weight $n-1$ in $H^{n-1}(Z_f)$. 
As application, we give another proof of some 
results of Danilov and Knovansi\^i \cite{dan.hov}. 
We  calculate the nondegenerate pairing 
between components of the pure Hodge structure in terms of 
the algebraic multiplication law,  give also the precise 
description of the Hodge-Deligne space $H^{n-2,1}(Z_f)$.   
\medskip

In sect. 10 we study the geometry of the moduli spaces 
${\cal M}({\Delta})$ of affine hypersurfaces $Z_f$ 
where $f$ has a fixed Newton polyhedron $\Delta$. The 
moduli spaces ${\cal M}({\Delta})$ are again projective toric varieties, they  
can be constructed   by means of geometric invariant theory \cite{mam}. 
On the other hand, in the category of projective toric varieties 
there exists the   
{\em universal moduli space} $\tilde{\cal M}(\Delta)$ such that the 
principal $A$-discriminat ${\rm Disc}_{\Delta}(f)$ of  Laurent 
polynomials $f$ is  represented 
by a section of an ample line bundle on $\tilde{\cal M}(\Delta)$.   
The toric variety $\tilde{\cal M}(\Delta)$ is associated with so called 
{\em secondary polyhedron}.  
\medskip

In sect. 11 we calculate the Gauss-Manin connection on the groups  
$H^{n}({\bf T} \setminus Z_f)$ which can be considered as a local system 
on ${\cal M}({\Delta})$. In particular, this gives an obtain algorithm for 
calculation of  Picard-Fuchs differential equations. We  
describe also the infinitesimal mixed period mapping corresponding to the 
variation  of the Hodge structure on $H^{n-1}(Z_f)$. 
\medskip

In sect. 12 we consider several applications of our results  to variations 
of complex structure of {\em Calabi-Yau hypersurfaces in toric varieties} 
defined by  Laurent polynomials with a reflexive Newton polyhedron $\Delta$.   
First, using the method of Griffiths,  we prove the infinitesimal mixed 
Torelli theorem for Calabi-Yau hypersurfaces. Then we show  (under some 
restrictions) that the Jacobian ring $R_f$ can be "deformed" 
to the Chow ring of the toric Fano variety ${\bf P}_{\Delta^*}$ 
corresponding to the dual polyhedron $\Delta^*.$ 
Since the dual reflexive polyhedron $\Delta^*$ defines the mirror  
family of Calabi-Yau hypersurfaces in ${\bf P}_{\Delta^*}$ \cite{bat.mir}, 
this shows some properties of the mirror mapping. 
\medskip

In sect. 13 we study the geometry of the moduli space ${\cal M}(\Delta)$ 
of the family of affine Calabi-Yau hypersurfaces defined by a reflexive 
polyhedron $\Delta$. For that purpose we use the Oda-Park conic version of 
the combinatorial construction of the secondary polyhedra 
\cite{oda.park}. We show that geometry of toric singularities arising 
at the boundary of the moduli space ${\cal M}(\Delta)$ can be described 
in terms of subcones in the cone of effective divisors of mirrors.  
This leds us to the conjecture about a relation 
between the geometry of the cone of effective divisors on Calabi-Yau 
varieties and singularities at  cusp points  on the moduli space of their 
mirrors.  
\medskip
 
In sect. 14 we show that the variation of mixed Hodge structure of 
affine hypersurfaces $Z_f \subset {\bf T}$ can be described also by means of 
hypergeometric functions introduced by Gelfand, Kapranov and Zelevinsky 
\cite{gel0}. Considering affine hypersurfaces $Z_f$ over fields of positive 
characteristic, we come to $p$-adic hypergeometric functions investigated 
by Dwork in the connection with  the Weil's  conjectures.   
\medskip

We hope that the $p$-adic hypergeometric functions can be also applied to  
non-arhimedian physics. On the other hand, the famous 
"trick of S. Mori in positive characteristic" shows that a  
consideration of families of 
algebraic varieties over fields of  characteristic $p >0$ 
is  very powerful method, and it  might happens that it  can be 
succesfully applied  not 
only to the study of rational curves on Fano varieties \cite{mori}, 
but also to the study of rational curves on Calabi-Yau manifolds.   
\bigskip

{\bf Acknowledgements.}    
I wish to thank Y. Kawamata, J. Stienstra and E. Viehweg for their 
stimulating advices   
during my preparation of the manuscript. 
It is a pleasure to acknowledge also very 
helpful discussions on the  mirror symmetry with Yu. Manin and D. Morrison.  
\bigskip

\section{Ehrhart polynomials and toric varieties}

\hspace*{\parindent}


Let $\Delta$ be a convex $n$-dimensional polyhedron in ${\bf R}^n$ such that 
all vertices of $\Delta$ belong to the integral lattice ${\bf Z}^n \subset 
{\bf R}^n$. For any positive  integer $k$, denote by $l(k\Delta)$ the number 
of integral points in $k\Delta \cap {\bf Z}^n$. Put $l(0\Delta) =1$. The  
following properties of the numbers $l(k\Delta)$ were proved by  E. Ehrhart  
\cite{ehrhart1}: 
 
\begin{theo}
There exists  a polynomial $\Lambda(t)$ of degree $n$ such that

{\rm (i)} for any non-negative integer $k$, one has 
\[ \Lambda(k) = l(k\Delta) ;\]

{\rm (ii)} for any positive integer $k$, one has   
\[ (-1)^n \Lambda(-k) = l^*(k\Delta) , \]
where $l^*(k\Delta)$ is the number of integral points in the interior of 
$k\Delta$; 

\label{ehr}
\end{theo}

\begin{rem}
{\rm  It is easy to see that the leading coefficient  of $\Lambda(t)$ equals 
${\rm vol}\, \Delta$, the volume of $\Delta$ relative to the standard Haare 
measure on ${\bf R}^n$. }
\end{rem}

The polynomial $\Lambda(t)$ is called the {\em Ehrhart polynomial of 
$\Delta$}. This polynomial 
 can be studied  using the technique of commutative algebra and 
 algebraic geometry \cite{stan1,hibi1}. The main idea of these  
 algebraic methods is  to compare the Ehrhart polynomial 
 with the  Hilbert polynomial of the 
graded 
commutative algebra $S_{\Delta}$ associated with $\Delta$. The graded 
algebra $S_{\Delta}$ is constructed as follows.

\begin{opr}
{\rm Let ${\bf L} = {\bf C} \lbrack X_1 ^{\pm 1}, \ldots,  
X_n^{\pm 1} \rbrack$ be the Laurent polynomial 
ring in $n$ variables over the complex number field ${\bf C}$. 
Consider the polynomial extension 
\[ {\bf L} \lbrack X_0 \rbrack = {\bf C} \lbrack X_0, X_1 ^{\pm 1}, \ldots,  
X_n^{\pm 1} \rbrack \]
of ${\bf L}$, and define the ring $S_{\Delta}$ as the 
subalgebra of ${\bf L}\lbrack X_0 \rbrack$ generated as 
${\bf C}$-vector space by elements 
of ${\bf C}$ and all monomials $X^k_0 X^m = X_0^kX_1^{m_1} \cdots X_n^{m_n}$ 
such that the rational point 
$m/k = (m_1/k, \ldots , m_n/k )$ belongs to $\Delta$. }
\end{opr}

The standard grading of the polynomial ring ${\bf L}\lbrack X_0 \rbrack$ 
induces  the grading of $S_{\Delta}$:  
\[ {\rm deg}\, X^k_0X^m =k.  \]
We now can associate with the polyhedron $\Delta$ the projective algebraic 
variety ${\bf P}_{\Delta}$. 

\begin{opr}
{\rm The $n$-dimensional projective algebraic variety 
\[ {\bf P}_{\Delta} = {\rm Proj}\, S_{\Delta} \]
we call {\em the toric variety associated with the polyhedron} $\Delta$. 
Let  $S_{\Delta}(-1)$ denotes the free graded $S_{\Delta}$-module with a 
homogeneous generator of degree $-1$. We denote by ${\cal O}_{\Delta}(1)$ the 
corresponding ample invertible sheaf on ${\bf P}_{\Delta}$. 
Put ${\cal O}_{\Delta}(k) = {\cal O}_{\Delta}(1)^{\otimes k}$. }
\end{opr}

For $k \geq 0$, the number $l(k\Delta)$ is  equal to the dimension of the 
$k$-th homogeneous component $S^k_{\Delta}$ of the ring $S_{\Delta}$, or 
the dimension of the space $H^0({\bf P}_{\Delta}, {\cal O}_{\Delta}(k))$ 
of global sections  of ${\cal O}_{\Delta}(k)$.  
By general theorem on Hilbert functions of graded commutative algebras 
\cite{atiyah}, there are integers 
$\lambda_i$ such that  for $k \gg 0$
\[ \Lambda (k)  = \lambda_0 {k \choose n} + 
\lambda_1 {k \choose n-1}  + \cdots + \lambda_n {k \choose 0}, \]
where
\[ { z \choose r } = \frac{1}{r!} z(z-1) \cdots (z - r + 1) . \]
\bigskip

The number 
\[ \lambda_0 = n! \, ({\rm the \; leading\;  
coefficient\; of } \; \Lambda(t)) = 
n! {\rm vol} \, \Delta \]
is called the {\em degree} of the projective variety ${\bf P}_{\Delta}$ 
\cite{hartsh}.

\begin{opr}
{\rm The positive integer  $n!\, {\rm vol}\, \Delta$ 
we call the {\em degree} of the 
polyhedron $\Delta$ and denote  by $d(\Delta)$.}
\end{opr}

The statement \ref{ehr}(i) follows from the following 
two facts.   

\begin{theo} {\rm (Serre)} 
The Euler characteristic 
\[ \chi({\bf P}_{\Delta}, {\cal O}_{\Delta}(k)) \]
is always a polynomial in $k$ for any $k \in {\bf Z}$. 
\end{theo}

\begin{theo} 
{\rm \cite{dan2}} 
Cohomology groups of invertible sheaves ${\cal O}_{\Delta}(k)$ on the 
toric variety ${\bf P}_{\Delta}$ satysfies the following conditions 
\[ H^n({\bf P}_{\Delta}, {\cal O}_{\Delta}(k) ) = 0\;\;{\rm  for}\;\;  
k \geq 0, \]
 and 
\[ H^i({\bf P}_{\Delta}, {\cal O}_{\Delta}(k) ) = 0, \;\; 
{\rm for}\;  0 < i < n, \; k \in {\bf Z}. \]
\label{vanish}
\end{theo}

Define now  special  ideals in $S_{\Delta}$. 

\begin{opr}
{\rm Denote by $I^{(l)}_{\Delta}$ $(0 \leq l \leq n+1)$ the 
homogeneous ideal in 
$S_{\Delta}$ generated as the ${\bf C}$-subspace by all 
monomials $X^k_0 X^m$ such 
that the rational point  $m/k = (m_1/k, \ldots , m_n/k )$ belongs to $\Delta$, 
but does not belong to any face $\Delta' \subset \Delta$ of codimension $l$. 
 }
\end{opr}

So we obtain  the increasing chain of homogeneous ideals in $S_{\Delta}$ 
\begin{equation}
 0 = I^{(0)}_{\Delta} \subset I^{(1)}_{\Delta} \subset \cdots \subset 
I^{(n-1)}_{\Delta} \subset I^{(n)}_{\Delta} \subset I^{(n+1)}_{\Delta} 
= S_{\Delta}^+ ,  
\label{chain}
\end{equation}
where $S_{\Delta}^+$ is the maximal 
homogeneous ideal in $S_{\Delta}$.  

\begin{opr}
{\rm Denote by  ${\omega}_{\Delta}$  the coherent 
sheaf on ${\bf P}_{\Delta}$ corresponding to  the 
homogeneous ideal $I^{(1)}_{\Delta}$ considered as a graded 
$S_{\Delta}$-module.}
\end{opr}

The statement \ref{ehr}(ii) follows from the 
Serre duality and  the 
following theorem contained in \cite{dan2}:

\begin{theo}
A toric variety ${\bf P}_{\Delta}$ is an  arithmetically 
Cohen-Macaulay projective algebraic variety of dimension $n$, i.e.,  
the homogeneous coordinate ring $S_{\Delta}$ is a 
Cohen-Macaulay ring of dimension $n+1$. Moreover, $\omega_{\Delta}$ is the 
dualizing $($or canonical$)$ sheaf  on ${\bf P}_{\Delta}$, and the 
graded module 
\[ \oplus_{k \geq 0} H^0({\bf P}_{\Delta}, \omega_{\Delta}\otimes 
{\cal O}_{\Delta}(k)) \]
is isomorphic to the ideal $I^{(1)}_{\Delta}$. 
\label{C.-M.}
\end{theo}

The above duality can be reformulated also via the Hilbert-Poincare functions 
of the graded modules $S_{\Delta}$ and $I^{(1)}_{\Delta}$ as follows. 

\begin{theo} 
Let 
\[ P_{\Delta}(t) = \sum_{k \geq 0} l(k\Delta) t^k  \] 
 the  Hilbert-Poincare series of the graded algebra $S_{\Delta}$,  
\[ Q_{\Delta} (t) =  \sum_{k > 0} l^* (k\Delta) t^k   \]
the Hilbert-Poincare series of the graded $S_{\Delta}$-module 
$I^{(1)}_{\Delta}$. Define  
two functions 
\[ \Psi_{\Delta} (t) = \sum_{i \geq 0} \psi_i (\Delta) t^i =  
(1 -t )^{n+1}  P_{\Delta} (t)  \]
and 
\[ \Phi_{\Delta} (t) = \sum_{i \geq 0} \varphi_i (\Delta) t^i =  
(1 -t )^{n+1}  Q_{\Delta} (t).  \]
Then $\Psi_{\Delta}(t)$ and $\Phi_{\Delta}(t)$ are polynomials satisfying 
the condition 
\[ t^{n+1} \Psi_{\Delta}(t^{-1}) = \Phi_{\Delta}(t), \] 
i.e., $\varphi_{i}(\Delta)  = \psi_{n+1 -i} (\Delta)$ $( 0 \leq i \leq n+1)$. 
Moreover,   
\[ \Psi_{\Delta} (t)  = \lambda_0 t^n + \lambda_1 t^{n-1}(1 -t) + 
\cdots + \lambda_n (1 - t)^n  \] 
is a polynomial  of degree $\leq n$ with non-negative 
integral coefficients such that 
\[ \Psi_{\Delta} (1)  = \sum_{i \geq 0} \psi_i(\Delta) = \lambda_0 = 
d(\Delta), \;\; \;\;\Psi_{\Delta}(0) = \psi_0(\Delta) = \lambda_n = 1,  \] 
and 
\[ \psi_1(\Delta) = l(\Delta) - n -1. \]
\label{hilb}
\end{theo}

\begin{exam}
{\rm  Let $\Delta$ be an arbitrary  $2$-dimensional integral polyhedron, then  
$\psi_0 (\Delta) = \varphi_3(\Delta)= 1$, 
$\psi_1 (\Delta) = \varphi_2(\Delta) = l(\Delta) -3$, 
$\psi_2 (\Delta) = \varphi_1(\Delta)= l^*(\Delta)$. The relation 
\[ d(\Delta) = 2! \,{\rm vol} \, \Delta = \psi_0(\Delta) + 
\psi_1(\Delta) + \psi_2(\Delta) = l(\Delta) + l^*(\Delta) -2 \]
is exactly the Pick's formula for  $\Delta$.}
\end{exam}

\begin{rem}
{\rm The integral vector $(\psi_0(\Delta), \psi_1(\Delta), \ldots , 
\psi_n(\Delta) )$ is called also $\delta$-vector of the integral polyhedron 
$\Delta$ (cf. \cite{stan0}). The properties of the $\delta$-vector  
in the theorem \ref{hilb} follow 
from the fact that the polynomial $\Psi_{\Delta}(t)$  is  
the Hilbert-Poincare series of an  Artinian ring 
$S_{\Delta}/ (s_0,s_1, \ldots, s_n)S_{\Delta}$, where  
$s_0,s_1, \ldots, s_n$ is a  regular sequence of elements of degree 1 
in $S_{\Delta}$ \cite{stan00}. The polynomial $\Phi_{\Delta}(t)$ is 
the Hilbert-Poincare series of the Artinian $R_f$-module 
\[ I^{(1)}_{\Delta}/ (s_0,s_1, \ldots, s_n)I^{(1)}_{\Delta}. \]
The existence of the regular sequence 
of global sections 
\[ s_0, s_1, \ldots , s_n \in H^0 ({\bf P}_{\Delta}, {\cal O}_{\Delta}(1)) 
\cong S^1_{\Delta} \]
follows from \ref{C.-M.} and Bertini theorem.}
\label{reg}
\end{rem}

\section{Hodge-Deligne numbers of hypersurfaces in ${\bf T}$}

\hspace*{\parindent}


Recall that ${\bf L}$ denotes  the Laurent polynomial ring 
${\bf C}\lbrack X_1 ^{\pm 1}, \ldots,  X_n^{\pm 1} \rbrack$ 
over the complex number field ${\bf C}$. 
Let $M$ be  the free abelian group of rank $n$. We identify elements 
$m \in M$ with exponents of monomials $X^m \in {\bf L}$.  The Newton 
polyhedron $\Delta(f)$ of a Laurent polynomial $f = f(X) = \sum c_m X^m \in 
{\bf L}$ is the convex hull of the integral 
points $m \in M \subset M \otimes {\bf R}$ such that $c_m \neq 0$. 
We denote by $L(\Delta)$ the space of all Laurent polynomials with the 
Newton polyhedron $\Delta$.  The dimension of $L(\Delta)$ equals $l(\Delta)$. 
Let 
\[{\bf T}: = {\rm Spec}\, {\bf L} \cong ({\bf C}^*)^n. \] 
Every Laurent polynomial $f$ defines  the affine hypersurface 
\[ Z_f = \{ X \in {\bf T}   \mid f(X) = 0 \}. \]

\begin{opr}
{\rm Let $f(X) = \sum c_m X^m  \in L(\Delta)$, $\Delta' \subset \Delta$ an 
$l$-dimensional face of $\Delta$. Define also the Laurent polynomial with 
the Newton polyhedron $\Delta'$. 
\[ f^{\Delta'}(X) = \sum_{m'\in \Delta'} c_{m'} X^{m'}. \]
} 
\end{opr}

\begin{opr}
Let $g = g(X)$ be any Laurent polynomial. Denote by $g_i$ $( 1 \leq i \leq n)$
the logarithmic derivatives of $g$:
\[ g_i(X) = X_i \frac{\partial}{\partial X_i} g(X)\;\; (1 \leq i \leq n). \]
\end{opr}

\begin{opr}
A Laurent polynomial $f \in L(\Delta)$ and the corresponding 
affine hypersurface $Z_f$ are called $\Delta$-regular if  
$\Delta(f) = \Delta$, and for every $l$-dimensional edge 
$\Delta' \subset \Delta$ $(l >0)$ the polynomial equations  
\[ f^{\Delta'}(X) = f_1^{\Delta'}(X) =  \cdots = f_n^{\Delta'}(X) = 0 \] 
have no common solutions in ${\bf  T}$. 
\label{regul2}
\end{opr}

Let us  recall  some already known results on   
the cohomology groups  of an affine hypersurface $Z_f \subset {\bf T}$, 
where $f$ is always assumed 
to be a $\Delta$-regular Laurent polynomial. 
\medskip

It follows from  general properties of smooth affine 
algebraic varieties that 
$H^i (Z_f) = 0$ for all $i > {\rm dim}\, Z_f = n -1$. 
On the other hand, there exists the following toric 
analog of the  Lefschetz theorem for affine hypersurfaces 
$Z_f \subset {\bf T}$.

\begin{theo} 
{\rm (Bernstein-Danilov-Khovanski\^i \cite{dan.hov})} 
The natural homomorphism
\[ H^i( {\bf T}) \rightarrow H^i ( Z_f ) \]
is bijective for $i < n -1$ and injective for $i = n - 1$. 
\label{Lef}
\end{theo}

Let $M_{\bf C} = M \otimes {\bf C}$. Since $H^*({\bf T})$ is 
canonically isomorphic to the 
exterior algebra $\Lambda^* M_{\bf C}$, 
theorem \ref{Lef} allows us to calculate  all cohomology groups  
$H^i(Z_f)$ for   $i < n-1$:

\begin{coro}
For $i < n-1$, one has  
\[ H^{i}(Z_f) \cong \Lambda^i M_{\bf C}. \] 
In particular, 
$ {\rm dim }\, H^{i}(Z_f)  = { n \choose i }$. 
\label{i.dim}
\end{coro} 

In order to compute the dimension of the remaining cohomology 
group $H^{n-1}(Z_f)$,  we need  to know only the 
Euler characteristic of $Z_f$. 

\begin{theo}
{\rm (Bernstein-Khovanski\^i-Kushnirenko \cite{kuch,hov.genus})}
The Euler characteristic of $Z_f$ equals 
\[  e(Z_f) = \sum_{i \geq 0} (-1)^i {\rm dim}\, H^i(Z_f)  = 
(-1)^{n-1} d(\Delta). \]
\end{theo}

Since $e({\bf T}) = 0$,  we obtain.  

\begin{coro}
$ {\rm dim }\, H^{n-1}(Z_f)  = d(\Delta) + n -1$, 
\label{n.dim}
\end{coro} 
\medskip

By general theorem of Deligne \cite{deligne}, the cohomology groups 
$H^m(V)$ $( 0 \leq m \leq {\rm dim}\, V)$ of 
any  smooth affine algebraic variety $V$ have  the mixed 
Hodge structure consisting of the Hodge filtration 
\[  H^m(V) = {\cal F}^{0}H^m(V) \supset {\cal F}^{1} H^m(V) \supset 
\cdots \supset {\cal F}^{m+1}H^m(V) = 0 \]
and the weight filtration
\[ 0 = {\cal W}_{m-1}H^m(V) \subset {\cal W}_{m} H^m(V) \subset 
\cdots \subset {\cal W}_{2m}H^m(V) =  H^m(V). \]
Moreover, for any $i \geq m$ the Hodge filtration induces on the quotient 
space 
\[ {\rm Gr}^{\cal W}_i H^m(V) = {\cal W}_{i}H^m(V) / {\cal W}_{i-1}H^m(V) \] 
a pure Hodge structure of weight $i$. This means that the space 
${\rm Gr}^{\cal W}_i H^m(V)$ decomposes into the direct sum 
\[ {\rm Gr}^{\cal W}_i H^m(V) = \bigoplus_{p+q = i} H^{p,q}(H^m(V)), \]
such that 
\[ {\cal F}^p {\rm Gr}^{\cal W}_i H^m(V) = \bigoplus_{s \geq p} 
H^{s,i - s}(H^m(V)).\]

\begin{opr}
{\rm The numbers 
\[ h^{p,q}(H^m(V)) = {\rm dim}_{\bf C}\, H^{p,q}(H^m(V)) \] 
are called 
{\em the Hodge-Deligne numbers}.}
\end{opr}
  
\begin{exam}
{\rm Consider the 
case $V = {\bf T}$. Then the Hodge and weight filtrations on the 
space 
$H^{m}({\bf T})  \cong \Lambda^m M_{\bf C}$ have the following properties 

\[ {\cal F}^{0}H^m({\bf T}) = {\cal F}^{1} H^m({\bf T} ) =  
\cdots = {\cal F}^{m}H^m({\bf T}) = H^m({\bf T}),   \]

\[ 0 = {\cal W}_{m-1}H^m({\bf T}) = {\cal W}_{m} H^m({\bf T}) =  
\cdots = {\cal W}_{2m-1}H^m({\bf T}). \]
In other words, all non-zero elements of $H^{m}({\bf T})$ have the 
maximal possible  weight $2m$ and these elements 
generate the minimal nonzero 
Hodge subspace ${\cal F}^{m}H^m({\bf T}) \subset H^m({\bf T})$. 
In particular, 
\[ h^{m,m}(H^m({\bf T})) =  { n \choose m }\;\; (0 \leq m \leq n) \]
are the only possible nonzero Hodge-Deligne numbers of  the cohomology 
groups of ${\bf T}$. }
\label{h.tori}
\end{exam}

Since the natural restriction mapping 
$H^{i}({\bf T}) \rightarrow H^{i}( Z_f )$ 
respects the Hodge filtrations and the weight filtrations, by 
toric version of the Lefschetz theorem (cf. \ref{Lef}), we obtain:

\begin{coro}
For $i < n-1$, one has 
\[H^i(Z_f) \cong \Lambda^i M_{\bf C}\]
where $H^i(Z_f) = {\cal F}^iH^i(Z_f)$ and all nonzero elements of 
$H^i(Z_f)$ have weight $2i$.  The nonzero elements of  the 
image of $H^{n-1}({\bf T})$ in $H^{n-1}(Z_f)$ have the maximal 
possible weight $2n -2$, and these elements are contained in the minimal 
Hodge subspace ${\cal F}^{n-1}H^{n-1}(Z_f)$. 
\label{min.hodge}
\end{coro}

Let us show how one can calculate using the results of 
Danilov and Khovanski\^i \cite{dan.hov} 
the dimensions of quotients 
\[ {\cal F}^{i}H^{n-1}(Z_f)/ {\cal F}^{i+1}H^{n-1}(Z_f) \]
of the Hodge filtration on the "interesting" cohomology group 
$H^{n-1}(Z_f)$. 
\medskip

First recall that for arbitrary complex variety $V$ the 
cohomology groups $H^m_c(V)$ with compact supports also 
have the natural mixed Hodge structure. If $V$ is smooth and 
irreducible, there exists 
the Poincare duality 
\[ H^i(V) \otimes H^{2r  -i}_c(V) \rightarrow  H_c^{2r}(V) \]
compatible with the Hodge structures, 
where $r = {\rm dim}_{\bf C} \,V$.  In particular,  
\begin{equation}
 h^{p,q}(H^i(V)) = h^{r-p,r-q}(H^{2r -i}_c(V)). 
\label{poincare} 
\end{equation} 

\begin{theo}
{\rm (Danilov-Khovanski\^i, \cite{dan.hov} \S4)}  
For any $\Delta$-regular affine hypersurface $Z_f \subset {\bf T}$, one has  
\[ (-1)^{n-1} \sum_{ q \geq 0} \sum_{k \geq 0}(-1)^k h^{p,q}(H^{k}_c(Z_f))  = 
(-1)^p {n \choose p+1 } + \psi_{p+1}(\Delta). \]
\end{theo}

Applying the Poincare duality (\ref{poincare}),  
and the calculations in \ref{h.tori}, 
we obtain: 

\begin{coro}
The dimension  
\[ {\rm dim}\, {\cal F}^{i}H^{n-1}(Z_f)/ {\cal F}^{i+1}H^{n-1}(Z_f) = 
\sum_{ q \geq 0  } h^{i,q}(H^{n-1}(Z_f)) \] 
equals $\psi_{n-i}(\Delta)$ for $ i < n-1$ and $\psi_1(\Delta) + n = 
l(\Delta) -1$ for $i = n-1$. 
\label{n.hodge}
\end{coro}

In this paper we will be  interested in the dependence of the mixed 
Hodge structure on the cohomology groups  of $Z_f$ on the choice of 
the Laurent polynomial $f$. 
It is clear that the cohomology classes coming from the cohomology 
groups of the algebraic torus ${\bf T}$ do not depend on $f$. 
This observation motivates the following 
definition.

\begin{opr}
{\rm The cokernel of the injective mapping 
$H^{n-1}({\bf T}) \rightarrow H^{n-1}( Z_f )$ we call 
the {\em primitive part of the cohomology group} $H^{n-1}( Z_f )$ 
and denote it by 
$PH^{n-1}( Z_f )$.}
\end{opr}

Using \ref{min.hodge} and the equality 
$\varphi_{n}(\Delta) = \psi_1 (\Delta) = l(\Delta) - n - 1$, 
we can reformulate 
\ref{n.hodge} as follows. 

\begin{coro}
For $i \leq n-1$, one has 
\[ {\rm dim}\, {\cal F}^{i}PH^{n-1}(Z_f)/ {\cal F}^{i+1}PH^{n-1}(Z_f)) = 
\sum_{ q \geq 0  } h^{i,q}(PH^{n-1}_c(Z_f)) = \psi_{n-i}(\Delta). \] 
\label{n.hodge1}
\end{coro}

\begin{rem}
{\em We  see now  that the dimensions of quotients 
of the Hodge filtration on $PH^{n-1}(Z_f)$ coincide with the dimensions of 
the homogeneous components of the quotient of the graded algebra $S_{\Delta}$ 
by a regular sequence $s_0, s_1, \ldots , s_n$ 
consisting of $n+1$ elements of degree $1$  (cf. \ref{reg}).  
This motivates the 
existence of the canonical isomorphisms between  the spaces   
\[ {\cal F}^{i}PH^{n-1}(Z_f)/ {\cal F}^{i+1}PH^{n-1}(Z_f)  \]
and  homogeneous components of an Artinian  quotient of $S_{\Delta}$. 
Of course, as we will see later, 
the choice of the regular sequence of length $n+1$ 
depends  on  the Laurent polynomial $f$.}
\label{artin}
\end{rem}

\section{Geometric and homological characterizations of $\Delta$-regularity}


\hspace*{\parindent}

Let $\overline{M} = {\bf Z} \oplus M$. We denote by $\overline{M}_{\bf R}$ the 
scalar extension  $\overline{M} \otimes {\bf R}$. We identify elements 
$\overline{m}= (m_0, m) \in \overline{M}$ with the exponents of Laurent 
monomials $X_0^{m_0} X^m$ in the Laurent polynomial 
ring 
\[ {\bf L} \lbrack X_0^{\pm 1} \rbrack = 
{\bf C}\lbrack X_0^{\pm 1}, X_1 ^{\pm 1}, \ldots,  X_n^{\pm 1} \rbrack.  \]

\begin{opr}
{\rm Let $f(X) \in {\bf C}\lbrack X^{\pm  } \rbrack$ 
be a Laurent polynomial, $\Delta = \Delta(f)$.  Define the new  polynomial 
$F(X_0, {X}) = X_0f(X) -1  \in 
{\bf L}\lbrack X_0 \rbrack$. The Newton polyhedron $\overline{\Delta}$ of $F$ 
in $\overline{M}_{\bf R}$ is the convex hull of 
$\{0 \} \in \overline{M}_{\bf R}$ and 
\[  \{ \overline{x} = (x_0, x)  \in 
\overline{M}_{\bf R} \mid x_0 =  1, \;\;  x \in \Delta  \}. \]
}
\end{opr}

The polynomial $F(X_0,{X})$ is an element of the subring $S_{\Delta} 
\subset {\bf L}\lbrack X_0, \rbrack$. 
One has the following  easy statement.

\begin{prop}
A Laurent  polynomial $f(X)$ is $\Delta$-regular if and only if 
$F(X_0, {X})$ is $\overline{\Delta}$-regular. 
\end{prop}
\label{f.reg}
\bigskip

For the polynomial $F(X_0, {X}) \in S_{\Delta}$,  define $n$ new 
polynomials  
 \[  F_i (X_0, {X}) =  X_i\frac{\partial }{\partial X_i} 
F(X_0, {X}) \;\;(0 \leq i \leq n). \]
All polynomials $F_i$ have the same Newton polyhedron 
\[ \Delta(F_i) = \{ \overline{x} = (x_0, x)  \in 
\overline{M}_{\bf R} \mid x_0 = 1, \;\;  x \in \Delta  \}, \]
i.e.,  $F_0, \ldots, F_n $ are elements of degree 1 in $S_{\Delta}$. 
One  has  $F_0 = X_0 f$, and $F_i = X_0 f_i$ ($1 \leq i \leq n$). 
We can now consider polynomials $F_0, F_1, \ldots, F_n $ 
as global sections of the ample invertible sheaf 
${\cal O}_{\Delta}(1)$ over the toric variety ${\bf P}_{\Delta}$. 
If we denote by $Z_i$ the zero set of $F_i$ in ${\bf P}_{\Delta}$, then  
we obtain the following  geometric formulation of 
the $\Delta$-regularity \ref{regul2}.

\begin{prop}
A Laurent polynomial $f \in L(\Delta)$ is $\Delta$-regular if and only if 
\[ \bigcap_{i=0}^n Z_i = \emptyset. \]
In other words, $f$ is $\Delta$-regular if and only if the open subsets  
\[ U_i = {\bf P}_{\Delta} \setminus Z_i \;\; (0 \leq i \leq n) \]
define an affine covering of the projective toric variety 
${\bf P}_{\Delta}$. 
\label{regul3}
\end{prop}

\begin{opr}
{\rm The covering 
\[ {\cal U} = \{ U_i \}_{ i =0, \ldots,n} \]
will be called a {\em Jacobian covering of ${\bf P}_{\Delta}$.} }
\end{opr}

\begin{rem} 
{\rm Toric variety ${\bf P}_{\Delta}$ is  a compactification of the 
$n$-dimensional torus ${\bf T}_{\Delta} \stackrel{\rm def}{=} {\bf T}$ 
by algebraic tori 
${\bf T}_{\Delta'}$ of smaller dimensions one-to-one corresponding to 
faces $\Delta' \subset \Delta$ of the polyhedron $\Delta$ \cite{dan2}. }
\end{rem}
\medskip
 
Let $\overline{Z}_f$  be the closure of $Z_f$ in ${\bf P}_{\Delta}$. For any 
face  $\Delta' \subset \Delta$ we obtain an affine  hypersurface 
$Z_{f, \Delta'} = \overline{Z}_f \cap  {\bf T}_{\Delta'}$ in the 
algebraic torus ${\bf T}_{\Delta'}$. Using the standard 
criterion of smoothness, we immediatelly obtain another  
geometric characterization of $\Delta$-regularity.

\begin{prop}
A Laurent polynomial $f \in L(\Delta)$ is $\Delta$-regular if and only if 
$Z_{f,\Delta'}$ is a smooth affine subvariety in ${\bf T}_{\Delta'}$ of 
codimension $1$ for any face $\Delta' \subset \Delta$ $($here we assume 
that empty algebraic set has dimension $-1$$).$ 
\label{regul4}
\end{prop}

Let us now discuss homological characterizations of $\Delta$-regularity. 

Put  $\overline{N} = {\rm Hom}(\overline{M}, {\bf Z})$.  
Consider the Koszul complex $K_{\cdot} (f)$ on elements $F_0, \ldots, 
F_n  \subset 
S_{\Delta}$:
\begin{equation}
K_{\cdot} (f)\; : \;  0 \rightarrow S_{\Delta} \otimes_{\bf Z} \Lambda^{n+1} 
\overline{N} \rightarrow \cdots 
\rightarrow S_{\Delta} \otimes_{\bf Z} \Lambda^1  \overline{N}  
{\rightarrow} S_{\Delta} \rightarrow 0. 
\label{koszul}
\end{equation}

Since $S_{\Delta}$ is a graded ring, we can define the grading of 
$S_{\Delta} \otimes_{\bf Z} \Lambda^i \overline{N}$   
by putting ${\rm deg}\, \overline{N} =1$. So we obtain  
$K_{\cdot}(f)$  is a complex 
of graded free $S_{\Delta}$-modules, and the differential of  $K_{\cdot}(f)$ 
has degree zero. We denote by $K_{\cdot}^{(l)}(f)$ the $l$-th homogeneous 
component of the complex $K_{\cdot}(f)$. 
\medskip

\begin{opr}
{\rm Let $J_{f,\Delta}$  be the ideal of $S_{\Delta}$ generated by 
$F_0, F_1 , \dots , F_n$. 
The ideal $J_{f, \Delta}$ we call {\em the Jacobian ideal}, and the 
quotient ring  $R_f  = S_{\Delta} / J_{f,\Delta}$ we call 
the {\em Jacobian ring}  corresponding to the  Laurent polynomial $f$.} 
\end{opr}

Consider the infinite filtration in $S_{\Delta}$
\begin{equation}
{\cal J} \;\;: \;\;  
J_{f,\Delta} \supset J_{f,\Delta}^2 \supset J_{f,\Delta}^3 \cdots 
\label{filt.jac} 
\end{equation} 
 by powers of the ideal $J_{f,\Delta}$. Denote by ${\rm Gr}_{\cal J} 
 S_{\Delta}$ 
the associated graded ring.

The following statement is very important in the sequel. 

\begin{theo}
Let $f \in L(\Delta)$ be a Laurent polynomial. 
Then the  following conditions are equivalent:

{\rm (i)} the elements $F_0,F_1, \cdots , F_n$ give rise to a 
regular sequence in $S_{\Delta}$$;$

{\rm (ii)} the elements $F_0,F_1, \cdots , F_n$ are algebraically 
independent$;$

{\rm (iii)} the homology groups  $H_i(f)$ of the Koszul complex $K_{\cdot}(f)$ 
are zero for $i >0$ $;$.

{\rm (iv)} the ${\bf C}$-space $H_0(f) =R_f$ has a finite dimension $;$

{\rm (v)} the dimension of ${\bf C}$-space $H_0(f) =R_f$ is  $d(\Delta)$$;$ 

{\rm (vi)} ${\rm Gr}_{\cal J} S_{\Delta}$ is isomorphic to the  polynomial
ring $ R_f \lbrack t_0, t_1, \ldots, t_n \rbrack$ over $R_f$ , where 
$t_i$ is the image of $F_i$ in ${\rm Gr}_{\cal J} S_{\Delta}$$;$

{\rm (vii)} The complex $K_{\cdot}^{(l)}(f)$ of finite dimensional 
${\bf C}$-spaces 
is acyclic for $l > n$$;$

{\rm (viii)} $f$ is $\Delta$-regular. 
\label{theo.reg}
\end{theo}  
 
\Proof The equivalences 
${\rm (i)}\Leftrightarrow{\rm (ii)}\Leftrightarrow{\rm (iii)}
\Leftrightarrow{\rm (iv)}\Leftrightarrow{\rm (vi)}$  are well-known properties of  
 Cohen-Macauley rings of the Krull 
dimension $n+1$  (cf. \cite{hartsh} Ch. II, Th. 8.21A). 
Obviously, ${\rm (v)}\Rightarrow{\rm (iv)}$. On the other hand, the statement 
${\rm (v)}$ follows from ${\rm (i)}$ and theorem \ref{hilb}. 

By \ref{regul3}, 
$f$ is $\Delta$-regular if and only if the maximal 
homogeneous ideal $S_{\Delta}^+$ is the only  prime ideal in $S_{\Delta}$ 
associated with $J_{f,\Delta}$, i.e., $J_{f,\Delta}$ contains some 
power of the ideal $S_{\Delta}^+$. This proves 
${\rm (iv)}\Leftrightarrow{\rm (viii)}$.  

Obviously, (vii)$\Rightarrow$(iv). On the other hand, (vii) follows from 
(i),(iii) and from the fact that the Hilbert-Poincare series of the Artinian 
ring $R_f$ is a polynomial of degree $\leq n$ (cf. \ref{hilb}). 
\bigskip

Using standard facts on regular sequences, we obtain:

\begin{coro} 
$f\in L(\Delta)$ is a $\Delta$-regular Laurent polynomial, if and only if 
the elements $F_1, \ldots, F_n$ form a regular sequence in the quotient 
ring $S_{\Delta}/ F_0 S_{\Delta}$.
\label{f0}
\end{coro} 
\bigskip

Let us now define another Koszul complex depending on 
the Laurent polynomial $f$. 
First we recall a general definition of the Koszul complex 
$K^{\cdot}(\{H_1, \ldots, H_p \}, E)$ 
of $p$ commuting operators $H_1, \ldots, H_p$ acting on a 
${\bf C}$-space $E$.  
 
\begin{opr}
{\rm Let $I = \{ i_1, \ldots, i_q \}$ be an ordered subset of ${\cal Q}_p = 
\{1, \ldots, p\}$, i.e., we assume that $i_1 < \cdots < i_q$. Denote by 
$\{ I \}$ the set of elements in $I$, and by $\mid I \mid = q$ 
the cardinality of 
$I$. Let $A$ be a free abelian group of rank $p$ with a chosen ${\bf Z}$-basis 
$e_1, \ldots, e_p$, $\Lambda^q A$ the $q$-th exterior power of $A$ with the 
${\bf Z}$-basis $e_I = e_{i_1} \wedge \cdots \wedge e_{i_q}$, where 
$\mid I \mid = q$, $I \subset {\cal Q}_p$. Denote by 
$K^q (\{H_1, \ldots, H_p \}, E)$ the ${\bf C}$-space 
\[ E \otimes_{\bf Z} \Lambda^q A.\]
Define the differential   
\[ d_H \;:\;  E \otimes_{\bf Z} \Lambda^{q-1} A \rightarrow  
E \otimes_{\bf Z} \Lambda^{q} A\]
of the complex $K^{\cdot} (\{H_1, \ldots, H_p \}, E)$ as follows 
\[ d_H  (\sum_{\mid I \mid = q-1} w^I \otimes e_I ) = 
\sum_{\mid I' \mid = q} d_H(w)^{I'} \otimes e_{I'}, \]
where }
\[ d_H(w)^{I'} = \sum_{ \{i\} \cup \{I\} = \{I'\}} 
 {\rm sgn} { I' \choose {i \; I }} H_i (w^I). \]
\end{opr}

\begin{exam}
{\rm Let  $E = S_{\Delta}$, $A = \overline{M}$, $p = n+1$, 
and $\{ H_1, \ldots, H_{n+1} \}$ are  operators induced 
by multiplications on the 
elements $F_0, F_1, \ldots, F_n$. 
Then there exists the standard $*$-isomorphism between complexes $K^{\cdot} 
(\{F_0, F_1, \ldots, F_n \}, S_{\Delta})$ and 
$K_{\cdot} (f)$, i.e., one has the isomorphisms 
\[ K^{q} (\{F_0, F_1, \ldots, F_n \}, S_{\Delta}) \cong K_{n+1 -q}(f) \;\; 
(0 \leq q \leq n+1).\]  }
\end{exam}	 

Using the above $*$-isomorphism, the characterizations of 
$\Delta$-regularity in \ref{theo.reg} and in \ref{f0}  can be reformulated 
as follows. 

\begin{prop}
A Laurent polynomial $f \in L(\Delta)$ is $\Delta$-regular if and only if 
the cohomology groups $H^i(\{F_0,F_1, \ldots F_n \}, S_{\Delta})$ are zero 
for $i < n+1$, or if and only if 
the cohomology groups $H^i(\{F_1, \ldots F_n \}, 
S_{\Delta}/F_0 S_{\Delta})$ are zero 
for $i < n$. 
\label{ex.reg}
\end{prop}
\bigskip

For the calculation in the sequel, we  need the following theorem. 

\begin{theo}
Let $f$ be a $\Delta$-regular Laurent polynomial. Denote by 
\[ K^{\cdot}(\{F_1, \ldots, F_n\}, S^+_{\Delta}/F_0 S^+_{\Delta}) \] 
the Koszul 
complex of $n$ operators induced by multiplications by  
$F_1, \ldots, F_n$ on the graded $S_{\Delta}$-module 
$S^+_{\Delta}/F_0 S^+_{\Delta}$.  
There exist the following canonical isomorphisms 
\[ H^i(\{F_1, \ldots, F_n\}, S^+_{\Delta}/F_0 S^+_{\Delta}) 
\cong \Lambda^{i}\overline{M}_{\bf C}  \;\; ( 0 \leq i < n) \]
and the  exact sequence 
\[ 0 \rightarrow  \Lambda^n\overline{M}_{\bf C} \rightarrow  
H^n (\{F_1, \ldots, F_n\}, S^+_{\Delta}/F_0 S^+_{\Delta}) \rightarrow 
R_f  \rightarrow   
\Lambda^{n+1}\overline{M}_{\bf C} \rightarrow 0. \]  

In particular, the dimension of 
\[ H^i(\{F_1, \ldots, F_n\}, S^+_{\Delta}/F_0 S^+_{\Delta}) \]
equals  ${ n+1 \choose i }$ 
for $0 \leq i < n $, and $d(\Delta) + n$ for $i =n$.   
\label{grade}
\end{theo}

\Proof One has the following short sequence of complexes
\[ 0 \rightarrow 
K^{\cdot}(\{F_1, \ldots, F_n\}, S^+_{\Delta}) \rightarrow 
K^{\cdot}(\{F_1, \ldots, F_n\}, S_{\Delta}) \rightarrow \Lambda^{\cdot} 
M_{\bf C} \rightarrow 0, \]
where $\Lambda^{\cdot} M_{\bf C}$ is the exterior algebra of the 
$n$-dimensional ${\bf C}$-space $M_{\bf C}$ considered as 
a complex with trivial differential. 
By \ref{ex.reg}, for $0 \leq i < n$, we obtain 
\[ H^{i}(\{F_1, \ldots, F_n \}, S_{\Delta}) = 0. \]
Applying the long exact cohomology sequence, we obtain the exact sequence 
\[ 0 \rightarrow  \Lambda^{n-1} {M}_{\bf C} \rightarrow  
H^n (\{F_1, \ldots, F_n\}, S^+_{\Delta}) \rightarrow 
S_{\Delta}/ ( F_1, \ldots, F_n) S_{\Delta}  \rightarrow   
\Lambda^{n}{M}_{\bf C} \rightarrow 0 \]    
and the isomorphisms 
\[ H^i(\{F_1, \ldots, F_n\}, S^+_{\Delta}) 
\cong \Lambda^{i}{M}_{\bf C}  \;\; ( 0 < i < n),  \]
\[ H^0(\{F_1, \ldots, F_n\}, S^+_{\Delta}) \cong 0. \]

Consider now the short exact sequence of complexes
\[ 0 \rightarrow 
K^{\cdot}(\{F_1, \ldots, F_n\}, S^+_{\Delta}) \stackrel{F_0}{\rightarrow} 
K^{\cdot}(\{F_1, \ldots, F_n\}, S^+_{\Delta}) \rightarrow  
 K^{\cdot}(\{F_1, \ldots, F_n\}, S^+_{\Delta}/F_0 S^+_{\Delta}) 
\rightarrow 0. \]
Since $F_0 \subset S_{\Delta}^+$, 
the multiplication by $F_0$ induces trivial  homomorphism 
of the cohomology groups 
\[ H^i(\{F_1, \ldots, F_n\}, S^+_{\Delta}) 
\stackrel{F_0}{\rightarrow} H^i(\{F_1, \ldots, F_n\}, S^+_{\Delta}) 
\;\; ( 0 \leq i < n).  \]
Now the statement of the theorem follows immediately from the corresponding 
long exact cohomology sequence and the canonical isomorphisms:  
\[ \Lambda^i \overline{M}_{\bf C} = \Lambda^i {M}_{\bf C} 
\oplus \Lambda^{i-1} {M}_{\bf C}. \]
\bigskip

Finally, we give the characterization of $\Delta$-regularity via the notion 
of the 
principal discriminant due to Gelfand, Kapranov and Zelevinsky \cite{gel1}. 
\medskip

Let 
\[ f(\alpha, X) = \sum_{a_m \in \Delta \cap {\bf Z}^n} a_m X^m \]
be the Laurent polynomial with algebraically independent coefficients 
$(a_{m_1}, \ldots, a_{m_s}) = \alpha$,  where $s$ is the number of integral 
points in $\Delta$. Denote by ${\bf C} \lbrack \alpha \rbrack$ the polynomial 
ring ${\bf C} \lbrack a_{m_1}, \ldots, a_{m_s} \rbrack$. Consider 
the polynomial ring 
\[ S_{\Delta} \lbrack \alpha \rbrack  : = S_{\Delta} {\otimes}_{\bf C} 
{\bf C} \lbrack a_{m_1}, \ldots, a_{m_s} \rbrack. \] 
Then every  $l$-th homogeneous component $K^{(l)}_{\cdot} (f(\alpha, X))$ 
of the Koszul complex $K_{\cdot} (f(\alpha, X))$ can be considered as a finite 
complex of  ${\bf C} \lbrack \alpha \rbrack$-modules. 

\begin{opr}
{\rm \cite{gel1} Define the {\em principal $\Delta$-discriminant of 
the polynomial $f(\alpha,X)$ } as the determinant of the complex 
$K^{(l)}_{\cdot} (f(\alpha, X))$ for $l \gg 0$. }
\end{opr}

\begin{rem}
{\rm Here we use the notion of the determinant of a complex as in 
\cite{bism}. }
\end{rem}

It is proved in \cite{gel1} that for $l \gg 0$ the determinant of 
$K^{(l)}_{\cdot} (f(\alpha, X))$ does not depend on $l$,  and 
it is a polynomial in variables $a_{m_1}, \ldots, a_{m_s}$. We denote 
this polynomial by ${\rm Disc}_{\Delta}(f)$.   
\medskip

From \ref{theo.reg}(vii) we obtain: 

\begin{prop}
The Laurent polynomial $f$ is $\Delta$-regular if and only if 
${\rm Disc}_{\Delta}(f) \neq 0$. 
\end{prop}

\section{Poincare residue mappings}


\hspace*{\parindent}

There exists the long exact sequence of cohomology with compact 
supports 
\[ \cdots \rightarrow H_c^{i-1}(Z_f) \rightarrow H^i_c ({\bf T} \setminus 
Z_f) \rightarrow H^i_c ({\bf T}) \rightarrow H^i_c(Z_f) \rightarrow 
H^{i+1}_c ({\bf T} \setminus Z_f ) \rightarrow \cdots ,  \]
where all homomorphisms are morphisms of Hodge structures. We know 
that for $i \geq n$ the space $H^i_c({\bf T})$ is dual to 
$H^{2n-i}({\bf T})$.  
So the space $H^i_c({\bf T})$ has the Hodge type $(i-n,i-n)$. 
By the toric version of the Lefschetz theorem 
\ref{Lef}, the space $H^i_c(Z_f)$  
has the Hodge type $(i-n +1, i-n +1)$, this space is dual to 
$H^{2n -2 -i}(Z_f)$.  Hence the homomorphisms  
\[ H^i_c ({\bf T}) \rightarrow H^i_c(Z_f) \]
should be trivial for $i \geq n$. Therefore, the long exact sequence 
splits into 
short exact sequences
\[ 0 \rightarrow H_c^{i-1}(Z_f) \rightarrow H^i_c ({\bf T} \setminus 
Z_f) \rightarrow H^i_c ({\bf T}) \rightarrow 0.\]
By the Poincare duality, for $j \leq n$ we obtain the short 
exact sequences  of usial cohomology groups
\[ 0  \rightarrow H^j ({\bf T}) 
\rightarrow H^j ({\bf T}\setminus Z_f) \rightarrow H^{j-1}(Z_f) \rightarrow 
0.  \]

\begin{opr}
{\rm The surjective homomorphisms 
\[ {\rm res}\;: \; H^j({\bf T}\setminus Z_f ) \rightarrow H^{j-1} (Z_f) \]
obtained in above short exact 
sequences are  called {\em Poincare residue mappings}. }
\end{opr}

 The Poincare residue mapping  establishes the connection between 
the cohomology  of $Z_f$ and the cohomology of the complementary 
open set ${\bf T} \setminus Z_f$. 
\medskip 

Our main purpose is  to give the explicit representation of the 
spaces  $H^{j}({\bf T} \setminus Z_f)$ by rational differential $j$-forms 
on ${\bf T}$ with poles on  $Z_f$, and then to describe $H^{j-1}(Z_f)$ 
and $PH^{n-1}( Z_f )$ by 
residues of these rational differential $j$-forms. 
\bigskip

First note  that the open subset ${\bf T} \setminus Z_f \subset {\bf T}$ 
is  isomorphic to  the  affine  hypersurface $Z_F$ in the 
$(n+1)$-dimensional torus 
\[ \overline{\bf T} = {\rm Spec}\, {\bf L} \lbrack X_0^{\pm1 } \rbrack. \] 
defined by 
the polynomial $F(X_0, X) = X_0f(X) -1$.  The Newton polyhedron 
$\overline{\Delta}$ of $F$ satisfies the condition 
$d(\overline{\Delta}) = d(\Delta)$. Hence, repeating 
the arguments for $n$-dimensional $\overline{\Delta}$-regular 
hypersurface $Z_F \subset \overline{\bf T}$ (cf. \ref{i.dim}, and 
\ref{n.dim}), we obtain:   

\begin{prop}
$ {\rm dim }\, H^{n}(Z_F)  = d(\Delta) + n$, $ {\rm dim }\, H^{i}(Z_F)  = 
{ n +1 \choose i }$ for $i < n$. 
\end{prop}

By theorem \ref{Lef},  for $i <n$ the residue mapping $H^i(Z_F) 
\stackrel{\rm res}{\rightarrow}  H^{i-1}(Z_f)$ reduces to the residue mapping 
$H^i(\overline{\bf T}) \stackrel{\rm res}{\rightarrow}  H^{i-1}({\bf T})$.
Hence  for $i < n$ we obtain the following isomorphisms 
of short exact sequences 

\begin{center}
\begin{tabular}{ccccccccc}
  &  & $0$  &  & $0$ & & 0  & & \\
  &  & $\downarrow$  &  & $\downarrow$ & & $\downarrow$  & & \\
$0$ & $\rightarrow$ & $H^{i}({\bf T})$ & ${\rightarrow}$ &
 $H^{i}(\overline{\bf T})$& $\stackrel{\rm res}{\rightarrow}$ & 
 $H^{i-1}({\bf T})$ & & \\

  &  & $\cong \downarrow$  &  & $\cong \downarrow$ & & $\cong \downarrow$ 
  & & \\

$0$ & $\rightarrow$ & $H^{i}({\bf T})$  & 
${\rightarrow}$ & $H^{i}({Z_F})$ & $\stackrel{\rm res}{\rightarrow}$ 
& $H^{i-1}({Z}_f)$ 
& $\rightarrow$ & $0$ \\

  &  & $ \downarrow$$\;\;$  &  & $\downarrow$$\;\;$ & 
  & $\downarrow$ & & \\
    &  & $0$  &  & $0$ & & $0$ & & \\
  \end{tabular}
\end{center}

For $i \leq  n$, the exterior multiplication 
\[ \frac{dX_0}{X_0} \wedge * \;\; : \;\; 
 H^{i-1}({\bf T}) \rightarrow 
H^{i}(\overline{\bf T}).\] 
canonically splits the short exact sequences 
\[ 0 \rightarrow H^i({\bf T}) \rightarrow  H^{i}(\overline{\bf T}) 
\stackrel{\rm res}{\rightarrow}
 H^{i-1}({\bf T}) \rightarrow 0.  \]
Moreover, there exists the  commutative diagram with exact rows and columns 
\begin{center}
\begin{tabular}{ccccccccc}
  &  & $0$  &  & $0$ & & 0  & & \\
  &  & $\downarrow$  &  & $\downarrow$ & & $\downarrow$  & & \\
$0$ & $\rightarrow$ & $H^{n}({\bf T})$ & ${\rightarrow}$ &
 $H^{n}(\overline{\bf T})$& $\stackrel{\rm res}{\rightarrow}$ & 
 $H^{n-1}({\bf T})$ & & \\

  &  & $\cong \downarrow$  &  & $\downarrow$ & & $\downarrow$ & & \\

$0$ & $\rightarrow$ & $H^{n}({\bf T})$  & 
${\rightarrow}$ & $H^{n}({Z_F})$ & $\stackrel{\rm res}{\rightarrow}$ 
& $H^{n-1}({Z}_f)$ 
& $\rightarrow$ & $0$ \\

  &  & $ \downarrow$$\;\;$  &  & $\downarrow$$\;\;$ & 
  & $\downarrow$ & & \\

 &  & 0 & 
${\rightarrow}$ & $PH^{n}({Z_F})$ & $\stackrel{\cong}{\rightarrow}$ & 
$PH^{n-1}({Z_f})$ 
& $\rightarrow$ & $0$ \\ 
  &  &   &  & $\downarrow$ & & $\downarrow$ & & \\
    &  &   &  & $0$ & & $0$ & & \\
  \end{tabular}
\end{center}

One should remark that the Poincare residue mapping is a morphism of mixed 
Hodge structures of the Hodge type $(-1, -1)$, i.e., it is a 
morphism of  weight $-2$. 
This property of the residue mapping puts some conditions 
on the Hodge and weight 
 filtrations in $PH^{n}(Z_F)$. 

\begin{prop}
The Poincare residue mapping defines the isomorphism between the primitive 
cohomology groups $PH^{n}(Z_F)$ and $PH^{n-1}(Z_f)$ such that 
\[ {\rm res}\,  ({\cal F}^i PH^{n}(Z_F)) = {\cal F}^{i-1}PH^{n-1}(Z_f) \]
and 
\[ {\rm res}\,  ({\cal W}_j PH^{n}(Z_F)) = {\cal W}_{j-2}PH^{n-1}(Z_f).  \] 
In particular, we obtain 
\[ PH^n(Z_F) = {\cal F}^1 PH^{n}(Z_F) = {\cal F}^0 PH^{n}(Z_F) \] 
and 
\[ {\cal W}_{n}PH^{n}(Z_F) = 0 . \]
\end{prop}
\bigskip

For $i <n $ the description of the spaces $H^{i}({\bf T} \setminus Z_f)$ 
by rational differential $i$-forms is much more easy than for the case 
$i = n$. 
\medskip

In this section we consider only the case  $i <n$. 
The space $H^{i}({\bf T} \setminus Z_f)$ is isomorphic to $\Lambda^i 
\overline{M}_{\bf C}$. We can identify $\Lambda^i 
\overline{M}_{\bf C}$ with the cohomology $H^i (\overline{\bf T})$. Choose 
the basis in $H^1(\overline{\bf T})$ consisting of 
$\overline{\bf T}$-invariant 
differential $1$-forms 
\[ \frac{dX_0}{X_0},  \frac{dX_{1}}{X_{1}}, \ldots ,  \frac{dX_n}{X_n}. \]
We obtain also the ${\bf T}$-invariant basis in $H^i(\overline{\bf T})$ 
\[ \frac{dX_I}{X_I}, \;\; I \subset \{ 0, 1, \ldots, n\}, \; \mid I \mid = i. \] 
In order to get a basis in $H^{i}({\bf T} \setminus Z_f) = H^i(Z_F)$, we 
have to consider restriction of the $i$-forms  $dX_I/X_I$ on the hypersurface 
$Z_F \subset \overline{\bf T}$. We note that the restriction of the 
differential $1$-form ${dX_0}/{X_0}$ on $Z_F   
\subset \overline{\bf T}$ coincides with $-df/f$. So we obtain: 

\begin{theo}
For $i <n$ the space $H^{i}({\bf T} \setminus Z_f)$ has the canonical basis 
consisting of the rational differential forms  
\[ \frac{dX_I}{X_I}, \;\; I \subset \{ 1, \ldots, n\}, \; \mid I \mid = i,  \]
and 
\[ \frac{df}{f} \wedge 
\frac{dX_{I'}}{X_{I'}}, \;\; I' \subset \{1, \ldots, n\}, \; 
\mid I' \mid = i-1. \]
Moreover, the Poincare residue mapping vanishes exactly on the $i$-forms 
${dX_I}/{X_I}$, and 
\[ {\rm res}\,( \frac{df}{f} \wedge 
\frac{dX_{I'}}{X_{I'}}) = \frac{dX_{I'}}{X_{I'}} \in H^{i-1}(Z_f). \]
\end{theo}

As a consequence, we obtain:

\begin{coro}
The image of $H^n(\overline{\bf T})$ in  $H^{n} (Z_F)$ is spanned by $n+1$ 
differential $n$-forms  
\[ \omega_0 = \frac{dX_1}{X_1} \wedge \cdots \wedge \frac{dX_n}{X_n} \]
and  
\[ \frac{f_i}{f} \omega_0  \;\; (1 \leq i \leq n). \]
The kernel of the residue mapping 
\[ {\rm res}\, \;:\; H^n(Z_F) \rightarrow H^{n-1}(Z_f) \]
is spanned by $\omega_0$.  
\end{coro}

\section{Logarithmic de Rham cohomology }  


\hspace*{\parindent}  

Fix an $n$-dimensional integral polyhedron $\Delta$.    
We denote the corresponding projective toric variety ${\bf P}_{\Delta}$ 
simply by  ${\bf P}$. 
 
Let $D$ be the union of ${\bf D} = {\bf P} 
 \setminus {\bf T}$ and $\overline{Z}_f$. 
The pair $({\bf P}, D)$ is toroidal (cf. \cite{dan2}), i.e., this pair is 
locally  isomorphic to affine toric variety with an invariant toric 
divisor. 

There exists a degenerating at $E_1$-term  spectral sequence 
\[ E_1^{p,q} = H^q ({\bf P}, \Omega^p_{({\bf P},D)}) \Rightarrow 
H^{p+q}_c ({\bf T} \setminus Z_f),  \]
where $\Omega^p_{({\bf P},D)}$ denotes the sheaf of the differential 
$p$-forms with zeros along $D$ (cf. \cite{dan3}). The limit filtration 
of the spectral sequence coincides with the Hodge filtration on 
$H^{\cdot}_c ({\bf T} \setminus Z_f)$
Let  
$\Omega^{\cdot}_{\bf P}({\rm log}\,D)$ be the complex 
of sheaves of differential forms on 
projective toric variety ${\bf P}$ with logarithmic poles along 
$D$. 
Since $\Omega^p_{({\bf P},D)} \cong Hom\,( 
\Omega^{n-p}_{\bf P}({\rm log}\, D),  \Omega^n_{\bf P})$, using the Poincare 
duality and the Serre's duality, we obtain that for the logarithmic 
de Rham complex 
\[ \Omega^{\cdot}({\rm log}\,D) = \{ \Omega^0({\rm log}\,D) 
\stackrel{d}{\rightarrow}  \Omega^1({\rm log}\,D)  \stackrel{d}{\rightarrow}
\cdots \} \]
there exists the spectral sequence 
\begin{equation}
 E_1^{pq} = H^q ({\bf P}, \Omega^p_{{\bf P}}({\rm log}\,D)) 
\Rightarrow H^{p+q} ({\bf T} \setminus Z_f).
\label{spec.s}
\end{equation}
which {\em  degenerates at the $E_1$-term and converges to the Hodge 
filtration on} $H^{\cdot}({\bf T} \setminus Z_f)$.
\medskip

In our situation we want to give a direct proof of the degeneration of 
the spectral sequence (\ref{spec.s}).   
\medskip

First we note that if $({\bf P}, D)$ were a pair consisting of a 
normal crossing divisor $D$ on a smooth 
projective variety ${\bf P}$, then the deneration of the spectral 
sequence would follow from the general theorem of Deligne \cite{deligne}.  
If ${\bf P}$ is not smooth,  we can use  a toric desingularization 
$ \pi\,:\, {\bf P}'\rightarrow {\bf P}$. Let ${\bf D}' = {\bf P}' \setminus  
{\bf T}$ and $D'$ be ${\bf D}' \cup \overline{Z}'_f$, where 
$\overline{Z}'_f$ is the closure  of $Z_f$ in ${\bf P}'$. Since 
\[ {\bf P}' \setminus D' = {\bf P} \setminus D = {\bf T} \setminus Z_f, \]
the degeneration of the spectral sequence $E_1^{pq}$ in our case 
reduces  to the Deligne's theorem by  the following lemma.  

\begin{lem}
\[ R^0\pi_*  \Omega^p_{{\bf P}'}({\rm log}\,D') =  
\Omega^p_{{\bf P}}({\rm log}\,D) \]
and 
\[ R^j\pi_*   \Omega^p_{{\bf P}'}({\rm log}\,D') = 0\;\; 
{\rm for}\;\; j > 0. \] 
In particular, applying the Leray spectral sequence,  one has 
\[ H^q ({\bf P}', \Omega^p_{{\bf P}'}({\rm log}\,D' )) \cong 
H^q ({\bf P}, \Omega^p_{{\bf P}}({\rm log}\,D))  \]
for $p, q \geq 0$. 
\label{reduce}
\end{lem}

The  proof of lemma \ref{reduce} will be given  later, since it uses 
a  general theorem \ref{log.s} on special  resolutions 
of the sheaves $\Omega^p_{{\bf P}}({\rm log}\,D)$.  
\bigskip

Let  $\Omega^i_{{\bf P, D}} \cong  {\cal O}_{{\bf P}} 
\otimes_{\bf Z} \Lambda^i M$ denotes  the ${\bf T}$-invariant sheaf of 
differential $i$-forms on ${\bf P}$ with logarithmic poles along 
${\bf D}$ (cf. \cite{dan2}). On has the following 
statement. 

\begin{theo}
Let $d$ be the exterior derivative. Then the sequence of sheaves on ${\bf P}$ 
\[ 0 \rightarrow \Omega^p_{\bf P}({\rm log}\,D) \rightarrow 
\Omega^p_{{\bf P,D}}(1) \stackrel{d}{\rightarrow} 
\Omega^{p+1}_{{\bf P,D}}(2)/ \Omega^{p+1}_{{\bf P,D}}(1)
\stackrel{d}{\rightarrow} \cdots \]
\[ \cdots \stackrel{d}{\rightarrow}  
\Omega^{p+k}_{{\bf P,D}}(k +1)/ \Omega^{p+k}_{{\bf P,D}}(k) 
\stackrel{d}{\rightarrow} \cdots 
\stackrel{d}{\rightarrow}  
\Omega^{n}_{{\bf P,D}}(n -p  +1)/ \Omega^{n}_{{\bf P,D}}(n-p) 
{\rightarrow} 0\]
is exact.
\label{log.s}
\end{theo}

\Proof Since the pair $({\bf P}, D)$ is toroidal, it is sufficient to prove  
the following  local statement  for  affine toric models:
\medskip

\begin{theo}
Let $W = {\rm Spec}\, A_{\sigma}$ be an $(n-1)$-dimensional affine 
toric variety defined by a convex rational polyhedral 
$(n-1)$-dimensional cone $\sigma$ {\rm \cite{dan2}}. Let $D_0$ be a toric 
divisor on $W$ corresponding to a set of facets of $\sigma$. 
Consider the $n$-dimensional toric variety $W'= W \times {\bf A}_1 = 
{\rm Spec}\, B$, where $B = A_{\sigma}\lbrack t \rbrack$.   Let $\wp \;:\; 
W' \rightarrow W$ be the projection. 
Denote by $\Omega^i_{W'}({\rm log}\, D)$ the $B$-module of differential 
$i$-forms on $W'$ having  only logarithmic poles along $D \subset W'$, where 
$D = W\times \{0\} \cup {\wp}^*D_0$. Let  
\[ \Omega^i_{W', {\wp}^*D_0} \]
be the $B$-module of differential $i$-forms on $W'$ with only logarithmic poles 
along ${\wp}^*D_{0} \subset W'$. Then there exists the 
following exact sequence of $B$-modules 
\begin{equation}
 0 \rightarrow \Omega^p_{W'}({\rm log}\, D) \rightarrow 
\Omega^p_{W', {\wp}^*D_0}(1) \stackrel{d}{\rightarrow} 
\Omega^{p+1}_{W', {\wp}^*D_0}(2)/ 
\Omega^{p+1}_{W', {\wp}^*D_0}(1)
\stackrel{d}{\rightarrow} \cdots 
\label{long}
\end{equation}
\[ \cdots \stackrel{d}{\rightarrow}  
\Omega^{p+k}_{W', {\wp}^*D_0}(k +1)/ 
\Omega^{p+k}_{W', {\wp}^*D_0}(k) 
\stackrel{d}{\rightarrow} \cdots \]
\[ \cdots \stackrel{d}{\rightarrow}  
\Omega^{n}_{W', {\wp}^*D_0}(n -p  +1)/ 
\Omega^{n}_{W', {\wp}^*D_0}(n-p) 
{\rightarrow} 0, \]
where for any graded $B$-module $M = \oplus_i M^i$ we denote by $M(k) = 
\oplus_i M^i(k)$ the graded $B$-module with the homogeneous components 
$M^i(k) = M^{i+k}$.  
\label{local}
\end{theo}

\Proof  
Let $\Lambda^iQ$ be the $A_{\sigma}$-module of 
differential $i$-forms  on $W$ with only logarithmic poles along $D_0$. 
Note  that 
\[ \Omega^1_{W', {\wp}^*D_0} \cong B dt \oplus {\wp}^* \Omega_{W, D_0}^1 \cong 
Bdt \oplus 
\Lambda^1 Q\lbrack t \rbrack \cong B(-1) \oplus 
\Lambda^1 Q \lbrack t \rbrack.   \]
On the other hand, 
\[ \Omega^1_{W'}({\rm log}\, D) \cong B \frac{dt}{t} \oplus  
\Lambda^1 Q \lbrack t \rbrack \cong B \oplus 
\Lambda^1 Q \lbrack t \rbrack.  \]
Therefore, 
\[ \Omega^i_{W', {\wp}^*D_0} \cong 
\Lambda^{i-1}Q \lbrack t \rbrack \wedge dt  
\oplus  \Lambda^i Q \lbrack t \rbrack \cong 
\Lambda^{i-1}Q \lbrack t \rbrack(-1) 
\oplus  \Lambda^i Q \lbrack t \rbrack,  \]
and 
\[ \Omega^i_{W'} ({\rm log}\, D) \cong \Lambda^{i-1}Q \lbrack t \rbrack  
\oplus  \Lambda^i Q \lbrack t \rbrack.  \]
Thus the cokernel of the injective homomorphism 
\[  \Omega^i_{W'}({\rm log}\, D) \rightarrow 
 \Omega^i_{W', {\wp}^*D_0}(1)  \]
 is isomorphic to  $\Lambda^i Q(1)$, and the quotient  
 \[ \Omega^{i+k}_{W', {\wp}^*D_0}(k +1)/ 
\Omega^{i+k}_{W', {\wp}^*D_0}(k) \] 
is isomorphic to $\Lambda^{i+k-1}Q(k) \oplus \Lambda^{i+k}Q(k+1)$. 
We can see that the long sequence of homomorphisms 
(\ref{long}) is obtained by glueing of the short 
exact sequences  
\[ 0 \rightarrow \Lambda^{p+1 -k}Q(k) 
\stackrel{\mu(-k)}{\rightarrow} 
\Lambda^{p+k-1}Q(k) \oplus \Lambda^{p+k}Q(k+1) 
\stackrel{\mu(-k-1)}{\rightarrow}  \Lambda^{p+k}Q(k+1) \rightarrow 0, \] 
where $\mu(l)$ denotes the multiplication by $l$. This completes proof. 
\bigskip

\begin{rem}
{\rm It follows from the proof of \ref{local} that the theorem \ref{log.s} 
holds also in the case when $\overline{Z}_f$ is not necessary an ample 
Cartier divisor on ${\bf P}$, i.e., we obtain also analogous resolutions of 
sheaves $\Omega^p_{\bf P'}({\rm log \, D'})$.}
\end{rem}

\begin{coro}
{\rm (i)} $H^q({\bf P}, \Omega^{p+k}_{\bf  P,D}(k +1)/ 
\Omega^{p+k}_{\bf P,D}(k)) = 0$ for $q >0$, $k \geq 0$. 

{\rm (ii)} $H^0({\bf P}, \Omega^{p+k}_{\bf P,D}(k +1)/ 
\Omega^{p+k}_{\bf P,D}(k)) = H^0({\bf P}, \Omega^{p+k}_{\bf P,D}(k +1))/ 
H^0({\bf P}, \Omega^{p+k}_{\bf P,D}(k))$ for $k \geq 0$. 
\label{nul}
\end{coro}

\Proof The statement follows from the vanishing theorems  
for ${\cal O}_{\Delta}(k)$ 
(cf. \ref{vanish}):    
\[ H^i({\bf P}, \Omega^p_{{\bf P,D}}(k))  = 0\;\;  {\rm for}\;\; 
i > 0, \; \; k \geq 0, \]
and from  the  short exact sequences  
\begin{equation}
 0 \rightarrow \Omega^p_{\bf P,D}(k) \rightarrow 
\Omega^p_{\bf P,D}(k +1) \rightarrow  \Omega^p_{\bf P,D}(k +1)/ 
\Omega^p_{\bf P,D}(k) \rightarrow 0. 
\label{short1}
\end{equation}
\bigskip

By theorem \ref{log.s} and corollary \ref{nul}, we obtain. 

\begin{theo}
The space  $H^q ({\bf P}, \Omega^p_{{\bf P}}({\rm log}\,D))$ is isomorphic 
to the $q$-th 
cohomology group of the complex 
\[ C^{p,\cdot}_{\rm log}(D)\;:\;  0 \rightarrow H^0({\bf P}, 
\Omega^p_{{\bf P,D}}(1)) \stackrel{d}{\rightarrow} 
H^0({\bf P},\Omega^{p+1}_{{\bf P, D}}(2))/ H^0({\bf P}, 
\Omega^{p+1}_{{\bf P,D}}(1)) 
\stackrel{d}{\rightarrow} \cdots \]
\[  \cdots \stackrel{d}{\rightarrow}  
H^0({\bf P}, \Omega^{p+k}_{{\bf P,D}}(k +1))/ 
H^0({\bf P}, \Omega^{p+k}_{{\bf P,D}}(k) )
\stackrel{d}{\rightarrow} \cdots \]
\[ \cdots \stackrel{d}{\rightarrow}  
H^0({\bf P}, \Omega^{n}_{{\bf P,D}}(n -p  +1))/ 
H^0({\bf P}, \Omega^{n}_{{\bf P,D}}(n-p)) {\rightarrow} 0 .\]

In particular, 
\[H^{n-p} ({\bf P}, \Omega^p_{{\bf P}}({\rm  log}\,D)) \cong 
\frac{H^0({\bf P}, \Omega^{n}_{{\bf P,D}}(n -p +1))} {H^0({\bf P}, 
\Omega^{n}_{{\bf P,D}}(n-p)) + d H^0({\bf P}, 
\Omega^{n-1}_{{\bf P,D}}(n-p))}\]
and 
\[ H^{0} ({\bf P}, \Omega^p_{{\bf P}}({\rm log}\,D)) \cong
{\rm Ker} \lbrack H^0({\bf P}, 
\Omega^p_{{\bf P,D}}(1)) \stackrel{d}{\rightarrow} 
H^0({\bf P},\Omega^{p+1}_{{\bf P,D}}(2))/ H^0({\bf P}, 
\Omega^{p+1}_{{\bf P,D}}(1)) \rbrack.\]
\end{theo}

{\bf Proof of lemma \ref{reduce}.} Since ${\bf P}$ is a normal variety,  
$R^0\pi_* {\cal O}_{{\bf P}'} = {\cal O}_{\bf P}$. Since toric varieties have 
only rational singularities, $R^j \pi_* {\cal O}_{{\bf P}'} = 0$ for $j >0$. 
Therefore, the  twisted sheaves ${\cal O}_{{\bf P}'}(k) = 
{\cal O}_{{\bf P}'} \otimes \pi^*{\cal O}_{\Delta}(k)$ 
for any $k \in {\bf Z}$ are $R\pi_*$-acyclic: 
 \[ R^j \pi_* {\cal O}_{{\bf P}'}(k)  = 0 \;\;{\rm for}\;\;  j >0,  \]
$R^0\pi_* {\cal O}_{{\bf P}'}(k) = {\cal O}_{\Delta}(k)$.  
Since $\Omega^p_{\bf P',D'}(k)$ are direct sums of ${ n \choose p}$ 
copies of  ${\cal O}_{{\bf P}'}(k)$, they are also  $R^j\pi_*$-acyclic. 
Moreover, using the short exact sequence (\ref{short1}), we again obtain the 
$R\pi_*$-acyclicity 
\[  R^j \pi_* (\Omega^p_{\bf P',D'}(k+1) /  
\Omega^p_{\bf P',D'}(k))  = 0 \;\;{\rm for}\;\;  j >0, \]
and 
\[ R^0 \pi_* (\Omega^p_{\bf P',D'}(k+1) /  
\Omega^p_{\bf P',D'}(k) ) = \Omega^p_{\bf P,D}(k+1) /  
\Omega^p_{{\bf P,D}}(k).\]  
Applying theorem \ref{log.s} simultaneously 
for ${\bf P}'$ and ${\bf P}$, we see that the direct image of the $R\pi_*$-acyclic 
resolution of  $\Omega^p_{{\bf P}'}({\rm log}\,D')$ is the analogous 
resolution of  $\Omega^p_{{\bf P}}({\rm log}\,D)$. 
Therefore, $\Omega^p_{{\bf P}'}({\rm log}\,D')$ are also $R\pi_*$-acyclic, 
i.e., $R^j\pi_*  \Omega^p_{{\bf P}'}({\rm log}\,D') = 0$ for $j >0$ 
and  $R^0\pi_*  \Omega^p_{{\bf P}'}({\rm log}\,D') =  
\Omega^p_{{\bf P}}({\rm log}\,D)$. 
\bigskip

It is easy to describe the complex $C_{\rm log}^{p,\cdot}(D)$ explicitly via  
the homogeneous components $S^i_{\Delta}$ of the graded homogeneous coordinate 
ring $S_{\Delta}$.

\begin{prop}
The complex $C_{\rm log}^{p,\cdot}(D)$ has the following representation 
\[ C^{p,\cdot}_{\rm log}(D)\;:\;  0 \rightarrow \frac{1}{F_0}S^1_{\Delta} 
\otimes_{\bf Z} \Lambda^p M  \stackrel{d_1}
{\rightarrow} \frac{1}{F_0^2}(S^2_{\Delta}/ F_0S^1_{\Delta}) \otimes_{\bf Z} 
\Lambda^{p+1}M 
\stackrel{d_2}{\rightarrow} \cdots \]

\[  \cdots \stackrel{d_k}{\rightarrow}  
\frac{1}{F_0^{k+1}}(S^{k+1}_{\Delta}/F_0S^k_{\Delta}) 
\otimes_{\bf Z} \Lambda^{p+k} 
M \stackrel{d_{k+1}}{\rightarrow} 
\cdots \stackrel{d_{n-p}}{\rightarrow}  
\frac{1}{F_0^{n-p+1}}(S_{\Delta}^{n-p+1}/ F_0 S_{\Delta}^{n-p}) 
\otimes_{\bf Z} \Lambda^{n}M  
{\rightarrow} 0, \]
where $d_k$ is the exterior multiplication by 
\[ -k\frac{X_0 df}{F_0}  = -k\frac{1}{F_0}( F_1 \frac{dX_1}{X_1} + 
\cdots + F_n \frac{dX_n}{X_n}). \]
\end{prop}

\begin{coro}
There exist  the following isomorphisms 

{\rm (i)}
\[H^{0} ({\bf P}, \Omega^n_{{\bf P}}({\rm log}\,D)) \cong 
H^{0} ({\bf P}, \Omega^n_{{\bf P,D}} (1) ) \cong  S_{\Delta}^1 
\cong L(\Delta); \]

{\rm (ii)} 
\[H^{n-p} ({\bf P}, \Omega^p_{{\bf P}}({\rm log}\,D)) \cong R_f^{n-p+1} 
\;\;for\;\; n-p >0; \]

{\rm (iii)}  
\[H^{0} ({\bf P}, \Omega^p_{{\bf P}}({\rm log}\,D)) \cong 
\Lambda^{p-1} M_{\bf C} \oplus \Lambda^p M_{\bf C} 
\cong \Lambda^p \overline{M}_{\bf C} \cong 
H^p(\overline{T})\;\; { for\; n > p \geq 0}, \]
where 
\[ \Lambda^{p-1} M_{\bf C} \cong \frac{X_0df}{F_0} \wedge 
\Lambda^{p-1}M_{\bf C} \subset 
\frac{1}{F_0} S_{\Delta}^1 \otimes \Lambda^p M, \] 
\[ \Lambda^{p} M_{\bf C} \cong \frac{F_0}{F_0} \otimes  \Lambda^{p}M_{\bf C} 
 \subset 
\frac{1}{F_0} S_{\Delta}^1 \otimes \Lambda^p M. \] 
\end{coro}

We come now to our first result on the Hodge filtration 
in  the primitive part of $H^n(Z_F)$. 

\begin{theo}
Let 
\[ PH^{n}(Z_F) = 
{\cal F}^0 = {\cal F}^1 \supset \cdots \supset {\cal F}^n \supset 
{\cal F}^{n+1} = {0} \] 
be the Hodge filtration in $PH^{n}(Z_F)$. Denote by 
\[ {\rm Gr}_{\cal F} \, PH^{n}(Z_F) = \bigoplus_{ a} {\cal F}^a / 
{\cal F}^{a+1} \]
the associated graded space. 
Then for any $i \geq 0$ there exists the  canonical isomorphism  
\[ {\rm Gr}_{\cal F}^i PH^{n}(Z_F) =  {\cal F}^{i} / {\cal F}^{i+1} 
\cong R^{n+1 -i}, \] 
i.e., ${\rm Gr}_{\cal F} \, PH^{n}(Z_F)$ is canonically isomorphic to 
the maximal homogeneous ideal $R^+_f$ in the Jacobian ring $R_f$. 
\label{d.mixed} 
\end{theo}

Now we can formulate the precise statement for our observation 
in \ref{artin}. 

\begin{coro}
For any $i \geq 0$ there exists the canonical isomorphism 
\[ {\cal F}^{i}PH^{n-1}(Z_f)/ {\cal F}^{i+1}PH^{n-1}(Z_f)  \cong R_f^{n-i}. \]
\end{coro}

\section{Global de Rham complexes: Dwork-Katz method}  


\hspace*{\parindent}

By the Grothendieck theorem \cite{groth},  we  can 
describe the spaces $H^{i}(Z_F)$ using the 
global de Rham  complex  
\[ \Gamma \Omega^*(Z_F)  = \{ \Gamma \Omega^0(Z_F) 
\stackrel{d}{\rightarrow}  \Gamma \Omega^1(Z_F)  
\stackrel{d}{\rightarrow} \cdots \} \]
of meromorphic differential forms on 
${\bf T}$  with poles of arbitrary order on  $Z_f$. 
Since $\Gamma \Omega^i({\bf T}) \cong {\bf L} \otimes_{\bf Z} 
\Lambda^i M$ is a free module over the affine coordinate ring 
${\bf L}$ of ${\bf T}$, one has      
\[ \Gamma \Omega^i(Z_F) \cong 
{\bf L}_f \otimes_{\bf Z} \Lambda^i M, \] 
i.e., $\Gamma \Omega^i(Z_F)$ is a free module over the affine coordinate ring 
\[ {\bf L}_f \cong {\bf L} \lbrack X_0 \rbrack / (X_0f(X) -1) \] 
of  the open subset ${\bf T} \setminus Z_f \subset {\bf T}$.
\medskip 

Our purpose now is to find the explicit description of the global 
de Rham complex. This description is a generalized version of the approach 
of Dwork \cite{dwork1,dwork2} and N.Katz \cite{katz} for 
hypersurfaces in projective spaces. First we need 
to describe the affine coordinate ring 
${\bf L} \lbrack  X_0 \rbrack / (X_0 f(X) -1)$ by means of so called 
{\em Laplace transform} ${\cal R}$ \cite{dwork.geom}.

\begin{opr}
Define the  ${\bf C}$-linear mapping 
\[ {\cal R}_0\;:\; X_0 {\bf L} \lbrack X_0 \rbrack 
\rightarrow {\bf L} \lbrack X_0 \rbrack \]
as follows 
\[ {\cal R}_0 \; : \; X_0^a X^m \mapsto (-1)^{a-1} (a-1)! X_0^a X^m, \;\; 
a \in {\bf Z}\; \;(a>0),\; m \in M. \]
Let 
\[ {\cal R} \; : \;  X_0 {\bf L} \lbrack X_0 \rbrack 
\rightarrow {\bf L} \lbrack X_0 \rbrack / (X_0f(X) -1) \] 
the composition of ${\cal R}_0$ with the homomorphism 
\[   {\bf L} \lbrack X_0 \rbrack \rightarrow {\bf L} 
\lbrack X_0 \rbrack / (X_0f(X) -1). \]
\end{opr}
\medskip

Recall that $F = F(X_0,X)$ denotes the polynomial $X_0f(X) -1$. 
\medskip

\begin{opr}
For any Laurent polynomial $f \in {\bf L}$ define 
the following commuting differential operators on the ring 
${\bf L} \lbrack X_0 \rbrack $ 
\[ {\cal D}_{X_0} = X_0 \frac{\partial}{\partial X_0} + F_0, \]
\[  {\cal D}_{X_i} = X_i \frac{\partial}{\partial X_i} + F_i\;\;  
(i =1, \ldots, n ). \] 
\end{opr}

\begin{prop}
The mapping ${\cal R}$ is surjective and the kernel of 
${\cal R}$ coincides with the image of the operator 
${\cal D}_{X_0}$ acting  on the ${\bf L}\lbrack X_0 \rbrack$-ideal 
$X_0 {\bf L} \lbrack X_0 \rbrack$. 
\label{R.map}
\end{prop}

\Proof The surjectivity of ${\cal R}$ follows from the fact that 
${\cal R}(X_0)$ is an invertible element in 
${\bf L}\lbrack X_0 \rbrack / (X_0f(X) -1)$.  

One has
\[ {\cal R} {\cal D}_{X_0} (X_0^a X^m) = a (-1)^{a-1}(a-1)! X_0^a X^m 
+ (-1)^a a! X_0^a X^m X_0f(X) \]
\[ = (-1)^a a! X_0^a X^m (X_0f(X) -1 ), \]
i.e., ${\rm Im}\,  {\cal D}_{X_0} \subset {\rm Ker}\,  {\cal R}$.  

Assume now that ${\cal R}(g(X_0, X) ) = 0$  for some element $g \in 
X_0 {\bf C} \lbrack X_0, X^{\pm1} \rbrack$. Consider the expansion of $g$ 
in powers of $X_0$ 
\[ g(X_0, X)= X_0^{a_1} g_1(X) + X_0^{a_2}g_2(X) + \cdots + 
X_0^{a_s}g_s(X), \]   
where $ a_1 > a_2 > \cdots a_s > 0$. Note that $a_1  \neq  1$. 
Indeed, assume that  $s =1$. Then $X_0g_1(X)$ is divided by $F(X_0,X)$, i.e., 
$X_0 g_1(X) = (X_0 f(X) -1) q(X)$. The latter is impossible, since 
$X_0$ can not divide $(X_0f(X) -1)q(X)$.  Thus, $a_1 > 1$. The polynomial 
\[ {\cal R}_0 (g(X_0, X)) = (-1)^{a_1 -1} (a_1 -1 )! X_0^{a_1}g_1(X) +  
\cdots + (-1)^{a_s -1} (a_s -1 )! X_0^{a_s}g_s(X) \]
is divided by $F$. Therefore, the leading  coefficient $g_1(X)$ is 
divided by $f(X)$. Let $g_1(X) = f(X) p_1(X)$. Clearly,  the maximal power 
of $X_0$ in the polynomial  
\[ g(X_0,X) - {\cal D}_{X_0} ( X_0^{a_1 -1} p_1(X) ) \]
is less than $a_1$. Now we can apply induction. 

\begin{coro}
The  mapping ${\cal R}$ induces  a ${\bf C}$-linear  isomorphism between 
the infinite dimensional ${\bf C}$-spaces
\[ V(f) =  X_0 {\bf L} \lbrack X_0 \rbrack  /  {\cal D}_{X_0} 
(X_0 {\bf L} \lbrack X_0 \rbrack)  \] 
and 
\[ \Omega^0 (Z_F) \cong {\bf L} \lbrack X_0 \rbrack / (F) \cong {\bf L}_f. \]
\label{R.iso1}
\end{coro}

\begin{theo}
For an arbitrary Laurent polynomial $f(X) \in {\bf L}$,  
the mapping ${\cal R}$ induces the isomorphism between the Koszul complex 
\[ K^{\cdot}( {\cal D}, V(f)) : = 
K^{\cdot}(\{ {\cal D}_{X_1}, \ldots, {\cal D}_{X_n} \}, V(f)) \] 
of the operators ${\cal D}_{X_i}$ $( 1 \leq i \leq n)$ 
acting on the ${\bf C}$-space $V(f)$ 
and the global de Rham complex 
$\Gamma \Omega^{\cdot} ({\bf T} \setminus Z_f)$. 
\label{derham.glob}
\end{theo}

\Proof By \ref{R.iso1}, for any $i \geq 0$ one has the isomorphisms  
\[ K^{i}( {\cal D}, V(f))  \cong 
\Gamma \Omega^{0} ({\bf T} \setminus Z_f) \otimes_{\bf Z} \Lambda^i M \cong 
\Gamma \Omega^{i} ({\bf T} \setminus Z_f). \]

The mapping ${\cal R}$ sends the element 
\[ X_0^kX^m \otimes e_I \in  
K^{i}({\cal D}, V(f)), \;\; 
\mid I \mid = i \]
to the rational differential $i$-form
\[ \frac{(-1)^{k-1} (k -1)! X^m}{f(X)^k}\frac{dX_I}{X_I}. \]
Direct calculations show that the differential of the complex 
$K^{i}( {\cal D}, V(f))$ gives rise to  
the exterior differential in $\Gamma \Omega^{\cdot} ({\bf T} \setminus Z_f)$. 
\bigskip

In this section, we did not assume so far  any regularity condition for the 
Laurent polynomial $f$. However, if $f$ is $\Delta$-regular polynomial, we can 
find in the Koszul complex 
$K^{\cdot}( {\cal D}, V(f))$  the  
more simple quasi-isomorphic subcomplex 
$K^{\cdot}( {\cal D} , V_{\Delta}(f))$, 
where $V_{\Delta}(f)$ is a $\{ {\cal D}_{X_i} \}$-invariant 
subspace in $V(f)$.   
\bigskip

From now on we always assume that $f$ is  $\Delta$-regular.
\bigskip 

Recall that the  graded ring $S_{\Delta}$ is a subring of 
${\bf L} \lbrack X_0\rbrack$.  The maximal homogeneous ideal 
$S^+_{\Delta}$ is the intersection of $S_{\Delta}$ and the principal ideal  
$X_0{\bf L} \lbrack  X_0 \rbrack$ in 
${\bf L} \lbrack  X_0 \rbrack$. It is clear that $S^+_{\Delta}$ 
is a ${{\cal D}_{X_i}}$-invariant subspace for  $ 0 \leq i \leq n$. 

\begin{opr}
Define the ${\bf C}$-linear space $V_{\Delta}(f)$ as follows  
\[ V_{\Delta}(f) = S^+_{\Delta} / {\cal D}_{X_0}(S^+_{\Delta}). \]
\end{opr}

The same arguments as in the proof  \ref{R.map} show the following.  

\begin{prop}
The restriction of the mapping ${\cal R}_0$ on the subspace  
$S^+_{\Delta}$ induces the isomorphism between $V_{\Delta}(f)$ and 
the affine coordinate ring $S_{\Delta} / F(X_0, X) S_{\Delta}$ of the 
affine variety $U_0 = {\bf P} \setminus \overline{Z}_f$. 
\end{prop}

Since $f$ is $\Delta$-regular,  
the pair $(U_0 , U_0 \cap {\bf D})$ is toroidal. Therefore,  the 
$ S_{\Delta} / (F)S_{\Delta}$-modules 
$\Omega^i_{U_0}({\rm log}\,U_0 \cap {\bf D})$ 
of  differential $i$-forms on $U_0$ 
with only logarithmic poles along  $U_0 \cap {\bf D}$ are 
well-defined. These modules are global sections of the globally 
free sheaves $\Omega^i_{\bf P,D} \cong {\cal O}_{\bf P} \otimes_{\bf Z} 
\Lambda^i M$ over the open affine chart $U_0 \subset {\bf P}$. 
So we obtain:    

\begin{prop}
There exist the natural isomorphisms of  
$S_{\Delta} /(F) S_{\Delta}$-modules  
\[ \Omega^i_{U_0}({\rm log}\,U_0 \cap {\bf D})  \cong  
S_{\Delta} / (F)S_{\Delta} \otimes_{\bf Z} \Lambda^i M. \]
\end{prop}

Using the same arguments as in \ref{derham.glob}, one proves the following.  

\begin{theo}
The mapping  ${\cal R}$ induces the isomorphism between the Koszul complex 
\[ K^{\cdot}( {\cal D}, V_{\Delta}(f))
: = K^{\cdot}(\{ {\cal D}_{X_1}, \ldots, {\cal D}_{X_n} \}, V_{\Delta}(f)) \]
of the operators ${\cal D}_{X_i}$ $( 1 \leq i \leq n)$ 
acting on the ${\bf C}$-space $V_{\Delta}(f)$ 
and the global logarithmic de Rham complex $\Gamma 
\Omega^{\cdot}_{U_0}({\rm log}\,U_0 \cap {\bf D})$. 
\label{kos.log}
\end{theo}

Now we come to the most important statement of this section.  

\begin{theo}
The natural injective morphism of  the  Koszul complexes  
\[ K^{\cdot}( {\cal D}, V_{\Delta}(f)) 
\rightarrow 
K^{\cdot}( {\cal D}, V(f))  \]
is a quasi-isomorphism.
\end{theo}

\Proof Let $\pi\;:\; {U_0'} \rightarrow {U_0}$ be a toric 
desingularization of the affine variety ${U_0}$, then 
$U_0' \cap {\bf D'}$ is a normal crossing divisor. By theorem of Deligne, 
the logarithmic de Rham complex 
\[ \Omega^{\cdot}_{U_0'}({\rm log} \,U_0'\cap {\bf D'}) \]
is quasi-isomorphic to the global Grothendieck's de Rham complex on 
${U_0'} \cap {\bf T} = {\bf T} \setminus Z_f$.  

Applying local version of lemma 6.1, we obtain that $\pi$-direct image 
of $\Omega^{\cdot}_{U_0'}({\rm log}\,U_0'\cap {\bf D'})$ is quasi-isomorphic 
to the complex $\Omega^{\cdot}_{U_0}({\rm log}\,U_0 \cap {\bf D})$ on $U_0$. 
On the other hand, since $U_0$ is an affine algebraic variety, the  
hypercohomology of  $\Omega^{\cdot}_{U_0}({\rm log}\,U_0 \cap {\bf D})$ can be 
compute by the global logarithmic de Rham complex 
\[ \Gamma \Omega^{\cdot}_{U_0}({\rm log}\,U_0 \cap {\bf D}),  \]
i.e., this complex is also quasi-isomorphic to the global Grothendieck's 
de Rham complex on ${\bf T} \setminus Z_f$. Now the statement follows from 
theorems \ref{derham.glob} and \ref{kos.log}. 
\bigskip

Let us calculate the cohomology groups of the complex 
$K^{\cdot}( {\cal D}, V_{\Delta}(f))$. 
\medskip 

Define on ${\bf L}\lbrack X_0^{\pm 1} \rbrack$ the infinite filtration 
\[ {\cal E} \;: \; 
\supset \cdots  {\cal E}^{-k} \cdots \supset {\cal E}^{-1} 
\supset {\cal E}^{0} \supset  
{\cal E}^{1} \supset \cdots,  \]
where ${\cal E}^{-k}$ is the subspace spanned by all monomials of 
the $X_0$-degree $\leq k$. 
This filtration induces  also the filtration on $S^+_{\Delta}$ 
and on the complex 
$K^{\cdot}( {\cal D}, V_{\Delta}(f))$. We denote the induced 
filtrations also by  ${\cal E}$.    

Since the action of the differential operator ${\cal D}_{X_i}$ on quotients 
${\rm Gr}_{\cal E}$ is the multiplication  by $F_i$, we obtain. 

\begin{prop}
The complex of graded ${\bf C}$-spaces  
\[ {\rm Gr}_{\cal E} \, K^{\cdot}( {\cal D}, V_{\Delta}(f)) \]
is canonically isomorphic to the Koszul complex 
\[ K^{\cdot}(\{ F_1, \ldots, F_n \}, S^+_{\Delta}/F_0 S^+_{\Delta})  \]
of graded $S_{\Delta}$-modules $($ cf. {\rm \ref{grade})}. 
\end{prop}

In \ref{grade} we have shown that  the cohomology groups of 
$K^{\cdot}(\{ F_1, \ldots, F_n \}, 
S^+_{\Delta}/F_0 S^+_{\Delta})$ have already 
the same dimensions as the cohomology groups of $Z_F$. 

\begin{coro}
The spectral sequence of the 
${\cal E}$-filtred complex 
$K^{\cdot}( {\cal D}, V_{\Delta}(f))$ degenerates,  
i.e.,  the cohomology groups of 
$K^{\cdot}(\{ F_1, \ldots, F_n \}, 
S^+_{\Delta}/F_0 S^+_{\Delta})$ are isomorphic 
to the limit cohomology groups 
\[ {\rm Gr}_{\cal E} \, H^{\cdot}( {\cal D}, V_{\Delta}(f)). \]
\end{coro}

Let us now resume our results on the $n$-th cohomology group of $Z_F =
{\bf T} \setminus Z_f$. 
\medskip

Denote also by ${\cal R}$  the ${\bf C}$-linear mapping 
\[ {\cal R}\; :\; 
S^+ \rightarrow \Omega^n({\bf T} \setminus  Z_f) \]
defined as   
\[ {\cal R}(X_0^k{X}^{m}) = \frac{ (-1)^{k-1} (k-1)!\, {X}^{m}}{
 f(X)^k} \frac{dx_1}{x_1} \wedge \cdots \wedge \frac{dx_n}{x_n}. \]

\begin{theo}
There exists the following commutative diagram
\begin{center}
\begin{tabular}{ccc}
$S^+$ & ${\rightarrow}$ & $S^+/ ({\cal D}_{X_0} S^+  + 
\sum_{i=1}^n {\cal D}_{X_i} S^+)$ \\
 ${\cal R} \downarrow$  &  & $\rho \downarrow$  \\
${\Omega}^n({\bf T} \setminus Z_f)$ & $\;\;\;{\rightarrow}$ & 
$H^{n}({\bf T} 
\setminus Z_f)$ 
  \end{tabular}
\end{center}
such that $\rho$ is an isomorphism. 

The filtration 
${\cal E}$ on $S^+_{\Delta}$ induces via ${\cal R}$ the filtration 
on $\Omega^n({\bf T} \setminus Z_f )$, where $ {\cal E}^{-k} 
\Omega^n({\bf T} \setminus Z_f )$ consists of $n$-forms with poles at most 
$k$ along $Z_f$. Moreover, there exists the canonical isomorphism 
\[ {\rm Gr}_{\cal E}\,H^{n}({\bf T} \setminus Z_f) 
\cong S^1_{\Delta} \oplus R^2_f \oplus 
\cdots \oplus R_f^n. \]
\label{com.diag}
\end{theo}

\section{The mixed Hodge structure}


\hspace*{\parindent}

In the previous section we have found the isomorphism $\rho$ between 
the space 
\[S^+/ ({\cal D}_{X_0} S^+  + \sum_{i=1}^n {\cal D}_{X_i} S^+) \]
and $ H^n(Z_F)$. Let us now describe the Hodge and weight filtrations on 
this cohomology group. 
\medskip

\begin{theo}
The ${\cal E}$-filtration on $H^n(Z_F)$ coincides with the Hodge filtration. 
\label{hodge.filt}
\end{theo}

\Proof Recall that if we consider the cohomology $H^{\cdot}(Z_F)$ 
as hypercohomology   
\[ {\bf H}^{\cdot}\Omega^{\cdot}_{\bf P} ({\rm log}\,D), \]
then the Hodge filtration is induced by the "stupid" filtration on 
$\Omega^{\cdot}_{\bf P}({\rm log}\,D)$. 

On the other hand, the same cohomology groups are computed by the 
hypercohomology  of the complex 
\[ \jmath_* \Omega^{\cdot}_{U_0} ({\rm log}\, U_0 \cap {\bf D}),\]
where $\jmath$ is the open embeddng $U_0 \hookrightarrow {\bf P}$. 

The natural injective homomorphism of complexes 
\[  \Omega^{\cdot}_{\bf P}({\rm log}\,D) \rightarrow  
j_* \Omega^{\cdot}_{U_0}({\rm log}\, U_0 \cap {\bf D}) \]
is a quasi-isomorphism. Since 
$H^i({\bf P}, j_* \Omega^{\cdot}_{U_0}({\rm log}\, U_0 
\cap {\bf D})) = 0$ for 
$i > 0$, one has 
\[ {\bf H}^n (j_* \Omega^{\cdot}_{U_0}({\rm log}\, U_0 \cap {\bf D})) \cong 
H^n (j_* \Omega^{\cdot}_{U_0}({\rm log}\, U_0 \cap {\bf D}))  \]
\[ \cong S^+/ ({\cal D}_{X_0} S^+  + \sum_{i=1}^n {\cal D}_{X_i} S^+). \]

The hypercohomology groups can be calculated using the \v{C}ech-de Rham 
complex corresponding  to the Jacobian covering 
\[ {\cal U}\;:\; U_i\;\; (0 \leq i \leq n)\]
of ${\bf P}_{\Delta}$. 
Thus we obtain two quasi-isomorphic  double complexes 
\[ {\cal C}^q({\cal U}, \Omega^{p}_{\bf P} ({\rm log}\,D)) \rightarrow 
{\cal C}^q({\cal U}, j_* \Omega^{p}_{U_0}({\rm log}\, U_0 \cap {\bf D})). \]

We can now apply the Carlson-Griffiths method \cite{carl}. Working inside 
the \v{C}ech-de Rham complex 
${\cal C}^{\cdot}({\cal U}, j_* \Omega^{\cdot}_{U_0}({\rm log}\, 
U_0 \cap {\bf D})) $, we represent a rational $n$-form with polar 
oder $k+1$ by a \v{C}ech-de Rham cocycle 
\[ \tilde{\omega}_0 + \cdots + \tilde{\omega}_k \in 
\bigoplus_{i \leq k} {\cal C}^i ({\cal U}, 
\Omega^{n+1 - i}_{\bf P}({\rm log}\, D) ).\]
By this way we obtain that 
the ${\cal E}$-filtration is compatible with the Hodge filtration. Since 
${\rm Gr}_{\cal E}\,H^{n}(Z_F) \cong {\rm Gr}_{\cal F}\,H^{n}( Z_F)$, 
the Hodge filtration and ${\cal E}$-filtration must coincide. This proves 
theorem. (One can also prove the compatibility of ${\cal E}$-filtration and 
the Hodge filtration by another method explained in \ref{e-filt}.) 
\bigskip

\begin{theo}
Let  
\[ 0 = {\cal W}_n \subset {\cal W}_{n+1} \subset \cdots \subset 
{\cal W}_{2n-1} 
\subset {\cal W}_{2n} = H^n({\bf T} \setminus Z_f) \]
be the weight filtration. Then for $0 < i < n$ the subspace 
${\cal W}_{n+i}$ equals $\rho({\cal I}^{(i)})$, where 
${\cal I}^{(i)}$ is the image of the ideal $I^{(i)}_{\Delta} 
\subset S^+_{\Delta}$ in the  space  $S^+/ ({\cal D}_{X_0} S^+  + 
\sum_{i=1}^n {\cal D}_{X_i} S^+)$.
\label{ideal.filt}
\end{theo}

\Proof In our situation all the pairs $({\bf P}, D)$ and 
$(U_0, U_0 \cap {\bf D})$ are toroidal. So we can use the local definitions 
of the weight filtrations on sheaves of logarithmic differentials on affine 
toric varieties \cite{dan2}. 

Recall that if we consider the cohomology 
$H^{\cdot}({\bf T} \setminus 
Z_f)$ 
as hypercohomology   
\[ {\bf H}^{\cdot}\Omega^{\cdot}_{\bf P}({\rm log}\,D), \]
then the weight  filtration ${\cal W}$ on the cohomology 
is induced by the weight filtration ${\cal V}$ of  
sheaves 
\[ {\cal V}_0 \Omega^{i}_{\bf P}({\rm log}\,D) \subset 
\cdots \subset {\cal V}_i \Omega^{i}_{\bf P}({\rm log}\,D), \]
where 
\[ {\cal V}_k \Omega^{i}_{\bf P}({\rm log}\,D) = 
\Omega^{i-k}_{\bf P} \wedge  \Omega^{k}_{\bf P}({\rm log}\,D). \]
By definition, one has 
\[ {\cal W}_{k+i} H^i({\bf T} \setminus Z_f) = 
{\cal V}_k H^i({\bf T} \setminus Z_f).\]

On the other hand, the complex 
$\Omega^{\cdot}_{U_0} ({\rm log}\, U_0 \cap {\bf D})$ also has the  
filtration ${\cal V}'$ 
coming from logarithmic singularities along $U_0 \cap {\bf D}$:
\[ {\cal V}_k' \Omega^{i}_{U_0}({\rm log}\,U_0 \cap {\bf D}) = 
\Omega^{i-k}_{U_0} \wedge  \Omega^{k}_{U_0}({\rm log}\,U_0 \cap {\bf D}). \]

Since the complex $\Omega^{\cdot}_{U_0} ({\rm log}\, U_0 \cap {\bf D})$ 
computes cohomology of ${\bf T} \setminus Z_f$, we obtain a 
${\cal W}'$-filtration on the cohomology defined as 
\[  {\cal W}_{k+i}' H^i({\bf T} \setminus Z_f) = 
{\cal V}_k' H^i({\bf T} \setminus Z_f).\] 
Althought, in general the ${\cal W}'$-filtration must not coincide with the 
"correct" ${\cal W}$-filtration, we can compare these two filtrations. 

\begin{lem}
Define ${\cal V}''$-filtration on  $\Omega^{\cdot}_{\bf P}({\rm log}\,D)$  
as follows
\[ {\cal V}_k'' \Omega^{i}_{\bf P}({\rm log}\,D) = 
\Omega^{i-k}_{\bf P} \wedge  \Omega^{k}_{\bf P,T} + 
({\cal V}_{k-1}'\,  
\jmath_* \Omega^{i}_{U_0}({\rm log}\,U_0 \cap {\bf D})) \cap (
\Omega^{i}_{\bf P}({\rm log}\,D)).  \]
Then 
\[ {\cal V}_k \Omega^{i}_{\bf P}({\rm log}\,D) \subset 
{\cal V}_k'' \Omega^{i}_{\bf P}({\rm log}\,D). \]
\end{lem}

\Proof The statement follows from the fact that locally any 
differential $i$-form $\omega \in {\cal V}_k \Omega^{i}_{\bf P}({\rm log}\,D)$ 
can be represented as a sum 
\[ \omega = \omega_1 + \omega_2, \]
where $\omega_1 \in  \Omega^{i-k}_{\bf P} \wedge  \Omega^{k}_{\bf P,T}$ and 
$g\cdot \omega_2 \in {\cal V}_{k-1} \Omega^{i}_{\bf P}({\rm log}\,D)$ 
($g$ is a local 
equation of $\overline{Z}_f$ in ${\bf P}$).  
We note that one has the embedding 
\[ \Omega^{i-k}_{\bf P} \wedge  \Omega^{k}_{\bf P}({\rm log}\,D) 
\hookrightarrow  
\jmath_* \Omega^{i-k}_{U_0} \wedge  
\jmath_* \Omega^{k}_{U_0}({\rm log}\,U_0 \cap {\bf D}).  \]
\medskip

\begin{coro}
Define ${\cal W}''$-filtration on cohomology as 
\[  {\cal W}_{k+i}'' H^i({\bf T} \setminus Z_f) = 
{\cal V}_k'' H^i({\bf T} \setminus Z_f).\] 
Then one has 
\[  {\cal W}_{k+i} H^i({\bf T} \setminus Z_f) \subset 
 {\cal W}_{k+i}'' H^i({\bf T} \setminus Z_f) \;\; {\rm for}\;\;  i \geq 0. \] 
\end{coro}
 
Since the image of the hypercohomology of $\Omega^{\cdot}_{\bf P,D}$ 
in the hypercohomology of   $\Omega^{\cdot}_{\bf P}({\rm log}\,D)$ coincides 
with the image of $H^{\cdot}({\bf T})$ in $H^{\cdot}({\bf T} \setminus Z_f)$, 
we obtain also. 

\begin{coro}
On the primitive parts of cohomology groups, the ${\cal W}''$-filtration 
coincides with the ${\cal W}'$-filtration 
and it is compatible with the standard ${\cal W}$-filtration: 
\[ {\cal W}_{k+i} PH^i({\bf T} \setminus Z_f) \subset 
{\cal W}_{k+i}'' H^i({\bf T} \setminus Z_f) =  
{\cal W}_{k+i-1}' H^i({\bf T} \setminus Z_f). \]
\label{compat}
\end{coro}

We want now to  calculate the ${\cal W}'$-filtration explicitly. The 
"interesting" case is  only the corresponding ${\cal V}'$-filtration on the 
$S_{\Delta}/(F)S_{\Delta}$-module 
$\Omega^{n}_{U_0}({\rm log}\,U_0 \cap {\bf D})$ whose elements represent 
classes in $H^n({\bf T} \setminus Z_f)$. 

The next statement immediately follows from results in \cite{dan2}. 

\begin{lem}
The weight filtration on the sheaf  $\Omega^n_{\bf P,T}$ of $n$-forms on 
${\bf P}$ with logarithmic poles along ${\bf D} = {\bf P} \setminus {\bf T}$ 
consists of subsheaves ${\cal S}_i$ $(0 \leq i \leq n)$  
associated with the homogeneous $S_{\Delta}$-ideals  $I_{\Delta}^{(i+1)}$  
$($cf. {\rm (\ref{chain}))}. 
\label{toric.log}
\end{lem}

\medskip

Since $S_{\Delta}/(F)S_{\Delta}$-module 
$\Omega^{n}_{U_0}({\rm log}\,U_0 \cap {\bf D})$ consists 
 of global sections of 
$\Omega^n_{\bf P,T}$ over the affine open subset $U_0$, we conclude. 

\begin{lem}
The subsheaf 
\[ {\cal V}_i' \Omega^{n}_{U_0}({\rm log}\,U_0 \cap {\bf D}) 
\subset 
\Omega^{n}_{U_0}({\rm log}\,U_0 \cap {\bf D}) \] 
is the restriction of the sheaf ${\cal S}_i$ on $U_0$.  Moreover, 
\[ {\cal V}_{n-1}' \Omega^{n}_{U_0}({\rm log}\,U_0 \cap {\bf D}) = 
{\cal V}_{n}' \Omega^{n}_{U_0}({\rm log}\,U_0 \cap {\bf D}) \]
\[ =  \Omega^{n}_{U_0}({\rm log}\,U_0 \cap {\bf D}). \]
\end{lem}

Applying \ref{toric.log} and  \ref{compat}, we obtain: 

\begin{coro}
The weight subspace ${\cal W}_{n +k}$ on $H^{n}({\bf T} \setminus Z_f)$ 
is contained in the $\rho({\cal I}^{(k)})$. 
\end{coro}

Thus  we have proved the compatibility of the weight filtration and the 
filtration induced by the chain of the homogeneous ideals (\ref{chain}). 
In order to prove theorem \ref{ideal.filt}, it is suffices to check that
these two filtrations on $H^n({\bf T} \setminus Z_f)$ has the same 
image in ${\rm Gr}_{\cal F}H^n({\bf T} \setminus Z_f)$. The latter follows 
from  the following statement. 

\begin{lem}
Define the ${\cal V}$-filtration on the sheaves $\Omega^p_{\bf P,T} (k)$ 
as follows 
\[ {\cal V}_i \Omega^p_{\bf P,T} (k) = 
\Omega^{p-i+1}_{\bf P} \wedge  \Omega^{i-1}_{\bf P,T}(k). \]
Then the iduced ${\cal V}$-filtration on the resolution 
$C^{p,\cdot}_{\rm log} (D)$ can be considered as a resolution of the 
${\cal V}$-filtered sheaf $\Omega^p_{\bf P} ({\rm log}\, D)$, 
and the filtration  on cohomology groups 
$H^{n-p}({\bf P}, \Omega^p_{\bf P} ({\rm log}\, D)$ 
coincides with the filtration 
defined by the chain of homogeneous 
ideals $(\ref{chain})$ on homogeneous components of the Jacobian ring $R_f$.   
\end{lem}

The proof of the lemma is based on the checking of the corresponding local statement 
for affine toric models. 
\bigskip

\begin{coro}
The weight filtrations on the space $L(\Delta)$ considered as 
the minimal nonzero Hodge subspace 
${\cal F}^n H^{n}({\bf T} \setminus Z_f)$ can be defined as follows. 
For $k < n$ the space 
\[ {\cal W}_{n + k}  L(\Delta) \]
is spanned by all monomials $X^m$ $(m \in \Delta \cap {\bf Z}^n)$ such 
that $m$ does not belong to any face $\Delta'$ of codimension $k$. 
\label{n-filt}
\end{coro}

\section{Hodge structures  of weight $n-1$}


\hspace*{\parindent}

\begin{opr}
{\rm Denote by $H_f$ be the image of the ideal $I^{(1)}_{\Delta}$ in $R_f$.}
\end{opr}

As a consequnce of our  results in \ref{hodge.filt} and 
\ref{ideal.filt}, we obtain. 

\begin{prop}
There exists the isomorphism 
\[ H_f \cong {\rm Gr}_{\cal F} {\cal W}_{n-1} H^{n-1}(Z_f). \]
The homogeneous components $H_f^i$ of the graded space 
\[ H_f = H^1_f \oplus \cdots \oplus H^n_f \]
are isomorphic to the Hodge-Deligne spaces $H^{n-i,i-1} (H^{n-1}(Z_f))$. 
In particular, one has 
\[ {\rm dim}\, H^i_f = {\rm dim}\, H^{n +1 -i}_f\;\; ( 1 \leq i \leq n).  \]
\end{prop}

We would like now to decribe the duality between the spaces 
$H_f^i$ and $H_f^{n+1 -i}$ more precisely.

\begin{opr}
{\rm Denote by $D_f$ the quotient of the homogeneous 
$S_{\Delta}$-ideal  $I^{(1)}_{\Delta}$ by the regular sequence 
$F_0, F_1, \ldots, F_n$. }
\end{opr}

By well-known properties of dualizing modules over Cohen-Macaulay rings, 
one has.  

\begin{prop}
The Hilbert-Poincare series of $D_f$ is the polynomial $\Phi_{\Delta}(t)$ 
of degree $n+1$ $($cf. \ref{hilb}$)$.  
The graded Artinian $R_f$-module  $D_f$ is the dualizing $R_f$-module.  
The socle of $D_f$ is its maximal nonzero homogeneous component 
$D_f^{n+1}$. The multiplicative $R_f$-module structure 
defines the nondegenerate pairing
\[ \langle *, * \rangle \; \; :\;\; 
R^i_f \times D_f^{n+1 -i} \rightarrow D_f^{n+1} \cong {\bf C}. \]
\label{duality}
\end{prop}

\begin{rem}
{\rm Let $PH^{n-1}_c(Z_f)$ denotes the kernel of the surjective 
homomorphism 
\[H^{n-1}_c(Z_f) \rightarrow H_c^{n-1}({\bf T}). \]
Using local duality theorem, we can describe the quotients of the 
Hodge filtration on $PH^{n-1}_c(Z_f)$ as follows
\[ {\cal F}^{i}PH^{n-1}_c(Z_f)/ {\cal F}^{i+1}PH^{n-1}_c(Z_f)  
\cong D_f^{n-i}. \] 
The Poincare duality
\[{\rm Gr}_{\cal F}\, PH^{n-1}(Z_f) 
\otimes {\rm Gr}_{\cal F}\, PH^{n-1}_c (Z_f) \rightarrow {\bf C} \]
is induced by pairings
\[ R_f^{i} \otimes D_f^{n+1-i} \rightarrow D_f^{n+1} \cong {\bf C}. \]}
\end{rem}
\medskip

By definition of the space $H_f$, one gets. 

\begin{prop}
The {\bf C}-space $H_f$ is the image 
of the homomorphism 
\[ h \;\; :\;\; D_f \rightarrow R_f \]
induced by the embedding $I_{\Delta}^{(1)} \hookrightarrow S_{\Delta}$. 
\end{prop}

We can now calculate the nondegenerate pairing between the 
Hodge-Deligne spaces $H^i_f$ and $H^{n +1 - i}_f$. 

\begin{prop}
Let $a \in H^i_f$, $b \in H^{n+1 -i}_f$. Take 
$a' \in D^i_f$, $b' \in D^{n+1 -i}_f$  such that $h(a' ) = a,$ 
$h(b') = b$.  Then two elements 
$\langle a, b' \rangle$ and $\langle b, a' \rangle$ coincide, and the pairing 
\[ \{ *, * \} \;\; : \;\; H^{i}_f \times H^{n+1-i}_f \rightarrow {\bf C} \]
\[ \{ a, b \} = \langle a, b' \rangle = \langle a', b \rangle \]
is nondegenerate. 
\end{prop}

\Proof Choose the elements  $x, y \in I_{\Delta}^{(1)} \subset S_{\Delta}$ 
representing $a',b'$ in $D_f$.  Let $z = xy \in S_{\Delta}$. Then $z \in 
I^{(1)}_{\Delta}$, and the image of $z$ in $D_f^{n+1}$ concides with 
$\langle a, b' \rangle$ and $\langle a', b \rangle$. The nondegeneracy 
of $\{*, *\}$ follows from the nondegeneracy of  $\langle *, * \rangle$ 
(cf. \ref{duality}). 
\bigskip

From \ref{n-filt} we obtain that  ${\bf C}$-space 
\[ H^{n-1, 0}(H^{n-1}(Z_f)) \cong  H^1_f \cong 
( H^{n}_f)^*\]  
is  isomorphic to the subspace $L^*(\Delta) \subset L(\Delta)$ 
spanned by all monomials $m \in \Delta \cap {\bf Z}^n$  such that  
$m$ is in the interior of $\Delta$. Therefore, the dimension of $L^*(\Delta)$ 
is $l^*(\Delta)$. 
\medskip

Our next purpose is to describe explicitly the  Hodge-Deligne subspace  
\[ H^{n-2,1}(H^{n-1}(Z_f)) \cong H^2_f \cong (H^{n-1}_f)^* .\]   
\medskip

\begin{theo}
There exists the following short exact sequence of ${\bf C}$-spaces 
\[  0 \rightarrow \bigoplus_{{\rm codim}\, \Delta'=1} L^*(\Delta')  
\rightarrow L^*(2\Delta)/(f, f_1, \ldots, f_n)L^*(\Delta) 
\rightarrow H^{n-2,1}(H^{n-1}(Z_f)) \rightarrow 0, \]
where $L^*(2\Delta)$ is the space spanned by all monomials $X^m \in 
L(2\Delta)$ such that $m$ belongs to the interior of the polyhedron 
$2\Delta$. 
\label{hodge.sub}
\end{theo}

\Proof Since the homogeneous component $R_f^{n-1}$ is isomorphic to 
the space 
\[ {\cal F}^1 H^{n-1}(Z_f)/{\cal F}^2 H^{n-1}(Z_f), \] 
the weight 
filtration on this space has the  following properties 
\[ {\cal W}_{n-1} R^{n-1}_f = H_f^{n-1} \subset 
{\cal W}_n R^{n-1}_f = R^{n-1}_f. \]
Hence we obtain the short exact sequence
\[ 0 \rightarrow H_f^{n-1} \rightarrow  R^{n-1}_f \rightarrow 
H^{0,n-2}(H^{n-1}(Z_f)) \rightarrow 0. \]

Take the dual exact sequence
\[ 0  \rightarrow H^{n-2,0}(H^{n-1}(Z_f))  \rightarrow  
D^{2}_f \rightarrow H_f^{2}  \rightarrow 0. \]   
The statement follows now from the isomorphisms 
\[ D_f^2 \cong L^*(2\Delta)/(f, f_1,  \ldots, f_n)L^*(\Delta), \]
and 
\[ H^{n-2,0}(H^{n-1}(Z_f)) \cong \bigoplus_{{\rm codim}\, 
\Delta'=1} L^*(\Delta') \]
(cf. \ref{n-filt}).

\section{Moduli spaces of affine hypersurfaces} 


\hspace*{\parindent}

Let $A = \Delta \cap {\bf Z}^n$ the set of integral points 
in an $n$-dimensional polyhedron $\Delta$. Then denote by ${\bf C}^{A}$ the 
affine complex space parametrizing coefficients of Laurent 
polynomials 
\[ f(X) = \sum_{m \in A}\, a_m X^m. \]
Let us ask the following question:

\begin{ques}
What is the moduli space of affine hypersurfaces $Z_f$ in ${\bf T}$ 
with the  fixed  Newton polyhedron $\Delta$ $?$ 
\label{moduli}
\end{ques}

There is another natural question connected with  \ref{moduli}:

\begin{ques}
When two hypersurfaces $Z_{f_1}, Z_{f_2} \subset {\bf T}$ are isomorphic$?$ 
\label{isom}
\end{ques}

Of course, complete answers to the questions \ref{moduli} and \ref{isom} are 
rather subtle problems. However, we will try to solve them using a naive 
approach, i.e., we assume that the our "moduli spaces" are quotients  of 
${\bf C}^A$ by "appropriate" group of authomorphisms.  In our case this 
group of authomorphisms will be $(n+1)$-dimensional algebraic torus 
$\overline{\bf T} \cong {\bf C}^* \times {\bf T}$. We note 
that affine hypersurfaces $Z_{f_j}$ $(j =1,2)$ defined by 
Laurent polynomials  
\[  f_j(X) = \sum_{m \in A} a_m^{(j)} X^m \]
are isomorphic in the following two cases:

{\em Case }I. There exists a nonzero complex number $\lambda$ such that  
for all $ m \in A$ one has $a_m^{(2)} = \lambda a_m^{(1)}$. 

{\em Case }II. There exists an element ${\bf t} = (t_1, \ldots, t_n) \in 
{\bf T}$ such that for all $m \in A$ one has  
$a_n^{(2)} = {\bf t}^m a_m^{(1)}$. 

These two cases define the  actions  of ${\bf C}^*$ and  
${\bf T}$ on ${\bf C}^A$. 

\begin{opr}
{\rm Denote by ${\bf P}(A)$ the projective space obtained from quotient 
of ${\bf C}^A \setminus \{0 \}$ by ${\bf C}^*$. }
\end{opr}

\begin{opr}
{\rm Define the moduli space ${\cal M}(\Delta)$ 
of affine hypersurfaces $Z_f \subset {\bf T}$ as projective spectrum of 
the graded ${\bf T}$-invariant subring ${\bf C}\lbrack a_m \rbrack^{\bf T}$  
in the polynomial ring ${\bf C}\lbrack a_m \rbrack$. In other words, 
${\cal M}(\Delta)$ is a quotient of ${\bf P}(A)$ by ${\bf T}$.}
\end{opr}

By geometric invariant theory \cite{mam}, the moduli space 
${\cal M}(\Delta)$ set-theoretically consists of the equivalence classes 
of semistable ${\bf T}$-orbits on ${\bf P}(A)$ relative to a {\em 
linearization } of the {\bf T}-action. Now we want to explain what is 
the ${\bf T}$-linearization in our case and why the {\em choice of 
the linearization}  is very important for the construction of the moduli 
spaces. 
\medskip

First note that we have "missed"  the following third  case:  

{\em Case} III. If there exists a monomial $X^{m'}$ such that 
$f_2(X) = X^{m'}f_1(X)$, then  two affine hypersurfaces $Z_{f_1}$ 
and $Z_{f_2}$ are isomorphic.
\medskip

One should remark that in the third  case the Newton polyhedra 
$\Delta_1 = \Delta(f_2)$ and  $\Delta_2 = \Delta(f_1)$ are not 
{\em the same}, but they are isomorphic. The polyhedron $\Delta_2$  
is a shift of  $\Delta_1$ by the integral vector $m \in M \cong {\bf Z}^n$. 
The crucial observation is that members of the family of affine hypersurfaces 
defined by polyhedra $\Delta_2$  are isomorphic to members of 
family defined  by $\Delta_1$ while the projective 
varieties ${\cal M}(\Delta_2)$ and ${\cal M}(\Delta_1)$ are {\em not 
isomophic if $m ' \neq 0$}. Moreover, these spaces might have even 
different dimensions. In this situation ${\cal M}(\Delta_2)$ and 
${\cal M}(\Delta_1)$ are GIT-quotients of ${\bf P}(A)$ by ${\bf T}$ 
corresponding to {\em two different ${\bf T}$-linearzation}. The first 
linearization  can be obtained from the  second one by twisting with the 
character of ${\bf T}$ corresponding to $m \in M$. In order to fix a 
linearization it is suffices to fix the location of the zero point 
$0 \in M_{\bf R}$ relative to a polyhedron $\Delta$.  
\medskip

Recall definitions of semistability and stability for our situation.  
\medskip

\begin{opr}
{\rm Let  $x$ be a point in ${\bf P}(A)$, $v$   a non-zero 
vector $v \in {\bf C}^A$ representing $x$. The point $x$ is called 
{\em semistable} if the ${\bf T}$-orbit ${\bf T}v$ does not contain 
$\{0 \} \in {\bf C}^A$ in its closure. The point $x$ is called {\em stable} if 
${\bf T}v$ is closed and has dimension $n$. }
\end{opr} 

\begin{prop} {\rm (\cite{kap} 3.5)} 
The set ${\cal M}^s(\Delta)$ of stable points is non-empty if and only if 
the zero $0$ of the lattice $M \cong {\bf Z}^n$ 
is in the interior of the polyhedron $\Delta$. 
\label{stable}
\end{prop}

Let ${\rm Sec}(\Delta)$ be the {\em secondary polyhedron} corresponding 
to $\Delta$ \cite{gel1}. The dimension of the polyhedron 
${\rm Sec}(\Delta)$ is equal to $l(\Delta) - n -1$. 
It was proved in \cite{kap} 
that the toric variety ${\bf P}_{{\rm Sec}(\Delta)}$  satisfies an universal 
property relative to all moduli spaces ${\cal M}(\Delta)$ constucted 
by means of the geometric invariant theory. Vertices of ${\rm Sec}(\Delta)$ 
one-to-one correspond to {\em regular triangulations}  of the polyhedron 
$\Delta$.  The secondary polyhedron is the Newton polyhedron of 
the principal discriminant ${\rm Disc}_{\Delta} (f)$. Thus we obtain. 

\begin{prop}
Let $\tilde{M}(\Delta)$ be the toric variety ${\bf P}_{{\rm Sec}(\Delta)}$ 
corresponding to the secondary polyhedron ${\rm Sec}(\Delta)$. Then 
the principal discriminant function of Laurent polynomials $f \in 
L(\Delta)$  is a global section of the ample sheaf 
${\cal O}_{{\rm Sec}(\Delta)}(1)$ on $\tilde{M}(\Delta)$. 
\end{prop}

\section{Variations  of Hodge structures}


\hspace*{\parindent}

\begin{opr}
{\rm Let $\Delta$ be an integral polyhedron. 
Denote by $U_{\rm reg}$ the Zariski open subset of the $L(\Delta)$ consisting 
of $\Delta$-regular Laurent polynomials with the Newton polyhedron $\Delta$. }
\end{opr}

Following the ideas of Griffiths \cite{griff1}, we can observe that 
the Jacobian ring $R_f = S_{\Delta}/ J_{f,\Delta}$ contains also a local 
information on the moduli 
space of affine $\Delta$-regular hypersurfaces $Z_f$ in ${\bf T}$. 

\begin{prop}
Consider the action of the $(n + 1)$-dimensional torus 
$\overline{\bf T}$ on $L(\Delta)$
\[ \overline{\bf T} \times L(\Delta) \rightarrow L(\Delta)\; :\; 
(\lambda, {\bf t}) \times f(X) \mapsto \lambda f({\bf t}X). \]
Let $\Theta_f$ be the tangent space to the orbit $\overline{\bf T}f$ 
at $f$. If we identify $L(\Delta)$ with $S_{\Delta}^1$, then 
\[ \Theta_f \cong J^1_{f}. \]
\end{prop}

\Proof The Lie algebra of the toris ${\bf T}$ is canonically isomorphic to 
$\overline{N}_{\bf C}$. If we identify $L(\Delta)$ with $S_{\Delta}^1$, then 
$f(X)$ corresponds to $X_0f(X)$. There is a  basis of 
$\overline{N}_{\bf C}$ consisting of  
\[ X_0\frac{\partial}{\partial X_0},  
X_1\frac{\partial}{\partial X_1}, \ldots , X_n\frac{\partial}{\partial X_n}. \]
Applying these vector fields to $X_0f(X)$,  we obtain the elements 
$F_0, F_1, \ldots, F_n$  spanning  the ${\bf C}$-space $J^{1}_f$.

\begin{coro}
Take  $f \in U_{\rm reg} \subset L(\Delta)$. Assume that 
$\lbrack f \rbrack$ is a smooth stable point on the moduli space 
${\cal M}(\Delta)$.  
Then the 
homogeneous component $R^1_f$ is canonically to the tangent space 
$\Theta_{\lbrack f \rbrack}$ 
at $\lbrack f \rbrack  \in {\cal M}(\Delta)$.  
\end{coro}

\Proof Since $F_0, F_1, \ldots, F_n$ are linearly independent, the dimension 
of the $\overline{T}$-orbit $\overline{T}f$ is $n +1$. The ${\bf C}$-space 
$S_{\Delta}^1$ can be identified with the tanget space to $L(\Delta)$. 
The dimension of $S^1_{\Delta} / J_f^1 \cong R_f^1$ is $l(\Delta) - n -1$. 
The space $R_f^1$ has canonical morphism to the tangent space 
$\Theta_{\lbrack f \rbrack}$. 
Since $\lbrack f \rbrack$ is a smooth point, its tangent space 
must be isomorphic to $R_f^1$. 
\bigskip 

We can describe the tangent space to $L(\Delta)$ by means of derivations. 
\medskip

Consider the universal affine hypersurface $Z_{f(\alpha, X)}$ 
in ${\bf T}$ defined by the Laurent polynomial 
\[ f(X) = \sum_{m \in A } a_m X^m \]
with algebraically independent coefficients $a_m$ 
$( A = \Delta \cap {\bf Z}^n)$. 
The universal hypersurface can be considered as a family of 
hypersurfaces over the affine space 
\[ {\bf C}^A : = {\rm Spec}\, {\bf C} 
\lbrack a_m \rbrack_{m \in A }. \]
 There exist $l(\Delta)$ linearly 
 independent derivations $\partial/ \partial a_m$ 
of ${\bf C} \lbrack \alpha \rbrack$ 
$({\bf C} \lbrack \alpha \rbrack = 
{\bf C} \lbrack a_m \rbrack_{m \in A} )$. 
\medskip

We can obviously extend the differential operators ${\cal D}_{X_0}, 
{\cal D}_{X_1}, \ldots, {\cal D}_{X_n}$ to universal operators acting 
on $S_{\Delta} \lbrack \alpha \rbrack$. We denote these extensions 
by the same letters:  
\[ {\cal D}_{X_0} := X_0 \frac{\partial}{\partial X_0} +  X_0f(\alpha,X), \]
\[  {\cal D}_{X_i} := X_i \frac{\partial}{\partial X_i} + X_0f_i(\alpha,X) 
\;\; (i =1, \ldots, n ). \] 

On the other hand, we define new $l(\Delta)$ differential operators.

\begin{opr}
{\rm Define the differential operators 
\[ {\cal D}_{a_m}, \;\;\; m \in A  \]
acting on 
$S_{\Delta}  \lbrack \alpha \rbrack$ as: 
\[ {\cal D}_{a_m} : = \frac{\partial}{\partial a_m} 
+ X_0X^{m}. \]
One can give also a formal representation of the operators ${\cal D}_{a_m}$ 
as follows
\[ {\cal D}_{X_{a_m}} = 
{\rm exp}(-X_0f(\alpha, X)) \frac{\partial}{\partial a_m} 
{\rm exp}(X_0f(\alpha,X)). \] }
\label{opr.diff}
\end{opr}

One gets immediately from the definition. 

\begin{prop}
Consider $S_{\Delta}  \lbrack \alpha \rbrack$ as an 
${\bf C} \lbrack \alpha \rbrack$-module of infinite rank with a basis 
$\{ X_0^kX^l\}$ $(l = (l_1, \ldots, l_n) \in k\Delta)$.  Then the action of 
a differential operator ${\cal D}_{a_m}$  on a  monomial $X_0^kX^l$ is 
the  multiplication by $X_0X^m$. 
\label{multi}
\end{prop}
\medskip

By \ref{com.diag}, the space 
\[H_{\rm DR}^{n-1}(\alpha) \cong S^+_{\Delta} \lbrack \alpha \rbrack / 
({\cal D}_{X_0} S_{\Delta}^+ \lbrack 
\alpha \rbrack  + \sum_{i=1}^n {\cal D}_{X_i} S^+_{\Delta}\lbrack \alpha \rbrack) \]
is a ${\bf C} \lbrack \alpha \rbrack$-module representing 
the $(n-1)$-th relative de Rham cohomology group 
$H^{n-1}(Z_{f(\alpha,X)})$ of the universal hypersurface $Z_{f(\alpha,X)}$ 
over ${\bf C}^A$.  
The action of the derivations $\partial/ \partial a_m$ is  the 
Gauss-Manin connection on $H_{\rm DR}^{n-1}(\alpha)$. 
\medskip

\begin{theo}
The action of a derivation $\partial/ \partial a_m$ on 
$H_{\rm DR}^{n-1}(\alpha)$ is induced by the action of the differential 
operator ${\cal D}_{a_m}$ on $S_{\Delta} \lbrack \alpha \rbrack$. 
\label{diff.act}
\end{theo}

\Proof The statement follows from the direct checking of the property 
${\cal R}{\cal D}_{a_m} = \partial/\partial a_m {\cal R}$. 
\medskip

\begin{rem}
{\rm Using our description of the Gauss-Manin connection, one gets another 
proof of the compatibility of the ${\cal E}$-filtration and the Hodge 
filtration on $H^n({\bf T} \setminus Z_f)$ in the case when 
$S_{\Delta}$ is generated by $S^1_{\Delta}$. Indeed, by results in sect. 6, 
the minimal Hodge subspace ${\cal F}^n \subset H^n({\bf T} \setminus Z_f)$ is 
represented by rational differential forms on ${\bf T}$ having poles of 
order at most $1$ along $Z_f$. A derivation  of a rational differential 
$n$-form  increases the order of its pole at most by $1$. 
On the other hand, by Griffiths transversality property, a derivative      
of a rational differential $n$-form representing an element in ${\cal F}^i$ 
should represent an element in ${\cal F}^{i-1}$.  
By \ref{multi} and \ref{diff.act}, all rational differential forms having 
poles of order at most $j$ along $Z_f$ can be obtained by $(j-1)$-times 
derivations form rational differential forms of order at most $1$ along 
$Z_f$. Therefore, all these forms belong to ${\cal F}^{n+1 -j}$. }
\label{e-filt}
\end{rem}
\medskip
 
Using the fact that the action of the differential 
operators ${\cal D}_{a_m}$  on 
${\rm Gr}_{\cal F}$ is the algebraic multiplication (cf. \ref{opr.diff}), 
one gets.  

\begin{prop}
The differential of the mixed period mapping  for $PH^{n-1}(Z_f)$ 
\[ v = \oplus_p v_p \; : \; \Theta_{\lbrack f \rbrack} \rightarrow 
\bigoplus_p {\rm Hom }( {\cal F}^{p+1}/{\cal F}^p, 
{\cal F}^{p-1} / {\cal F}^p ) \]
is induced by the $R_f$-module multiplication 
\[ R^1_f \otimes R^+_f \rightarrow R^{+}_f . \]
\label{infinitesimal}
\end{prop}

\section{Calabi-Yau hypersurfaces in toric varieties}


\hspace*{\parindent}

Recall the definition of the dual polyhedron.  

\begin{opr}
{\rm Let $\Delta$ be an $n$-dimensional convex polyhedron in ${\bf R}^n$. 
Assume that the zero point $(0, \ldots, 0) \in 
{\bf R}^n$ is contained in the interior of $\Delta$. Then the dual polyhedron 
$\Delta^* \subset {\bf R}^n$ is defined by 
\[ \Delta^* : = \{ (x_1, \ldots, x_n) \in {\bf R}^n \mid 
\sum_{i =1}^n x_iy_i \geq -1\;\; {\rm for}\;\;{\rm any}\;\;  (y_1, \ldots, y_n) 
\in \Delta \}. \]}
\label{dual.pol}
\end{opr}   

Obviously, one has $(\Delta^*)^* = \Delta$. Define the notion of reflexive 
polyhedron by the following theorem. 

\begin{theo}
Let $\Delta$ be a polyhedron as in {\rm \ref{dual.pol}}.
 Assume that $\Delta$ is 
integral, i.e., all vertices of $\Delta$ belong to the integral lattice 
${\bf Z}^n$. Let $f$ be a $\Delta$-regular Laurent polynomial 
with the Newton polyhedron $\Delta$. Then the following conditions are 
equivalent

{\rm (i)} the dual polyhedron $\Delta^*$ is integral$;$  

{\rm (ii)} for any hyperplane 
\[ {\cal H} = \{  (x_1, \ldots, x_n) \in {\bf R}^n \mid 
\sum_{i =1}^n a_iy_i = b \;\;a_i, b \in {\bf Z}, \; b >0 \; {\rm g.c.d.} 
\{a_1, \ldots, a_n, b\} = 1 \} \]
such that $\Delta \cap {\cal H}$ is a face of codimension $1$ of $\Delta$, 
one has $b =1$$;$

{\rm (iii)} the $\delta$-vector of the polyhedron $\Delta$ satisfies the 
symmetry $\psi_i(\Delta) = \psi_{n - i} (\Delta)$ $(0 \leq i \leq n)$$;$  

{\rm (iv)} ${\bf P}_{\Delta}$ is a toric Fano variety with Gorenstein 
singularities, and ${\cal O}_{\Delta}(1)$ is the 
anticanonical sheaf on ${\bf P}_{\Delta}$$;$

{\rm (v)} the closure $\overline{Z}_f$ of $Z_f$ in ${\bf P}_{\Delta}$ is 
a Calabi-Yau variety with only canonical singularities$;$ 

{\rm (vi)} ${\rm dim}_{\bf C} R^{n} = 1$ and the multiplication law on 
the Jacobian ring $R_f$ defines for any $i\; (0 \leq i \leq n)$ 
the perfect pairing 
\[ R_f^{i} \times R^{n-i}_f \rightarrow R_f^{n}. \]
\label{theo.reflex}
\end{theo}

\Proof The equivalence (i)$\Leftrightarrow$(ii)$\Leftrightarrow$(iii) 
was proved by Hibi in \cite{hibi2}. The equivalence 
(ii)$\Leftrightarrow$(iv)$\Leftrightarrow$(v)  was proved by author in 
\cite{bat.mir}. Since $\psi_i (\Delta) = {\rm dim}_{\bf C} R^{i}$, one has 
(vi)$\Rightarrow$(iii). 

Let us deduce (vi) from other conditions. 
It follows from (ii) that the ideal 
$I^{(1)}_{\Delta} \subset S_{\Delta}$ is principal (it is generated by 
$X_0$). Therefore, $S_{\Delta}$ is a Gorenstein graded ring. 
Since $R_f$ is a quotient of the Gorenstein $S_{\Delta}$ by a regular 
sequence,  $R_f$ is also a Gorenstein 
graded ring. It follows from (iii) that ${\rm dim}_{\bf C} R^{n} = 1$. 
Thus $R_f^{n}$ is the socle of $R_f$, or the annulet of the 
maximal ideal $R_f^+$ in $R_f$. The statement (vi) follows now 
from the well-known duality for Artinian Gorenstein rings  \cite{bass}.  
\bigskip

\begin{opr}
{\rm An integral polyhedron $\Delta$ satisfying the equivalent conditions 
in \ref{theo.reflex} is called {\em a reflexive polyhedron}. }
\end{opr}

We obtain now the first application of our main result. 
\medskip

\begin{theo}
Let $\Delta$ be a reflexive polyhedron. Then for the family of 
$\Delta$-regular affine Calabi-Yau hypersurfaces $Z_f \subset {\bf T}$ 
the infinitesimal 
mixed period mapping is injective, i.e. infinitesimal mixed 
Torelli theorem holds.
\end{theo}  

\Proof The statement immediatelly follows from \ref{infinitesimal} and 
\ref{theo.reflex}(vi). 
\bigskip

Consider now a special type of reflexive polyhedra, so called {\em Fano 
polyhedra}. 

\begin{opr}
{\rm An $n$-dimensional  reflexive polyhedron $\Delta$ is called a 
{\em Fano polyhedron} if it satisfies the following conditions 

{\rm (i)} $\Delta$ is simplicial, i.e., all faces  of $\Delta$  are 
simplices$;$ 

{\rm (ii)} vertices of any $(n-1)$-dimensional face $\Delta' \subset 
\Delta$ form a ${\bf Z}$-basis of the lattice ${\bf Z}^n$. }
\end{opr}

\begin{rem} 
{\rm Note that $\Delta$ is a Fano polyhedron if and only if the 
dual polyhedron $\Delta^*$ defines a {\em smooth} toric Fano variety 
${\bf P}_{\Delta^*}$. The classification of all $n$-dimensional Fano polyhedra 
is known for $n \leq 4$. In particular, there exist $123$ different types 
of $4$-dimensional Fano polyhedra \cite{bat11}.}
\end{rem}  

Recall the notion of the Stanley-Reisner ring $R_{\Sigma}$ 
of a simplicial complex 
$\Sigma$. 

\begin{opr}
{\rm Let $\Sigma$ be a simplicial complex with vertices $v_1, \ldots, v_p$. 
Consider the polynomial ring ${\bf C} \lbrack Y_1, \ldots, Y_p \rbrack$. 
Establish one-to-one correspondence $v_i \leftrightarrow Y_i$ between  
the vertices of $\Sigma$ and the variables in $P$. Let $I_{\Sigma}$ be 
the homogeneous ideal in $P$ generated by all  monomials 
$Y_{i_1} \cdots Y_{i_k}$ such that the set $\{v_{i_1}, \ldots, v_{i_k} \}$ 
does not belong to any simplex in $\Sigma$. The quotient ring 
$R_{\Sigma} = {\bf C} \lbrack Y_1, \ldots, Y_p \rbrack/I_{\Sigma}$ is called 
{\em the Stanley-Reisner ring associated with 
the simplicial complex} $\Sigma$.}
\end{opr}

Our next purpose is the following theorem. 

\begin{theo}
Let $\Delta$ be a Fano polyhedron. Then the quotient 
$S_{\Delta}/ I_{\Delta}^{(1)} = S_{\Delta}/(X_0)$  is isomorphic to the 
Stanley-Reisner ring $R_{\partial \Delta}$ 
of the simplicial complex $\partial \Delta$, where 
$\partial \Delta$ denotes the $(n-1)$-dimensional boundary of $\Delta$. 
\label{coin}
\end{theo}

\begin{lem}
If $\Delta$ is a Fano polyhedron, then $S_{\Delta}$ is generated by 
elements of degree $1$.  
\end{lem}

\Proof There exists the canonical triangulation of $\Delta$ into the union of 
$n$-dimensional simplices with volume $1/n!$. Every simplex $\sigma(\Delta')$ 
of this triangulation is the convex hull of $\{0\}$ and a 
$(n-1)$-dimensional face $\Delta' \subset \Delta$.  The graded rings 
$S_{\sigma(\Delta')}$ are polynomial rings in $n+1$ variables. 
The graded ring $S_{\Delta}$ is the sum of the subrings  
$S_{\sigma(\Delta')}$. Since the polynomial rings are generated by elements of 
degree $1$, $S_{\Delta}$ is also generated by elements of degree $1$. 
\medskip

\begin{coro}
Let $v_1, \ldots, v_p$ be vertices of a Fano polyhedron $\Delta$. Then there 
exists the canonical surjective homomorphism of the polynomial ring 
\[ h\;:\; 
{\bf C} \lbrack Y_0, Y_1, \ldots, Y_p \rbrack \rightarrow S_{\Delta} \]
such that 
\[ h(Y_0) = X_0, \;\; h(Y_i) = X_0X^{v_i}\;\; (1 \leq i \leq p).  \]
\end{coro}

Let ${\cal J}_{\Delta}$ denotes the kernel of $h$. We describe now minimal 
generators of the ideal ${\cal J}_{\Delta}$. Our description is based on 
methods in \cite{bat}.  

\begin{opr}
{\em Let $\Delta$ be a Fano polyhedron. 
A nonempty subset ${\cal P} = \{ v_{i_1}, \ldots, v_{i_k} \}$ of vertices of 
$\Delta$ is called a {\em primitive collection} if $ v_{i_1}, \ldots, v_{i_k}$ 
are not contained in any face of $\Delta$, while any proper subset 
${\cal P}' \subset {\cal P}$ is contained in a face of $\Delta$. }
\end{opr}

By defintion, for any primitive collection  ${\cal P} = 
\{ v_{i_1}, \ldots, v_{i_k} \}$ the integral vector 
\[ s({\cal P}) =   v_{i_1} + \ldots + v_{i_k}  \]
is  in the interior of the polyhedron $k\Delta$. Let $\kappa({\cal P})$ be the 
minimal nonnegative integer $k$ such that $s({\cal P}) \in k\Delta$. 
Then $s({\cal P})$ belongs to $\kappa({\cal P})\Delta'$, where $\Delta'$ is a 
proper face of $\Delta$. We can represent $s({\cal P})$ as an integral 
linear combination with positive integral coefficients 
\[ a_1 v_{j_1} + \cdots + a_l v_{j_l},\]
where $v_{j_1}, \ldots,  v_{j_l}$ are vertices of $\Delta'$, and 
$a_1  + \cdots + a_l = \kappa({\cal P})$. Thus we obtain 
the integral linear relation
\[ v_{i_1} + \ldots + v_{i_k} = a_1 v_{j_1} + \cdots + a_l v_{j_l} \]
which gives rise to an element 
\[ J({\cal P}) = Y_{i_1}  \ldots  Y_{i_k} - Y_0^{k - \kappa({\cal P})} 
Y_{j_1}^{a_1}  \cdots  Y_{j_l}^{a_l} \]
of the ideal ${\cal J}_{\Delta}$. 

\begin{opr}
{\rm The element $J({\cal P})$ is called the {\em primitive relation} 
corresponding to the primitive collection ${\cal P}$.}
\end{opr}

\begin{lem}
The elements $\{ J({\cal P}) \}$, where ${\cal P}$ runs over all primitive 
collections of vertices of $\Delta$, are generators of the 
homogeneous ideal  ${\cal J}_{\Delta}$. 
\label{gener}
\end{lem} 

\Proof The homogeneous ideal of a projective toric variety 
${\bf P}_{\Delta}$ has a ${\bf C}$-basis consisting of elements 
\[ J = Y_{j_1}^{b_1}  \cdots  Y_{j_r}^{b_r} - 
Y_{i_1}^{c_1}  \cdots  Y_{i_t}^{c_t},\; \] 
such that $b_1 + \cdots b_r = c_1 + \cdots + c_t = {\rm deg}\, J$, 
and 
\[ b_1v_{j_1}  + \cdots b_r v_{j_r} = c_1 v_{i_1} + \cdots + c_t v_{i_t}. \]

Let $\Delta'$ be the convex hull of vertices $v_{j_1},  \ldots,  v_{j_r}$, 
$\Delta''$ be the convex hull of vertices $v_{i_1},  \ldots,  v_{i_t}$. 
Assume that $\Delta'$  is a $(r-1)$-dimensional face  of $\Delta$. 
Since $b_1v_{j_1}  + \cdots b_r v_{j_r} $ is 
in the interior of $\Delta'$, and 
$c_1 v_{i_1} + \cdots + c_t v_{i_t}$ is in the interior of $\Delta''$, 
the face  $\Delta'$ must coincide with $\Delta''$. Contradiction. Therefore 
no one of two sets 
\[ V_1 = \{v_{j_1},  \ldots,  v_{j_r}\} ,\; V_2 = 
\{ v_{i_1}, \ldots, v_{i_t} \} \]
is not contained in any proper face of $\Delta$. Thus each of these two 
sets contains either a primitive collection, or the zero vector $v_0$.

We will use now induction on the degree of $J$. If the  both sets 
$V_1$ and $V_2$ contain $v_0$, then $J = Y_0^aJ'$ such that 
${\rm deg}\, J' < {\rm deg}\, J$. In this case we are done by induction. 

Assume that $V_1$ does not contain $v_0$. 
Using the existence of a primitive collection in $V_1$, 
we can find a monomial 
$Y_{j'_1}^{b'_1}  \cdots  Y_{j'_l}^{b'_l}$ such that 
\[ Y_{j_1}^{b_1}  \cdots  Y_{j_r}^{b_r} - 
Y_{j'_1}^{b'_1}  \cdots  Y_{j'_l}^{b'_l} \]
is divided by an element $J({\cal P})$ for some primitive collection 
${\cal P}$ in $V_1$, and the set $V_1' = 
\{v_{j'_1},  \ldots,  v_{j'_l}\}$ contains $v_0$. 

Repeating the same arguments for $V_2$, we find (if necessary) a monomial 
$Y_{i'_1}^{c'_1}  \cdots  Y_{i'_u}^{c'_u}$ such that  
\[ Y_{i'_1}^{c'_1}  \cdots  Y_{i'_u}^{c'_u} - 
Y_{i_1}^{c_1}  \cdots  Y_{i_t}^{c_t} \]
is divided by an element $J({\cal P}')$ for some primitive collection 
${\cal P}'$ in $V_2$, and the set $V_2' = 
\{v_{i'_1},  \ldots,  v_{i'_u}\}$ contains $v_0$. 

We have obtained a new element 
\[ J' = Y_{j'_1}^{b'_1}  \cdots  Y_{j'_l}^{b'_l} - 
Y_{i'_1}^{c'1}  \cdots  Y_{i'_u}^{c'_u} \in {\cal J}_{\Delta} \]
such that $J - J'$ belongs to the ideal generated by $J({\cal P})$, and 
$J' = Y_0^aJ''$ for some $a >0$. We can again apply induction. 
\bigskip

{\bf Proof of theorem \ref{coin}}. By definition, the ideal $I_{\partial \Delta}$ is 
generated by all monomials 
$Y_{i_1} \cdots Y_{i_k}$ such that the set $\{v_{i_1}, \ldots, v_{i_k} \}$ 
is a primitive collection. On the other hand, by lemma \ref{gener}, 
these monomials are also 
generators of ${\cal J}_{\Delta}$ modulo $Y_0$. Since $h(Y_0) = X_0$, we 
obtain the isomorphism $R_{\partial \Delta} \cong S_{\Delta}/(X_0)$.

\begin{coro}
Let $\Delta$ be a Fano polyhedron, $\{ v_1, \ldots, v_p \}$ 
vertices of $\Delta$. Consider the Laurent polynomial  
\[  f^0(X) =\sum_{i =1}^p X^{m_i}. \]
Then the elements 
\[ X_0, X_0f_1^0, \ldots, X_0f_n^0 \;\; 
(f_i^0 = X_i \frac{\partial f^0}{\partial X_i})\]
 form a regular sequence in $S_{\Delta}$, and the quotient ring 
\[ S_{\Delta} / (X_0, X_0f_1^0, \ldots, X_0f_n^0) \]
is isomorphic to the cohomology ring 
\[ H^*({\bf P}_{\Delta^*}, {\bf C}) = \oplus_{i=0}^n 
H^{2i}({\bf P}_{\Delta^*}, {\bf C}).\]
\label{chow}
\end{coro}

\Proof The statement follows from the description of the cohomology ring 
of smooth projective toric varieties \cite{dan2,oda1} 
and from the identification of Laurent polynomials 
$f_1^0, \ldots, f_n^0$ with $n$ independent  
Picard relations between $p$ irreducible divisors on ${\bf P}_{\Delta^*}$ 
corresponding to vertices of $\Delta$.

\begin{theo}
Let $\Delta$ be a Fano polyhedron, $\{ v_1, \ldots, v_p \}$ 
vertices of $\Delta$. Consider one-parameter family of Laurent polynomials 
\[ f({\lambda},X) =  f^0(X) - \lambda  , 
\;\; \lambda \in {\bf C}. 
\]
Then for almost all values of $\lambda$ the polynomial $f({\lambda},X)$ is 
$\Delta$-regular. 
\label{fano.regul}
\end{theo}

\Proof By \ref{chow}, the polynomials  
\[ X_i \frac{\partial}{\partial X_i}X_0 f(\lambda, X) = X_0f_i^0(X) \;\; 
(1 \leq i \leq n) \]
form a regular sequence in $S_{\Delta}$. Since the sections 
\[ X_0, X_0f_1^0, \ldots X_0f_n^0 \in 
H^0({\bf P}_{\Delta}, {\cal O}_{\Delta}) \]
have no common zeros in ${\bf P}_{\Delta}$, there are only finitely many 
common zeros $\{ z_j \}$ of the sections  
\[  X_0f_1^0, \ldots X_0f_n^0, \]
and all these zeros are contained in the torus 
${\bf T} = \{ x \in {\bf P}_{\Delta} \mid 
X_0 = 0\}$. 
Note that the number of the common zeros 
$\{ z_j \}$ with their multiplicities is equal to 
the degree $d(\Delta)$ of ${\bf P}_{\Delta}$.  If we choose 
$\lambda$ such that 
\[ \lambda \neq \sum_{i =1}^p  Y^{v_i},  \;\; Y \in \{z_j \}, \]
then the polynomial $f(\lambda, X)$ will be $\Delta$-regular. 

\begin{exam}
{\rm Let $\Delta$ be the convex hull of $e_1, \ldots, e_n, 
- (e_1 + \cdots + e_n)$, where $e_1, \ldots e_n$ is an integral basis 
of the lattice ${\bf Z}^n$. Consider the polynomial 
\[ f(\lambda, X) = X_1 + \cdots + X_n + (X_1 \cdots X_n)^{-1} - 
\lambda.\]
Applying arguments from  \ref{fano.regul}, we obtain that 
$f(\lambda, X)$ is not $\Delta$-regular if and only if 
\[ \lambda \neq (n+1) e^{{2\pi ik}/{n+1}}, \;\; k = 0, 1, \ldots, n. \]}
\label{simp}
\end{exam}  
\medskip

Let us investigate the situation  when $\lambda \rightarrow \infty$. Put $f = 
f(\lambda, X)$. 
\medskip

By  theorem \ref{ideal.filt}, the  subspace 
\[ H_f = {\rm Gr}_{\cal F}{\cal W}_{n-1}PH^{n-1}(Z_f) 
\subset {\rm Gr}_{\cal F}PH^{n-1}(Z_f) \] 
is generated by $X_0$ in $R_f$. Let 
\[{\rm Ann}\,(X_0) =   \{ r \in R_f \mid r \cdot X_0 = 0 \}.  \] 
Then $H_f \cong R_f / {\rm Ann}\,(X_0)$. 
Take  $X_0f^0(X)  = \sum_{i =1}^p X_0 X^{v_i}$. Define also 
\[ {\rm Ann}\,(X_0f^0(X)) = \{ r \in R_f \mid r \cdot X_0f^0(X) = 0 \}.  \] 
Denote by $G_f$ the ideal in $R_f$ generated by $X_0f^0(X)$. 
Since $X_0f(X) = X_0f^0(X) - \lambda X_0$ represents zero element 
in $R_f$, one has 
\[  {\rm Ann}\,(X_0f^0(X)) = 
{\rm Ann}\,(X_0)\;\; {\rm for}\;\; \lambda \neq 0, \]
i.e., $G_f = R_f/ {\rm Ann}\, (X_0f^0(X)) \cong H_f$.  
We have obtained that $G_f$ is isomorphic to 
the chiral ring of the infinitesiomal deformation of 
the pure Hodge structure in the middle cohomology group of affine 
Calabi-Yau hypersurfaces $Z_f$. 

On the other hand, the element $X_0f^0(X)$ represents anticanonical class 
in the cohomology ring $H^*({\bf P}_{\Delta^*}, {\bf C})$ of the smooth 
toric Fano variety ${\bf P}_{\Delta^*}$. 
Define 
\[ \tilde{\rm Ann}\,(X_0f^0(X)) = \{ s \in S_{\Delta} \mid   
s \cdot X_0f^0(X)  \in (X_0, X_0f_1^0, \ldots, X_0f_n^0) \}. \]
By \ref{chow},  elements of the ring 
\[ S_{\Delta} / \tilde{\rm Ann}\, (g) \]
are classes of  algebraic cycles on the mirror  
Calabi-Yau hypersurfaces $\overline{Z}$ 
in the toric Fano manifold ${\bf P}_{\Delta^*}$. When $\lambda \rightarrow 
\infty$, the closures $\overline{Z}_{f(\lambda,X)}$ of the one-parameter 
family of Calabi-Yau hypersurfaces in ${\bf P}_{\Delta}$ degenerate to the 
union ${\rm Proj}\, S_{\Delta}/(X_0)$ 
of $d(\Delta)$ copies of $(n-1)$-dimensional projective spaces. 
Thus we obtain.  

\begin{theo}
The cohomology ring $S_{\Delta} / \tilde{\rm Ann}\, (X_0f^0(X))$ 
generated by restrictions from ${\bf P}_{\Delta^*}$ of algebraic cycles to 
smooth Calabi-Yau hypersurfaces $\overline{Z}$ 
is a "limit" $(\lambda \rightarrow \infty)$ 
of the infinitesimal deformation ring 
for the mirror one-parameter family of Calabi-Yau hypersurfaces in 
${\bf P}_{\Delta}$ defined by the Laurent polynomials $f(\lambda, X)$. 
\label{l.limit}
\end{theo}
\bigskip

In \cite{bat.mir} we have proved that for $n \geq 4$  the dimension of 
the Hodge-Deligne space $H^{n-2,1}(Z_f)$ of affine Calabi-Yau hypersurfaces 
in ${\bf T}$ coincides with the dimension of the Hodge space $H^{1,1}$ 
of a maximal crepant resolution $\hat{Z}^*$ 
of the mirror Calabi-Yau hypersurfaces 
in ${\bf P}_{\Delta^*}$. Using the description of the Hodge-Deligne subspace 
$H^{n-2,1}(Z_f)$ in \ref{hodge.sub}, we get the following 
more precise information. 

\begin{theo}
For Calabi-Yau hypersurfaces $Z_f \subset {\bf T}$, one has the 
following short exact sequence
\[  0 \rightarrow \bigoplus_{{\rm codim}\, \Delta' =1} L^*(\Delta')  
\rightarrow L(\Delta)/{\bf C}\langle f, f_1, \ldots, f_n \rangle 
\rightarrow H^{n-2,1}(H^{n-1}(Z_f)) \rightarrow 0, \]
where ${\bf C}\langle f, f_1, \ldots, f_n \rangle$ is the {\bf C}-space 
spanned by  the polynomials $f, f_1, \ldots, f_n$. 
\end{theo}

On the other hand, for $n \geq 4$ 
it was proved  in \cite{bat.mir} the following.

\begin{theo}
For maximal crepant resolutions $\hat{Z}$ of mirror Calabi-Yau hypersurfaces, 
one has the following short exact sequence 
\[  0 \rightarrow \bigoplus_{{\rm codim}\, \Delta' =1} L^*(\Delta')  
\rightarrow L(\Delta)/{\bf C}\langle 1,  f_1^0, \ldots, f^0_n \rangle 
\rightarrow H^{1,1}(H^{n-1}(\hat{Z})) \rightarrow 0, \]
where 
\[ f_i^0(X) = \frac{\partial}{\partial X_i}f^0(X), \;\;(1 \leq i \leq n),  \]
\[ f^0(X) = \sum_{m \partial \Delta \cap {\bf Z}^n} X^m. \]
\end{theo}

\begin{coro}
The ${\bf C}$-space $H^{1,1}(H^{n-1}(\hat{Z}))$ can be considered as 
a "degeneration" of the ${\bf C}$-spaces  $H^{n-2,1}(H^{n-1}(Z_f))$.
\end{coro}

\section{Moduli spaces of Calabi-Yau hypersurfaces}


\hspace*{\parindent}

We begin with general remarks on properties of the  mirror symmetry 
for Calabi-Yau varieties. 
\medskip

Let ${\cal M}$ be a moduli space of the complex structure on $n$-dimensional 
Calabi-Yau variety $W$, ${\cal M}'$ a  moduli space of the complex structure 
on the $n$-dimensional mirror Calabi-Yau variety $W'$. The general phylosophy 
is that the infinitesimal approximation  of ${\cal M}'$ at generic point 
is the complexified Picard space  ${\rm Pic}(W)\otimes {\bf C}$ of $W$, or 
to the complexified N\'eron-Severi space $NS(W)_{\bf C}$.  

In general there might be rather complicated singularities at 
the boundary of the moduli space ${\cal M}'$. However, 
there are known {\em toriodal compactifications} 
(due to Mamford et. al.) of important classes of moduli spaces. 
The idea of such compactification bases on analitical  descriptions  of  
singularities 
at boundary points in terms of appropriate convex cones. 
This observation leads us to the following question:

\begin{ques}
Is it  possible also to find a natural projective model  of 
the moduli space ${\cal M}'$ 
in such a way that the singularities at boundary points of 
${\cal M}'$ are analitically isomorphic to toric singularities and 
therefore admit a description  by some naturally defined convex cones $?$   
\end{ques}

From viewpoint of the mirror symmetry phylosophy, it is important to know 
singularities arising at very special (cusp) points of the moduli 
space ${\cal M}'$. These cusp poins correspond to  "maximal degenerations"  of 
Calabi-Yau varieties, or zero-dimensional strata in the case when 
the compactified moduli space admits a natural stratification. 
Now we formulate our main conjecture:

\begin{conj}
Singularities at cusp points of "maximal degeneration" 
of a "good" moduli space the Calabi-Yau mirrors $W'$ should 
have  toroidal structure defined by "natural" subcones of the cone of 
effective divisors in $NS( W)_{\bf C}$. 
\label{conj.cone}
\end{conj}

Althougth our conjecture is not precisely formulated, the main purpose of 
this section is to show that this conjecture 
is true for some moduli spaces of Calabi-Yau hypersurfaces in toric varieties. 
\bigskip

First we consider the moduli space ${\cal M}(\Delta)$ 
of affine Calabi-Yau hypersurfaces defined by an $n$-dimensional 
reflexive polyhedron $\Delta$. 
We know that ${\cal M}(\Delta)$ is a projective toric variety of dimension 
$l(\Delta) - (n+1)$. Therefore, by definition, all singularities of 
${\cal M}(\Delta)$ are toric, i.e., these singularities are defined by 
some convex cones. The only problem is to find a geometric interpretation of 
these cones.  
\medskip

Recall the explicit construction of the toric variety ${\cal M}(\Delta)$.  
\medskip

Let $s = l(\Delta)$, and let  $m_0, m_1, \ldots, m_{s-1}$ be integral 
points in ${\Delta}$ such that $m_0$ is the zero of $M_{\bf R}$. 
Since the action of ${\bf T}$ on 
$L(\Delta)$ is diagonal, 
every monomial $a_{m_0}^{k_0} \cdots a_{m_{{s-1}}}^{k_{s-1}}$ generates a 
$1$-dimensional ${\bf T}$-invariant subspace in 
${\bf C} \lbrack a_{\Delta} \rbrack$  with the ${\bf T}$-action by 
the character $\chi_m\;:\; {\bf t} \rightarrow {\bf t}^m$, 
where $m = {k_0}{m_0} + \cdots + {k_{s-1}}{m_{s-1}} \in {\bf Z}^n$. 
So the ${\bf T}$-invariant subring ${\bf C} \lbrack \alpha  
\rbrack^{\bf T}$ in  ${\bf C} \lbrack \alpha \rbrack$ 
consists of monomials 
\[ \{ a_{m_0}^{k_0} \cdots a_{m_{s-1}}^{k_{s-1}} \; \mid \;  
{k_0}{m_0} + \cdots + {k_{s-1}}{m_{s-1}} = 0,\; k_i \geq 0, \; m_i \in 
\Delta \}. \]

Let $g_{\Delta}\, :\, {\bf R}^{s} \rightarrow M_{\bf R}$ be the 
linear map 
\[ g_{\Delta}(x_0, \ldots, x_{s-1}) =  {x_0}{m_0} + \cdots + 
{x_{s-1}}{m_{s-1}}. \]
The kernel 
$W_{\Delta}$ of $g_{\Delta}$ is a $(s-n)$-dimensional rational subspace in 
${\bf R}^{s}$ containing the $1$-dimensional subspace $W_{\Delta}^0 \subset 
{\bf R}^{s}$ generated by $(1,0, \ldots, 0)$. We can define a rational 
polyhedral cone ${\sigma}_{\Delta}$ in $W_{\Delta}$,  
\[ \sigma_{\Delta} = \{  (x_0, x_1, \ldots, x_{s-1} ) \in W_{\Delta} 
 \mid x_i \geq 0\;\;( 0 \leq i \leq s -1 ) \}. \] 
Monoid of integral points in ${\sigma}_{\Delta}$ is a submonoid of 
the free monoid ${\bf N}^{s}$. Thus, we can describe 
${\bf C} \lbrack \alpha  \rbrack^{\bf T}$  as 
a semigroup algebra of the monoid of integral points in the 
convex cone ${\sigma}_{\Delta}$, i.e., ${\cal M}(\Delta)$ is a projective 
toric variety such that affine cone over ${\cal M}(\Delta)$ is 
the affine toric variety of dimension $s-n$  corresponding to the 
rational polyhedral cone $\sigma_{\Delta}$. Let us study more  
the cone $\sigma_{\Delta}$. 

Note that powers of $a_0$ are always ${\bf T}$-invariant.  Thus, in order to 
describe integral points in  $\sigma_{\Delta}$,  
it is sufficient to consider only linear relations  with non-negative 
integral coefficients $k_i$ between $s-1$ integral points $m_i$ lying on the 
boundary $\partial \Delta$ of the polyhedron $\Delta$.

We can restrict $g_{\Delta}$ on 
\[ {\bf R}^{s-1} = \{  (x_0, x_1, \ldots, x_{s-1}) \in  {\bf R}^{s}  
\mid x_0 = 0 \} \] 
and obtain  
\[ g_{\partial \Delta}\; :\; {\bf R}^{s-1} \rightarrow M_{\bf R} \]
whose kernel we denote by $W_{\partial \Delta}$. Let 
${\sigma}_{\partial \Delta} = \sigma_{\Delta} \cap W_{\partial \Delta}$, 
then 
\[ {\sigma}_{\Delta} = {\bf R}_{\geq 0} \oplus {\sigma}_{\partial \Delta}. \]
Therefore, ${\bf C} \lbrack a_{\Delta} \rbrack^{\bf T}$  is 
a polynomial  ring in one variable $a_0$ over 
the semigroup algebra  $S({\partial \Delta})$ 
of the monoid of integral points in the cone ${\sigma}_{\partial \Delta}$. 
\bigskip

We come to the following.  

\begin{theo} 
The moduli space ${\cal M}(\Delta)$ of affine Calabi-Yau 
hypersurfaces is a projective cone 
over projective toric variety ${\rm Proj}\, S({\partial \Delta})$ with 
the vertex $p_0 \in {\cal M}(\Delta)$ defined by the conditions 
$a_0 =1, a_{m_i} = 0$ $(m_i \in \partial \Delta)$. In particular, 
for every point $p \in {\cal M}(\Delta)$ $(p \neq p_0)$, there exists 
a unique projective line joining $p$ and $p_0$.   
\end{theo}

\begin{coro}
The condition $a_0 \neq 0$ defines an affine toric chart 
in ${\cal M}(\Delta)$ containing $p_0$. This chart is 
isomorphic to ${\rm Spec}\, S({\partial \Delta}).$ 
\end{coro} 

The next statement is very important. 

\begin{theo}
Let $\Delta^*$ be the dual reflexive polyhedron, ${\bf P}_{\Delta^*}$ 
the mirror toric Fano variety. Then the cone ${\sigma}_{\partial \Delta}$ is 
dual to the cone of effective divisors 
of any maximal crepant partial desingularization of ${\bf P}_{\Delta^*}$. 
\end{theo} 

\Proof The integral points $m_1, \ldots m_{s-1}$ correspond to invariant 
toric divisors on a maximal crepant resolution of ${\bf P}_{\Delta^*}$. 
These divisors generate the cone of effective divisors in 
$NS({\bf P}_{\Delta^*})_{\bf R}.$ 
On the other hand,  the ${\bf R}$-space  
$NS({\bf P}_{\Delta^*})_{\bf R}$ is dual to the space of all 
${\bf R}$-linear relations among $m_1, \ldots m_{s-1}$. So the cone 
of the effective divisors is dual to the cone of all linear 
relations with  {\em non-negative} real coefficients among 
$m_1, \ldots m_{s-1}$. 
\medskip

There is a nice construction due to Oda and Park \cite{oda.park} 
of the  subdivision 
of the cone of effective divisors of the toric Fano variety 
${\bf P}_{\Delta^*}$ (in fact, the Oda-Park method works for an arbitrary 
projective toric variety). By standard toric arguments (cf. \cite{dan2}), 
any subdivision of the cone of effective divisors of ${\bf P}_{\Delta^*}$ 
defines a partial toric desingularization of ${\cal M}(\Delta)$. 
In the  Oda-Park construction the subcones of maximal dimension 
in $\sigma_{\partial \Delta}$ one-to-one correspond to 
projective ${\bf Q}$-factorial birational models of ${\bf P}_{\Delta^*}$. 
This means that these models define zero-dimensional strata on the 
partial resolution. So we come to the following.  

\begin{theo}
The Oda-Park construction of the subdivision of the cone of 
effective divisor $\sigma_{\partial \Delta}$ define a partial 
toric desingularization $\overline{\cal M}(\Delta)$ 
of the moduli space ${\cal M}(\Delta)$ such that 
maximal crepant desingularizations of ${\bf P}_{\Delta^*}$ one-to-one 
correspond to some zero-dimensional toric strata lying over $p_0$. 
\label{toric.o-p}
\end{theo}

From geometrical point of view, this partial desingularization is 
important to understand the local geometry of ${\cal M}(\Delta)$ at the 
point $p_0$. It was obsereved by Aspinwall and L\"utken in \cite{aspin3} that 
the notion of "large complex structure limit" depends not only on the 
choice of the special point (e.g., $p_0 \in {\cal M}(\Delta)$) 
on the moduli space, but also on the {\em direction} which we use to come to 
this point. This remark immediately leds us to a consideration  of 
a new  birational model  which looks like "blow up" at the special point of 
the original moduli space of Calabi-Yau varieties. 
\medskip

We can describe local properties of the moduli space 
${\cal M}(\Delta)$ at the point $p_0$ in terms of the 
principal dicriminant ${\rm Disc}_{\Delta}\,(f)$. From 
\cite{gel1} we obtain.

\begin{theo}
Let $\Delta$ be a reflexive polyhedron. Then 

{\rm (i) } the multiplicity of the principal 
discriminant  at $p_0$ equals 
$n d(\Delta)$ $($equivalently, the maximal power of $a_0$ in 
${\rm Disc}_{\Delta}\,(f)$ 
equals $d(\Delta))$ $;$

{\rm (ii)} the number of vertices of the Newton polyhedron of 
${\rm Disc}_{\Delta}\,(f)$ 
corresponding to monomials with the maximal power of $a_0$ equals 
the number of different maximal crepant desingularizations 
of the mirror  toric Fano variety ${\bf P}_{\Delta^*}$, 
where $\Delta^*$ is the dual reflexive polyhedron. 
\end{theo}
\medskip

Assume now that ${\bf P}_{\Delta}$ has already ${\bf Q}$-factorial 
terminal singularities, i.e., ${\bf P}_{\Delta}$ admits no nontrivial 
crepant partial resolutions. This is equivalent  to the condition that 
the dual polyhedron 
$\Delta^*$ has no faces of codimension $1$ with interior integral points 
(e.g., this is the case when $\Delta^*$ is a Fano polyhedron).  
In this case, by Lefschetz theorem,  any $\Delta$-regular 
Calabi-Yau hypersurface $\overline{Z} \in {\bf P}_{\Delta}$ 
has the same Picard group as the 
ambient toric variety ${\bf P}_{\Delta}$. Moreover, the 
cones of effective divisors of ${\bf P}_{\Delta}$ and $\overline{Z}$ are 
isomorphic. On the other hand, one can prove that the connected component 
of the authomorphism group of the toric variety ${\bf P}_{\Delta^*}$ is 
isomorphic to ${\bf T}$. Thus, up to a finite covering, the moduli 
space of the mirror {\em affine} Calabi-Yau hypersurfaces  and {\em their 
compactifications} in ${\bf P}_{\Delta^*}$ are the same.   As a result we 
come to the following. 

\begin{theo}
Assume that ${\bf P}_{\Delta}$ admits no nontrivial 
crepant partial resolutions. 
If $W'$ denotes a Calabi-Yau hypersurface in 
${\bf P}_{\Delta^*}$, and $W$ denotes a Calabi-Yau hypersurface in 
${\bf P}_{\Delta}$, then for the moduli spaces ${\cal M}(\Delta)$ and 
 $\overline{\cal M}(\Delta)$ 
the conjecture {\rm \ref{conj.cone}} is true. 
\end{theo}

It seems that the conjecture \ref{conj.cone} is also true without additional 
restrictions, but we do not know  a complete proof.   
\medskip

We want to finish this section by a short 
remark on the relation between the Minimal Model Program of Mori and 
properties of special subcones in the cone of effective divisors.  

\begin{rem}
{\rm Let $W$ and $W'$ be two birationally isomorphic 
projective algebraic varieties 
with ${\bf Q}$-factorial terminal singularities varieties. We assume 
that they are biregularly isomorphic at least in codimension $1$. Then 
we obtain an isomorphism of the N\'eron-Severi spaces $NS(W)_{\bf R} \cong 
NS(W')$. Moreover, we have also an isomorphisms  
between the cones of 
effective  divisors $NE^1(W) \subset NS(W)_{\bf R}$ and   
$NE^1(W') \subset NS(W')_{\bf R}$. But it is important to remark that 
the K\"ahler subcones $NK^1(W), NK^1(W')  \subset NE^1(W) \cong NE^1(W')$ 
{\em have a common interior point 
if and only if $W$ is biregularly isomorpic to $W'$.} 
In particular, this means that the K\"ahler subcones of different birational 
equivalent and biregularly isomorphic in codimension $1$ 
minimal models of Calabi-Yau varieties have a "fan-like"  structure: 
intersection   of any two cones is a common face of the both. This observation 
and \ref{toric.o-p} shows that the K\"aler cones of different minimal  
models  obtained by flopping transformations from a Calabi-Yau variety $W$ 
are good 
candidates for the "natural" subcones in the conjecture \ref{conj.cone}. }
\end{rem}

\section{Periods and generalized hypergeometric functions}


\hspace*{\parindent}

Recall the construction of the $A$-hypergeometic 
system due to Gelfand, Kapranov, and Zelevinsky \cite{gel0}.

Let 
\[ A = \{ v_1, \ldots, v_s \} \subset  {\bf Z}^n \]
be a finite  set of integral vectors  in  $n$-dimensional 
space ${\bf R}^n$. Using the standard basis of ${\bf Z}^n$,   one has 
\[ v_j = (v_{j1}, \ldots, v_{jn} ). \]

Consider $s$ integral 
vectors 
\[ \overline{A} = \{ \overline{v}_1, \ldots,  \overline{v}_s \} 
\subset {\bf Z}^{n+1},  \]
where $\overline{v}_j = (1, v_j)$. Define the lattice $L$ of 
integral relations 
among $\overline{v}_1, \ldots , \overline{v}_s$ as follows 
\[ L = \{ (l_1, \ldots, l_s ) \in {\bf Z}^s \mid \sum^s_{j=1} l_j 
\overline{v}_j = 0 \}. \]
We will assume that the integral vectors 
$\overline{v}_1, \ldots,  \overline{v}_s$ generate the lattice 
${\bf Z}^{n+1}$. In this case $L$ has rank $s - n -1$.

Let ${\bf C}^A$ be the affine complex space with complex 
coordinates $a_1, \ldots, a_s$.  
For any $l = (l_1, \ldots, l_s) \in L$ we define the differential operator 
\[ \Box_l = \prod_{l_j >0 } (\frac{\partial}{\partial a_j})^{l_j} - 
\prod_{l_j < 0} (\frac{\partial}{\partial a_j})^{-l_j}. \]
Take a complex vector $\beta = (\beta_0, \beta_1, \ldots, \beta_n) 
\in {\bf C}^{n+1}$.  

\begin{opr}
{\rm Let $\Phi(\alpha)$  be a complex function 
on ${\bf C}^A$, $\alpha = (a_1, \ldots, a_s)$. 
The differential system 
\[ \Box_l \Phi(\alpha) = 0, \;\; l \in L , \]
\[ \sum_{j =1}^{s} a_j \frac{\partial}{\partial a_j} \Phi(\alpha)  
= \beta_0 \Phi(\alpha) \]
\[ \sum_{j =1}^{s} v_{ji}a_j  \frac{\partial}{\partial a_j} \Phi(\alpha) 
= \beta_i \Phi(\alpha)\;\; 
(1 \leq i \leq n). \]
is called the  {\em $A$-hypergeometric system with a vector of 
exponents} $\beta$. }
\label{diff}
\end{opr}

It is known that $A$-hypergeometric system is holonomic and its solutions 
form a locally constant sheaf of finite rank outside a hypersurface in 
${\bf C}^A$. 
\bigskip

We will be interested in the case when $A = \Delta \cap {\bf Z}^n$ for 
some $n$-dimensional integral polyhedron $\Delta$. 
In this case  the differential system in \ref{diff} we call {\em 
$\Delta$-hypergeometric system}. 

With every point $\alpha = (a_1, \ldots, a_s) \in {\bf C}^{A}$ 
one can associate the Laurent polynomial 
\[ f(\alpha, X) = \sum_{j =1}^s a_j X^{v_j} \in {\bf L}. \]
By theorem of Gelfand-Kapranov-Zelevinskiy, if $\alpha \in {\bf C}^A$ 
corresponds to a $\Delta$-regular Laurent 
polynomial $f(\alpha, X)$, then the number of linearly independent solutions 
of the $\Delta$-hypergeometric system at the points $\alpha$ is equal to 
$d(\Delta)$. Therefore, outside the discriminant 
hypersurface 
\[ \{ \alpha \in {\bf C}^{A} \mid {\rm Disc}_{\Delta}\, (f(\alpha, X)) = 0 
\} \]
solutions of the $\Delta$-hypergeometric system form a locally constant 
sheaf of rank $d(\Delta)$. The following theorem shows  
the relation between 
  $\Delta$-hypergeometric functions and Euler integrals, solutions of 
$\Delta$-hypergeometric system,  
and  periods of rational forms on ${\bf T} \setminus Z_f$. 
It is a special case of the Euler integral representation in \cite{gel01}. 

\begin{theo}
Let  $\omega(\alpha, k, m_1, \ldots, m_n)$ be  a rational 
differential $n$-form 
\[  \frac{X_1^{m_1} \cdots X_n^{m_n}}
{f(\alpha, X)^k} \omega_0,\;\; (\omega_0 = 
 \frac{dX_1}{X_1} \wedge \cdots \wedge \frac{dX_n}{X_n}).  \] 
Then periods of this $n$-form 
\[ \Pi(\alpha) = \int_{\gamma} \omega(\alpha,k,m_1, \ldots, m_n), \]
where $\gamma$ is an $n$-cycle in $H_n({\bf T} \setminus Z_f)$, are solutions 
of the $\Delta$-hypergeometric system with the exponents  
\[ \beta = (-k, -m_1, \ldots, -m_n) .\]
\label{hyper.geom}
\end{theo}

\Proof Instead of applying differential operators $\partial /\partial a_j$ 
to the function  $\Pi({\alpha})$, we 
 apply these operators to the $n$-form $\omega(\alpha)$. 
Note that $\omega(\alpha) = {\cal R} (X_0^kX^m)$. Hence, applying the 
Laplace transform, it is suffices to check the action 
of ${\cal R}$-transforms of the differential operators  
\[ \Box_l \;\;(l \in L), \;\; \sum_{j =1}^{s} a_j 
\frac{\partial}{\partial a_j}, \]
and 
\[ \sum_{j =1}^{s} v_{ji}a_j  \frac{\partial}{\partial a_j} \;\; 
(1 \leq i \leq n) \]
on the class of the monomial $X_0^kX^m$ in the space 
\[ H^{n-1}({\bf T} \setminus Z_f) \cong 
S^+/ ({\cal D}_{X_0} S^+  + \sum_{i=1}^n {\cal D}_{X_i} S^+). \]

By  \ref{multi},  the differential operator 
$\partial/ \partial a_j$ corresponds to the 
operator of multiplication by $X_0X^{v_j}$ in $H^{n-1}({\bf T} \setminus Z_f)$.  
Thus, if $\sum^s_{j=1} l_j \overline{v}_j = 0$ is an integral linear 
relation among $\{ \overline{v}_j \}$, then 
\[ \sum_{l_j > 0} l_j \overline{v}_j = \sum_{l_j < 0} (-l_j) 
\overline{v}_j, \]
i.e., we obtain 
\[ \prod_{l_j > 0} (X_0X^{v_j})^{l_j} = 
\prod_{l_j < 0} (X_0X^{v_j})^{-l_j}. \]
Hence $\Box_l \Pi(\alpha) = 0$, $l \in L$.  

The differential operator 
\[ \sum_{j =1}^{s} a_j \frac{\partial}{\partial a_j} \]
corresponds to the multiplication by $F_0$ in $H^{n-1}({\bf T} \setminus Z_f)$. 
On the other hand, one has 
\[ {\cal D}_{X_0}(X_0^kX^m) = kX_0^kX^m + F_0X_0^kX^m. \]
Therefore, the monomial $X_0^kX^m$ represents 
in $H^{n-1}({\bf T} \setminus Z_f)$ 
an eigenvector with the eigenvalue $-k$. 

The differential operators  
\[ \sum_{j =1}^{s} v_{ji}a_j  \frac{\partial}{\partial a_j} \;\; 
(1 \leq i \leq n) \]
correspond to multiplications by $F_i$. The same  arguments show that the 
monomial $X_0^kX^m$ represents in $H^{n-1}({\bf T} \setminus Z_f)$ an
eigenvector with the eigenvalue $-m_i$. 
\bigskip

Consider several examples. 
\medskip

\begin{exam}
{\rm Let $\Delta$ be the $n$-dimensional reflexive polyhedron from 
the example \ref{simp}. It is shown by Gelfand et. al. that the solutions of 
the corresponding $\Delta$-hypergeometric system are well-known 
generalized hypergeometric functions ${\,}_{n+1}F_n$ (\cite{gel0}, 3.1). }
\end{exam}

\begin{exam}
{\rm Let $\Delta$ be the $2$-dimensional reflexive polyhedron with vertices 
$(\pm1, 0)$, $(0, \pm1)$. Every $\Delta$-regular Laurent polynomial 
$f \in {\bf C} \lbrack  X_1^{\pm1}, X_2^{\pm1} \rbrack$ defines an affine 
elliptic curve whose smooth compactification has $4$ more points. 
It was shown by Gelfand et.al. that the $\Delta$-hypergeometric system 
has as a solution the Appell's hypergeometric function $F_4$ \cite{gel0}. 
A variation of mixed Hodge structure for $F_4$ based on the family of the 
affine elliptic curves  was constructed explicitly  by J. Stienstra 
\cite{stien3}. }
\end{exam}

\begin{exam}
{\rm Let $\Delta$ be an $n$-dimensional  reflexive polyhedron, 
$f({\alpha},X)$ a $\Delta$-regular Laurent polynomial of the form 
\[ f({\alpha}, X) = a_0 - \sum_{m \in {\bf Z}^n \cap \partial \Delta} 
a_m X^m. \]
Consider the integral 
\[ \Pi_0 ({\alpha}) = \frac{1}{(2\pi i)^n} \int_{\gamma_0} 
\frac{a_0}{f({\alpha}, X)} \omega_0,\]
where $\gamma_0$ is the restriction of the element  
\[ \gamma_0 = \{ X = (X_1, \ldots, X_n ) \in {\bf T} : \mid X_1 \mid = \cdots 
= \mid X_n \mid = 1 \} \]
of the group $H^n({\bf T})$ to $H^n({\bf T} \setminus Z_f)$.

By theorem \ref{hyper.geom}, $\Pi_0({\alpha})$ is a hypergeometric function 
with the exponents  $\beta = (-1, 0, \ldots, 0)$. 

One can easy calculate the expansion of $\Pi_0({\alpha})$ in a 
power series for small values of $a_0^{-1}$ and 
$a_m$ ($m \in \partial \Delta$). 
Indeed, let $f_0({\alpha}, X) = 1 - f({\alpha}, X)/a_0$, then we obtain 
\[ \frac{a_0}{f(\alpha,X)} = (1 -f_0(\alpha,X))^{-1} = 
\sum_{i = 0}^{\infty} (f_0(\alpha,X))^i. \]
By the Cauchi residue formula, 
\[ B_i({\alpha}) = \frac{1}{(2\pi i)^n} \int_{\gamma_0} 
(f_0({\alpha}, X))^{i} \omega_0  \]
equals the constant coefficient in the Laurent polynomial 
\[ (\frac{1}{a_0})^i (\sum_{m \in {\bf Z}^n \cap \partial \Delta} a_m X^m)^i. \]
The function $\Pi_0({\alpha})$ has the following 
properties:

(i) $\Pi_0({\alpha})$ has {\em integral expansion coefficients} as power 
series  in $a_0^{-1}$  and $a_m$ $(m \in \partial \Delta)$; 

(ii) $\Pi_0({\alpha})$ is  {\em regular} at the boundary 
point $p_0 \in {\cal M}(\Delta)$ corresponding to values $a_0 = 1$, $a_m = 0$ 
("large complex structure limit"). 

Note that this type of hypergeometric function was used 
in the calculation of P. Candelas et. al. in \cite{cand2}.    
The hypergeometric function $\Pi_0$ for the one-paramenter family 
of Laurent polynomials 
\[ f (\psi, X) = 5\psi - X_1 -X_2 - X_3 - X_4 - (X_1 X_2 X_3 X_4)^{-1} \]
is as follows 
\[ \Pi_0(\psi) = \sum_{k =0}^{\infty} \frac{(5k)!}{(k!)^5 (5\psi)^{5k}}. \]  }
\end{exam}
\medskip

One can calculate  the power expansion of the hypergeometric 
function $\Pi_0$ by counting  integral points  in the cone $\sigma_{\Delta}$ 
dual to the cone of effective divisors of the toric Fano variety $\Delta^*$ 
as follows. 

\begin{prop}
Let $\Delta$ be a polyhedron with $r$ integral points $m_1, \ldots, 
m_r$ on its boundary. 
Then 
\[ \Pi_0 (\alpha) = 1 + \sum \frac{k!}{k_1! \cdots k_r!} a_0^{-k} 
a_{m_1}^{k_1}  \cdots a_{m_r}^{k_r}, \]
where $k = k_1 + \cdots + k_r$, and $(k_1, \ldots , k_r)$ runs over all 
non-negative integral solutions of the equation
\[ k_1 m_1 + \cdots + k_r m_r = 0. \]
\end{prop}
\medskip

\begin{rem}
{\rm The hypergeometric function $\Pi_0$ is closely related 
to the Hasse invariant of Calabi-Yau hypersurfaces $Z_f$ and the $p$-adic 
analytic representation of the unit root of the Weil zeta function of 
Calabi-Yau varieties \cite{katz1,stien2}. 
The coefficients $B_{\alpha}$ define also 
a logarithm of the formal Artin-Mazur group of the family of Calabi-Yau 
varieties \cite{stien1}. The main idea to prove this is that  for Calabi-Yau 
hypersurfaces in toric varieties there exists a canonical representative  
of the one-dimensional space $H^n ({\bf P}_{\Delta}, {\cal O}_{\Delta}(-1))$,  
a ${\bf T}$-invariant rational differential $n$-form  $\omega_0$ 
(recall that  ${\cal O}_{\Delta}(-1) \cong \Omega^n_{\bf P}$ is the 
canonical line bundle). In particular, we obtain:}
\end{rem}

\begin{prop} 
Let $\Delta$ be an $n$-dimensional reflexive polyhedron, $f \in L(\Delta)$ a 
$\Delta$-regular Laurent polynomial with coefficients in a field of 
characteristic $p$. Then 
the Hasse invariant ${\rm Hass}(f)$ is zero if and only if 
the constant coefficient in the Laurent polynomial 
$f^{p-1}$ is zero. 
\end{prop}

The proof of this theorem repeats the arguments from the book of Hartshorn 
(\cite{hartsh}, Ch.4, Prop. 4.21) based on the commutative diagram 
 
\begin{center}
\begin{tabular}{ccc}
 $H^{n-1}(\overline{Z}_f, {\cal O}_{\overline{Z}_f})$ & 
 $\stackrel{\cong}{\rightarrow}$ &
$H^{n}({\bf P}_{\Delta}, {\cal O}_{\Delta}(-1))$  \\
 ${\rm Fr}^* \downarrow$  & & ${\rm Fr}^* \downarrow$  \\
 $H^{n-1}(\overline{Z}_f^p, {\cal O}_{\overline{Z}_f^p})$ & 
 $\stackrel{\cong}{\rightarrow}$ &
$H^{n}({\bf P}_{\Delta}, {\cal O}_{\Delta}(-p))$  \\
 $ \downarrow$  & & $f^{p-1} \downarrow$  \\
 $H^{n-1}(\overline{Z}_f, {\cal O}_{\overline{Z}_f})$ & 
 $\stackrel{\cong}{\rightarrow}$ &
$H^{n}({\bf P}_{\Delta}, {\cal O}_{\Delta}(-1))$  
  \end{tabular}
  \end{center}
where ${\rm Fr}$  is the Frobenius morphism on the toric variety 
${\bf P}_{\Delta}$, $\overline{Z}_f^p$ is the subscheme 
defined by the equation $f^p =0$. We need only to note 
that the canonical generator of the space  
$H^{n}({\bf P}_{\Delta}, {\cal O}_{\Delta}(-1))$ corresponds to the unique 
integral point in the interior of $\Delta$.
\bigskip

\begin{rem}
{\rm Note that $p$-adic representation for roots of the Weil 
zeta-function has links with the  $p$-adic string theory \cite{gross}.}
\end{rem}  

\begin{rem}
{\rm A remarkable idea of S. Cecotti  \cite{ceco} that all problems 
in $N=2$ superconformal theory are Hodge-theoretical problems as consequence 
shows why  $p$-adic  analysis was recently  so 
succesfully applied to  modern theoretical physics. 
The reason is that the notion of Hodge stucture and periods 
admits $p$-adic realizations (e.g. in the crystalline 
cohomology theory, or in the  Dwork's theory of $p$-adic hypergeometric 
functions). In particular, this gives an evidence that there should be a 
general $p$-adic theory of periods  of  modular forms connected with moduli 
spaces of higher-dimensional  Calabi-Yau varieties,  a generalization of the 
classical case of modular families of  elliptic curves.}
\end{rem}